\documentstyle[11pt]{article}

\newtheorem{theorem}{Theorem}[section]
\newtheorem{lemma}[theorem]{Lemma}
\newtheorem{corollary}[theorem]{Corollary}

\newtheorem{proposition}[theorem]{Proposition}

\newtheorem{scholium}[theorem]{Scholium}

\makeindex \makeglossary
\begin{document}
%
%

\long\def\ig#1{\relax}
\ig{Thanks to Roberto Minio for this def'n.  Compare the def'n of
\comment in AMSTeX.}

\newcount \coefa
\newcount \coefb
\newcount \coefc
\newcount\tempcounta
\newcount\tempcountb
\newcount\tempcountc
\newcount\tempcountd
\newcount\xext
\newcount\yext
\newcount\xoff
\newcount\yoff
\newcount\gap%
\newcount\arrowtypea
\newcount\arrowtypeb
\newcount\arrowtypec
\newcount\arrowtyped
\newcount\arrowtypee
\newcount\height
\newcount\width
\newcount\xpos
\newcount\ypos
\newcount\run
\newcount\rise
\newcount\arrowlength
\newcount\halflength
\newcount\arrowtype
\newdimen\tempdimen
\newdimen\xlen
\newdimen\ylen
\newsavebox{\tempboxa}%
\newsavebox{\tempboxb}%
\newsavebox{\tempboxc}%

\makeatletter
\setlength{\unitlength}{.01em}%
\def\settypes(#1,#2,#3){\arrowtypea#1 \arrowtypeb#2 \arrowtypec#3}
\def\settoheight#1#2{\setbox\@tempboxa\hbox{#2}#1\ht\@tempboxa\relax}%
\def\settodepth#1#2{\setbox\@tempboxa\hbox{#2}#1\dp\@tempboxa\relax}%
\def\settokens[#1`#2`#3`#4]{%
     \def\tokena{#1}\def\tokenb{#2}\def\tokenc{#3}\def\tokend{#4}}
\def\setsqparms[#1`#2`#3`#4;#5`#6]{%
\arrowtypea #1
\arrowtypeb #2
\arrowtypec #3
\arrowtyped #4
\width #5
\height #6
}
\def\setpos(#1,#2){\xpos=#1 \ypos#2}

\def\bfig{\begin{picture}(\xext,\yext)(\xoff,\yoff)}
\def\efig{\end{picture}}

\def\putbox(#1,#2)#3{\put(#1,#2){\makebox(0,0){$#3$}}}

\def\settriparms[#1`#2`#3;#4]{\settripairparms[#1`#2`#3`1`1;#4]}%

\def\settripairparms[#1`#2`#3`#4`#5;#6]{%
\arrowtypea #1
\arrowtypeb #2
\arrowtypec #3
\arrowtyped #4
\arrowtypee #5
\width #6
\height #6
}

\def\resetparms{\settripairparms[1`1`1`1`1;500]\width 500}

\resetparms

\def\mvector(#1,#2)#3{
\put(0,0){\vector(#1,#2){#3}}%
\put(0,0){\vector(#1,#2){30}}%
}
\def\evector(#1,#2)#3{{
\arrowlength #3
\put(0,0){\vector(#1,#2){\arrowlength}}%
\advance \arrowlength by-30
\put(0,0){\vector(#1,#2){\arrowlength}}%
}}

\def\horsize#1#2{%
\settowidth{\tempdimen}{$#2$}%
#1=\tempdimen
\divide #1 by\unitlength
}

\def\vertsize#1#2{%
\settoheight{\tempdimen}{$#2$}%
#1=\tempdimen
\settodepth{\tempdimen}{$#2$}%
\advance #1 by\tempdimen
\divide #1 by\unitlength
}

\def\vertadjust[#1`#2`#3]{%
\vertsize{\tempcounta}{#1}%
\vertsize{\tempcountb}{#2}%
\ifnum \tempcounta<\tempcountb \tempcounta=\tempcountb \fi
\divide\tempcounta by2
\vertsize{\tempcountb}{#3}%
\ifnum \tempcountb>0 \advance \tempcountb by20 \fi
\ifnum \tempcounta<\tempcountb \tempcounta=\tempcountb \fi
}

\def\horadjust[#1`#2`#3]{%
\horsize{\tempcounta}{#1}%
\horsize{\tempcountb}{#2}%
\ifnum \tempcounta<\tempcountb \tempcounta=\tempcountb \fi
\divide\tempcounta by20
\horsize{\tempcountb}{#3}%
\ifnum \tempcountb>0 \advance \tempcountb by60 \fi
\ifnum \tempcounta<\tempcountb \tempcounta=\tempcountb \fi
}

\ig{ In this procedure, #1 is the paramater that sticks out all the way,
#2 sticks out the least and #3 is a label sticking out half way.  #4 is
the amount of the offset.}

\def\sladjust[#1`#2`#3]#4{%
\tempcountc=#4
\horsize{\tempcounta}{#1}%
\divide \tempcounta by2
\horsize{\tempcountb}{#2}%
\divide \tempcountb by2
\advance \tempcountb by-\tempcountc
\ifnum \tempcounta<\tempcountb \tempcounta=\tempcountb\fi
\divide \tempcountc by2
\horsize{\tempcountb}{#3}%
\advance \tempcountb by-\tempcountc
\ifnum \tempcountb>0 \advance \tempcountb by80\fi
\ifnum \tempcounta<\tempcountb \tempcounta=\tempcountb\fi
\advance\tempcounta by20
}

\def\putvector(#1,#2)(#3,#4)#5#6{{%
\xpos=#1
\ypos=#2
\run=#3
\rise=#4
\arrowlength=#5
\arrowtype=#6
\ifnum \arrowtype<0
    \ifnum \run=0
        \advance \ypos by-\arrowlength
    \else
        \tempcounta \arrowlength
        \multiply \tempcounta by\rise
        \divide \tempcounta by\run
        \ifnum\run>0
            \advance \xpos by\arrowlength
            \advance \ypos by\tempcounta
        \else
            \advance \xpos by-\arrowlength
            \advance \ypos by-\tempcounta
        \fi
    \fi
    \multiply \arrowtype by-1
    \multiply \rise by-1
    \multiply \run by-1
\fi
\ifnum \arrowtype=1
    \put(\xpos,\ypos){\vector(\run,\rise){\arrowlength}}%
\else\ifnum \arrowtype=2
    \put(\xpos,\ypos){\mvector(\run,\rise)\arrowlength}%
\else\ifnum\arrowtype=3
    \put(\xpos,\ypos){\evector(\run,\rise){\arrowlength}}%
\fi\fi\fi
}}

\def\putsplitvector(#1,#2)#3#4{
\xpos #1
\ypos #2
\arrowtype #4
\halflength #3
\arrowlength #3
\gap 140
\advance \halflength by-\gap
\divide \halflength by2
\ifnum \arrowtype=1
    \put(\xpos,\ypos){\line(0,-1){\halflength}}%
    \advance\ypos by-\halflength
    \advance\ypos by-\gap
    \put(\xpos,\ypos){\vector(0,-1){\halflength}}%
\else\ifnum \arrowtype=2
    \put(\xpos,\ypos){\line(0,-1)\halflength}%
    \put(\xpos,\ypos){\vector(0,-1)3}%
    \advance\ypos by-\halflength
    \advance\ypos by-\gap
    \put(\xpos,\ypos){\vector(0,-1){\halflength}}%
\else\ifnum\arrowtype=3
    \put(\xpos,\ypos){\line(0,-1)\halflength}%
    \advance\ypos by-\halflength
    \advance\ypos by-\gap
    \put(\xpos,\ypos){\evector(0,-1){\halflength}}%
\else\ifnum \arrowtype=-1
    \advance \ypos by-\arrowlength
    \put(\xpos,\ypos){\line(0,1){\halflength}}%
    \advance\ypos by\halflength
    \advance\ypos by\gap
    \put(\xpos,\ypos){\vector(0,1){\halflength}}%
\else\ifnum \arrowtype=-2
    \advance \ypos by-\arrowlength
    \put(\xpos,\ypos){\line(0,1)\halflength}%
    \put(\xpos,\ypos){\vector(0,1)3}%
    \advance\ypos by\halflength
    \advance\ypos by\gap
    \put(\xpos,\ypos){\vector(0,1){\halflength}}%
\else\ifnum\arrowtype=-3
    \advance \ypos by-\arrowlength
    \put(\xpos,\ypos){\line(0,1)\halflength}%
    \advance\ypos by\halflength
    \advance\ypos by\gap
    \put(\xpos,\ypos){\evector(0,1){\halflength}}%
\fi\fi\fi\fi\fi\fi
}

\def\putmorphism(#1)(#2,#3)[#4`#5`#6]#7#8#9{{%
\run #2
\rise #3
\ifnum\rise=0
  \puthmorphism(#1)[#4`#5`#6]{#7}{#8}{#9}%
\else\ifnum\run=0
  \putvmorphism(#1)[#4`#5`#6]{#7}{#8}{#9}%
\else
\setpos(#1)%
\arrowlength #7
\arrowtype #8
\ifnum\run=0
\else\ifnum\rise=0
\else
\ifnum\run>0
    \coefa=1
\else
   \coefa=-1
\fi
\ifnum\arrowtype>0
   \coefb=0
   \coefc=-1
\else
   \coefb=\coefa
   \coefc=1
   \arrowtype=-\arrowtype
\fi
\width=2
\multiply \width by\run
\divide \width by\rise
\ifnum \width<0  \width=-\width\fi
\advance\width by60
\if l#9 \width=-\width\fi
\putbox(\xpos,\ypos){#4}
{\multiply \coefa by\arrowlength
\advance\xpos by\coefa
\multiply \coefa by\rise
\divide \coefa by\run
\advance \ypos by\coefa
\putbox(\xpos,\ypos){#5} }%
{\multiply \coefa by\arrowlength
\divide \coefa by2
\advance \xpos by\coefa
\advance \xpos by\width
\multiply \coefa by\rise
\divide \coefa by\run
\advance \ypos by\coefa
\if l#9%
   \put(\xpos,\ypos){\makebox(0,0)[r]{$#6$}}%
\else\if r#9%
   \put(\xpos,\ypos){\makebox(0,0)[l]{$#6$}}%
\fi\fi }%
{\multiply \rise by-\coefc
\multiply \run by-\coefc
\multiply \coefb by\arrowlength
\advance \xpos by\coefb
\multiply \coefb by\rise
\divide \coefb by\run
\advance \ypos by\coefb
\multiply \coefc by70
\advance \ypos by\coefc
\multiply \coefc by\run
\divide \coefc by\rise
\advance \xpos by\coefc
\multiply \coefa by140
\multiply \coefa by\run
\divide \coefa by\rise
\advance \arrowlength by\coefa
\ifnum \arrowtype=1
   \put(\xpos,\ypos){\vector(\run,\rise){\arrowlength}}%
\else\ifnum\arrowtype=2
   \put(\xpos,\ypos){\mvector(\run,\rise){\arrowlength}}%
\else\ifnum\arrowtype=3
   \put(\xpos,\ypos){\evector(\run,\rise){\arrowlength}}%
\fi\fi\fi}\fi\fi\fi\fi}}

\def\puthmorphism(#1,#2)[#3`#4`#5]#6#7#8{{%
\xpos #1
\ypos #2
\width #6
\arrowlength #6
\putbox(\xpos,\ypos){#3\vphantom{#4}}%
{\advance \xpos by\arrowlength
\putbox(\xpos,\ypos){\vphantom{#3}#4}}%
\horsize{\tempcounta}{#3}%
\horsize{\tempcountb}{#4}%
\divide \tempcounta by2
\divide \tempcountb by2
\advance \tempcounta by30
\advance \tempcountb by30
\advance \xpos by\tempcounta
\advance \arrowlength by-\tempcounta
\advance \arrowlength by-\tempcountb
\putvector(\xpos,\ypos)(1,0){\arrowlength}{#7}%
\divide \arrowlength by2
\advance \xpos by\arrowlength
\vertsize{\tempcounta}{#5}%
\divide\tempcounta by2
\advance \tempcounta by20
\if a#8 %
   \advance \ypos by\tempcounta
   \putbox(\xpos,\ypos){#5}%
\else
   \advance \ypos by-\tempcounta
   \putbox(\xpos,\ypos){#5}%
\fi}}

\def\putvmorphism(#1,#2)[#3`#4`#5]#6#7#8{{%
\xpos #1
\ypos #2
\arrowlength #6
\arrowtype #7
\settowidth{\xlen}{$#5$}%
\putbox(\xpos,\ypos){#3}%
{\advance \ypos by-\arrowlength
\putbox(\xpos,\ypos){#4}}%
{\advance\arrowlength by-140
\advance \ypos by-70
\ifdim\xlen>0pt
   \if m#8%
      \putsplitvector(\xpos,\ypos){\arrowlength}{\arrowtype}%
   \else
      \putvector(\xpos,\ypos)(0,-1){\arrowlength}{\arrowtype}%
   \fi
\else
   \putvector(\xpos,\ypos)(0,-1){\arrowlength}{\arrowtype}%
\fi}%
\ifdim\xlen>0pt
   \divide \arrowlength by2
   \advance\ypos by-\arrowlength
   \if l#8%
      \advance \xpos by-40
      \put(\xpos,\ypos){\makebox(0,0)[r]{$#5$}}%
   \else\if r#8%
      \advance \xpos by40
      \put(\xpos,\ypos){\makebox(0,0)[l]{$#5$}}%
   \else
      \putbox(\xpos,\ypos){#5}%
   \fi\fi
\fi
}}

\def\topadjust[#1`#2`#3]{%
\yoff=10
\vertadjust[#1`#2`{#3}]%
\advance \yext by\tempcounta
\advance \yext by 10
}
\def\botadjust[#1`#2`#3]{%
\vertadjust[#1`#2`{#3}]%
\advance \yext by\tempcounta
\advance \yoff by-\tempcounta
}
\def\leftadjust[#1`#2`#3]{%
\xoff=0
\horadjust[#1`#2`{#3}]%
\advance \xext by\tempcounta
\advance \xoff by-\tempcounta
}
\def\rightadjust[#1`#2`#3]{%
\horadjust[#1`#2`{#3}]%
\advance \xext by\tempcounta
}
\def\rightsladjust[#1`#2`#3]{%
\sladjust[#1`#2`{#3}]{\width}%
\advance \xext by\tempcounta
}
\def\leftsladjust[#1`#2`#3]{%
\xoff=0
\sladjust[#1`#2`{#3}]{\width}%
\advance \xext by\tempcounta
\advance \xoff by-\tempcounta
}
\def\adjust[#1`#2;#3`#4;#5`#6;#7`#8]{%
\topadjust[#1``{#2}]
\leftadjust[#3``{#4}]
\rightadjust[#5``{#6}]
\botadjust[#7``{#8}]}

\def\putsquarep<#1>(#2)[#3;#4`#5`#6`#7]{{%
\setsqparms[#1]%
\setpos(#2)%
\settokens[#3]%
\puthmorphism(\xpos,\ypos)[\tokenc`\tokend`{#7}]{\width}{\arrowtyped}b%
\advance\ypos by \height
\puthmorphism(\xpos,\ypos)[\tokena`\tokenb`{#4}]{\width}{\arrowtypea}a%
\putvmorphism(\xpos,\ypos)[``{#5}]{\height}{\arrowtypeb}l%
\advance\xpos by \width
\putvmorphism(\xpos,\ypos)[``{#6}]{\height}{\arrowtypec}r%
}}

\def\putsquare{\@ifnextchar <{\putsquarep}{\putsquarep%
   <\arrowtypea`\arrowtypeb`\arrowtypec`\arrowtyped;\width`\height>}}
\def\square{\@ifnextchar< {\squarep}{\squarep
   <\arrowtypea`\arrowtypeb`\arrowtypec`\arrowtyped;\width`\height>}}
\def\squarep<#1>[#2`#3`#4`#5;#6`#7`#8`#9]{{
\setsqparms[#1]
\xext=\width                                          
\yext=\height                                         
\topadjust[#2`#3`{#6}]
\botadjust[#4`#5`{#9}]
\leftadjust[#2`#4`{#7}]
\rightadjust[#3`#5`{#8}]
\begin{picture}(\xext,\yext)(\xoff,\yoff)
\putsquarep<\arrowtypea`\arrowtypeb`\arrowtypec`\arrowtyped;\width`\height>%
(0,0)[#2`#3`#4`#5;#6`#7`#8`{#9}]%
\end{picture}%
}}

\def\putptrianglep<#1>(#2,#3)[#4`#5`#6;#7`#8`#9]{{%
\settriparms[#1]%
\xpos=#2 \ypos=#3
\advance\ypos by \height
\puthmorphism(\xpos,\ypos)[#4`#5`{#7}]{\height}{\arrowtypea}a%
\putvmorphism(\xpos,\ypos)[`#6`{#8}]{\height}{\arrowtypeb}l%
\advance\xpos by\height
\putmorphism(\xpos,\ypos)(-1,-1)[``{#9}]{\height}{\arrowtypec}r%
}}

\def\putptriangle{\@ifnextchar <{\putptrianglep}{\putptrianglep
   <\arrowtypea`\arrowtypeb`\arrowtypec;\height>}}
\def\ptriangle{\@ifnextchar <{\ptrianglep}{\ptrianglep
   <\arrowtypea`\arrowtypeb`\arrowtypec;\height>}}

\def\ptrianglep<#1>[#2`#3`#4;#5`#6`#7]{{
\settriparms[#1]%
\width=\height                         
\xext=\width                           
\yext=\width                           
\topadjust[#2`#3`{#5}]
\botadjust[#3``]
\leftadjust[#2`#4`{#6}]
\rightsladjust[#3`#4`{#7}]
\begin{picture}(\xext,\yext)(\xoff,\yoff)
\putptrianglep<\arrowtypea`\arrowtypeb`\arrowtypec;\height>%
(0,0)[#2`#3`#4;#5`#6`{#7}]%
\end{picture}%
}}

\def\putqtrianglep<#1>(#2,#3)[#4`#5`#6;#7`#8`#9]{{%
\settriparms[#1]%
\xpos=#2 \ypos=#3
\advance\ypos by\height
\puthmorphism(\xpos,\ypos)[#4`#5`{#7}]{\height}{\arrowtypea}a%
\putmorphism(\xpos,\ypos)(1,-1)[``{#8}]{\height}{\arrowtypeb}l%
\advance\xpos by\height
\putvmorphism(\xpos,\ypos)[`#6`{#9}]{\height}{\arrowtypec}r%
}}

\def\putqtriangle{\@ifnextchar <{\putqtrianglep}{\putqtrianglep
   <\arrowtypea`\arrowtypeb`\arrowtypec;\height>}}
\def\qtriangle{\@ifnextchar <{\qtrianglep}{\qtrianglep
   <\arrowtypea`\arrowtypeb`\arrowtypec;\height>}}

\def\qtrianglep<#1>[#2`#3`#4;#5`#6`#7]{{
\settriparms[#1]
\width=\height                         
\xext=\width                           
\yext=\height                          
\topadjust[#2`#3`{#5}]
\botadjust[#4``]
\leftsladjust[#2`#4`{#6}]
\rightadjust[#3`#4`{#7}]
\begin{picture}(\xext,\yext)(\xoff,\yoff)
\putqtrianglep<\arrowtypea`\arrowtypeb`\arrowtypec;\height>%
(0,0)[#2`#3`#4;#5`#6`{#7}]%
\end{picture}%
}}

\def\putdtrianglep<#1>(#2,#3)[#4`#5`#6;#7`#8`#9]{{%
\settriparms[#1]%
\xpos=#2 \ypos=#3
\puthmorphism(\xpos,\ypos)[#5`#6`{#9}]{\height}{\arrowtypec}b%
\advance\xpos by \height \advance\ypos by\height
\putmorphism(\xpos,\ypos)(-1,-1)[``{#7}]{\height}{\arrowtypea}l%
\putvmorphism(\xpos,\ypos)[#4``{#8}]{\height}{\arrowtypeb}r%
}}

\def\putdtriangle{\@ifnextchar <{\putdtrianglep}{\putdtrianglep
   <\arrowtypea`\arrowtypeb`\arrowtypec;\height>}}
\def\dtriangle{\@ifnextchar <{\dtrianglep}{\dtrianglep
   <\arrowtypea`\arrowtypeb`\arrowtypec;\height>}}

\def\dtrianglep<#1>[#2`#3`#4;#5`#6`#7]{{
\settriparms[#1]
\width=\height                         
\xext=\width                           
\yext=\height                          
\topadjust[#2``]
\botadjust[#3`#4`{#7}]
\leftsladjust[#3`#2`{#5}]
\rightadjust[#2`#4`{#6}]
\begin{picture}(\xext,\yext)(\xoff,\yoff)
\putdtrianglep<\arrowtypea`\arrowtypeb`\arrowtypec;\height>%
(0,0)[#2`#3`#4;#5`#6`{#7}]%
\end{picture}%
}}

\def\putbtrianglep<#1>(#2,#3)[#4`#5`#6;#7`#8`#9]{{%
\settriparms[#1]%
\xpos=#2 \ypos=#3
\puthmorphism(\xpos,\ypos)[#5`#6`{#9}]{\height}{\arrowtypec}b%
\advance\ypos by\height
\putmorphism(\xpos,\ypos)(1,-1)[``{#8}]{\height}{\arrowtypeb}r%
\putvmorphism(\xpos,\ypos)[#4``{#7}]{\height}{\arrowtypea}l%
}}

\def\putbtriangle{\@ifnextchar <{\putbtrianglep}{\putbtrianglep
   <\arrowtypea`\arrowtypeb`\arrowtypec;\height>}}
\def\btriangle{\@ifnextchar <{\btrianglep}{\btrianglep
   <\arrowtypea`\arrowtypeb`\arrowtypec;\height>}}

\def\btrianglep<#1>[#2`#3`#4;#5`#6`#7]{{
\settriparms[#1]
\width=\height                         
\xext=\width                           
\yext=\height                          
\topadjust[#2``]
\botadjust[#3`#4`{#7}]
\leftadjust[#2`#3`{#5}]
\rightsladjust[#4`#2`{#6}]
\begin{picture}(\xext,\yext)(\xoff,\yoff)
\putbtrianglep<\arrowtypea`\arrowtypeb`\arrowtypec;\height>%
(0,0)[#2`#3`#4;#5`#6`{#7}]%
\end{picture}%
}}

\def\putAtrianglep<#1>(#2,#3)[#4`#5`#6;#7`#8`#9]{{%
\settriparms[#1]%
\xpos=#2 \ypos=#3
{\multiply \height by2
\puthmorphism(\xpos,\ypos)[#5`#6`{#9}]{\height}{\arrowtypec}b}%
\advance\xpos by\height \advance\ypos by\height
\putmorphism(\xpos,\ypos)(-1,-1)[#4``{#7}]{\height}{\arrowtypea}l%
\putmorphism(\xpos,\ypos)(1,-1)[``{#8}]{\height}{\arrowtypeb}r%
}}

\def\putAtriangle{\@ifnextchar <{\putAtrianglep}{\putAtrianglep
   <\arrowtypea`\arrowtypeb`\arrowtypec;\height>}}
\def\Atriangle{\@ifnextchar <{\Atrianglep}{\Atrianglep
   <\arrowtypea`\arrowtypeb`\arrowtypec;\height>}}

\def\Atrianglep<#1>[#2`#3`#4;#5`#6`#7]{{
\settriparms[#1]
\width=\height                         
\xext=\width                           
\yext=\height                          
\topadjust[#2``]
\botadjust[#3`#4`{#7}]
\multiply \xext by2 
\leftsladjust[#3`#2`{#5}]
\rightsladjust[#4`#2`{#6}]
\begin{picture}(\xext,\yext)(\xoff,\yoff)%
\putAtrianglep<\arrowtypea`\arrowtypeb`\arrowtypec;\height>%
(0,0)[#2`#3`#4;#5`#6`{#7}]%
\end{picture}%
}}

\def\putAtrianglepairp<#1>(#2)[#3;#4`#5`#6`#7`#8]{{
\settripairparms[#1]%
\setpos(#2)%
\settokens[#3]%
\puthmorphism(\xpos,\ypos)[\tokenb`\tokenc`{#7}]{\height}{\arrowtyped}b%
\advance\xpos by\height
\advance\ypos by\height
\putmorphism(\xpos,\ypos)(-1,-1)[\tokena``{#4}]{\height}{\arrowtypea}l%
\putvmorphism(\xpos,\ypos)[``{#5}]{\height}{\arrowtypeb}m%
\putmorphism(\xpos,\ypos)(1,-1)[``{#6}]{\height}{\arrowtypec}r%
}}

\def\putAtrianglepair{\@ifnextchar <{\putAtrianglepairp}{\putAtrianglepairp%
   <\arrowtypea`\arrowtypeb`\arrowtypec`\arrowtyped`\arrowtypee;\height>}}
\def\Atrianglepair{\@ifnextchar <{\Atrianglepairp}{\Atrianglepairp%
   <\arrowtypea`\arrowtypeb`\arrowtypec`\arrowtyped`\arrowtypee;\height>}}

\def\Atrianglepairp<#1>[#2;#3`#4`#5`#6`#7]{{%
\settripairparms[#1]%
\settokens[#2]%
\width=\height
\xext=\width
\yext=\height
\topadjust[\tokena``]%
\vertadjust[\tokenb`\tokenc`{#6}]
\tempcountd=\tempcounta                       
\vertadjust[\tokenc`\tokend`{#7}]
\ifnum\tempcounta<\tempcountd                 
\tempcounta=\tempcountd\fi                    
\advance \yext by\tempcounta                  
\advance \yoff by-\tempcounta                 %
\multiply \xext by2 
\leftsladjust[\tokenb`\tokena`{#3}]
\rightsladjust[\tokend`\tokena`{#5}]%
\begin{picture}(\xext,\yext)(\xoff,\yoff)%
\putAtrianglepairp
<\arrowtypea`\arrowtypeb`\arrowtypec`\arrowtyped`\arrowtypee;\height>%
(0,0)[#2;#3`#4`#5`#6`{#7}]%
\end{picture}%
}}

\def\putVtrianglep<#1>(#2,#3)[#4`#5`#6;#7`#8`#9]{{%
\settriparms[#1]%
\xpos=#2 \ypos=#3
\advance\ypos by\height
{\multiply\height by2
\puthmorphism(\xpos,\ypos)[#4`#5`{#7}]{\height}{\arrowtypea}a}%
\putmorphism(\xpos,\ypos)(1,-1)[`#6`{#8}]{\height}{\arrowtypeb}l%
\advance\xpos by\height
\advance\xpos by\height
\putmorphism(\xpos,\ypos)(-1,-1)[``{#9}]{\height}{\arrowtypec}r%
}}

\def\putVtriangle{\@ifnextchar <{\putVtrianglep}{\putVtrianglep
   <\arrowtypea`\arrowtypeb`\arrowtypec;\height>}}
\def\Vtriangle{\@ifnextchar <{\Vtrianglep}{\Vtrianglep
   <\arrowtypea`\arrowtypeb`\arrowtypec;\height>}}

\def\Vtrianglep<#1>[#2`#3`#4;#5`#6`#7]{{
\settriparms[#1]
\width=\height                         
\xext=\width                           
\yext=\height                          
\topadjust[#2`#3`{#5}]
\botadjust[#4``]
\multiply \xext by2 
\leftsladjust[#2`#3`{#6}]
\rightsladjust[#3`#4`{#7}]
\begin{picture}(\xext,\yext)(\xoff,\yoff)%
\putVtrianglep<\arrowtypea`\arrowtypeb`\arrowtypec;\height>%
(0,0)[#2`#3`#4;#5`#6`{#7}]%
\end{picture}%
}}

\def\putVtrianglepairp<#1>(#2)[#3;#4`#5`#6`#7`#8]{{
\settripairparms[#1]%
\setpos(#2)%
\settokens[#3]%
\advance\ypos by\height
\putmorphism(\xpos,\ypos)(1,-1)[`\tokend`{#6}]{\height}{\arrowtypec}l%
\puthmorphism(\xpos,\ypos)[\tokena`\tokenb`{#4}]{\height}{\arrowtypea}a%
\advance\xpos by\height
\putvmorphism(\xpos,\ypos)[``{#7}]{\height}{\arrowtyped}m%
\advance\xpos by\height
\putmorphism(\xpos,\ypos)(-1,-1)[``{#8}]{\height}{\arrowtypee}r%
}}

\def\putVtrianglepair{\@ifnextchar <{\putVtrianglepairp}{\putVtrianglepairp%
    <\arrowtypea`\arrowtypeb`\arrowtypec`\arrowtyped`\arrowtypee;\height>}}
\def\Vtrianglepair{\@ifnextchar <{\Vtrianglepairp}{\Vtrianglepairp%
    <\arrowtypea`\arrowtypeb`\arrowtypec`\arrowtyped`\arrowtypee;\height>}}

\def\Vtrianglepairp<#1>[#2;#3`#4`#5`#6`#7]{{%
\settripairparms[#1]%
\settokens[#2]
\xext=\height                  
\width=\height                 
\yext=\height                  
\vertadjust[\tokena`\tokenb`{#4}]
\tempcountd=\tempcounta        
\vertadjust[\tokenb`\tokenc`{#5}]
\ifnum\tempcounta<\tempcountd%
\tempcounta=\tempcountd\fi
\advance \yext by\tempcounta
\botadjust[\tokend``]%
\multiply \xext by2
\leftsladjust[\tokena`\tokend`{#6}]%
\rightsladjust[\tokenc`\tokend`{#7}]%
\begin{picture}(\xext,\yext)(\xoff,\yoff)%
\putVtrianglepairp
<\arrowtypea`\arrowtypeb`\arrowtypec`\arrowtyped`\arrowtypee;\height>%
(0,0)[#2;#3`#4`#5`#6`{#7}]%
\end{picture}%
}}

\def\putCtrianglep<#1>(#2,#3)[#4`#5`#6;#7`#8`#9]{{%
\settriparms[#1]%
\xpos=#2 \ypos=#3
\advance\ypos by\height
\putmorphism(\xpos,\ypos)(1,-1)[``{#9}]{\height}{\arrowtypec}l%
\advance\xpos by\height
\advance\ypos by\height
\putmorphism(\xpos,\ypos)(-1,-1)[#4`#5`{#7}]{\height}{\arrowtypea}l%
{\multiply\height by 2
\putvmorphism(\xpos,\ypos)[`#6`{#8}]{\height}{\arrowtypeb}r}%
}}

\def\putCtriangle{\@ifnextchar <{\putCtrianglep}{\putCtrianglep
    <\arrowtypea`\arrowtypeb`\arrowtypec;\height>}}
\def\Ctriangle{\@ifnextchar <{\Ctrianglep}{\Ctrianglep
    <\arrowtypea`\arrowtypeb`\arrowtypec;\height>}}

\def\Ctrianglep<#1>[#2`#3`#4;#5`#6`#7]{{
\settriparms[#1]
\width=\height                          
\xext=\width                            
\yext=\height                           
\multiply \yext by2 
\topadjust[#2``]
\botadjust[#4``]
\sladjust[#3`#2`{#5}]{\width}
\tempcountd=\tempcounta                 
\sladjust[#3`#4`{#7}]{\width}
\ifnum \tempcounta<\tempcountd          
\tempcounta=\tempcountd\fi              
\advance \xext by\tempcounta            
\advance \xoff by-\tempcounta           %
\rightadjust[#2`#4`{#6}]
\begin{picture}(\xext,\yext)(\xoff,\yoff)%
\putCtrianglep<\arrowtypea`\arrowtypeb`\arrowtypec;\height>%
(0,0)[#2`#3`#4;#5`#6`{#7}]%
\end{picture}%
}}

\def\putDtrianglep<#1>(#2,#3)[#4`#5`#6;#7`#8`#9]{{%
\settriparms[#1]%
\xpos=#2 \ypos=#3
\advance\xpos by\height \advance\ypos by\height
\putmorphism(\xpos,\ypos)(-1,-1)[``{#9}]{\height}{\arrowtypec}r%
\advance\xpos by-\height \advance\ypos by\height
\putmorphism(\xpos,\ypos)(1,-1)[`#5`{#8}]{\height}{\arrowtypeb}r%
{\multiply\height by 2
\putvmorphism(\xpos,\ypos)[#4`#6`{#7}]{\height}{\arrowtypea}l}%
}}

\def\putDtriangle{\@ifnextchar <{\putDtrianglep}{\putDtrianglep
    <\arrowtypea`\arrowtypeb`\arrowtypec;\height>}}
\def\Dtriangle{\@ifnextchar <{\Dtrianglep}{\Dtrianglep
   <\arrowtypea`\arrowtypeb`\arrowtypec;\height>}}

\def\Dtrianglep<#1>[#2`#3`#4;#5`#6`#7]{{
\settriparms[#1]
\width=\height                         
\xext=\height                          
\yext=\height                          
\multiply \yext by2 
\topadjust[#2``]
\botadjust[#4``]
\leftadjust[#2`#4`{#5}]
\sladjust[#3`#2`{#5}]{\height}
\tempcountd=\tempcountd                
\sladjust[#3`#4`{#7}]{\height}
\ifnum \tempcounta<\tempcountd         
\tempcounta=\tempcountd\fi             
\advance \xext by\tempcounta           %
\begin{picture}(\xext,\yext)(\xoff,\yoff)
\putDtrianglep<\arrowtypea`\arrowtypeb`\arrowtypec;\height>%
(0,0)[#2`#3`#4;#5`#6`{#7}]%
\end{picture}%
}}

\def\setrecparms[#1`#2]{\width=#1 \height=#2}%
%

\def\recursep<#1`#2>[#3;#4`#5`#6`#7`#8]{{%
\width=#1 \height=#2
\settokens[#3]
\settowidth{\tempdimen}{$\tokena$}
\ifdim\tempdimen=0pt
  \savebox{\tempboxa}{\hbox{$\tokenb$}}%
  \savebox{\tempboxb}{\hbox{$\tokend$}}%
  \savebox{\tempboxc}{\hbox{$#6$}}%
\else
  \savebox{\tempboxa}{\hbox{$\hbox{$\tokena$}\times\hbox{$\tokenb$}$}}%
  \savebox{\tempboxb}{\hbox{$\hbox{$\tokena$}\times\hbox{$\tokend$}$}}%
  \savebox{\tempboxc}{\hbox{$\hbox{$\tokena$}\times\hbox{$#6$}$}}%
\fi
\ypos=\height
\divide\ypos by 2
\xpos=\ypos
\advance\xpos by \width
\xext=\xpos \yext=\height
\topadjust[#3`\usebox{\tempboxa}`{#4}]%
\botadjust[#5`\usebox{\tempboxb}`{#8}]%
\sladjust[\tokenc`\tokenb`{#5}]{\ypos}%
\tempcountd=\tempcounta
\sladjust[\tokenc`\tokend`{#5}]{\ypos}%
\ifnum \tempcounta<\tempcountd
\tempcounta=\tempcountd\fi
\advance \xext by\tempcounta
\advance \xoff by-\tempcounta
\rightadjust[\usebox{\tempboxa}`\usebox{\tempboxb}`\usebox{\tempboxc}]%
\bfig
\putCtrianglep<-1`1`1;\ypos>(0,0)[`\tokenc`;#5`#6`{#7}]%
\puthmorphism(\ypos,0)[\tokend`\usebox{\tempboxb}`{#8}]{\width}{-1}b%
\puthmorphism(\ypos,\height)[\tokenb`\usebox{\tempboxa}`{#4}]{\width}{-1}a%
\advance\ypos by \width
\putvmorphism(\ypos,\height)[``\usebox{\tempboxc}]{\height}1r%
\efig
}}

\def\recurse{\@ifnextchar <{\recursep}{\recursep<\width`\height>}}

\def\puttwohmorphisms(#1,#2)[#3`#4;#5`#6]#7#8#9{{%
%
\puthmorphism(#1,#2)[#3`#4`]{#7}0a
\ypos=#2
\advance\ypos by 20
\puthmorphism(#1,\ypos)[\phantom{#3}`\phantom{#4}`#5]{#7}{#8}a
\advance\ypos by -40
\puthmorphism(#1,\ypos)[\phantom{#3}`\phantom{#4}`#6]{#7}{#9}b
}}

\def\puttwovmorphisms(#1,#2)[#3`#4;#5`#6]#7#8#9{{%
%
%
%
\putvmorphism(#1,#2)[#3`#4`]{#7}0a
\xpos=#1
\advance\xpos by -20
\putvmorphism(\xpos,#2)[\phantom{#3}`\phantom{#4}`#5]{#7}{#8}l
\advance\xpos by 40
\putvmorphism(\xpos,#2)[\phantom{#3}`\phantom{#4}`#6]{#7}{#9}r
}}

\def\puthcoequalizer(#1)[#2`#3`#4;#5`#6`#7]#8#9{{%
%
\setpos(#1)%
\puttwohmorphisms(\xpos,\ypos)[#2`#3;#5`#6]{#8}11%
\advance\xpos by #8
\puthmorphism(\xpos,\ypos)[\phantom{#3}`#4`#7]{#8}1{#9}
}}

\def\putvcoequalizer(#1)[#2`#3`#4;#5`#6`#7]#8#9{{%
%
%
%
%
\setpos(#1)%
\puttwovmorphisms(\xpos,\ypos)[#2`#3;#5`#6]{#8}11%
\advance\ypos by -#8
\putvmorphism(\xpos,\ypos)[\phantom{#3}`#4`#7]{#8}1{#9}
}}

\def\putthreehmorphisms(#1)[#2`#3;#4`#5`#6]#7(#8)#9{{%
\setpos(#1) \settypes(#8)
\if a#9 %
     \vertsize{\tempcounta}{#5}%
     \vertsize{\tempcountb}{#6}%
     \ifnum \tempcounta<\tempcountb \tempcounta=\tempcountb \fi
\else
     \vertsize{\tempcounta}{#4}%
     \vertsize{\tempcountb}{#5}%
     \ifnum \tempcounta<\tempcountb \tempcounta=\tempcountb \fi
\fi
\advance \tempcounta by 60
\puthmorphism(\xpos,\ypos)[#2`#3`#5]{#7}{\arrowtypeb}{#9}
\advance\ypos by \tempcounta
\puthmorphism(\xpos,\ypos)[\phantom{#2}`\phantom{#3}`#4]{#7}{\arrowtypea}{#9}
\advance\ypos by -\tempcounta \advance\ypos by -\tempcounta
\puthmorphism(\xpos,\ypos)[\phantom{#2}`\phantom{#3}`#6]{#7}{\arrowtypec}{#9}
}}

\def\putarc(#1,#2)[#3`#4`#5]#6#7#8{{%
\xpos #1
\ypos #2
\width #6
\arrowlength #6
\putbox(\xpos,\ypos){#3\vphantom{#4}}%
{\advance \xpos by\arrowlength
\putbox(\xpos,\ypos){\vphantom{#3}#4}}%
\horsize{\tempcounta}{#3}%
\horsize{\tempcountb}{#4}%
\divide \tempcounta by2
\divide \tempcountb by2
\advance \tempcounta by30
\advance \tempcountb by30
\advance \xpos by\tempcounta
\advance \arrowlength by-\tempcounta
\advance \arrowlength by-\tempcountb
\halflength=\arrowlength \divide\halflength by 2
\divide\arrowlength by 5
\put(\xpos,\ypos){\bezier{\arrowlength}(0,0)(50,50)(\halflength,50)}
\ifnum #7=-1 \put(\xpos,\ypos){\vector(-3,-2)0} \fi
\advance\xpos by \halflength
\put(\xpos,\ypos){\xpos=\halflength \advance\xpos by -50
   \bezier{\arrowlength}(0,50)(\xpos,50)(\halflength,0)}
\ifnum #7=1 {\advance \xpos by
   \halflength \put(\xpos,\ypos){\vector(3,-2)0}} \fi
\advance\ypos by 50
\vertsize{\tempcounta}{#5}%
\divide\tempcounta by2
\advance \tempcounta by20
\if a#8 %
   \advance \ypos by\tempcounta
   \putbox(\xpos,\ypos){#5}%
\else
   \advance \ypos by-\tempcounta
   \putbox(\xpos,\ypos){#5}%
\fi
}}

\makeatother

\sloppy

\newcommand{\nl}{\hspace{2cm}\\ }

\def\nec{\Box}
\def\pos{\Diamond}
\def\diam{{\tiny\Diamond}}

\def\lc{\lceil}
\def\rc{\rceil}
\def\lf{\lfloor}
\def\rf{\rfloor}
\def\lk{\langle}
\def\rk{\rangle}
\def\blk{\dot{\langle\!\!\langle}}
\def\brk{\dot{\rangle\!\!\rangle}}

\def\lse{[\!|}
\def\rse{|\!]}
\def\le{(\!|}
\def\re{|\!)}

\def\homo{{\approx\!\! >}}
\def\inn{\in\!\!\!\!\ra}
\def\dsum{\stackrel{\cdot}{\sqcup}}
\def\dsum{\sqcup\!\!\!\!\cdot\;}

\def\tl{\triangleleft}
\def\tr{\triangleright}

\def\lhb{\lhd \hspace {-1mm}\bullet}

\newcommand{\pa}{\parallel}
\newcommand{\lra}{\longrightarrow}
\newcommand{\hra}{\hookrightarrow}
\newcommand{\hla}{\hookleftarrow}
\newcommand{\ra}{\rightarrow}
\newcommand{\la}{\leftarrow}
\newcommand{\lla}{\longleftarrow}
\newcommand{\da}{\downarrow}
\newcommand{\ua}{\uparrow}
\newcommand{\dA}{\downarrow\!\!\!^\bullet}
\newcommand{\uA}{\uparrow\!\!\!_\bullet}
\newcommand{\Da}{\Downarrow}
\newcommand{\DA}{\Downarrow\!\!\!^\bullet}
\newcommand{\UA}{\Uparrow\!\!\!_\bullet}
\newcommand{\Ua}{\Uparrow}
\newcommand{\Lra}{\Longrightarrow}
\newcommand{\Ra}{\Rightarrow}
\newcommand{\Lla}{\Longleftarrow}
\newcommand{\La}{\Leftarrow}
\newcommand{\nperp}{\perp\!\!\!\!\!\setminus\;\;}
\newcommand{\pq}{\preceq}

\newcommand{\lms}{\longmapsto}
\newcommand{\ms}{\mapsto}
\newcommand{\subseteqnot}{\subseteq\hskip-4 mm_\not\hskip3 mm}

\newcommand{\bth}{\begin{theorem}}
\newcommand{\eth}{\end{theorem}}

\def\phi{\varphi}
\def\ve{\varepsilon}
\def\o{{\omega}}

\def\bA{{\bf A}}
\def\bM{{\bf M}}
\def\bN{{\bf N}}
\def\bC{{\bf C}}
\def\bI{{\bf I}}
\def\bL{{\bf L}}
\def\bT{{\bf T}}
\def\bS{{\bf S}}
\def\bD{{\bf D}}
\def\bB{{\bf B}}
\def\bW{{\bf W}}
\def\bP{{\bf P}}
\def\bX{{\bf X}}
\def\bY{{\bf Y}}
\def\ba{{\bf a}}
\def\bb{{\bf b}}
\def\bc{{\bf c}}
\def\bd{{\bf d}}
\def\bh{{\bf h}}
\def\bi{{\bf i}}
\def\bj{{\bf j}}
\def\bk{{\bf k}}
\def\bm{{\bf m}}
\def\bn{{\bf n}}
\def\bp{{\bf p}}
\def\bq{{\bf q}}
\def\be{{\bf e}}
\def\br{{\bf r}}
\def\bi{{\bf i}}
\def\bs{{\bf s}}
\def\bt{{\bf t}}
\def\jeden{{\bf 1}}
\def\dwa{{\bf 2}}
\def\trzy{{\bf 3}}

\def\cA{{\cal A}}
\def\cB{{\cal B}}
\def\cC{{\cal C}}
\def\cD{{\cal D}}
\def\cE{{\cal E}}
\def\cF{{\cal F}}
\def\cG{{\cal G}}
\def\cI{{\cal I}}
\def\cJ{{\cal J}}
\def\cK{{\cal K}}
\def\cL{{\cal L}}
\def\cN{{\cal N}}
\def\cM{{\cal M}}
\def\cO{{\cal O}}
\def\cP{{\cal P}}
\def\cQ{{\cal Q}}
\def\cR{{\cal R}}
\def\cS{{\cal S}}
\def\cT{{\cal T}}
\def\cU{{\cal U}}
\def\cV{{\cal V}}
\def\cW{{\cal W}}
\def\cX{{\cal X}}
\def\cY{{\cal Y}}

\def\bCat{{{\bf Cat}}}

\pagenumbering{arabic} \setcounter{page}{1}

\title{Lax Monoidal Fibrations}
\author{Marek Zawadowski
\\
Instytut Matematyki, Uniwersytet Warszawski\\
ul. S.Banacha 2, 00-913 Warszawa, Poland\\
zawado@mimuw.edu.pl\\}
\date{December 22, 2009
}
\maketitle
\begin{flushright}
\em Dedicated to Mihaly Makkai\\ on the occasion of his
70th birthday.
\end{flushright}

\begin{abstract}  We introduce the notion of a lax monoidal fibration
and we show how it can be conveniently used to deal with various
algebraic structures that play an important role in some
definitions, cf. \cite{BaezDolan}, \cite{HMP}, \cite{SZ},
\cite{S} of the opetopic sets. We present the 'standard' such
structures, the exponential fibrations of basic
fibrations and three areas of applications. The first area is related
to the $T$-categories of A. Burroni.  The monoids in the Burroni
lax monoidal fibrations form the fibration of $T$-categories.
The construction of the relative Burroni fibrations and
free $T$-categories in this context, allow us to extend the
definition of the set of opetopes given in \cite{Le} to the
category of opetopic sets (internally to any Grothendieck topos,
if needed). We also show that the fibration of (1-level)
multicategories considered in \cite{HMP} is equivalent to the
fibration of (finitary, cartesian) polynomial monads. This
equivalence is induced by the equivalence of lax monoidal
fibrations of amalgamated signatures, polynomial diagrams, and
polynomial (finitary, endo) functors. Finally, we develop a
similar theory for symmetric signatures, analytic diagrams (a
notion introduced here), and (finitary, multivariable) analytic
(endo)functors, cf. \cite{J2}. Among other things we show that the
fibrations of symmetric multicategories is equivalent to the
fibration of analytic monads. We also give a characterization
(Corollary \ref{analytic monads}) of such a fibration of analytic
monads. An object of this fibration is a weakly cartesian monad on a
slice of $Set$ whose functor part is a finitary functor weakly
preserving wide pullbacks. A morphism of this fibration is  a weakly cartesian
morphism of monads whose functor part is a pullback functor.
\end{abstract}

MS Classification 18D10, 18D30, 18D50, 18C15   (AMS 2010).

\newpage
\tableofcontents
\newpage
\section{Introduction}
The notion of a lax monoidal fibration studied here,
is designed to help understand connections between
various definitions of opetopic sets. More specifically, various algebraic/categorical mechanisms
used to define them. The primary goal was to understand the relation between the notion of
opetopic set of Baez-Dolan, cf. \cite{BaezDolan}  and multitopic set of Hermida-Makkai-Power, cf. \cite{HMP}.
However, it turned out that many other notions connected with the development of higher categories can be
successfully organized into a lax monoidal fibration. This is not a mere encoding for its sake but
in the context of lax monoidal fibration many notions can be conveniently compared, characterized, and developed
beyond what was previously known. This paper provides many examples of such applications of lax monoidal fibrations
but the comparison of the mentioned definitions of opetopic sets,
as well as a yet another definition of opetopic sets(!),
will be presented in the forthcoming paper \cite{SZ}.

The lax monoidal fibrations provide a convenient tool to deal with many-level structures, like categories that have
objects and morphisms, multicategories\footnote{We follow mostly the terminology from \cite{Le}, in particular  for
us (various) multicategories are the same things as (various) colored set operads.} of various kinds (that have objects-types
and multiarrows-function symbols) or $T$-categories of Burroni, cf. \cite{B}. These are examples of structures that have
just two levels but by building a tower of fibrations, see Subsections \ref{tower-opetopic-sets} and \ref{tower-n-cat},
or iterating a construction inside a single fibration, cf. \cite{SZ}, \cite{S}, we can deal with many-level structures like
opetopic sets, polygraphs, $n$-categories, $\o$-categories and others. To define monoids of interest in this setting
we are not doing it in 'one big step' but we divide it into three smaller steps. First we define a fibration, then we define
the monoidal structure in this fibration and finally we define a fibration of monoids over the same base as the fibration
we started with. In that way if we want, as we do in Section \ref{symm-sig},  to compare multicategories with non-standard amalgamation
with symmetric multicategories we can compare the fibrations of amalgamated and symmetric signatures first, then compare
the tensors (there is more than one possibility) and finally we get a comparison of suitable multicategories.

A lax monoidal fibration is a fibration $p:\cE\ra \bB$ equipped with two functors $\otimes :\cE\times_\bB\cE\ra \cE$ and $I:\bB\ra \cE$
commuting over the base but not required to be morphisms of fibrations\footnote{This means that we do not require that $\otimes$ or $I$ send
prone (formerly cartesian) morphisms to prone morphisms.}. We call such morphism lax morphism of fibrations (or fibred morphisms),
as the fact that a morphism of fibrations commute over the base already forces some lax preservation of prone morphisms.
There are also coherence morphisms $\alpha$, $\lambda$, $\rho$
satisfying the usual conditions but they are not required to be isomorphisms, as in many examples they are not.
The direction of these morphism are so chosen to cover all our examples. The fibres of such fibration are monoidal categories
and reindexing functors are monoidal functors. It is in fact often the case, that the fibres are strong monoidal categories but
reindexing functors are almost never strong even in the lax monoidal fibration whose monoids are small categories!
This makes the whole context unavoidably lax.
The morphism of lax monoidal fibrations are lax morphisms of fibrations that are monoidal in the only reasonable sense.
The 2-cells are also defined in the only reasonable way. Then in the analogy with the non-fibred situation, cf. \cite{BaezDolan},
a lax monoidal fibration may act on arbitrary fibrations. So we have a 2-category of actions of lax monoidal fibrations, as well.
It is quite surprising how many things can be explained in terms of actions and their exponential adjoints. This will be carefully explained in
Section \ref{sec-exp-fib}  and used many times in the following sections.

The paper is organized as follows. In sections \ref{sec-lmf} and \ref{sec-actions of lmf},
we introduce the main notions of the paper of a lax monoidal fibration,
an action of a lax monoidal fibrations and $2$-categories of these structures. The examples presented there are very basic.

In Section \ref{sec-exp-fib}, we discuss the lax monoidal fibrations $\cE^\cE\ra \bB$ that
arises as exponential fibrations of bifibrations.
It turns out that the exponentiation in $Cat_{/\bB}$, the slice of $Cat$ over the base, is much
more interesting than the exponentiation in $Fib(\bB)$,
the category of fibrations over $\bB$. Among such fibrations there are even more special ones,
the exponential fibrations of the basic fibrations $cod:\bB^\ra\lra \bB$.
If $\bB$ has pullbacks such a fibration, denoted $Exp(\bB)\ra \bB$, always exists and if a lax
monoidal fibration $\cE\ra \bB$ acts on the basic fibrations $cod:\bB^\ra\lra \bB$ we have
a representation morphism of lax monoidal fibrations
  \begin{center} \xext=700 \yext=400
\begin{picture}(\xext,\yext)(\xoff,\yoff)
\settriparms[1`1`1;350]
  \putVtriangle(0,0)[\cE`Exp(\bB)`\bB;``]
 \end{picture}
\end{center}
that compares an arbitrary fibration with a standard one.

In Section \ref{burroni}, we show how one can split the definition of a $T$-category of
Burroni, cf. \cite{B} p. 225-227, into three parts.
The fibration of $T$-graphs, denoted $p_T:Gph(T)\ra \cC$, the monoidal part and finally
the monoids in such a lax monoidal fibrations. We call such fibrations Burroni
fibrations to honor A.Burroni who was the first to consider them, cf. \cite{B} p. 262.
The monoids in such a fibrations are exactly the $T$-categories
of Burroni. Since it is not necessary to have a cartesian
monad\footnote{The only requirement is that the category have pullbacks.}
$T$ to build such a fibration we can recover that way all the $T$-categories
that were considered in \cite{B}. The fibres of such a fibration $p_T$
are not necessarily strong monoidal unless $T$ is cartesian. However,
as we already mentioned, the reindexing functors are almost never strong monoidal
functors. The Burroni fibrations always acts on basic fibrations and hence they
have representation morphisms of lax monoidal fibration into the standard ones
   \begin{center} \xext=700 \yext=400
\begin{picture}(\xext,\yext)(\xoff,\yoff)
\settriparms[1`1`1;350]
  \putVtriangle(0,0)[Gph(T)`Exp(\cC)`\cC;rep_T`p_T`]
 \end{picture}
\end{center}
If $T$ is cartesian then this morphism is a morphism of bifibrations,
Proposition \ref{T-cart-repT}. The construction of $T$-categories can
be made relative with respect to a fibration if the monad $T$
is already fibred. Moreover, in this relative context the construction
due to M. Kelly, cf. \cite{Ke}
p.69, see also \cite{BJT}, together with the characterization of
T. Leinster, cf. \cite{Le} p. 334, gives a characterization of those
fibred cartesian monads for which the free relative $T$-categories
exists.  This allows us to extend the definition of the set of opetopes
given by T. Leinster, cf. \cite{Le} p. p.179, to the definition of
the whole category of opetopic sets, and this category can be build
internally in any Grothendieck topos not only in $Set$. We simply
iterate $\o$ times the construction of relative $T$-graph fibration
starting from the identity monad.

In Section \ref{am-sig}, we show that two seemingly different
languages used to define opetopes and opetopic sets, cf.
\cite{HMP} and \cite{KJBM}, are in fact equivalent. We show
that the lax monoidal fibrations of amalgamated signatures
$p_a:Sig_a\ra Set$ and of polynomial diagrams $p_{pd}:\cP oly \cD
iag\ra Set$ are equivalent. The difference is rather in style that
can be easily explained in the context of lax monoidal fibrations.
The amalgamated signatures are 'more concrete' and come naturally
equipped with an action on the basic fibration $cod:Set^\ra\lra Set$
whereas polynomial diagrams come equipped with a representation into
the exponential fibration $Exp(Set)\ra Set$. This representation
is the exponential adjoint of the action
and its essential image is the lax monoidal fibration of (finitary)
polynomial (endo)functors $p_{poly}:\cP oly\ra Set$. As these three
fibrations are equivalent as lax monoidal fibrations we obtain in
particular that the fibration of (1-level) multicategories with
non-standard amalgamations is equivalent to the fibration of
polynomial monads (i.e. cartesian monads on slices of $Set$ whose
functor parts are  polynomial functors) and as morphisms cartesian
morphism of monads whose functor parts are pullback functors
(counted as morphism in the opposite direction), see Corollary
\ref{equiv-mona-polymonads}. It is possible to give a natural
definition of opetopic sets in this context, see \cite{SZ},
\cite{S}. We end this section by showing how to deal with the so
called single tensor and 2-level objects, the original setting for
the definition in \cite{HMP}.

In Section \ref{symm-sig}, we develop a parallel theory to the one from
the previous section but this time we start with the lax monoidal fibration
of symmetric signatures $p_s:Sig_s\ra Set$ instead of amalgamated signatures,
whose monoids form exactly the fibration symmetric multicategories.
This fibration is also naturally equipped with an action on the basic fibration
and taking an adjoint we get again representation morphism
   \begin{center} \xext=700 \yext=400
\begin{picture}(\xext,\yext)(\xoff,\yoff)
\settriparms[1`1`1;350]
  \putVtriangle(0,0)[Sig_s`Exp(Set)`Set;rep_s`p_s`]
 \end{picture}
\end{center}
As in the previous case, this morphism is faithful and full on isomorphisms.
Its essential image, denoted by $p_{an}:\cA n\ra Set$, is the lax monoidal fibration
of multivariable analytic (endo)functors, cf. \cite{J2}, and analytic natural transformations between them.
As a consequence, the fibration of symmetric multicategories in $Set$ is equivalent to
the fibration of analytic monads. We also provide an abstract characterization of the fibration
of analytic functors extending the one from \cite{J2}.
We show, see Theorem \ref{analytic char}, that this fibration of (multivariable) analytic (endo)functors, which is
a lax monoidal subfibration of the exponential fibration $Exp(Set)\ra Set$,
consists of finitary functors on slices of $Set$ that preserve weakly wide pullbacks and has as morphisms weakly cartesian
natural transformations. The proof of this characterization is based on ideas from \cite{J2} and \cite{AV}.
As a consequence, we obtain Corollary \ref{analytic monads} saying that the notions of a symmetric multicategory and
of an analytic monad are equivalent. So analytic functors is yet another tool that could be used to
define the category of opetopic sets. In Subsection \ref{an-vs-poly} we introduce an intermediate notion of an analytic
diagrams that is related to symmetric signatures and analytic functors as polynomial diagrams are related to
amalgamated signatures and polynomial functors. These diagrams are polynomial diagrams of a special kind
in the category of symmetric sets $\sigma Set$, i.e. the category of presheaves on the coproduct (in $Cat$)
of finite symmetric groups.  However the representation is given via a composition of five functors not three
as in the case of usual polynomial diagrams.
In the last Subsection of the paper  we  compare the notions studied in  Sections \ref{am-sig} and \ref{symm-sig}.
We describe the following diagram of lax monoidal fibrations and their morphisms
 \begin{center} \xext=800 \yext=1680
\begin{picture}(\xext,\yext)(\xoff,\yoff)
\setsqparms[1`1`1`0;800`400]
\putsquare(0,1200)[ Sig_a`Sig_s`\phantom{\cP oly\cD iag}`\phantom{\cA n\cD iag};K_{sig}`\iota_{a}`\iota_{s}`]
\setsqparms[1`1`1`0;800`400]
\putsquare(0,800)[\cP oly\cD iag`\cA n\cD iag`\phantom{\cP oly}`\phantom{\cA n};K_{diag}`rep_{pd}`rep_{and}`]
\setsqparms[1`1`1`0;800`400]
\putsquare(0,400)[\cP oly`\cA n`\phantom{Cart(Set)}`\phantom{wCart(Set)};K_{fu}```]
\settriparms[1`1`1;400]
  \putVtriangle(0,0)[Cart(Set)`wCart(Set)`Exp(Set);``]
\put(540,930){\makebox(200,100){$\Phi$}}
\put(560,850){\makebox(200,100){$\Ra$}}
\put(540,1330){\makebox(200,100){$\Psi$}}
\put(560,1250){\makebox(200,100){$\Ra$}}
 \end{picture}
\end{center}
All the arrows are strong morphisms of lax monoidal fibrations.
We show among other things that the horizontal morphisms are full and faithful.
The three named horizontal arrows are morphisms comparing signatures, diagrams,
and functors, respectively. The four named vertical arrows are equivalences of
lax monoidal fibrations. The five unnamed arrows are inclusions.
Many more interesting connections between these lax monoidal fibrations will be explained in \cite{SZ}.
We finish with an observation that a weakly cartesian natural transformations between
polynomial functors are cartesian.

There are two possible notions of an analytic functor on $Set_{/O}$.
The species  and the analytic functors of one  variable $Set\ra Set$
and of many variables $Set_{/O}\ra Set$, for a finite set $O$, were introduced
by A. Joyal in \cite{J1} and \cite{J2} to study enumerative combinatorics.
Clearly, an $O$-tuple of multivariable analytic functors
taken together form an endofunctor $Set_{/O}\ra Set_{/O}$,
that should be considered as analytic, as well.
The concept of analytic functor was studied in the category $Set$,
in category of vector spaces $Vect$, cf. \cite{J2} but
also in an arbitrary monoidal category, cf. \cite{AV}. In that way,
we have two kinds of analytic functors on slices of $Set$
(and powers of other monoidal categories). An analytic functor from $Set_{/O}$
to $Set_{/O}$ can be defined as the left Kan extension of a functor
$f:\bB \ra  Set_{/O}$, where $\bB$ is the category of finite sets and
bijections, cf. \cite{AV}, or as an $O$-tuple of multivariable analytic
functors $Set_{/O}\ra Set$, cf. \cite{J2}. The first notion does not
allow functors that are not coproducts of functors between fibres.
In this paper we consider only the second notion.

The idea of equipping fibrations with some kind of a monoidal
structure goes back to N.S. Rivano \cite{Sa} and M.F. Gouzou-R. Grunig, cf. \cite{GG}.
It was taken up later by M. Shulman  in \cite{Sh}, p. 698. These
notions `$\bB-\otimes-$cat\'egories fibre\'e' in \cite{Sa},
`cat\'egorie fibr\'ee sur $\bB$ multiplicative' in \cite{GG},
and `monoidal fibration' in \cite{Sh} are different
than the notion of a lax monoidal fibrations presented here. Also
the motivations in each case are different than ours. The total category of
a lax monoidal fibration is not a monoidal category and in this
sense the notion is closer to the notions considered in \cite{Sa} and \cite{GG}.
On the other hand, we do not require our tensor or
unit to be morphisms of fibrations (i.e. preserve the prone
morphisms) as it would eliminate most of our examples. This causes
that our reindexing morphisms are not necessarily strong monoidal
functors. In our applications the actions of lax monoidal
fibrations play an important role. This does not have an analog in
the other approaches.

I would like to thank George Janelidze and Thomas Streicher for
the conversations related to the matters contained in
this paper, Andre Joyal for explaining to me some aspects of his
theory of analytic functors. Special thanks are due to Krzysztof Kapulkin,
Magdalena K\c{e}dziorek, Karol Szumi{\l}o and Stanis{\l}aw Szawiel,
the members of an informal Category Seminar held in Spring 2008 at
Warsaw University, for giving me an opportunity to present the essential
notions introduced and studied in this paper. I would like also to thank the anonymous
referee for the very thorough report that helped considerably to improve
the presentation of the paper. Last but not least I would
like to thank, our jubilee, Mihaly Makkai for introducing me to the subject
of higher-dimensional categories many years ago and to him and Victor Harnik
for countless discussions of the related matters.

The diagrams for this paper were prepared with the help of  {\em
catmac1} of Michael Barr.

\section{Lax monoidal fibrations}\label{sec-lmf}

\subsection{Preliminaries, fibrations}

Our standard reference no fibrations (opfibrations and
bifibrations) is \cite{St}. However the terminology used here
follows more the one used by P. Tayor and P.T. Johnston. We call
{\em prone} and {\em supine} morphisms what \cite{St} would call
{\em cartesian} and {\em cocartesian}.  The fibre of a
(bi)fibration $p:\cE\ra \bB$ over $B\in\bB$ will be denoted
 $\cE_B$. If $p$ is a fibration,  $u:B\ra B'$ is a morphism in
$\bB$ then we have (using axiom of choice for classes) a {\em
reindexing functor} $u^*: \cE_{B'}\ra \cE_B$ defined with the help of
prone morphisms; if $p$ is an opfibration, we have a {\em
coreindexing functor} $u_!: \cE_{B}\ra \cE_{B'}$ defined with the
help of supine morphisms. In a bifibration both functors exist
and are adjoint $u_!\dashv u^*$. The unit and counit of this
adjunction will be denoted by $\eta^u$ and $\varepsilon^u$,
respectively. We call a bifibration $P:\cE\ra \bB$ {\em cartesian}
if the fibres of $p$ have pullbacks, $u_!$ preserves them and both
$\eta^u$ and $\varepsilon^u$ are cartesian natural
transformations, for all morphisms $u$ in $\bB$. Note that this
notation suppresses the fact that these functors and natural
transformations are related to a specific (bi,op)fibration.
The fibration we have in mind should be always read from the
context. A fibration has fibred (co)limits of type $K$ if and only if
each fibre has (co)limits of type $K$ and reindexing functors
preserve them.
In the paper, we consider (bi)fibrations that are
equipped additionally with a monoidal structure. It is not
always the case that the morphisms we want to consider between (bi)fibrations
preserves all the structure involved (prone morphisms,
supine morphisms, and/or tensor). Therefore as the basic morphisms between fibration
we consider lax morphisms that only make the square below commute.
A {\em lax morphism} from a fibration $p$ to $p'$ is a pair of functors $(G,G')$ making the square
  \begin{center} \xext=00 \yext=450
\begin{picture}(\xext,\yext)(\xoff,\yoff)
 \setsqparms[1`1`1`1;500`400]
 \putsquare(0,0)[\cE`\cE'`\bB`\bB';G`p`p'`G']
\end{picture}
\end{center}
commute. If if $p$ and $p'$ are (op)fibrations and $G$ preserves prone
(supine) morphisms, $(G,G')$ will be called  a {\em  morphism
of fibrations (opfibrations)}. If $p$ and $p'$ are bifibrations and $G$ preserves
both prone and supine morphisms, $(G,G')$ will be called a {\em morphism of bifibrations}. The {\em
fibred natural transformation} $(\tau,\tau'): (G,G')\ra (H,H')$  is a pair of natural transformations
$\tau: G\ra H$ and $\tau': G'\ra H'$ such that $p'(\tau)=\tau'_p$. If  $g=G'=id_\bB$ and $\tau'=id_{id_\bB}$
then a fibered natural transformation is natural transformation $\tau: G\ra H$ whose
components are vertical morphisms. A {\em fibred left adjoint} to
$(G,G')$ is a fibred morphism $(F,F')$ such that $F\dashv G$ and
$F'\dashv G'$ are adjunctions with units and counits $(\eta,\varepsilon)$ and
$(\eta',\varepsilon')$, respectively, so that $p'(\eta)=\eta'$ and $p'(\varepsilon)=\varepsilon'$.

\subsection{The basic definition}
A {\em lax monoidal fibration} $(p:\cE\ra \bB,I,\otimes, \alpha,\lambda,\varrho)$ is
\begin{enumerate}
  \item a fibration $p:\cE \ra \bB$,
  \item equipped with two lax morphisms of fibrations $\otimes$ and $I$
 \begin{center} \xext=1000 \yext=500
\begin{picture}(\xext,\yext)(\xoff,\yoff)
 \settriparms[1`1`1;500]
 \putqtriangle(0,0)[\cE\times_\bB\cE`\cE`\bB;\otimes``p]
  \settriparms[-1`0`1;500]
 \putptriangle(500,0)[\phantom{\cE}`\bB`\bB;I``1_{\bB}]
\end{picture}
\end{center}
  \item three fibred natural transformations $\alpha$, $\lambda$, $\varrho$
   \begin{center}
\xext=1800 \yext=650
\begin{picture}(\xext,\yext)(\xoff,\yoff)
 \setsqparms[1`1`1`1;1800`550]
 \putsquare(0,50)[\cE\times_\bB(\cE\times_\bB\cE)\cong(\cE\times_\bB\cE)\times_\bB\cE`\cE\times_\bB\cE `\cE\times_\bB\cE`\cE;
 \otimes\times_\bB\cE`\cE\times_\bB\otimes`\otimes`\otimes]
 \put(1500,150){\makebox(50,50){$\Ra$}}
 \put(1500,220){\makebox(50,50){$\alpha$}}
 \end{picture}
\end{center}
(where the unnamed iso $\cong$ is the canonical one between pullbacks) i.e. there are morphisms
\begin{center} \xext=1000 \yext=0
\begin{picture}(\xext,\yext)(\xoff,\yoff)
  \putmorphism(0,0)(1,0)[A\otimes(B\otimes C)`(A\otimes B)\otimes C`\alpha_{A,B,C}]{1000}{1}a
\end{picture}
\end{center}
for any $O\in \bB$ and any $A,B,C \in \cE_O$ so that
$p(\alpha_{A,B,C})=1_O$ and these morphisms are natural in $A$,
$B$, and $C$, in the obvious sense. Moreover
 \begin{center} \xext=1200 \yext=800
\begin{picture}(\xext,\yext)(\xoff,\yoff)
 \settriparms[1`1`1;700]
 \putqtriangle(0,0)[\bB\times_\bB\cE`\cE\times_\bB\cE`\cE;I\times 1_\cE`\pi_2`\otimes]
  \settriparms[-1`0`1;700]
 \putptriangle(700,0)[\phantom{\cE\times_\bB\cE}`\cE\times_\bB\bB`\cE;1_\cE\times I``\pi_1]
  \put(500,250){\makebox(50,50){$\Ra$}}
 \put(500,320){\makebox(50,50){$\lambda$}}
  \put(850,250){\makebox(50,50){$\Ra$}}
 \put(850,320){\makebox(50,50){$\varrho$}}
\end{picture}
\end{center}
i.e. there are morphisms
\begin{center} \xext=1600 \yext=0
\begin{picture}(\xext,\yext)(\xoff,\yoff)
  \putmorphism(0,0)(1,0)[A\otimes I_O`A`\rho_A]{600}{1}a
  \putmorphism(600,0)(1,0)[A`I_O\otimes A`\lambda_A]{600}{1}a
\end{picture}
\end{center}
for any $O\in \bB$ and any $A\in \cE_O$ so that $p(\rho_A)=1_O=p(\lambda_A)$ and these morphisms are natural in $A$.
  \item The diagrams
  \begin{center} \xext=2300 \yext=1000
\begin{picture}(\xext,\yext)(\xoff,\yoff)
\put(800,900){$A\otimes(B\otimes (C\otimes D))$}
      \put(900,850){\vector(-3,-2){400}}
      \put(100,700){$1_A\otimes \alpha_{B,C,D}$}

      \put(1300,850){\vector(3,-2){400}}
      \put(1560,700){$\alpha_{A,B,C\otimes D}$}

\put(0,500){$A\otimes((B\otimes C)\otimes D$}
\put(1600,500){$(A\otimes B)\otimes (C\otimes D)$}
      \put(500,450){\vector(1,-2){180}}
      \put(180,200){$\alpha_{A,B\otimes C,D}$}

      \put(1700,450){\vector(-1,-2){180}}
      \put(1650,200){$\alpha_{A\otimes B,C,D}$}

\put(200,0){$(A\otimes(B\otimes C))\otimes D$}
\put(1400,0){$((A\otimes B)\otimes C)\otimes D$}

  \put(1050,0){\vector(1,0){300}}
      \put(950,-120){$\alpha_{A,B,C}\otimes 1_D$}
\end{picture}
\end{center}
and
  \begin{center} \xext=1200 \yext=550
\begin{picture}(\xext,\yext)(\xoff,\yoff)
 \setsqparms[1`1`-1`1;1200`500]
 \putsquare(0,0)[A\otimes_O B`A\otimes_O B`A\otimes_O(I_O\otimes B)`(A\otimes_OI_O)\otimes B;
 1_{A\otimes B}`1_A\otimes \lambda_B`\rho_A\otimes 1_B`\alpha_{A,I_O,B}]
\end{picture}
\end{center}
commute, and finally $\rho_{I_O}$ and $\lambda_{I_O}$ are isomorphisms
and
\[ \rho_{I_O}=\lambda_{I_O}^{-1}: I_O\otimes I_O\lra I_O \]
%
where $O\in \bB$ and  $A,B,C,D\in \cE_O$. \\  End of the definition of a lax monoidal fibration.
\end{enumerate}

{\bf Remarks}
\begin{enumerate}
  \item The tensor operation can be applied to objects in the same fibre of \mbox{$p:\cE\ra \bB$} only, and to morphisms that
lie over the same map in the base. Sometimes we emphasize this
by writing $a\otimes_O b$ and $f\otimes_u g$ to indicate that the
tensor is in the fibre over $O$ or over a morphism $u$. So the fibres
are (lax) monoidal categories and reindexing 'functors' are lax
monoidal.  But the total category $\cE$ is not monoidal. Both
facts are important for the examples we have in mind.

  \item  For any
$u:O\ra Q \in \bB$ and $A,B\in \cE_Q$ we have (unique) morphisms $\psi^0_u:I_O\ra u^*(I_Q)$ and $\psi^2_{u,A,B}:u^*(A)\otimes_Ou^*(B)\ra u^*(A\otimes_Q B)$
so that the triangles
\begin{center} \xext=1400 \yext=700
\begin{picture}(\xext,\yext)(\xoff,\yoff)
 \settriparms[-1`1`1;600]
 \putAtriangle(0,0)[u^*(A\otimes B)`u^*(A)\otimes u^*(B)`A\otimes B;\psi^2_{u,A,B}`pr_{u,A\otimes B}`pr_{u,A}\otimes pr_{u,B}]
\end{picture}
\end{center}
and
\begin{center} \xext=1000 \yext=500
\begin{picture}(\xext,\yext)(\xoff,\yoff)
 \settriparms[-1`1`1;500]
 \putAtriangle(0,0)[u^*(I_Q)`I_O`I_Q;\psi^0_u`pr_{u,I_Q}`I_u]
\end{picture}
\end{center}
commute, where $pr_{u,A}$ is a prone morphism over $u$ with codomain $A$.

Due to the fact that we deal with fibrations, lax morphisms preserve prone morphisms in the lax sense.
Thus even if we do not require the tensor and the unit to be morphisms of fibrations, we still have that the
'reindexing' functors are lax monoidal, i.e. they 'respect' the monoidal structure (to some extent).

There are many more diagrams involving $\psi$'s, $\alpha$'s, $\lambda$'s and $\varrho$'s that commute.
  \item The directions of the natural transformations $\alpha$'s, $\lambda$'s and $\varrho$ in the definition of a lax monoidal fibration
  are so chosen to cover all the examples we have in mind. But it is sometimes convenient to consider natural transformations $\lambda$ or
  $\varrho$ that go in the other direction.
\end{enumerate}

\subsection{Monoids in a lax monoidal fibration}

A {\em monoid} in a fibre over $O$ is a triple $(M,m:M\otimes M\ra M, e:I_O\ra M)$ where
$M$ is an object in $\cE_O$, $m$, $e$ are morphisms in $\cE_O$ making the diagrams
 \begin{center}
\xext=1800 \yext=750
\begin{picture}(\xext,\yext)(\xoff,\yoff)
 \setsqparms[1`1`1`1;1800`600]
 \putsquare(0,50)[\hskip 3cm M\otimes(M\otimes M)\stackrel{\alpha}{\lra}(M\otimes M)\otimes M`M\otimes M` M \otimes M`M;
 m\otimes 1_M`1_M\otimes m`m`m]
 \end{picture}
\end{center}
and
 \begin{center}
\xext=1600 \yext=550
\begin{picture}(\xext,\yext)(\xoff,\yoff)
 \setsqparms[1`-1`1`1;800`400]
 \putsquare(0,50)[I_O\otimes M`M\otimes M` M`M;
 e\otimes 1_M`\lambda_M`m`1_M]
  \setsqparms[-1`0`1`-1;800`400]
 \putsquare(800,50)[\phantom{M\otimes M,}` M\otimes I_O`\phantom{M}`M;
 1_M\otimes e ``\rho_M`1_M]
 \end{picture}
\end{center}
commute.

A {\em morphism of monoids} $f:(M,m,e)\lra (M',m',e')$ over $u :O\ra Q\in \bB$ is a morphism $f:M\lra M'$ in $\cE$ over $u$
such that the squares
 \begin{center}
\xext=800 \yext=950
\begin{picture}(\xext,\yext)(\xoff,\yoff)
 \setsqparms[1`1`1`0;800`400]
 \putsquare(0,450)[M\otimes_O M` M' \otimes_Q M'`\phantom{M}`\phantom{M'};
 f\otimes_uf`m`m'`]
 \setsqparms[1`-1`-1`1;800`400]
 \putsquare(0,50)[M`M'`I_O`I_Q;f`e`e'`I_u]
 \end{picture}
\end{center}
commute.

Then the category of monoids is again fibred over $\bB$ and the
forgetful functor is a morphism of fibrations
\begin{center} \xext=600 \yext=580
\begin{picture}(\xext,\yext)(\xoff,\yoff)
 \settriparms[-1`1`1;500] \putVtriangle(0,0)[\cE`Mon(\cE,\otimes, I)`Set;{\cal U} `p`q]
\end{picture}
\end{center}
In the interesting cases it should have a (fibred) left adjoint
which is not likely to be a morphism of fibrations.

\vskip 2mm
{\bf Remark}  As we will see it is not always true that the category of
all monoids is of real interest. If the coherence transformations are indeed not isomorphisms
  it may happen that we may want to consider only those monoids that
  satisfy some additional conditions, see \ref{single-tensor}.

\subsection{The 2-category of lax monoidal fibrations}

A {\em morphism of lax monoidal fibrations}
\[ (F,K,\varphi_0,\varphi_2):(p:\cE\ra \bB,\otimes,I,\alpha,\lambda,\varrho)\lra (p':\cE'\ra \bB',\otimes',I',\alpha',\lambda',\varrho') \]
is data 1-3 subject to conditions 4-6 below ($O\in \bB$, $A,B,C\in \cE_O$):
\begin{enumerate}
  \item $(F,K):(\cE,p)\ra (\cE',p')$ a lax morphism of fibrations,
  \item $\varphi_0: I_K \lra F\circ I$ a fibred natural transformation
  (i.e. for any $O\in \bB$ we have a morphisms $(\varphi_0)_O: I'_{K(O)} \lra F(I_O)$
  in $\cE'_O$ which is natural in $O$),
  \item $\varphi_2: F(-)\otimes' F(=) \lra F((-)\otimes(=))$ a fibred natural transformation (i.e. for $A,B\in \cE_O$ we have a morphism
  $(\varphi_2)_{A,B}: F(A)\otimes' F(B) \lra F(A\otimes B)$), in $\cE'_O$ which is natural in $A$ and $B$)
  \item the square
\begin{center} \xext=1500 \yext=600
\begin{picture}(\xext,\yext)(\xoff,\yoff)
 \setsqparms[1`-1`1`1;1500`400]
  \putsquare(0,50)[I'_{K(O)}\otimes F(A)`F(I_O)\otimes F(A)`F(A)`F(I_O\otimes A);
  \varphi_0\otimes 1_{F(A)}`\lambda'_{F(A)}`\varphi_2`F(\lambda_A)]
\end{picture}
\end{center}
commutes, where $O\in \bB$ and  $A\in \cE_O$;
  \item  the square \begin{center} \xext=1500 \yext=600
\begin{picture}(\xext,\yext)(\xoff,\yoff)
 \setsqparms[1`1`1`-1;1500`400]
  \putsquare(0,50)[F(A)\otimes I'_{K(O)} ` F(A)\otimes F(I_O)`F(A)`F( A\otimes I_O);
1_{F(A)}\otimes  \varphi_0`\varrho'_{F(A)}`\varphi_2`F(\varrho_A)]

\end{picture}
\end{center}
commutes, where $O\in \bB$ and  $A\in \cE_O$;
  \item the diagram
  \begin{center} \xext=2000 \yext=1050
\begin{picture}(\xext,\yext)(\xoff,\yoff)
 \setsqparms[1`1`1`0;2000`450]
  \putsquare(0,500)[F(A)\otimes' (F(B)\otimes' F(C))`(F(A)\otimes' F(B))\otimes' F(C)`F(A)\otimes' F(B\otimes C)`F(A\otimes B)\otimes' F(C);
  \alpha'_{F(A),F(B),F(C)}`1_{F(A)}\otimes'\varphi_2`\varphi_2\otimes' 1_{F(C)}`]
  \setsqparms[0`1`1`1;2000`450]
  \putsquare(0,50)[\phantom{F(A)\otimes' F(B\otimes C)}`\phantom{F(A\otimes B)\otimes' F(C)}`F(A\otimes (B\otimes C))`F((A\otimes B)\otimes C);
  `\varphi_2`\varphi_2`F(\alpha_{A,B,C})]
\end{picture}
\end{center}
  commutes, where $O\in \bB$ and  $A,B,C\in \cE_O$.
\end{enumerate}
End of the definition of a morphism of lax monoidal fibrations.

\vskip 3mm
A morphism of lax monoidal fibrations is called {\em strong} if the transition
morphisms $\varphi_0$, $\varphi_2$ are isomorphisms and $(F,K)$ is a morphism of fibrations.
\vskip 3mm

A {\em transformation} between two morphisms of lax monoidal fibrations is a pair of natural transformations
\[ (\tau,\sigma): (F,K,\varphi_0,\varphi_2)\lra (F',K',\varphi'_0,\varphi'_2)\]
such that
\begin{enumerate}
  \item $\sigma : K\lra K'$ and $\tau : F\lra F'$ are natural transformations,
  \item $p'(\tau)=\sigma_p$, i.e. $\tau$ is fibred over $\sigma$,
  \item the diagrams
\begin{center} \xext=800 \yext=500
\begin{picture}(\xext,\yext)(\xoff,\yoff)
\setsqparms[1`1`1`1;800`450]
 \putsquare(0,0)[I'_{K(O)}`I'_{K'(O)}`F(I_O)`F'(I_O);I'_{\sigma_O}`(\phi_0)_O`(\phi'_0)_O`\tau_{I_O}]
\end{picture}
\end{center}
and
\begin{center}
\xext=1600 \yext=650
\begin{picture}(\xext,\yext)(\xoff,\yoff)
 \setsqparms[1`1`1`1;1600`500]
 \putsquare(0,0)[F(A)\otimes_{K(O)}F(B)`F'(A)\otimes_{K(O)}F'(B)`F(A\otimes_O B)`F'(A\otimes_O B);
 \;\;\;\;\tau_A\otimes_{K(O)}\tau_B`(\phi_2)_{A,B}`(\phi'_2)_{A,B}`\tau_{A\otimes_O B}]
\end{picture}
\end{center}
commute, for $O\in\bB$.
\end{enumerate}
End of the definition of a transformation between two morphisms of lax monoidal fibrations.

\vskip 3mm

\begin{proposition}\label{mor-lmf-induce}
The morphisms of lax monoidal fibrations induce morphisms between the fibrations of monoids.
The transformations between morphisms of lax monoidal fibrations induce natural transformations between the induced functors.
\end{proposition}

{\it Proof.}~ Exercise. $~\Box$

\subsection{Simple examples} 

1. {\em Categories.} Probably the simplest non-trivial example of a lax monoidal
fibration (in the above sense) is the fibration of graphs over sets, say $p:Gph\lra Set$,
  where  $p$ sends the parallel pair of arrows to their common codomain. The tensor
  \[(A, d^A,c^A : A\ra O)\otimes_O (B, d^B,c^B : B\ra O)= (A\times_O B, c^A\circ \pi_1, d^B\circ \pi_2 : A\times_O B\lra O)\]
where $A\times_O B$ denotes the pullback of the following pair of morphisms
   \begin{center}
   \xext=500 \yext=450
\begin{picture}(\xext,\yext)(\xoff,\yoff)
 \setsqparms[1`1`1`1;500`400]
 \putsquare(0,0)[A\times_O B`B`A`O;
 \pi_2`\pi_1`c^B`d^A]
 \end{picture}
\end{center}
The unit for the tensor in the fibre over $O$ is a pair of identities on $O$, $(O, 1_O,1_O:O\ra O)$.
The total category of the fibration of monoids $q:Mon(Gph)\ra Set$ in this fibration is
the category of small categories and functors. The monoids in a fibre  $Mon(Gph)_O$
are categories with the set of objects $O$.

2. {\em Lambek's multicategories.} Let $((-)^*,\eta,\mu)$ be the monad for monoids on the category $Set$.
Then, we can define a fibration of multisorted signatures $p_m:Sig_m\lra Set$ as follows.
An object of $Sig_m$ in the fibre over the set $O$ is a triple $(A,\partial,O)$ such that $A$ is a set and
$\partial:A\lra O\times O^*$  function.
$(f,u): (A,\partial^A,O)\lra (A',\partial^{A'},O')$ is a morphism  in $Sig_m$ over a function $u:O\ra O'$
if $f:A\ra A'$ is a function making the square
\begin{center} \xext=800 \yext=500
\begin{picture}(\xext,\yext)(\xoff,\yoff)
 \setsqparms[1`1`1`1;800`400]
  \putsquare(0,50)[A`A'`O\times O^*`O'\times O'^*;f`\partial^A`\partial^{A'}`u\times u^*]
\end{picture}
\end{center}
commute. The tensor $(A,\partial^A,O)\otimes(B,\partial^B,O) =(A\times_{O^*}B^*,\partial,O)$ of
two object in the same fibre is given by the pullback and multiplication $\mu$ in the monad $(-)^*$
and the unit for this tensor is $(I_O,\partial,O)=(O,\lk 1_O, \eta \rk,O)$.
The category of monoids in this fibration is the category of Lambek's multicategories.
As this construction will be described in Section \ref{burroni} in a more general case of
arbitrary monad over a category with pullbacks we don't go into the details here.

\section{Actions of lax monoidal fibrations}\label{sec-actions of lmf}

\subsection{The basic definition}
An {\em action} $(\star, \psi_2, \psi_0)$ of a lax monoidal fibration
$(\cE,p,I,\otimes, \alpha,\lambda,\varrho)$ on a fibration $\pi:\cX\ra \bB$ is
\begin{enumerate}
  \item a lax morphism of fibrations
\begin{center} \xext=1000 \yext=500
\begin{picture}(\xext,\yext)(\xoff,\yoff)
 \settriparms[1`1`1;500]
 \putVtriangle(0,0)[\cE\times_\bB\cX`\cX`\bB;\star`p`\pi]
\end{picture}
\end{center}

  \item a fibred natural transformation $\psi_0$
     \begin{center} \xext=1000 \yext=550
\begin{picture}(\xext,\yext)(\xoff,\yoff)
 \settriparms[1`1`1;500]
 \putVtriangle(0,0)[\cX\cong\bB\times_{\bB}\cX`\cE\times_{\bB}\cX`\cX;I\times_\bB 1_\cX`1_\cX`\star]
    \put(500,300){$\psi_0$}
    \put(500,200){$\Rightarrow$}
\end{picture}
\end{center}
i.e. for $X\in\cX_O$, we have a morphism \[ (\psi_0)_X : X\ra I_O\star X \]
in the fibre $\cX_O$ which is natural in $O$;
  \item a fibred natural transformation $\psi_2$
     \begin{center}
\xext=1800 \yext=650
\begin{picture}(\xext,\yext)(\xoff,\yoff)
 \setsqparms[1`1`1`1;1800`550]
 \putsquare(0,50)[\hskip 15mm\cE\times_\bB(\cE\times_\bB\cX)\cong(\cE\times_\bB\cE)\times_\bB\cX`\cE\times_\bB\cX `\cE\times_\bB\cX`\cX;
\otimes\times_\bB\cX`\cE\times_\bB\star`\star`\star]
 \put(1500,150){\makebox(50,50){$\Ra$}}
 \put(1500,220){\makebox(50,50){$\psi_2$}}
 \end{picture}
\end{center}
i.e. for $O\in \bB$, $A,B\in \cE_O$ and $X\in \cX_O$  we have a morphism
\[ (\psi_2)_{A,B,X} : A\star(B\star X)\lra (A\otimes B)\star X  \]
in the fibre $\cX_O$, which is natural in $A$,  $B$ and $X$; the unnamed iso $\cong$
is the canonical one between pullbacks;
  \item making the pentagon
  \begin{center} \xext=2300 \yext=1000
\begin{picture}(\xext,\yext)(\xoff,\yoff)
\put(800,900){$A\star(B\star (C\star X))$}
      \put(900,850){\vector(-3,-2){400}}
      \put(120,740){$(\psi_2)_{A, B, C\star X}$} 

      \put(1300,850){\vector(3,-2){400}}
      \put(1520,740){$1_A\star(\psi_2)_{B,C,X}$} 

\put(0,500){$(A\otimes B)\star (C\star X)$} 
\put(1600,500){$A\star ((B\otimes C)\star X)$}
      \put(500,450){\vector(1,-2){180}}
      \put(40,230){$(\psi_2)_{A\otimes B, C,X}$} 

      \put(1700,450){\vector(-1,-2){180}}
      \put(1640,230){$(\psi_2)_{A,B\otimes C, X}$} 

\put(200,0){$((A\otimes B)\otimes C)\star X$} 
\put(1400,0){$(A\otimes (B\otimes C))\star X$}

  \put(1380,20){\vector(-1,0){300}}
      \put(950,150){$\alpha_{A,B,C}\star 1_X$} 
\end{picture}
\end{center}
\item  and two squares
\begin{center} \xext=2800 \yext=510
\begin{picture}(\xext,\yext)(\xoff,\yoff)
\setsqparms[1`1`1`1;900`450]
 \putsquare(0,0)[A\star X `I\star(A\star X) `A\star X`(I\otimes A)\star X ;(\psi_0)_{A\star X}`1_{A\star X}`(\psi_2)_{I,A,X}`\lambda_A\star 1_X]
 \setsqparms[1`1`1`-1;900`450]
 \putsquare(1800,0)[A\star X `A\star(I\star X) `A\star X `(A\otimes I)\star X;1_A\star (\psi_0)_X`1_{A\star X}`(\psi_2)_{A,I,X}`\varrho_A\star 1_X]
 \end{picture}
\end{center}
commute,  where $O\in \bB$, $A,B,C\in \cE_O$ and $X\in \cX_O$. \\  End of the definition of an action of a lax monoidal fibration on a fibration.
\end{enumerate}

\vskip 3mm
An action of a lax monoidal fibration is called {\em strong}
if the transition morphisms $\psi_0$, $\psi_2$ are isomorphisms and $\star$ is a morphism of fibrations.
\vskip 3mm

\subsection{Actions of monoids along an action of a lax monoidal fibration}

Let $(\star, \psi_2, \psi_0)$ be an action of a lax monoidal fibration $(p,I,\otimes, \alpha,\lambda,\varrho)$ on a fibration $\pi:\cX\ra \bB$, $O$ an object of $\bB$,
$(M,m,e)$ a monoid in $Mon(\cE)_O$, $X$ an object in $\cX_O$ and $\nu : M\star X \ra X$ a morphism in $\cX_O$.
The pair $(X,\nu)$ is an {\em action} of $(M,m,e)$ on $X$ along the action $(\star, \psi_2, \psi_0)$ (or just $\star$, for short)
if the following diagrams
  \begin{center}
\xext=2100 \yext=600
\begin{picture}(\xext,\yext)(\xoff,\yoff)
 \setsqparms[0`1`1`1;2100`450]
 \putsquare(0,0)[M\star(M\star X)`M\star X `M\star X `X;`1_M\star\nu`\nu`\nu]
 \putmorphism(0,450)(1,0)[\phantom{M\star(M\star X)}`(M\otimes M)\star X`(\psi_2)_{M,M,X}]{1100}{1}a
 \putmorphism(1100,450)(1,0)[\phantom{M\star(M\star X)}`\phantom{M\star X}`m\star 1_X]{1000}{1}a
 \end{picture}
\end{center}
and
   \begin{center}
\xext=800 \yext=600
\begin{picture}(\xext,\yext)(\xoff,\yoff)
 \setsqparms[1`-1`1`1;800`450]
 \putsquare(0,50)[I\star X`M\star X`X`X ;e\star 1_X`(\psi_0)_{X}`\nu`1_X]
 \end{picture}
\end{center}
commute. A {\em morphism of actions} $(f,g,u) : (M,X,\nu)\lra (M',X',\mu')$
is a triple of morphisms $u:O\ra O'$ in $\bB$, and $f : M \ra M'$ in $Mon(\cE)$,
$g:X\ra X'$ in $\cX$ both over $u$ so that the square in $\cX$
 \begin{center} \xext=800 \yext=450
\begin{picture}(\xext,\yext)(\xoff,\yoff)
 \setsqparms[1`1`1`1;800`400]
 \putsquare(0,0)[M\star X`M'\star X'`X`X';f\star g`\nu`\nu'`g]
\end{picture}
\end{center}
commutes.

The category of actions $Act(\cE,\cX,\star)$ is fibred over $Mon(\cE)$.
It might happen that monoids in the fibre over $O$ can be interpreted as algebras for
a single monoid $\cM_O$ in a fibre over $K(O)$.
If this association  is functorial we have commuting squares
 \begin{center} \xext=800 \yext=1100
\begin{picture}(\xext,\yext)(\xoff,\yoff)
 \setsqparms[1`1`1`1;800`500]
 \putsquare(0,500)[Mon(\cE)`Act(\cE,\cX,\star)`\bB`Mon(\cE);R`q`\pi_\mu`]
  \setsqparms[1`1`1`1;800`500]
 \putsquare(0,0)[\phantom{\bB}`\phantom{Mon(\cE)}`\bB`\bB;\cM`1_\bB`q`K]
\end{picture}
\end{center}
The functor $R$ is the representing functor that interprets monoids as algebras.
If the upper square is a pullback then we say that $\cM$ is the functor of
{\em metamonoid}\footnote{This is what seems to be the intension of the notion of operad
for (colored) operads introduced by J.Baez and J.Dolan in \cite{BaezDolan}.}
and the triple  $(R,\cM,K)$ {\em strongly represents}
monoids in $Mon(\cE)$ as algebras. This means in particular that the category of $O$-monoids is equivalent to the category
of $\cM_O$-algebras. If the functor $R$ is an embedding (faithful and reflects isomorphisms) on fibres then we say that the triple
 $(R,\cM,K)$ {\em weakly represents} monoids in $Mon(\cE)$.

\subsection{The 2-category of actions of lax monoidal fibrations}

We define below the morphisms of actions of lax monoidal fibrations and
transformations of such morphisms. In that way we shall define
the 2-category $ACTION$  of actions of lax monoidal fibrations on fibrations.

A {\em morphism of actions of lax monoidal fibrations}
\[ (F,H,K,\phi_0,\phi_2,\tau):  (\cE,p,\cX,\pi,\star,\psi_0,\psi_2)\lra  (\cE',p',\cX',\pi',\star',\psi'_0,\psi'_2)  \]
consists of the data 1-4 subject to the conditions 5-6 below:
\begin{enumerate}
  \item functors
  \[ F:\cE\lra \cE',\hskip 5mm H:\cX\lra \cX',\hskip 5mm K:\bB\lra \bB' \]
  \item a morphism of lax monoidal fibrations
  \[ (F,K,\phi_0,\phi_2):  (\cE,p,\otimes, \alpha,\lambda,\varrho)\lra  (\cE',p',\otimes', \alpha',\lambda',\varrho') \]
  \item a lax morphism of fibrations
  \[ (H,K):  (\cX,\pi)\lra  (\cX',\pi')  \]
  \item a natural transformation
  \[ \tau : \star'\circ (F\times_KH)  \lra H \circ \star \]
  i.e. we have a morphism
  \[ \tau_{A,X}:F(A)\star'H(X)\lra H(A\star X)\]
  which is natural in  $A\in \cE_O$, $X\in\cX_O$ and $O\in \bB$.

  So we have a diagram
  \begin{center} \xext=1400 \yext=950
\begin{picture}(\xext,\yext)(\xoff,\yoff)
 \setsqparms[1`1`1`1;1400`900]
  \putsquare(0,0)[\cE\times_\bB\cX`\cX`\cE'\times_\bB\cX'`\cX';\star`F\times_KH`H`\star'  ]

  \putmorphism(700,600)(0,1)[\bB`\bB'`K]{300}{1}r

  \putmorphism(300,130)(3,1)[\phantom{\cE\times_\bB\cX}`\phantom{\bB'}`]{140}{1}r
  \putmorphism(1100,130)(-3,1)[\phantom{\bB'}`\phantom{\cX'}`\pi']{140}{1}l

  \putmorphism(1000,700)(3,1)[\phantom{\cE\times_\bB\cX}`\phantom{\bB'}`\pi]{140}{-1}l
  \putmorphism(400,700)(-3,1)[\phantom{\bB'}`\phantom{\cX'}`]{140}{-1}l
\end{picture}
\end{center}
  in which the triangles and internal squares commute, but external square commutes up to the natural transformation $\tau$.
  \item The square
  \begin{center} \xext=1000 \yext=650
\begin{picture}(\xext,\yext)(\xoff,\yoff)
 \setsqparms[1`1`1`-1;1000`500]
\putsquare(0,0)[H(X)`I'_{K(O)}\star' H(X)`H(I_O\star
X)`F(I_O)\star'
H(X);(\psi'_0)_{H(X)}`H((\psi_0)_X)`(\phi_0)_O\star'1_{H(X)}`\tau_{I_O,X}]
\end{picture}
\end{center}
  commutes, for $X\in \cX_O$, and $O\in\bB$.
  \item The hexagon
  \begin{center} \xext=2000 \yext=1500
\begin{picture}(\xext,\yext)(\xoff,\yoff)
\put(600,1350){$F(A)\star'(F(B)\star' H(X))$}
      \put(900,1300){\vector(-3,-2){400}}
      \put(-120,1180){$(\psi'_2)_{F(A),F(B),H(X)}$}

     \put(1300,1300){\vector(3,-2){400}}
    \put(1520,1180){$1_{F(A)}\star'\tau_{B,X}$}

  \put(0,900){$(F(A)\otimes'F(B))\star' H(X)$}
  \put(1600,900){$F(A)\star' H(B\star X)$}
     \put(450,850){\vector(0,-1){350}}
     \put(-270,650){$(\phi_2)_{A,B}\star' 1_{H(X)}$}

    \put(1750,850){\vector(0,-1){350}}
     \put(1800,650){$\tau_{A, B\star X}$}

\put(150,400){$F(A\otimes B)\star'H(X)$}
\put(380,130){$\tau_{A\otimes B,X}$}

\put(1500,400){$H(A\star (B\star X))$}

\put(800,0){$H((A\otimes B)\star X)$}
\put(1660,350){\vector(-3,-2){400}}

     \put(600,350){\vector(3,-2){400}}

\put(1500,130){$H((\psi_2)_{A,B,X})$}

\end{picture}
\end{center}
  commutes, for $A,B\in\cE_O$,  $X\in \cX_O$, and $O\in\bB$.
\end{enumerate}
End of the definition of a morphism of actions of lax monoidal fibrations.

\vskip 2mm

Let \[ (F,H,K,\phi_0,\phi_2,\tau), (F',H',K',\phi'_0,\phi'_2,\tau):  (\cE,p,\cX,\pi,\star,\psi_0,\psi_2)\lra  (\cE',p',\cX',\pi',\star',\psi'_0,\psi'_2)  \]
be two morphisms of actions of lax monoidal fibrations.
A {\em transformation of morphisms of actions of lax monoidal fibrations}
\[ (\zeta_0,\zeta_1,\zeta_2): (F,H,K,\phi_0,\phi_2,\tau)\lra (F',H',K',\phi'_0,\phi'_2,\tau) \]
consists of data 1-3 subject to the condition 4 below:
\begin{enumerate}
  \item natural transformations $\zeta_2:F\lra F'$, $\zeta_1:H\lra H'$ and $\zeta_0:K\lra K'$;
  \item a transformation of lax monoidal fibrations
  \[  (\zeta_2,\zeta_0): (F,K,\phi_0,\phi_2)\lra (F',K',\phi'_0,\phi'_2); \]
  \item a fibred natural transformation of morphisms of fibrations
  \[ (\zeta_1,\zeta_0): (H,K)\lra (H',K') \]
  between fibrations $(\cX,\pi)$ and $(\cX',\pi')$;
  \item so that the square
\begin{center} \xext=1200 \yext=650
\begin{picture}(\xext,\yext)(\xoff,\yoff)
 \setsqparms[1`1`1`1;1200`550]
\putsquare(0,0)[F(A)\star'H(X)`H(A\star
X)`F'(A)\star'H'(X)`H'(A\star
X);\tau_{A,X}`(\zeta_2\star\zeta_1)_{A,X}`(\zeta_1)_{A\star X
}`\tau'_{A\star X}]
\end{picture}
\end{center}
  commutes, for $A\in \cE_O$ and $X\in \cX_O$.
\end{enumerate}

\begin{proposition}
The morphisms of actions of lax monoidal fibrations induce morphisms
between fibrations of actions of monoids along actions of
monoidal fibrations.
The transformations between morphisms of  actions of lax monoidal
fibrations induce natural transformations between the induced functors.
\end{proposition}

{\it Proof.}~ Exercise. $~\Box$

\subsection{Simple examples}

\begin{enumerate}
  \item The lax monoidal fibration of graphs $p:Gph\ra Set$ acts naturally
  on the basic fibration $cod: Set^\ra \lra Set$. The action of a graph
  $(d,c:A\ra O)$ on a function $F:X\ra O$ is defined as the composition of
  the horizontal arrows on the top of the following diagram
      \begin{center} \xext=1100 \yext=450
\begin{picture}(\xext,\yext)(\xoff,\yoff)
 \setsqparms[-1`1`1`-1;600`400]
  \putsquare(500,0)[A`A\star X`O`X^*;`c``\xi]
  \putmorphism(0,400)(1,0)[O`\phantom{A}`d]{500}{-1}a
\end{picture}
\end{center}
in which the square is a pullback. Then the action of monoids along this
action are all presheaves on all small categories.
\item If we replace in the previous example $Set$ by any category $\cC$
with pullbacks we get all internal presheaves on all internal categories in $\cC$.

  \item The fibration of multisorted signatures also acts on the basic fibration.
  But this example will be described in Section \ref{burroni}.
 \end{enumerate}

\section{The exponential fibrations}\label{sec-exp-fib}
Most of the material of this section belongs to folklore. We present it here as we need it later in this form.
$\bCat$ is the category of large categories,
so $Set$ and $Cat$ are objects of $\bCat$.

If $X$ is an object of a cartesian closed category $\cC$ then $X^X$ carries a natural structure of a monoid.
If $p_\cX:\cX\ra \bB$ is a fibration then we can form an exponential fibration $p:[\cX\Ra\cX]\ra \bB$ in $Fib(\bB)$, the category of fibrations over $\bB$,
which also carries a  natural structure of a lax monoidal fibration. Then any strong action
of a lax monoidal fibration $p_\cE :\cE \ra \bB$ on $p_\cX:\cX\ra \bB$ gives rise to a representation of
$p_\cE :\cE \ra \bB$ in the lax monoidal fibration of the internal endomorphisms of $p_\cX:\cX\ra \bB$, i.e.
a strong morphism of lax monoidal fibrations from $p_\cE :\cE \ra \bB$ to $p:[\cX\Ra\cX]\ra \bB$ in $Fib(\bB)$.
However, the examples of actions we have in mind, are almost never strong. But even in this case
we can still find reasonable representations if we will consider the exponentiation of $p_\cX:\cX\ra \bB$ in $\bCat/\bB$ instead of $Fib(\bB)$.
To distinguish these two kinds of exponentiation we denote the exponential object $p_\cX:\cX\ra \bB$ to $p_\cY:\cY\ra \bB$
in $\bCat/\bB$ as $p:\cY^\cX\ra \bB$. It is well known, cf. \cite{G}, for $p:\cY^\cX\ra \bB$ to be a well defined object of $\bCat/\bB$ it is
necessary and sufficient for $\cX$ to be a so called Conduch\'{e} fibration. But as we want $p:\cY^\cX\ra \bB$ to be a fibration we shall assume that $p_\cX$
is a bifibration, i.e. both fibration and opfibration. In fact, as we are mainly interested in the case where $\cX=\cY$,
in order to get a better description $p:\cY^\cX\ra \bB$, it won't be a big restriction when we shall
assume that both $\cX$ and $\cY$ are bifibrations.

\subsection{The exponential bifibrations in $\bCat/\bB$ }

For any bifibration  $p_\cX :\cX\ra \bB$ the exponential fibration
$p_{exp} : \cX^\cX\ra \bB$ in $\bCat/\bB$  is a lax monoidal fibration
with tensor being the composition of functors in fibres.  The monoids
in a fibre $\cX^\cX$ over $B$ are monads on $\cX_B$, and a morphism of
monoids over $u$ is a usual morphism of monads whose functor part is $u^*$.
The counit of the exponential adjunction, the evaluation
$ev_\cX : \cX^\cX\times_\bB\cX\lra \cX$, is an action of the lax monoidal fibration
$p_{exp} : \cX^\cX\ra \bB$ on  $p_\cX :\cX\ra \bB$. Finally, the algebras for this
action are Eilenberg-Moore algebras for all the monads taken together.
As we shall need it later, we shall describe all of this below in detail.
The case of real interest in this paper is when $p_\cX:\cX\lra \bB$
is a basic bifibration $cod:\cC^\ra\lra \cC$ of the category $\cC$
with pullbacks and very likely being just $Set$.

Let $\jeden$, $\dwa$, $\trzy$ be the obvious categories generated by the graphs
$\{ \bullet \}$, $\{ \bullet\ra \bullet \}$, $\{ \bullet\ra \bullet\ra \bullet \}$, respectively.
For an object $B\in\bB$, $r_B:\jeden \ra \bB$ is the functor picking the object $B$.
Similarly, $r_u:\dwa\ra \bB$ is a functor picking the morphism $u:B'\ra B$ in $\bB$,
and $r_{u,v}:\trzy \ra \bB$ is a morphism picking a composable pair $u\circ v$ in $\bB$.
Let $p_{\cX}:\cX\ra \bB$ and  $p_{\cY}:\cY\ra \bB$ be two bifibration. We can form pullbacks
\begin{center} \xext=2500 \yext=400
\begin{picture}(\xext,\yext)(\xoff,\yoff)
\setsqparms[1`1`1`1;500`350]
 \putsquare(0,0)[\cX_B`\cX`\jeden`\bB;``p_\cX`r_B]
 \putsquare(1000,0)[\cX_u`\cX`\dwa`\bB;``p_\cX`r_u]
 \putsquare(2000,0)[\cX_{u,v}`\cX`\trzy`\bB;``p_\cX`r_{u,v}]
 \end{picture}
\end{center}
in the category $\bCat$, i.e. products in $\bCat/\bB$. If the exponential object
$p_{exp}:\cY^\cX\lra \bB$ exists in $\bCat/\bB$ then the objects of $\cY^\cX_B$ correspond
to morphisms from $r_B$ to $p_{exp}$ in $\bCat/\bB$ and we have a sequence of correspondences
\begin{center} \xext=1000 \yext=850
\begin{picture}(\xext,\yext)(\xoff,\yoff)

       \putmorphism(200,800)(1,0)[\jeden`\cY^\cX`]{600}{1}a
       \put(0,500){\line(1,0){1000}}
         \put(200,740){\vector(2,-1){250}}
       \put(800,740){\vector(-2,-1){250}}
       \put(460,540){$\bB$}
       \put(200,640){$r_B$}
       \put(700,640){$p_{exp}$}

       \putmorphism(200,400)(1,0)[\cX_B`\cY`]{600}{1}a
       \put(200,340){\vector(2,-1){250}}
       \put(800,340){\vector(-2,-1){250}}
       \put(460,140){$\bB$}
       \put(700,240){$p_\cY$}

        \put(0,100){\line(1,0){1000}}
        \putmorphism(200,0)(1,0)[\cX_B`\cY_B`]{600}{1}a
 \end{picture}
\end{center}
showing that we can (and we will) identify objects of $\cY^\cX_B$ with functors from $\cX_B$ to $\cY_B$.
Similarly, using $r_u$ and $r_{u,v}$ we see that morphisms in $\cY^\cX$ over a morphism $u$ are functors
from $\cX_u$ to $\cY_u$ commuting over $\dwa$, and the composable pairs are morphisms from $\cX_{u,v}$
to $\cY_{u,v}$ commuting over $\trzy$. The Conduch\'{e} condition is saying in elementary
terms\footnote{This condition will never be used in the explicit form but for the interested reader
we recall it here, cf. \cite{G}. The functor $p_\cX :\cX\ra \bB$ is called a {\em Conduch\'{e} fibration},
if for any morphism $f$ in $\cX$ and a pair of morphisms $u$, $v$ in $\bB$ such that $p_\cX(f)=u\circ v$,
there are morphisms $g$ and $h$ in $\cX$ such that $f=g\circ h$, $p_\cX(g)=u$ and $p_\cX(h)=v$. Moreover,
such a factorization of $f$ is unique up to a zigzag of morphisms in $\cX$ that belong to the fibre over
the domain of $u$.} that for any composable pair of morphisms
$v:B''\ra B'$, $u:B'\ra B$ in $\bB$ the square of the obvious embeddings
\begin{center} \xext=500 \yext=400
\begin{picture}(\xext,\yext)(\xoff,\yoff)
\setsqparms[1`1`1`1;500`350]
 \putsquare(0,0)[\cX_{B'}`\cX_u`\cX_v`\cX_{u,v};d_u`c_v`\kappa_u`\kappa_v]
 \end{picture}
\end{center}
is a pushout in $\bCat$. Then the composition of morphisms $F:\cX_u\ra \cY$ over
$u$ and $G:\cX_v\ra \cY$ over $v$ such that $F\circ d_u=G\circ c_v$ is
the unique functor $[F,G]:\cX_{u,v}\ra \cY$ such that $[F,G]\circ \kappa_u=F$
and $[F,G]\circ \kappa_v=G$ composed with the embedding $\cX_{u\circ v}\ra \cX_{u,v}$.

Recall that for $u:B'\ra B\in \bB$ we have the reindexing functor $u^*:\cX_{B'}\ra \cX_B$
and the coreindexing functor $u_!:\cX_B\ra \cX_{B'}$ defined with the use
of prone and supine morphisms in $\cX$. We denote such functors in different fibrations by the same symbols.
The following Lemma describes morphisms in $\cY^\cX$ more conveniently in five different ways.

\begin{lemma}\label{morphXX} Let $p_{\cX}:\cX\ra \bB$ and  $p_{\cY}:\cY\ra \bB$ be two bifibrations.
Let $u:B'\ra B$ be a morphism  in $\bB$, $Q$ an object in $\cY^\cX_{B'}$ i.e. a functor from
$\cX_{B'}$ to $\cY_{B'}$,  $P$ an object in $\cY^\cX_B$ i.e. a functor from $\cX_{B}$ to $\cY_{B}$.
There is a natural correspondence between
\begin{enumerate}
  \item functors from $F:\cX_u\lra \cY_u$ over $\dwa$ such that $F\circ d_u = Q$ and $F\circ c_u = P$;
   \item\label{morXX} natural transformations $\tau : Qu^*\lra u^*P$ in $\bCat(\cX_B,\cY_{B'})$;
   \item natural transformations $\sigma : u_! Q\lra Pu_!$ in $\bCat(\cX_{B'},\cY_B)$.%
  \item natural transformations $\overline{\tau} : u_!Qu^*\lra P$ in $\bCat(\cX(B),\cY(B))$;
  \item natural transformations $\overline{\sigma} : Q\lra u^*Pu_!$ in $\bCat(\cX(B'),\cY(B'))$.
\end{enumerate}
Moreover, if both $p_\cX$ and $p_\cY$ are cartesian bifibrations and both $P$ and $Q$ (weakly)
preserve pullbacks in fibres, then under the above correspondences the (weakly) cartesian natural transformations
correspond to the (weakly) cartesian natural transformations.
\end{lemma}

Note that in the above Lemma, the two occurrences of the symbols $u^*$ in 2., and $u_!$ in 3.,  do NOT denote the same functors!
In each of the conditions 2. to 5. one of the functors $u^*$, $u_!$ refers to the bifibration $p_{\cX}:\cX\ra \bB$ and
one of the functors $u^*$, $u_!$ to the other  bifibration   $p_{\cY}:\cY\ra \bB$.

\vskip 2mm

{\it Proof.}~  As the conditions 4. and 5. are easily seen to be equivalent
to 2. and 3., respectively, we shall concentrate on equivalence 1., 2., 3.

Fix $u:B'\ra B$ in $\bB$, a functor
\begin{center} \xext=1000 \yext=320
\begin{picture}(\xext,\yext)(\xoff,\yoff)
       \putmorphism(200,260)(1,0)[\cX_u`\cY_u`F]{600}{1}a
       \put(200,200){\vector(2,-1){250}}
       \put(800,200){\vector(-2,-1){250}}
       \put(460,0){$\bB$}
 \end{picture}
\end{center}
as in 1., and two natural transformations
\[  \tau : Qu^* \lra u^*P, \hskip 2cm \sigma : u_!Q\lra Pu_!  \]
from the diagram
\begin{center} \xext=800 \yext=580
\begin{picture}(\xext,\yext)(\xoff,\yoff)
\setsqparms[1`0`0`1;800`500]
 \putsquare(80,0)[\cX_{B'}`\cY_{B'}`\cX_B`\cY_B;Q```P]
 \putmorphism(0,480)(0,-1)[``u^*]{440}{-1}l
 \putmorphism(100,480)(0,-1)[``u_!]{440}{1}r

  \putmorphism(800,480)(0,-1)[``u^*]{440}{-1}l
 \putmorphism(900,480)(0,-1)[``u_!]{440}{1}r
 \end{picture}
\end{center}
as in 2. and 3.

Let $f:X'\ra X$ be a morphism in $\cX$ over $u:B'\ra B$ with a factorization
via prone $pr_{u,X}$ and supine $su_{u,X'}$ morphisms in $\cX$ as follows
\begin{center} \xext=1000 \yext=580
\begin{picture}(\xext,\yext)(\xoff,\yoff)
\setsqparms[1`1`1`1;850`500]
 \putsquare(80,0)[X'`u_!(X')`u^*(X)`X;su_{u,X'}`\hat{f}`\check{f}`pr_{u,X}]
  \put(150,430){\vector(2,-1){700}}
       \put(460,320){$f$}
 \end{picture}
\end{center}
The mutual correspondence between the functor $F$ and transformations  $\tau$ and
$\sigma$ can be read off from the following diagram
\begin{center} \xext=1400 \yext=1550
\begin{picture}(\xext,\yext)(\xoff,\yoff)
\setsqparms[1`1`1`1;1150`500]
\putsquare(80,500)[F(X')`F(u_!(X'))`F(u^*(X))`F(X);F(su_{u,X'})`F(\hat{f})`F(\check{f})`F(pr_{u,X})]
  \put(150,930){\vector(3,-1){1000}}
       \put(600,820){$F(f)$}

      \putmorphism(1230,1500)(0,-1)[u_!(F(X'))`\phantom{F(u_!(X'))}`\sigma_{X'}]{500}{1}r
      \put(150,1080){\vector(3,1){1000}}
       \put(350,1350){$su_{u,F(X')}$}

      \putmorphism(80,500)(0,-1)[\phantom{F(u^*(X))}`u^*(F(X))`\tau_X]{500}{1}l
      \put(150,70){\vector(3,1){1000}}
       \put(600,120){$pr_{u,F(X)}$}
 \end{picture}
\end{center}
whose part is just $F$ applied to the previous diagram, having in mind that $F$
restricted to the fibre $\cX_{B'}$ is $Q$ and  to the fibre $\cX_{B}$ is $P$.

To see the remaining part of the Lemma, we recall the direct relation between $\tau$ and $\sigma$.
From $\tau$ we get $\sigma$ as follows
\begin{center} \xext=2400 \yext=120
\begin{picture}(\xext,\yext)(\xoff,\yoff)
 \putmorphism(0,0)(1,0)[u_!Q`u_!Qu^*u_!`u_!Q(\eta^u)]{800}{1}a
 \putmorphism(800,0)(1,0)[\phantom{u_!Q(\eta^u)}`u_!u^*Pu_!`u_!(\tau_{u_!})]{800}{1}a
 \putmorphism(1600,0)(1,0)[\phantom{u_!u^*Pu_!}`Pu_!`\varepsilon^u_{Pu_!}]{800}{1}a
 \end{picture}
\end{center}
and we get back $\tau$ from $\sigma$  as follows
\begin{center} \xext=2400 \yext=120
\begin{picture}(\xext,\yext)(\xoff,\yoff)
 \putmorphism(0,0)(1,0)[Qu^*`u^*u_!Qu^*`\eta_{Qu^*}]{800}{1}a
 \putmorphism(800,0)(1,0)[\phantom{u^*u_!Qu^*}`u^*Pu_!u^*`u^*(\sigma_{u^*})]{800}{1}a
 \putmorphism(1600,0)(1,0)[\phantom{u^*Pu_!u^*}`u^*P`u^*P(\varepsilon^u)]{800}{1}a
 \end{picture}
\end{center}
where, as usual, $\eta^u$ and $\varepsilon^u$ are the unit and the counit of the adjunction $u_!\dashv u^*$. From this description
 it is easy to see that with the assumptions of the Lemma, $\tau$ is (weakly) cartesian if and only if $\sigma$ is. $~\Box$

\vskip 2mm

The second of the above five above descriptions of morphisms in $\cY^\cX$ seems to be the most convenient for us,
and  from now on we shall assume that the morphisms in $\cY^\cX$ are given in that form.
The composition in $\cY^\cX$ is defined as follows.
For morphisms  $\sigma : Q\ra P$ and $\tau : R\ra Q$  in $\cY^\cX$  over $u:B'\ra B$ and $v:B''\ra B'$, respectively, we have
\begin{center} \xext=1700 \yext=600
\begin{picture}(\xext,\yext)(\xoff,\yoff)

      \put(150,380){\vector(3,1){510}}
      \put(250,480){$\tau_{u^*}$}
      \put(1000,550){\vector(3,-1){510}}
      \put(1250,480){$v^*(\sigma)$}
       \put(700,550){$v^*Qu^*$}

      \putmorphism(0,300)(1,0)[Rv^*u^*`v^*u^*P`]{1600}{0}a
       \put(0,100){\vector(0,1){100}}
       \put(-120,110){$\cong$}
        \put(1600,200){\vector(0,-1){100}}
        \put(1650,110){$\cong$}
      \putmorphism(0,0)(1,0)[R(u\circ v)^*`(u\circ v)^*P`\sigma\circ\tau]{1600}{1}a

 \end{picture}
\end{center}
where unnamed isomorphisms come from canonical isomorphisms between functors $(u\circ v)^*$ and $u^*\circ v^*$.
The prone morphism over $u:B'\ra B$ with the codomain $P$ in $\cY^\cX_B$
\[ pr_{u,P} :  u^*Pu_! \lra P  \]
is the natural transformation in $\bCat (\cX_B,\cY_{B'})$ defined with the help of the counit $\varepsilon^u$
\[ u^*P(\varepsilon^u) : u^*Pu_!u^* \lra u^*P \]
Then, for any morphism $v:B''\ra B'$ in $\bB$ and any morphisms $\tau :Q\ra P$ in $\cY^\cX$  over $u\circ v$ we have a (unique!) morphism
$\hat{\tau} :Q\ra u^*Pu_!$ in $\cY^\cX$  over $v$ defined as a natural transformation
 \begin{center} \xext=1800 \yext=150
\begin{picture}(\xext,\yext)(\xoff,\yoff)
\putmorphism(0,0)(1,0)[Qv^*`\phantom{Qv^*u^*u_!}`Qv^*(\eta^u)]{1000}{1}a
\putmorphism(1000,0)(1,0)[Qv^*u^*u_!`v^*u^*Pu_!`\tau_{u_!}]{800}{1}a
 \end{picture}
\end{center}
so that $\tau=pr_{u,P}\circ \hat{\tau}$ in $\cY^\cX$, i.e. the triangle of natural transformations
 \begin{center} \xext=2000 \yext=600
\begin{picture}(\xext,\yext)(\xoff,\yoff)
\putmorphism(500,550)(3,-1)[Qv^*u^*\cong Q(uv)^*`v^*u^*P\cong (uv)^*P\;\;\;\;\;\;\;\;\;\;\;\;\;\;\;\;\;\;\;\;\;\;\;\;\;\;`]{1500}{1}r
\putmorphism(500,50)(1,0)[v^*u^*Pu_!u^*``v^*(pr_{u,P})]{800}{1}b
 \put(300,500){\vector(1,-3){120}}
    \put(230,250){\makebox(50,80){$\hat{\tau}_{u^*}$}}
     \put(1400,350){\makebox(50,80){$\tau$}}
 \end{picture}
\end{center}
commutes.

Similarly, the supine morphism over $u:B'\ra B$ with the domain $Q$ in $\cY^\cX_{B'}$
\[ su_{u,Q} :  Q\lra u_!Qu^*  \]
is the natural transformation in $\bCat (\cX_B,\cY_{B'})$ defined with the help of the unit $\eta^u$
\[ \eta^u_{ Qu^*} : Qu^*\lra u^*u_!Qu^* \]
Then, for any $w:B'\ra B''$ in $\bB$ and any morphism $\sigma:Q\ra P$ over $w\circ u$
we have a (unique) morphism $\check{\sigma}:u_!Qu^*\lra P$ in $\cY^\cX$ over $w$ defined as a natural transformation
 \begin{center} \xext=1800 \yext=150
\begin{picture}(\xext,\yext)(\xoff,\yoff)
\putmorphism(0,0)(1,0)[u_!Qu^*w^*`\phantom{u_!u^*w^*P}`u_!(\sigma)]{800}{1}a
\putmorphism(800,0)(1,0)[u_!u^*w^*P`w^*P`\varepsilon^u_{w^*P}]{1000}{1}a
 \end{picture}
\end{center}
so that the triangle
 \begin{center} \xext=2000 \yext=600
\begin{picture}(\xext,\yext)(\xoff,\yoff)
\putmorphism(150,550)(1,0)[Q(vu)^*\cong Qu^*w^*`u^*u_!Qu^*w^*`(su_{u,Q})_{w^*}]{1350}{1}b
\put(0,480){\vector(3,-1){1130}}
\putmorphism(1550,550)(1,-3)[``u^*(\check{\sigma})]{160}{1}r
\put(1500,0){\makebox(50,80){$(vu)^*P\cong u^*w^*P$}}
\put(560,200){\makebox(50,80){$\sigma$}}
 \end{picture}
\end{center}
commutes in $\cY^\cX$.

\begin{proposition}For any bifibration $p_{\cX}:\cX\ra \bB$ the exponential object
$p_{exp}:\cX^\cX\ra \bB$ in ${\bf Cat}/\bB$ is a bifibration and it
has the structure of a lax monoidal fibration, whose fibres are strict monoidal categories.\end{proposition}
{\it Proof.}~ The fact that $p_{exp}$ is a bifibration we have already seen.
We describe the monoidal structure in $p_{exp}:\cX^\cX\ra \bB$, and leave the reader to verify the axioms.
We have an obvious isomorphism of fibrations
\[ \bB\times_\bB\cX\lra \cX \]
whose exponential adjoint in $\bCat/\bB$
\[ I : \bB\lra \cX^\cX \]
is the unit for the tensor. Thus, for $B\in \bB$, the unit $I_B$ in fibre $\cX^\cX_B$
is the identity functor on the fibre $\cX_B$. The tensor functor
 \begin{center} \xext=1400 \yext=150
\begin{picture}(\xext,\yext)(\xoff,\yoff)
  \putmorphism(0,50)(1,0)[\cX^\cX\times_\bB\cX^\cX`\cX^\cX`\otimes]{1400}{1}a
\end{picture}
\end{center}
is the exponential adjoint in $\bCat/\bB$ to the morphism
 \begin{center} \xext=2400 \yext=150
\begin{picture}(\xext,\yext)(\xoff,\yoff)
  \putmorphism(0,50)(1,0)[\cX\times_\bB\cX^\cX\times_\bB\cX^\cX`\cX\times_\bB\cX^\cX`ev\times 1_{\cX^\cX}]{1400}{1}a
  \putmorphism(1400,50)(1,0)[\phantom{\cX\times_\bB\cX^\cX}`\cX` ev]{1000}{1}a
\end{picture}
\end{center}
The tensor on objects is the composition of functors.
As we will use it later, we describe explicitly the action of the tensor on morphisms.
Let $\sigma : P_1\ra P_0$ and $\tau : Q_1\ra Q_0$ be two morphisms in $p_{exp}:\cX^\cX\ra \bB$ over a morphism $u:B_1\lra B_0$,
 i.e. they are natural transformations $\sigma : P_1u^*\ra u^*P_0$ and $\tau : Q_1u^*\ra u^*Q_0$ in $\bCat(\cX(B_0),\cX(B_1))$.
 Then their tensor $\sigma\otimes_u\tau : P_1\otimes_{B_1} Q_1 = P_1\circ Q_1 \lra P_0\circ Q_0=P_0\otimes_{B_0} Q_0$ is defined from
 the commutative diagram below
\begin{center} \xext=2400 \yext=750
\begin{picture}(\xext,\yext)(\xoff,\yoff)
 \setsqparms[1`-1`1`1;2000`500]
  \putsquare(200,100)[P_1u^*u_!Q_1u^*`u^*P_0u_!u^*Q_0`P_1Q_1u^*`u^*P_0Q_0;\sigma_{u_!}\ast\tau=\sigma\ast u_!(\tau)`P_1(\eta^u_{Q_1u^*})`u^*P_0(\varepsilon^u_{Q_0})`\sigma\otimes_u\tau]
      \put(2800,-50){$~\Box$}
\end{picture}
\end{center}

The monoids in $(p_{exp}:\cX^\cX\ra \bB,\otimes,I)$ are monads over fibres of $p$.
A morphism of monoids $(f,u):(M',m',e')\lra (M,m,e)$ over $u:B'\ra B$
is a morphism of monads from  $(M,m,e)$ to $(M',m',e')$ whose functor part is
$u^*: \cX_B\lra \cX_{B'}$ and $f: M'u^*\lra u^*M$ is a natural transformation
satisfying the usual conditions (see Section \ref{burroni}).
 The evaluation morphism \[ ev : \cX^\cX \times \cX \lra \cX \] in $\bCat/\bB$
 is the action of the lax monoidal fibration $p_{exp}$ on the fibration $p_\cX$.
 The algebras for this action are all algebras for all monads in
\mbox{$Mon(p_{exp}:\cX^\cX\ra \bB,\otimes,I)$} taken together, i.e.
organized into a single category fibred over $\bB$ as well as
over \mbox{$Mon(p_{exp}:\cX^\cX\ra \bB,\otimes,I)$}.

We will apply this construction mainly to the basic fibration
\mbox{$cod:\cC^\ra\lra\cC$}  of a category $\cC$ with pullbacks (usually $\cC=Set$).
Such a fibration is always a cartesian bifibration.
We write \mbox{$p_{exp}:Exp(\cC)\lra \cC$} (or $p_{exp,\cC}$ if we want to indicate
the category $\cC$)  for \mbox{$p_{exp}:{\cC^\ra}^{\cC^\ra}\lra \cC$}.
The monoids in the exponential fibration \mbox{$p_{exp}:Exp(\cC)\lra \cC$} are monads in slices of $\cC$.
The morphism between two monoids over $u:c\ra c'\in \cC$ is a morphism of monads (in the opposite direction!)
whose functor part is the pullback functor $u^* : \cC/c' \lra \cC/c$.

Let $\cC art(\cC)$ denote the subcategory of $Exp(\cC)$ whose objects are pullback preserving
functors and cartesian natural transformations between them.
Moreover, let $w\cC art(\cC)$ denote the subcategory of $Exp(\cC)$ whose objects are
functors weakly preserving pullbacks and weakly cartesian natural transformations
between them. Restricting $p_{exp}$ to $\cC art(\cC)$ and $w\cC art(\cC)$ we get functors
$p_{ca,\cC}: \cC art(\cC)\lra \cC$ and $p_{wca,\cC}: w\cC art(\cC)\lra \cC$, respectively. We have

\begin{proposition}\label{subfib of exp}
The functors $p_{ca,\cC}$ and $p_{wca,\cC}$ described above are 
lax monoidal bifibrations with all that structure inherited from $p_{exp}:Exp(\cC)\ra \cC$.
In particular, the embeddings
\begin{center} \xext=1700 \yext=480
\begin{picture}(\xext,\yext)(\xoff,\yoff)
\putmorphism(0,400)(1,0)[ \cC art(C)`w\cC art(\cC)`]{800}{1}a
\putmorphism(800,400)(1,0)[\phantom{w\cC art(\cC)}`Exp(\cC)`]{800}{1}a

      \put(0,320){\vector(3,-1){700}}
      \put(1600,320){\vector(-3,-1){700}}
      \put(800,320){\vector(0,-1){200}}
      \put(780,0){$\cC$}

       \put(100,160){$p_{ca,\cC}$}
       \put(530,200){$p_{wca,\cC}$}
        \put(1300,160){$p_{exp}$}
 \end{picture}
\end{center}
are morphisms of lax monoidal fibrations and of bifibrations, that are faithful and full on isomorphisms.
The monoids  in $p_{ca,\cC}$ ($p_{wca,\cC}$) are (weakly) cartesian monads on slices of $\cC$.
\end{proposition}

{\it Proof.}~ To see that $p_{ca,\cC}$ is a fibration one has to notice that the prone morphism over
$u:c\ra c'$ with the codomain $P:\cC_{/c'}\ra \cC_{/c'}$, being  a pullback preserving functor, is
a cartesian natural transformation. Moreover to see that a factorization via a prone morphism of a
morphism $\tau$ in $Exp(\cC)$, being a  cartesian natural transformation between pullback preserving
functors, is also a cartesian natural transformation. All this follows directly from the explicit formulas
given above for the prone morphisms and the factorization via a prone morphism in $p_{exp}:\cY^\cX\lra \bB$
and the fact that both indexing and reindexing functors in $p_{exp}$ preserve pullbacks.
All the above remains true if we replace pullbacks by weak pullbacks. Thus both $p_{ca,\cC}$ and $p_{wca,\cC}$ are
fibrations and subfibrations of $p_{exp}$. The argument that these functors are subopfibrations of  $p_{exp}$ is similar.

Clearly, the composition of functors (weakly) preserving pullbacks (weakly) preserves pullbacks.
From this it is easy to see that the whole lax monoidal structure of $p_{ca,\cC}$ and $p_{wca,\cC}$ is inherited from $p_{exp}$.

The remaining part of the proposition is obvious.
$~\Box$

\vskip 2mm

{\bf Remark} There are many more interesting subfibrations of $p_{exp}:Exp(\cC)\ra \cC$.
The fibrations $p_{ca,\cC}$ and $p_{wca,\cC}$ have also their 'wide pullback versions'.
If slices of $\cC$ are sufficiently cocomplete (e.g if $\cC$ is $Set$) then finitary or even accessible
functors form full subfibrations of $p_{exp}$. However, the functors preserving finite limits (or just the terminal object)
do not constitute a subfibration of $p_{exp}$, as the functor $u_!:\cC_{/c}\lra \cC_{/c'}$ does not preserve the terminal
object, in general.

\vskip 2mm

The following result is the main reason we consider exponential fibrations. It will be used later many times.

\begin{proposition}
Let $(p_\cE:\cE\ra \bB,I,\otimes, \alpha,\lambda,\varrho)$ be a lax monoidal fibration and
$p_\cX:\cX\ra \bB$ be a bifibration. Then the bijective
correspondence given by the exponential adjunction in $\bCat /\bB$ between morphisms
   \begin{center} \xext=800 \yext=450
\begin{picture}(\xext,\yext)(\xoff,\yoff)
 \settriparms[1`1`1;400]
 \putVtriangle(0,0)[\cE\times_\bB\cX`\cX`\bB;\star``p_\cX]
\end{picture}
\end{center}
and morphisms
   \begin{center} \xext=800 \yext=450
\begin{picture}(\xext,\yext)(\xoff,\yoff)
 \settriparms[1`1`1;400]
 \putVtriangle(0,0)[\cE`\cX^\cX`\bB;\check{\star}`p_\cE`p_{exp}]
\end{picture}
\end{center}
induces a bijective correspondence between actions $(\star,\psi_0,\psi_2)$ of $p_\cE$ on $p_\cX$
and morphisms of lax monoidal fibrations $(\check{\star},\phi_0,\phi_2)$ from  $p_\cE$ to  $p_{exp}$.
The correspondence relates the coherence natural transformations as follows. The transformation
\[ \phi_0 :I\lra \check{\star}\circ I \]
for $O\in \bB$, is a natural transformation of functors
$(\phi_0)_O : id_{\cX_O}\lra I_O\star (-)$, i.e. for $X\in \cX_O$
we have
\[ ((\phi_0)_O)_X = (\psi_0)_{O,X} : X \lra I_O\star X \]
Moreover, for $A,B\in \cE_O$, we have
\begin{center} \xext=1200 \yext=450
\begin{picture}(\xext,\yext)(\xoff,\yoff)
 \setsqparms[1`0`0`1;1200`300]
  \putsquare(0,50)[\check{\star}(A)\circ \check{\star}(B)`\check{\star}(A\otimes B)`A\star (B(-))`(A\otimes B)\star(-);(\phi_2)_{A,B}```(\psi_2)_{A,B,-}]
  \put(0,160){\makebox(50,50){$\parallel$}}
  \put(1100,160){\makebox(50,50){$\parallel$}}
\end{picture}
\end{center}
i.e. $X\in \cX_O$,
\[ ((\phi_2)_{A,B})_X=(\psi_2)_{A,B,X} : A\star(B\star X) \lra (A\otimes B)\star X \]
This correspondence is natural in $\cE$.
\end{proposition}

{\it Proof.}~ Exercise. $~\Box$

From the above Proposition follows that if we have an action of the lax monoidal fibration $p_\cE:\cE\ra \bB$ on a bifibration $p_\cX:\cX\ra \bB$,
we get a morphism from the lax monoidal fibration $p_\cE$ into a lax monoidal fibration whose fibres are strict monoidal categories. If
$\cX$ is sufficiently concrete (like $Set^\ra$) and this morphism is an embedding we can view this kind of phenomena as representation theorems.
We represent objects of $\cE$ as endofunctors of fibres of $p_\cX:\cX\ra \bB$, and monoids in $(p_\cE:\cE\ra \bB,\otimes,I)$
as monads over fibres of $p_\cX:\cX\ra \bB$. Similar things can be said about morphisms. We will see many examples
of such representations later.

\subsection{The exponential fibrations  in $Fib/\bB$ } $Fib/\bB$ is a cartesian closed category.
For any fibration $p_\cX :\cX\ra \bB$ the exponential fibration $p_{fiexp} : [\cX \Ra\cX]\ra \bB$
is lax monoidal with tensor being (again) the internal composition.
Monoids in $p_{fexp}$ are compatible families of monads and cartesian morphisms between them.

Having a strong action $\star : \cE\times_\bB \cX \lra \cX$ we could also represent $\cE$ in $[\cX \Ra\cX]$.
But strong actions are less common, and such (non-trivial) representations are more difficult to achieve in practice.
This is why we are not going consider this kind of exponential fibrations in the following.

\section{The Burroni fibrations and opetopic sets}\label{burroni}

\subsection{The Burroni fibrations and $T$-categories}

Let $\cC$ be a category with pullbacks,
$\lk T,\eta,\mu \rk$ a monad on $\cC$. The category
$Gph(T)$ is the category of $T$-graphs. An object $\lk A,O,\gamma,\delta\rk $ of $Gph(T)$ is
a span
\begin{center} \xext=800 \yext=430
\begin{picture}(\xext,\yext)(\xoff,\yoff)
\settriparms[1`1`0;400]
  \putAtriangle(0,0)[A`O`T(O);\gamma`\delta`]
\end{picture}
\end{center}
in $\cC$. The morphisms $\gamma$ and $\delta$ are called {\em codomains} and {\em domains} of the $T$-graph
$\lk A,O,\gamma,\delta\rk$, respectively.
Sometimes we write $A$ instead of $\lk A,O,\gamma,\delta\rk$,  for short,  when it does not lead to a confusion.

A morphism of $T$-graphs $\lk f,u\rk :\lk A,O,\gamma,\delta\rk \lra \lk A',O',\gamma',\delta'\rk $
is a pair of morphisms $f:A\ra A'$ and $u:O\ra O'$ in $\cC$ making
the squares
\begin{center} \xext=1900 \yext=550
\begin{picture}(\xext,\yext)(\xoff,\yoff)
 \setsqparms[1`1`1`1;500`400]
  \putsquare(0,50)[A`A'`O`O';f`\gamma`\gamma'`u]
  \setsqparms[1`1`1`1;600`400]
  \putsquare(1300,50)[A`A'`T(O)`T(O');
  f`\delta`\delta'`T(u)]
\end{picture}
\end{center}
commute. Let $Gph(T)$ denotes the category of $T$-graphs and $T$-graph morphisms. We have a projection functor
\[ p_T:Gph(T)\lra \cC \]
 sending the morphism $\lk f,u\rk :\lk A,O,\gamma,\delta\rk \lra \lk A',O',\gamma',\delta'\rk $ to the morphism $u:O\ra O'$
which is easily seen to be a fibration, cf. \cite{B} p. 235. The lax monoidal structure in $p_T$ is defined as follows.
Let $\lk A, O, \gamma_A,\delta_A\rk$ and $\lk B, O,
\gamma_B,\delta_B\rk$ be two objects in the fibre over $O$, i.e.
in $Gph(T)_O$. Then the tensor
\[ \lk A, O,\gamma_A,\delta_A\rk\otimes_O\lk B, O, \gamma_B,\delta_B\rk =
\lk A\otimes B, O, \gamma_\otimes, \delta_\otimes \rk\]
 is defined from the following diagram
\begin{center} \xext=1800 \yext=1230
\begin{picture}(\xext,\yext)(\xoff,\yoff)
\settriparms[1`1`0;400]
  \putAtriangle(0,400)[A`O`T(O);\gamma_A``]
    \put(500,530){$\delta_A$}
  \putAtriangle(800,400)[T(B)`\phantom{T(O)}`T^2(O);`T(\delta_B)`]
  \put(1020,530){$T(\gamma_B)$}
  \settriparms[0`1`0;400]
  \putAtriangle(1200,0)[\phantom{T^2(O)}``T(O);`\mu_O`]
  \settriparms[1`1`0;400]
   \putAtriangle(400,800)[A\otimes B`\phantom{A}`\phantom{T(B)};\pi_1`\pi_2`]
\end{picture}
\end{center}
in which the square is a pullback and
\[\gamma_\otimes =  \gamma_A\circ\pi_1,\;\;\; \delta_\otimes = \mu_O\circ T(\delta_B)\circ\pi_2.\]
The unit in the fibre over $O$ is
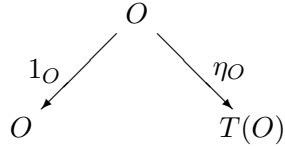
\begin{center} \xext=800 \yext=430
\begin{picture}(\xext,\yext)(\xoff,\yoff)
\settriparms[1`1`0;400]
  \putAtriangle(0,0)[O`O`T(O);1_O`\eta_O`]
\end{picture}
\end{center}

The coherence morphisms are defined using the universal properties of pullbacks.
For an object $\lk A, O, \gamma_A,\delta_A\rk$ on the fibre over $O$
the left unit morphism is
\[ \lambda_{\lk A, O, \gamma_A,\delta_A\rk}=\lk 1_O,\eta_A)  , 1_O \rk : \lk A, O, \gamma_A,\delta_A\rk \lra  \lk O\otimes A, O, \gamma_\otimes,\delta_\otimes\rk \]
the right unit morphism is
\[ \varrho_{\lk A, O, \gamma_A,\delta_A\rk}=  \lk 1_A  , 1_O \rk :
\lk A\otimes O, O, \gamma_\otimes,\delta_\otimes\rk \lra\lk A, O, \gamma_A,\delta_A\rk \]
(the right unit morphism is always an isomorphism and in fact we can assume that it is an identity as $A\otimes O$
is a pullback of $\delta_A$ along the identity).  The associativity morphism
\[ \alpha_{A,B,C}: A \otimes(B \otimes C) \lra (A \otimes B)\otimes C  \]
 is also defined similarly, using universal properties of pullbacks. We leave the details to the reader.

\begin{proposition}Let $\cC$ be a category with pullbacks.
The functor $p_T: Gph(T)\ra \cC$ is a bifibration and together with the
monoidal structure $(\otimes, I, \alpha,\lambda,\varrho)$ described above is a lax monoidal fibration.
The total category of the fibration $q_T:Mon(T)\lra \cC$ of monoids in $(Gph(T),p_T, \otimes, I, \alpha,\lambda,\varrho)$
is equivalent to the category of $T$-categories of Burroni.
If moreover,  the monad $(T,\eta,\mu)$ is cartesian then the fibres of $p_T$ are (strong) monoidal categories, i.e.
the coherence morphisms, $\lambda$ and $\alpha$, are isomorphisms.
\end{proposition}

{\it Proof.}~ A simple tedious check. $~\Box$

\vskip 2mm
{\bf Remark}
If the monad $(T,\eta,\mu)$ is cartesian then the fibres of $p_T: Gph(T)\ra \cC$ are strong monoidal categories
but the reindexing functors are still only lax monoidal. This is already so for the identity monad $(1_\cC, 1_{1_\cC}, 1_{1_\cC})$ on $\cC$.
The category $Mon(1_\cC)$ is the category of internal categories in $\cC$.

\subsection{Tautologous actions of Burroni fibrations}\label{taut-act-Burr}

If $(T,\eta,\mu)$ is a monad on a category $\cC$ with pullbacks then the lax monoidal fibration
$p_T:Gph(T)\lra \cC$ has a natural action on the basic fibration $cod:\cC^\ra \lra \cC$. The functor part
\begin{center} \xext=800 \yext=520
\begin{picture}(\xext,\yext)(\xoff,\yoff)
 \settriparms[1`1`1;400]
  \putVtriangle(0,0)[Gph(T)\times_\cC\cC^\ra`\cC^\ra`\cC;\star^T``cod]
\end{picture}
\end{center}
is defined on objects by
\begin{center} \xext=1800 \yext=430
\begin{picture}(\xext,\yext)(\xoff,\yoff)
\settriparms[1`1`0;400]
  \putAtriangle(0,0)[A`O`T(O);\gamma`\delta`]
  \putmorphism(1100,400)(0,-1)[X`O`d]{400}{1}r
    \putmorphism(1800,400)(0,-1)[A\star^T X`O`]{400}{1}r
\put(1350,250){\vector(1,0){250}}
\put(1350,200){\line(0,1){100}}
\end{picture}
\end{center}
where the right vertical arrow in the above diagram is the composite of the upper horizontal arrows in the following diagram
\begin{center} \xext=1200 \yext=600
\begin{picture}(\xext,\yext)(\xoff,\yoff)
 \setsqparms[-1`1`1`-1;700`400]
  \putsquare(500,50)[A`A\star^T X`T(O)`T(X);`\delta``T(d)]
  \putmorphism(0,450)(1,0)[O`\phantom{A}`\gamma]{500}{-1}a
\end{picture}
\end{center}
in which the square is a pullback. By adjunction, we get a morphism of lax monoidal fibrations
\begin{center} \xext=800 \yext=520
\begin{picture}(\xext,\yext)(\xoff,\yoff)
 \settriparms[1`1`1;400]
  \putVtriangle(0,30)[Gph(T)`Exp(\cC)`\cC;rep_T`p_T`p_{exp}]
\end{picture}
\end{center}
that represents $T$-graphs as endofunctors on slices of $\cC$. Under this representation the $T$-categories
correspond to (some) monads on  slices of $\cC$.

\vskip 2mm
{\bf Example} The actions $\star^T$ of the lax monoidal fibration defined above do not preserve prone
morphisms, in general. Even if $T$ is the free monoid monad on $Set$, the action $\star^T$ is only
a lax morphism of fibrations. To see this, consider a morphism $u:[1]\ra [0]$, from two element set
$[1]=\{ 0, 1 \}$ to one element set $[0]=\{ 0 \}$. Let $(A,\gamma,\delta)$ be an object in $Gph(T)$
over $[1]$, such that $A=  \{ a \}$, and $\partial(a)=00$, the word of length $2$ of zero's,
$\gamma(a)=0$. The identity $1_{[0]}$ on $[0]$ is a morphism in $Set_{/[0]}$.
Then the object $(A,\gamma,\delta)\star ([0],1_{[0]})$  (in the basic fibration over $Set$) has one element in the domain, operation $a$ with
inputs and outputs in $[0]$, i.e. $\lk a; 0,00\rk$.
The domain of the prone morphism $u^*((A,\gamma,\delta)\star ([0],1_{[0]})) \lra (A,\gamma,\delta)\star ([0],1_{[0]})$ over $u$ has two elements,
in the domain, the operation $a$ with the output either $0$ or $1$ and inputs as before, i.e.  $\{ \lk a; 0,00\rk , \lk a; 1,00\rk \}$.
On the other hand, the image under $\star^T$ of the prone morphisms over $u$, whose codomains are $(A,\gamma,\delta)$ and $1_{[0]}$, respectively,
$pr_{u,A}\star^T pr_{u,1_{[0]}} : u^*(A)\star^T u^*([0])\lra A\star^T [0]$ has in the domain of its domain eight elements, i.e. the operation $a$
with both inputs and outputs either $0$ or $1$, i.e.  $\{ \lk a; 0,00\rk , \lk a; 1,00\rk, \ldots, \lk a; 1,10\rk \lk a; 1,11\rk \}$. Thus the domains
of those morphisms are not isomorphic and hence the prone morphisms are not preserved.
\vskip 2mm
As the morphism $rep_T$ is the exponential transpose of $\star^T$ in $CAT_{/\cC}$ (not in $Fib(\cC)$) one of these morphisms can be
a morphism of fibrations even if the other one is not. We have

\begin{proposition}\label{T-cart-repT}Let $(T,\eta,\mu)$ be a cartesian monad on a category with  pullbacks $\cC$.
Then the functor $rep_T$ defined above is a strong morphism of lax monoidal fibrations and of bifibrations.
The image of $rep_T$ is in $p_{ca} : \cC art(\cC)\ra \cC$.
\end{proposition}

{\it Proof.}~First, we describe the functor $rep_T$ in details. For an object $A=(A,\gamma,\delta)$ in $Gph(T)_O$, we have a functor
\[ rep_T(A)= A\star^T(-) : \cC_{/O}\lra \cC_{/O} \]
In the following, we omit the superscript $T$. For $u:O\ra Q$ in $\cC$ and $(h,u):A\ra B$, a morphism in $Gph(T)$ over $u$, we have
a natural transformation in $\bCat(\cC_{/Q},\cC_{/O})$
\[ rep_T(h,u): A\star u^*(-)\lra u^*(B\star(-)) \]
so that for $d^Y:Y\ra Q$ in $\cC$, the value $rep_T(h,u)_Y$ is defined from the following diagram
  \begin{center}
\xext=2600 \yext=1780
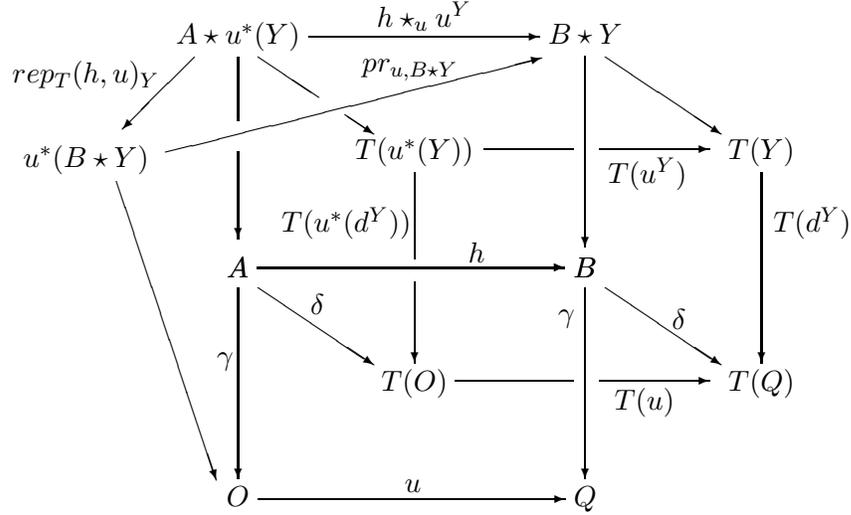
\begin{picture}(\xext,\yext)(\xoff,\yoff)
 \setsqparms[1`1`1`1;1200`800]
 \putsquare(550,50)[A`B`O`Q; ```]
 \setsqparms[1`0`1`1;1200`800]
 \putsquare(550,850)[A\star u^*(Y)`B\star Y`A`B;h\star_uu^Y```]
 \setsqparms[0`0`1`0;1200`800]
 \putsquare(1160,450)[T(u^*(Y))`T(Y)`T(O)`T(Q); ```]

 \put(0,1200){\makebox(50,50){$u^*(B\star Y)$}}

  \put(300,1250){\vector(4,1){1300}}
  \put(620,780){\vector(3,-2){400}}
   \put(620,1580){\line(3,-2){210}}
   \put(920,1380){\vector(3,-2){100}}
  \put(1820,780){\vector(3,-2){400}}
  \put(1820,1580){\vector(3,-2){400}}

    \put(400,1580){\vector(-1,-1){250}}

    \put(130,1150){\vector(1,-3){345}}

   \put(550,1250){\vector(0,-1){300}}
   \put(550,1580){\line(0,-1){220}}

 \put(1160,800){\vector(0,-1){280}}
   \put(1160,1180){\line(0,-1){300}}

  \put(1800,460){\vector(1,0){380}}
   \put(1300,460){\line(1,0){410}}
   \put(1800,1260){\vector(1,0){380}}
   \put(1400,1260){\line(1,0){310}}

    \put(1350,880){\makebox(50,50){$h$}}
    \put(1130,70){\makebox(50,50){$u$}}
    \put(1930,360){\makebox(50,50){$T(u)$}}
    \put(1940,1150){\makebox(50,50){$T(u^Y)$}}

    \put(400,980){\makebox(1050,50){$T(u^*(d^Y))$}}
    \put(2510,980){\makebox(50,50){$T(d^Y)$}}


     \put(480,500){\makebox(50,50){$\gamma$}}
      \put(800,700){\makebox(50,50){$\delta$}}
     \put(1660,650){\makebox(50,50){$\gamma$}}
     \put(2050,650){\makebox(50,50){$\delta$}}

\put(0,1490){\makebox(50,50){$rep_T(h,u)_Y$}}
\put(1120,1510){\makebox(50,50){$pr_{u,B\star Y}$}}
 \end{picture}
\end{center}
where
 \begin{center} \xext=500 \yext=430
\begin{picture}(\xext,\yext)(\xoff,\yoff)
\setsqparms[1`1`1`1;500`400]
\putsquare(0,0)[u^*(Y)`Y`O`Q;u^Y`u^*(d^Y)`d^Y`u]
\end{picture}
\end{center}
is a pullback. As three sides of the cube are pullbacks ($T$ preserves pullbacks), so is the front square.
In particular, if $h$ is an iso, so is $h\star_uu^Y$, for any $d^Y:Y\ra Q$ in $\cC$.

Note that the morphisms $pr_{u,B\star Y}: u^*(B\star Y)\ra (B\star Y$ and $u^Y,u):u^*(Y)\ra Y$, in the above two diagrams,
are prone morphisms (over $u$) in the basic fibration over $\cC$.
In the following, we will deal with prone and supine morphisms in two other fibrations $p_T:Gph(T)\ra \cC$ and $p_{exp}:Exp(\cC)\ra \cC$.
Thus in total, we have three different sorts of prone morphisms.

The codomain of a supine morphism over $u:O\ra Q$ whose domain is $A=(A,\gamma,\delta)$ is $u_!(A)=(A,u\circ\gamma,T(u)\circ\delta)$.
The supine morphism is $su_{u,A}=(1_A,u):A\ra u_!(A)$. We have a diagram in $Gph(T)$ over $u$:
 \begin{center} \xext=2800 \yext=1000
\begin{picture}(\xext,\yext)(\xoff,\yoff)
\setsqparms[1`-1`1`1;1200`700]
\putsquare(1600,200)[u^*u_!(A\star u^*(Y))`u_!(A\star u^*(Y))`A\star u^*(Y)`u_!(A)\star Y;``\xi^u_{A,Y}`1_A\star_u u^Y]

\putmorphism(0,200)(1,0)[u^*(u_!(A)\star Y)`\phantom{A\star u^*(Y)}`rep_T(1_A,u)_Y]{1600}{-1}a
\put(940,500){\makebox(200,100){$(su_{u,rep_T(A)})_Y =\eta^u_{A\star u^*(Y)}$}}

 \put(1700,300){\vector(2,1){1000}}
 \put(2000,650){\makebox(50,50){$su_u=1_{A\star u^*(Y)}$}}

 \put(0,600){\vector(0,-1){300}}
   \put(1200,900){\line(-4,-1){1200}}
   \put(350,800){\makebox(250,80){$u^*(\xi^u_{A,Y})$}}

   \put(0,0){\line(0,1){100}}
   \put(0,0){\line(1,0){2800}}
    \put(2800,0){\vector(0,1){100}}
    \put(1000,20){\makebox(250,80){$pr_{u,u_!(A)\star Y}$}}
 \end{picture}
\end{center}
The morphism $\xi^u_{A,Y}$ is the second part of the factorization
of $1_A\star_uu^Y$ via  a supine morphism.  Note that the morphisms
$1_A\star_uu^Y$ and $\xi^u_{A,Y}$ considered as morphisms in $\cC$
are equal but the first is a part of a morphism in $Gph(T)$ over
$u$, and the second is in the fibre over $Q$. By a remark below
the previous diagram, $1_A\star_uu^Y$ is an isomorphism and hence,
so are $\xi^u_{A,Y}$ and $u^*(\xi^u_{A,Y})$, as well.
One can verify that the left hand triangle commutes, as it is a triangle in the
fibre over $O$ in the basic fibration over $\cC$ and commutes when composed
with the prone morphism $pr_{u,u_!(A)\star Y}$.

Thus $rep_T$ preserves the supine morphisms.

One can verify that the prone morphism in $Gph(T)$ over $u:O\ra Q$ with codomain $B$ is $(u^B,u):u^*(B)\ra B$ where
\begin{center} \xext=2000 \yext=650
\begin{picture}(\xext,\yext)(\xoff,\yoff)
\settriparms[1`1`0;400]
  \putAtriangle(0,200)[u^*(B)`O`T(O);\gamma`\delta`]

  \putAtriangle(1200,200)[B`Q`T(Q);\gamma`\delta`]
  \putmorphism(400,600)(1,0)[\phantom{u^*(B)}`\phantom{B}`u^Y]{1200}{1}a

    \put(0,0){\line(0,1){100}}
   \put(0,0){\line(1,0){1200}}
    \put(1200,0){\vector(0,1){100}}
    \put(400,20){\makebox(50,80){$u$}}

     \put(800,30){\line(0,1){100}}
   \put(800,30){\line(1,0){350}}
   \put(1250,30){\line(1,0){750}}
    \put(2000,30){\vector(0,1){100}}
    \put(1600,50){\makebox(50,80){$T(u)$}}

\end{picture}
\end{center}
is a limiting cone. For $d^Y:Y\ra Q$ in $\cC$ with
 \begin{center} \xext=800 \yext=500
\begin{picture}(\xext,\yext)(\xoff,\yoff)
\setsqparms[0`0`1`1;800`400]
\putsquare(0,50)[`u_! u^*(Y)`u^*(Y)`Y;``\varepsilon^u_Y`u^Y]

 \put(200,150){\vector(2,1){500}}
 \put(300,330){\makebox(50,50){$pr_{u,Y}=1$}}
 \end{picture}
\end{center}
(i.e. $u^Y=\varepsilon^u_Y$) we can form a diagram
 \begin{center} \xext=2800 \yext=1050
\begin{picture}(\xext,\yext)(\xoff,\yoff)
\setsqparms[1`1`1`1;1200`600]
\putsquare(1600,50)[u^*(B)\star u^*(Y)`B\star u_!u^*(Y)`u^*(B\star Y)`B\star Y;u^B\star 1_{u^*(Y)}`rep_T(u^B,u)_Y`B\star \varepsilon^u_Y`pr_{u,B\star Y}]

\putmorphism(0,650)(1,0)[u^*(B\star u_!u^*(Y))`\phantom{u^*(B)\star u^*(Y)}`rep_T(1u^Y,u)_Y]{1600}{-1}a

    \put(0,380){\makebox(50,80){$(pr_{u,rep_T(B)})_Y$}}
     \put(0,580){\vector(3,-1){1400}}
   \put(0,750){\line(0,1){150}}
   \put(0,900){\line(1,0){2800}}
    \put(2800,900){\vector(0,-1){150}}
    \put(1200,920){\makebox(250,80){$pr_{u,B\star u_!u^*(Y)}$}}
 \end{picture}
\end{center}
($(pr_{u,rep_T(B)})_Y=u^*(B\star \varepsilon^u_Y)$)  in which one can verify, using properties of pullbacks,
 that $u^B\star_u1_{u^*(Y)}$ is prone in the basic fibration over $\cC$.
Thus $\zeta^u_{B,Y}$, the first part of the factorization of $u^B\star_u1_{u^*(Y)}$ via a prone morphism, is an iso.
As the left triangle commutes, $rep_T$ preserves prone morphisms, as well.
$~\Box$

\vskip 2mm

{\bf Remark} From the above proof, it follows that for any monad $T$ on a category $\cC$ with pullbacks, we have morphisms
\begin{center} \xext=1000 \yext=420
\begin{picture}(\xext,\yext)(\xoff,\yoff)
\putmorphism(0,300)(1,0)[u_!(A\star u^*(Y))`u_!(A)\star Y`\xi^u_{A,Y}]{1000}{1}a
\putmorphism(0,20)(1,0)[u^*(B)\star u^*(Y)`u^*(B\star u_!u^*(Y))`\zeta^u_{B,Y}]{1000}{1}a
 \end{picture}
\end{center}
natural in $A$, $B$ and $Y$, that are isomorphisms if the monad $T$ is cartesian.
Note that these isomorphisms express a kind of Beck-Chevalley
condition for actions of lax monoidal fibrations.

\vskip 2mm

{\bf Example.}  If $T$ is the identity monad on a category with pullbacks then $rep_T$ sends the internal category
$\bC=(C_1,C_0,m,i,d,c)$ in $\cC$ to a monad $rep_T(\bC)$ on the slice category $\cC_{/C_0}$ whose algebras are
internal presheaves on $\bC$.

\subsection{Multisorted signatures vs monotone polynomial diagrams}

In this section we shall examine the considerations from the previous section on a specific example
of the free monoid monad $T$ on the category $Set$. Note that $Gph(T)$ can be thought of as a category
of multisorted signatures.  An object $\lk A,O,\gamma,\delta\rk $ of $Gph(T)$ can be seen as a set
of operations $A$ a set of types $O$,  functions $\gamma$ and $\delta$ associating to operations
in $A$ their types of codomains in $O$ and their lists of types of their domains in $T(O)$.
To emphasize this, we shall denote the fibration $p_T:Gph(T)\ra Set$ for this particular monad
$T$ as $p_m:Sig_m\ra Set$. As we already mentioned, cf. \cite{B}, the category of monoids
in $p_m$ is equivalent to the category of Lambek's multicategories.
The action of $p_m:Sig_m\ra Set$ on $cod:Set^\ra\lra Set$ is as defined above.
Thus, by adjunction, we have a representation morphism
\begin{center} \xext=800 \yext=520
\begin{picture}(\xext,\yext)(\xoff,\yoff)
 \settriparms[1`1`1;400]
  \putVtriangle(0,30)[Sig_m`Exp(Set)`Set;rep_m`p_m`p_{exp}]
\end{picture}
\end{center}
We shall describe the image of this representation in a different way.
A {\em monotone polynomial diagram}\footnote{The name is so chosen to indicated the
 obvious relation with the notion of a polynomial diagram that will be considered in the next section.}
 over the set $O$ is a diagram of the following form
\begin{center} \xext=1500 \yext=200
\begin{picture}(\xext,\yext)(\xoff,\yoff)
  \putmorphism(0,50)(1,0)[O`E`s]{500}{-1}a
  \putmorphism(500,50)(1,0)[\phantom{E}`B`p]{500}{1}a
  \putmorphism(1000,50)(1,0)[\phantom{B}`O`t]{500}{1}a
\end{picture}
\end{center}
of sets and functions, moreover the fibres of the morphism $p$
are finite and linearly ordered. We write $(t,p,s)$ to denote such
a diagram. A {\em morphism of monotone diagrams}  $(f,g,u):
(t,p,s)\lra (t',p',s')$ over a function $u:O\ra O'$ is a triple of
functions with $g:E\ra E'$ and $f:B\ra B'$ so that the diagram
\begin{center} \xext=1800 \yext=650
\begin{picture}(\xext,\yext)(\xoff,\yoff)
 \setsqparms[-1`1`1`-1;600`500]
\putsquare(0,50)[O`E`O'`E';s`u`g`s']
 \setsqparms[1`0`1`1;600`500]
\putsquare(600,50)[\phantom{E}`B`\phantom{E}`B';p``f`p']
\putsquare(1200,50)[\phantom{B}`O`\phantom{B'}`O';t``u`t']
\end{picture}
\end{center}
 commutes, and the middle square is a pullback  in the category of posets
(i.e. $g_{\lceil p^{-1}(b)} :p^{-1}(b)\ra p'^{-1}(f(b))$ is an order isomorphism, for $b\in B$).
We compose morphisms of monotone polynomial diagrams in the obvious way, by placing one on top of the other.
In this way we
defined the category $\cM\cP oly\cD iag$ of monotone polynomial diagrams.
The category $\cM\cP oly\cD iag$ is fibred over $Set$,
where the projection functor
\[ p_{mpd} : \cM\cP oly\cD iag \lra Set \]
is given by
\[   (f,g,u): (t,p,s)\lra (t',p',s') \;\;\; \mapsto \;\;\;  u: O\lra O' \]
This is a lax monoidal fibration\footnote{The definition of the lax monoidal structure is left to be defined by the reader.
It is close to the structure on the fibration of polynomial functors defined in Section \ref{am-sig}.}
which also acts on the basic fibration $cod:Set^\ra\lra Set$.
The action
\begin{center} \xext=800 \yext=520
\begin{picture}(\xext,\yext)(\xoff,\yoff)
 \settriparms[1`1`1;400]
  \putVtriangle(0,0)[\cM\cP oly\cD iag\times_{Set}{Set}^\ra`{Set}^\ra`{Set};\star``cod]
\end{picture}
\end{center}
is given by the well known formula defining polynomial functors (see Section \ref{am-sig}), i.e. for $(t,p,s)$ in $\cM\cP oly\cD iag_O$,
and $d^X: X\ra O$ a function, we have
\[  (t,p,s)\star d^X = t_!p_*s^*(d^X) \]
Thus, by adjointness, we have a morphism of lax monoidal fibrations
\begin{center} \xext=800 \yext=520
\begin{picture}(\xext,\yext)(\xoff,\yoff)
 \settriparms[1`1`1;400]
  \putVtriangle(0,30)[\cM\cP oly\cD iag`Exp(Set)`Set;rep_{mpd}`p_{mpd}`p_{exp}]
\end{picture}
\end{center}
The class of functors in the image $rep_{mpd}$ coincides with the class of finitary polynomial endofunctors. However
the linear structure in the fibres of monotone polynomial diagrams restricts the class of natural transformations between them.
For every polynomial transformation $\tau:P\ra Q$ between  polynomial
endofunctors on $Set^I$, there is a monotone morphisms between monotone diagrams $(f,g,u): (t,p,s)\lra (t',p',s')$ so that
$rep_{mpd}(f,g,u)$ is isomorphic to $\tau$ (just order fibres in the polynomial diagrams defining $P$ and $Q$ in a compatible way).
This observation  says that the essential image of $rep_{mpd}$ consists of polynomial functors and
polynomial natural transformations, see Section \ref{am-sig}.  However this is not saying that the monotone polynomial monads on polynomial functors are the same as
polynomial monads. For a monad $(T,\eta, \mu)$ on a polynomial functor to be linear means\footnote{Here by a monotone monad we mean
a monad that is an image of a monoid in $p_{mpd} : \cM\cP oly\cD iag \lra Set$.}, that we can find one ordering
of the fibres of the polynomial diagram defining $T$ so that both morphisms $\eta:1_\cC \ra T$ and $\mu :T^2 \ra T$ are defined
by the morphisms of diagrams respecting these orderings (the order of
the fibres of the diagram defining $T^2$ is determined by the order of the diagram defining $T$). As we shall see later, this might be not possible.
We note for the record

\begin{proposition}\label{rep-monotone} 
The representations $rep_m$ and $rep_{mpd}$ are faithful and they are equivalent
as morphisms of lax monoidal fibrations into  $Exp(Set)\lra Set$.
As a consequence, $rep_{mpd}$ is a morphism of bifibrations and the category of Lambek's multicategories
is equivalent to the category of monoids in
$\cM\cP oly\cD iag$. Moreover, the monads in the image of the morphism of fibrations of monoids
\begin{center} \xext=900 \yext=490
\begin{picture}(\xext,\yext)(\xoff,\yoff)
 \settriparms[1`1`1;450]
  \putVtriangle(0,0)[Mon(\cM\cP oly\cD iag)`Mon(Exp(Set))`Set;``]
\end{picture}
\end{center}
induced by $rep_{mpd}$ are exactly monotone monads on polynomial functors.$~\Box$
\end{proposition}

{\bf Remark} A bad thing about the representations  $rep_m$ and $rep_{mpd}$ is that they are not full,  even on isomorphism.
As a consequence the monotone monads do not determine the monotone diagrams defining them uniquely (up to isomorphism).
The lack of fullness on isomorphisms is due to the fact that fibres in monotone diagrams
are linearly ordered.  As we shall see in the next two sections, similar representations of
both (finitary) polynomial (endo)functors and (finitary multivariable) analytic (endo)functors are full on isomorphisms.

\subsection{Morphisms of monads}

Morphisms of monads induce morphisms of Burroni fibrations
and morphisms of tautologous actions of Burroni fibrations.
In details, it looks as follows. Let $(S,\eta^S,\mu^S)$ be
a monad on $\cC$ and $(T,\eta^T,\mu^T)$ be a monad on $\cD$,
$F:\cC\ra \cD$ a functor preserving pullbacks,
and $\xi:TF\ra FS$ be a natural transformation so that
$(F,\xi): (S,\eta^S,\mu^S)\lra (T,\eta^T,\mu^T)$ is a monad morphism, i.e. the diagram
\begin{center} \xext=1300 \yext=680
\begin{picture}(\xext,\yext)(\xoff,\yoff)
\setsqparms[-1`1`0`-1;800`600]
\putsquare(460,20)[TF`T^2F`FS`FS^2;\mu^T_F`\xi``F(\mu^S)]
\putmorphism(1260,620)(0,-1)[\phantom{T^2F}`TFS`T(\xi)]{300}{1}r
\putmorphism(1260,320)(0,-1)[\phantom{TFS}`\phantom{FS^2}`\xi_S]{300}{1}r
\put(20,290){$F$}
\put(110,400){\vector(3,2){250}}
\put(100,490){$\eta^T_F$}
\put(110,240){\vector(3,-2){250}}
\put(0,70){$F(\eta^S)$}
 \end{picture}
\end{center}
commutes. Then we can define the functor from $S$-graphs to $T$-graphs
\begin{center} \xext=1000 \yext=550
\begin{picture}(\xext,\yext)(\xoff,\yoff)
 \setsqparms[1`1`1`1;1000`400]
  \putsquare(0,30)[Gph(S)`Gph(T)`\cC`\cD;Gph(F,\xi)`p_S`p_T`F]
\end{picture}
\end{center}
as
\begin{center} \xext=3200 \yext=430
\begin{picture}(\xext,\yext)(\xoff,\yoff)
\settriparms[1`1`0;400]
  \putAtriangle(0,0)[A`O`S(O);\gamma`\delta`]

 \putAtriangle(1400,0)[F(A)`F(O)`FS(O);F(\gamma)`F(\delta)`]
 \putmorphism(2200,0)(1,0)[\phantom{FS(O)}`TF(O)`\xi_O]{800}{-1}a
 \putmorphism(1800,400)(1,0)[\phantom{F(A)}`Gph(F,\xi)(A)`]{800}{-1}a
 \settriparms[0`1`0;400]
 \putAtriangle(2200,0)[\phantom{Gph(F,\xi)(A)}``\phantom{TF(O)};``]

\put(950,250){\vector(1,0){250}}
\put(950,200){\line(0,1){100}}

\end{picture}
\end{center}
where the square on the right in a pullback. This functor has an obvious structure
($\phi_0$ and $\phi_2$) of a morphism of lax monoidal fibrations.

In particular, as any monad $(T,\eta,\mu)$ on a category $\cC$ with pullbacks
has a monad morphism to the identity monad $1_\cC$, any Burroni fibration on
$\cC$ has a morphism into the $1_\cC$-fibration.
This is another way of saying that the category of $T$-categories has a forgetful
functors into the category of internal categories in $\cC$.

Now a routine verification will show that such a morphism of lax monoidal fibrations of graphs
together with a fibred morphism of basic fibrations
\begin{center} \xext=500 \yext=350
\begin{picture}(\xext,\yext)(\xoff,\yoff)
\setsqparms[1`1`1`1;500`300]
\putsquare(0,0)[\cC^\ra`\cD^\ra`\cC`\cD;F^\ra`cod`cod`F]
 \end{picture}
\end{center}
gives rise to a morphism of tautologous actions.

\subsection{Relative Burroni fibrations and relative $T$-categories}

The construction of a lax monoidal fibration of $T$-graphs can be
performed even on a fibred monad on a fibration.
Suppose $p:\cE\ra \bB$ is a fibration such that the fibres of $p$ have pullbacks.
Moreover $(T,\eta,\mu)$ is a monad on the category $\cE$ so that $T$ is a lax morphism of fibrations
\begin{center} \xext=600 \yext=360
\begin{picture}(\xext,\yext)(\xoff,\yoff)
 \settriparms[1`1`1;300]
  \putVtriangle(0,0)[\cE`\cE`\bB;T`p`p]
\end{picture}
\end{center}
and $\eta$, $\mu$ are fibred natural transformations (i.e. their components lie in the fibres of $p$).
Having such data we can repeat the construction of the category of $T$-graphs but restricting the objects
 to such spans
\begin{center} \xext=600 \yext=340
\begin{picture}(\xext,\yext)(\xoff,\yoff)
\settriparms[1`1`0;300]
  \putAtriangle(0,0)[A`O`T(O);\gamma`\delta`]
\end{picture}
\end{center}
that are in fibres of $p$ (i.e. $p(\gamma)=p(\delta)=1_{p(O)}$). The morphisms are defined as before. In this way,
we get a relative Burroni fibration $p_T:Gph(T,p)\ra \cE$ of $T$-graphs over $p$. Clearly, $p_T$ is a lax monoidal fibration with
the tensor structure defined as before. Thus we have a fibration of monoids with a forgetful to $Gph(T,p)$ as in the diagram
\begin{center} \xext=1000 \yext=1000
\begin{picture}(\xext,\yext)(\xoff,\yoff)
\putmorphism(0,900)(1,0)[Gph(T,p)`Mon(T,p)`\cU_T]{1000}{-1}a
\putmorphism(950,920)(-2,-1)[``]{1050}{1}a
\putmorphism(0,900)(0,-1)[\phantom{Gph(T,p)}`\phantom{\cE}`p_T]{500}{1}l
\putmorphism(0,400)(0,-1)[\cE`\bB`p]{400}{1}l
\put(550,600){\makebox(50,50){$q_T$}}
\end{picture}
\end{center}
of functors and categories. As for any category $\cC$, the functor $! : \cC \lra {\bf 1}$  into the terminal category $\bf 1$ is a fibration, this construction is a generalization of the previous one.

{\bf Remark} We can also define a basic fibration $cod : \cE^{\ra,p}\ra \cE$ relative to a fibration $p:\cE\ra \bB$, so that the objects of $\cE^{\ra,p}$
are morphisms of $\cE$ in fibres of $p$ and morphisms are commuting squares. Then, as previously for the Burroni fibrations, we have a tautologous action the lax monoidal fibration $p_T:Gph(T,p)\ra \cE$   on a fibration $cod : \cE^{\ra,p}\ra \cE$
\begin{center} \xext=800 \yext=580
\begin{picture}(\xext,\yext)(\xoff,\yoff)
 \settriparms[1`1`1;500]
  \putVtriangle(0,0)[Gph(T,p)\times_\cE\cE^{\ra,p}`\cE^{\ra,p}`\cE;\star^{(T,p)}``cod]
\end{picture}
\end{center}
If we take the exponential adjoint of this morphism, as in \ref{taut-act-Burr},  we obtain a (relative) representation of relative $T$-graphs
and relative $T$-categories.

\subsection{Free relative $T$-categories}

The full characterization of those monads $T$ for which the forgetful functor $\cU_T$ defined above has a left adjoint
seem to be unknown. However there are various reasonable sufficient conditions, cf. \cite{B}, \cite{Ke}, \cite{BJT}, \cite{Le}  in case
the monad $(T,\eta,\mu)$ is cartesian. Recall that a monad $(T,\eta,\mu)$ on a category $\cC$ with finite products is cartesian
if $T$ preserves pullbacks and both $\eta$, $\mu$  are cartesian natural transformations. A. Burroni in \cite{B} (pp. 267-269) provided one
such characterization and he noticed that if such an adjoint exists $\cU_T$ is automatically monadic \cite{B} (p. 304).
He also noticed that in certain cases one can iterate the $T$-category construction \cite{B} (p. 269). However the condition for the iteration
in \cite{B} is too strong\footnote{One of the requirement is that the monad $T$ commutes with coproducts.} to be used for our construction below.
T. Leinster in \cite{Le} used a weaker condition for iteration but
he was interested in iteration in particular fibres rather than of the whole fibration.  With the help
of this kind of iteration he defined the set of opetopes \cite{Le} (p. 179 and Appendix D).
The construction of the free monoids described in \cite{Le} is
the same as the earlier and more detailed, yet compact, construction described in \cite{BJT} in Appendix B.
The inductive formula defining the free monoids given in both \cite{BJT} and \cite{Le} seem to appear first in \cite{A} (p. 591)
to describe free algebras for a functor and then in a long comprehensive study
\cite{Ke} (p. 69) that extends and unifies some earlier developments of this and related subjects.
The prerequisites for the construction of the free monoids as well as the final goals differ in \cite{BJT} and \cite{Le}.
In \cite{BJT} the prerequisites are given directly in terms of the properties of the category and the tensor involved to get
a left adjoint to the forgetful functor from the monoids to the monoidal category. In \cite{Le}
the prerequisites are given also in terms of the properties of the category however the property of the tensor is not
specified directly but through the property of the monad the tensor is coming from. Moreover, in \cite{Le} the
aim is not only to get a left adjoint but also to make sure that a monad (and a category it is defined on) deduced
from the new adjunction satisfies the same properties, so that one can iterate the construction, as in \cite{B}.

Below we give a characterization of those fibrations $p$ and fibred
monads $T$ on them for which one can iterate
the process of taking $T$-graphs over a fibration $p$. In the exposition
we use ideas from all the mentioned papers.
The notions of a suitable fibrations and a fibrewise suitable monad are
very much inspired by the notions of a suitable
category and a suitable monad, respectively, cf. \cite{Le} Appendix D.
The main difference of our approach with respect to \cite{Le} is that
we iterate whole fibrations over fibrations and get as a final result
the category of opetopic sets, whereas in \cite{Le}
the construction is done fibre by fibre and gives the set of opetopes
as a result. From the perspective of our construction
this set of opetopes is the set of cells in the terminal opetopic set.

We say that a fibration $p:\cE\ra \bB$ is {\em suitable}
if and only if
\begin{enumerate}
  \item $p$ has fibred pullbacks, finite coproducts, and filtered colimits,
  \item finite coproducts and  filtered colimits are universal in fibres of $p$,
  \item filtered colimits commutes with pullbacks in fibres of $p$.
\end{enumerate}

Let $p:\cE\ra \bB$ be a fibration with fibred pullbacks. A monad $(T,\eta,\mu)$ on $\cE$ is {\em cartesian relative to} $p$ if and only if
$(T,\eta,\mu)$ is a fibred monad over $p$ (i.e. $p\circ T=p$, $p(\eta)=1_p=p(\mu)$) and the restriction of the monad
$(T,\eta,\mu)$ to every fibre of $p$ is a cartesian monad on this fibre.

Let $p:\cE\ra \bB$ be a suitable fibration. We say that a monad $(T,\eta,\mu)$ on $\cE$ is {\em suitable relative to} $p$ if and only if
$(T,\eta,\mu)$ is cartesian relative to $p$ and $T$ preserves  filtered colimits in the fibres of $p$.

The following theorem is the key to the definition of the tower of fibrations that defines the category of opetopic sets.

\begin{theorem} \label{suitable}
Let  $(T,\eta,\mu)$ be a suitable monad relative to a suitable fibration $p:\cE \ra \bB$.
Then
\begin{enumerate}
  \item the fibration $p_T$ over $p$ is again suitable;
  \item the forgetful functor $\cU_T$ is monadic;
  \item the monad $(\widetilde{T},\widetilde{\eta},\widetilde{\mu})$
induced by the adjunction $\cF_T\dashv\cU_T$ is suitable relative to $p_T$.
\end{enumerate}
\begin{center} \xext=1600 \yext=1000
\begin{picture}(\xext,\yext)(\xoff,\yoff)
\putmorphism(600,900)(1,0)[Gph(T,p)`Mon(T,p)`\cU_T]{1000}{-1}a
\putmorphism(1550,920)(-2,-1)[``]{1050}{1}a
\putmorphism(600,900)(0,-1)[\phantom{Gph(T,p)}`\phantom{\cE}`p_T]{500}{1}l
\putmorphism(600,400)(0,-1)[\cE`\bB`p]{400}{1}l
\put(1150,600){\makebox(50,50){$q_T$}}
\put(100,400){\oval(60,60)[l]}
\put(100,430){\vector(1,0){450}}
\put(100,370){\line(1,0){450}}
\put(0,370){\makebox(50,50){$T$}}

\put(100,900){\oval(60,60)[l]}
\put(100,930){\vector(1,0){250}}
\put(100,870){\line(1,0){250}}
\put(0,870){\makebox(50,50){$\widetilde{T}$}}
\end{picture}
\end{center}
\end{theorem}

In the proof of this theorem we shall use the following easy lemma.

\begin{lemma} \label{adj-fib}
Suppose $p$ and $q$ are fibrations and we have two lax morphisms of fibrations $\cU$ and $\cF$, as in the diagram
\begin{center} \xext=1000 \yext=600
\begin{picture}(\xext,\yext)(\xoff,\yoff)
\putmorphism(0,495)(1,0)[\phantom{\cE}`\phantom{\cM}`\cU]{900}{-1}b
\putmorphism(0,520)(1,0)[\cE`\cM`]{1000}{0}b
\putmorphism(0,545)(1,0)[\phantom{\cE}`\phantom{\cM}`\cF]{1000}{1}a
\putmorphism(1000,540)(-2,-1)[``]{1050}{1}a
\putmorphism(0,520)(0,-1)[\phantom{\cE}`\phantom{\bB}`p]{500}{1}l
\put(0,0){\makebox(50,50){$\bB$}}
\put(550,220){\makebox(50,50){$q$}}
\end{picture}
\end{center}
If $\cU$ is a morphism of fibrations and $\cF$ is a left adjoint to $\cU$ when restricted to each fibre then
$\cF$ is a left adjoint to $\cU$.
\end{lemma}

{\bf Remark}
This Lemma could be compared with Lemma 1.8.9 of \cite{Ja}. However we don't require the Beck-Chevalley condition
as we don't expect $\cF$ to be a morphism of fibrations, as in our application it won't be.

{\it Proof of Theorem \ref{suitable}.}~  The functor $Gph(T,p) \lra \cE$ sending $T$-graph $(A,\gamma,\delta)$ to $A$ creates
pullbacks, finite coproducts and  filtered colimits. Thus those limits and colimits have the same exactness properties
in the fibres of $p_T$ as they had in fibres of $p$.  The fact that they are fibred in $p_T$ follows from the fact that they are fibred
in $p$ and that finite coproducts and  filtered colimits are universal.
Thus, that $p_T$ is a suitable fibration.

Recall the construction of the free monoid from \cite{Ke},
\cite{BJT}\footnote{The assumptions that we have on the monad $(T,\eta,\mu)$ obviously sufficient to
for this construction to work.}, \cite{Le}.
For an object $(A,\gamma,\delta)$ in a fibre $\cE_O$ we construct a filtered diagram. We write $A$ for $(A,\gamma,\delta)$ and $O$ for the unit of the tensor  $(O,1_O,\eta_O)$, for short.
\begin{center} \xext=10 \yext=2050
\begin{picture}(\xext,\yext)(\xoff,\yoff)
\putmorphism(0,2000)(0,-1)[A^0=O`\phantom{A^1=O+(A,\gamma,\delta)\otimes O}`e_0]{500}{1}r
\putmorphism(0,1500)(0,-1)[A^1=O+A\otimes O `
\phantom{A^2=O+A\otimes (O   +A \otimes O )}`e_1=1+(1\otimes e_0)]{500}{1}r
\putmorphism(0,1000)(0,-1)[A^2=O+A\otimes (O   + A \otimes O)`
\phantom{A^3=O+A\otimes A^2}`e_2=1+(1\otimes e_1)]{500}{1}r
\putmorphism(0,500)(0,-1)[A^3=O+A\otimes A^2`\ldots`e_3=1+(1\otimes e_2)]{500}{1}r
\end{picture}
\end{center}
with the help of binary coproducts and tensors. The colimit of this diagram in $\cE_O$ is the universe of $\cF_T(A,\gamma,\delta)$.
To see the definition of multiplication for the monoid $\cF_T(A,\gamma,\delta)$ and unit see \cite{BJT}.
As all the operations involved are functorial in the whole fibration, $\cF_T$ is functorial, as well. Thus, by Lemma \ref{adj-fib},
to show that $\cF_T$ is a left adjoint to $\cU_T$  we need to verify that they are adjoint when restricted to each fibre. But this is clear
from \cite{BJT}, \cite{Le}. As $T$ preserves  filtered colimits in fibres of $p$ so does $\cU_T$ in the fibres of $p_T$
and hence $\widetilde{T}$ preserves them, as well.

The monadicity of $\cU_T$ follows from Lemme 1 page 304 of \cite{B} or can be proved directly
using the above explicit construction of the free functor $\cF_T$.
If $(M,m,e)$ is a monoid in  $Gph(T,p) \lra \cE$
then using $m$ and $e$ we can construct inductively an algebra $(M, \alpha: \widetilde{T}(M) \ra M)$ and having
a $\widetilde{T}$-algebra  $(M,\alpha)$ we define
a monoid by putting $e$ equal to
$I_O \ra \widetilde{T}(M) \stackrel{\alpha}{\lra} M$ and $m$ equal to
$M\otimes M\lra I+M\otimes(I+M\otimes I) \ra \widetilde{T}(M) \stackrel{\alpha}{\lra} M$.
The remaining details are left for the readers.

The fact that the induced monad $(\widetilde{T},\widetilde{\eta},\widetilde{\mu})$ is cartesian relative to $p_T$ is also easy.
$~\Box$

\vskip 2mm

{\em Remark.} The monadicity of $\cU_T$ was already noticed in \cite{B} Proposition II.1.19 for a monad $T$ satisfying slightly stronger conditions.

\subsection{A tower of fibrations for opetopic sets}\label{tower-opetopic-sets}

Using the above Theorem \ref{suitable}, and starting with any fibrewise
suitable monad $T_0$ on a fibrewise suitable fibration $p:\cE_0 \ra \bB$,
we can build a tower of (fibrewise suitable) lax monoidal fibrations and fibrewise suitable monads as
in the diagram below:
\begin{center} \xext=1600 \yext=2150
\begin{picture}(\xext,\yext)(\xoff,\yoff)

\put(530,2100){\makebox(50,50){$\vdots$}}
\put(1800,2100){\makebox(50,50){$\vdots$}}

\putmorphism(900,1900)(1,0)[\cE_3=Gph(T_2,p_{T_1})`Mon(T_2,p_{T_1})`\cU_{T_2}]{1000}{-1}a
\putmorphism(1800,1920)(-2,-1)[``]{1050}{1}a
\putmorphism(750,1900)(0,-1)[\phantom{\cE_3=Gph(T_2,p_{T_1})}`\phantom{\cE_1=Gph(T,p)}`p_{T_2}]{500}{1}l
\put(1400,1600){\makebox(50,50){$q_{T_2}$}}
\put(250,1900){\oval(60,60)[l]}
\put(250,1930){\vector(1,0){230}}
\put(250,1870){\line(1,0){230}}
\put(0,1870){\makebox(50,50){$T_3=\widetilde{T_2}$}}


\putmorphism(900,1400)(1,0)[\cE_2=Gph(T_1,p_{T_0})`Mon(T_1,p_{T_0})`\cU_{T_1}]{1000}{-1}a
\putmorphism(1800,1420)(-2,-1)[``]{1050}{1}a
\putmorphism(750,1400)(0,-1)[\phantom{\cE_2=Gph(T_1,p_{T_0})}`\phantom{\cE_1=Gph(T,p)}`p_{T_1}]{500}{1}l
\put(1400,1100){\makebox(50,50){$q_{T_1}$}}
\put(250,1400){\oval(60,60)[l]}
\put(250,1430){\vector(1,0){230}}
\put(250,1370){\line(1,0){230}}
\put(0,1370){\makebox(50,50){$T_2=\widetilde{T_1}$}}

\putmorphism(900,900)(1,0)[\cE_1=Gph(T_0,p)`Mon(T_0,p)`\cU_{T_0}]{1000}{-1}a
\putmorphism(1800,920)(-2,-1)[``]{1100}{1}a
\putmorphism(750,900)(0,-1)[\phantom{\cE_1=Gph(T_0,p)}`\phantom{\cE_0}`p_{T_0}]{500}{1}l
\putmorphism(750,400)(0,-1)[\cE_0`\bB`p]{400}{1}l
\put(1400,600){\makebox(50,50){$q_{T_0}$}}
\put(250,400){\oval(60,60)[l]}
\put(250,430){\vector(1,0){420}}
\put(250,370){\line(1,0){420}}
\put(130,370){\makebox(50,50){$T_0$}}

\put(250,900){\oval(60,60)[l]}
\put(250,930){\vector(1,0){250}}
\put(250,870){\line(1,0){250}}
\put(0,870){\makebox(50,50){$T_1=\widetilde{T_0}$}}
\end{picture}
\end{center}
So as $p$ and $T_0$ are fibrewise suitable, we have the monad $T_1=\widetilde{T}$ on $p_{T_0}$.
By Theorem \ref{suitable} $p_{T_0}$ and $T_1$ are again suitable and hence we can repeat
the construction again. The identity monad $1_{Set}$ on $Set$ is of course a fibrewise suitable
on the fibrewise suitable fibration $!:Set \ra {\bf 1}$, where ${\bf 1}$ is the terminal category.
Thus we can build a tower of fibrations, as above, starting form this fibration.  We obtain
\begin{center} \xext=1600 \yext=2100
\begin{picture}(\xext,\yext)(\xoff,\yoff)

\put(500,2050){\makebox(50,50){$\vdots$}}
\put(1800,2050){\makebox(50,50){$\vdots$}}

\putmorphism(900,1900)(1,0)[\cO_3=Gph(T_2,p_{T_1})`Mon(T_2,p_{T_1})`\cU_{T_2}]{1000}{-1}a
\putmorphism(1700,1920)(-2,-1)[``]{1050}{1}a
\putmorphism(750,1900)(0,-1)[\phantom{\cO_3=Gph(T_2,p_{T_1})}`\phantom{\cO_1=Gph(T,p)}`p_{T_2}]{500}{1}l
\put(1300,1600){\makebox(50,50){$q_{T_2}$}}
\put(250,1900){\oval(60,60)[l]}
\put(250,1930){\vector(1,0){230}}
\put(250,1870){\line(1,0){230}}
\put(0,1870){\makebox(50,50){$T_3=\widetilde{T_2}$}}


\putmorphism(900,1400)(1,0)[\cO_2=Gph(T_1,p_{T_0})`Mon(T_1,p_{T_0})`\cU_{T_1}]{1000}{-1}a
\putmorphism(1700,1420)(-2,-1)[``]{1050}{1}a
\putmorphism(750,1400)(0,-1)[\phantom{\cO_2=Gph(T_1,p_{T_0})}`\phantom{\cO_1=Gph(T,p)}`p_{T_1}]{500}{1}l
\put(1300,1100){\makebox(50,50){$q_{T_1}$}}
\put(250,1400){\oval(60,60)[l]}
\put(250,1430){\vector(1,0){230}}
\put(250,1370){\line(1,0){230}}
\put(0,1370){\makebox(50,50){$T_2=\widetilde{T_1}$}}

\putmorphism(900,900)(1,0)[\cO_1=Gph(T_0,!)`Mon(T_0,!)`\cU_{T_0}]{1000}{-1}a
\putmorphism(1700,920)(-2,-1)[``]{1050}{1}a
\putmorphism(750,900)(0,-1)[\phantom{\cO_1=Gph(T_0,p)}`\phantom{\cO_0}`p_{T_0}]{500}{1}l
\putmorphism(750,400)(0,-1)[\cO_0=Set`{\bf 1}`!]{400}{1}l
\put(1300,600){\makebox(50,50){$q_{T_0}$}}
\put(250,400){\oval(60,60)[l]}
\put(250,430){\vector(1,0){230}}
\put(250,370){\line(1,0){230}}
\put(-20,370){\makebox(50,50){$T_0=1_{Set}$}}

\put(250,900){\oval(60,60)[l]}
\put(250,930){\vector(1,0){250}}
\put(250,870){\line(1,0){250}}
\put(0,870){\makebox(50,50){$T_1=\widetilde{T_0}$}}
\end{picture}
\end{center}
An {\em opetopic set} is an infinite sequence of objects $\{A_n\}_{n\in\o}$ such that
\begin{enumerate}
  \item $A_n$ is an object in $\cO_n$,
  \item $A_{n+1}$ lies in the fibre over $A_n$, i.e. $p_{T_n}(A_{n+1})=A_n$,
\end{enumerate}
for $n\in \o$. A {\em morphism of opetopic sets} $\{f_n\}_{n\in\o} : \{A_n\}_{n\in\o} \lra \{B_n\}_{n\in\o}$ is
a family of morphisms such that
\begin{enumerate}
  \item $f_n : A_n \lra B_n$ is a morphism in $\cO_n$
  \item $f_{n+1}$ lies in the fibre over $f_n$, i.e. $p_{T_n}(f_{n+1})=f_n$,
\end{enumerate}
for $n\in \o$.

Unraveling this definition, we see that an opetopic set (in the above sense) is an
$\infty$-span as the diagram below:
\begin{center} \xext=600 \yext=2050
\begin{picture}(\xext,\yext)(\xoff,\yoff)
 \setsqparms[0`1`1`0;600`600]

\put(0,1950){\makebox(50,50){$\vdots$}}
\put(550,1950){\makebox(50,50){$\vdots$}}

  \putsquare(0,1200)[A_3`T_3(A_3)`\phantom{A_2}`\phantom{T_2(A_2)};`\gamma_2`\overline{\delta}_2`]
  \putmorphism(0,1800)(1,-1)[\phantom{A_3}`\phantom{T(A_2)}`]{600}{1}a
   \putmorphism(300,1500)(-1,-1)[`\phantom{A_2}`]{300}{1}a
   \put(530,1730){\line(-1,-1){200}}
   \put(100,1400){\makebox(50,50){$\overline{\gamma}_2$}}
   \put(450,1400){\makebox(50,50){$\delta_2$}}

  \putsquare(0,600)[A_2`T_2(A_2)`\phantom{A_1}`\phantom{T_1(A_1)};`\gamma_1`\overline{\delta}_1`]
  \putmorphism(0,1200)(1,-1)[\phantom{A_2}`\phantom{T(A_1)}`]{600}{1}a
   \putmorphism(300,900)(-1,-1)[`\phantom{A_1}`]{300}{1}a
   \put(530,1130){\line(-1,-1){200}}
   \put(100,800){\makebox(50,50){$\overline{\gamma}_1$}}
   \put(450,800){\makebox(50,50){$\delta_1$}}

  \putsquare(0,0)[A_1`T_1(A_1)`A_0`T_0(A_0);`\gamma_0`\overline{\delta}_0`]
  \putmorphism(0,600)(1,-1)[\phantom{A_1}`\phantom{T(A_0)}`]{600}{1}a
   \putmorphism(300,300)(-1,-1)[`\phantom{A_0}`]{300}{1}a
   \put(530,530){\line(-1,-1){200}}
   \put(100,200){\makebox(50,50){$\overline{\gamma}_0$}}
   \put(450,200){\makebox(50,50){$\delta_0$}}
\end{picture}
\end{center}
with

\[   \gamma_n\circ\gamma_{n+1} =\overline{ \gamma}_n\circ\delta_{n+1}, \hskip 10mm  \delta_n\circ\gamma_{n+1} =\overline{ \delta}_n\circ\delta_{n+1} \]
\[   \gamma_n\circ\overline{\gamma}_{n+1} =\overline{ \gamma}_n\circ\overline{\delta}_{n+1}, \hskip 10mm
 \delta_n\circ\overline{\gamma}_{n+1} =\overline{\delta}_n\circ\overline{\delta}_{n+1}  \]
for $n\in \o$. To describe the terminal opetopic set $A$, we need to start with $A_0=1$ the terminal object in $Set$. And then choose $A_{n+1}$ as the terminal object in the fibre of $p_{T_n}$ over $A_n$. Thus $A_1$ is $1$ and $A_{n+1}$ for $n>0$ can be taken as the limit in the following diagram:
\begin{center} \xext=600 \yext=1250
\begin{picture}(\xext,\yext)(\xoff,\yoff)
 \setsqparms[0`1`0`0;600`600]

  \putsquare(0,600)[A_{n+1}``\phantom{A_n}`\phantom{T_n(A_n)};`\gamma_n``]
  \putmorphism(0,1200)(1,-1)[\phantom{A_{n+1}}`\phantom{T(A_n)}`]{600}{1}a

   \put(380,870){\makebox(50,50){$\delta_n$}}

  \setsqparms[0`1`1`0;600`600]
  \putsquare(0,0)[A_n`T_n(A_n)`A_{n-1}`T_{n-1}(A_{n-1});`\gamma_{n-1}`\overline{\delta}_{n-1}`]
  \putmorphism(0,600)(1,-1)[\phantom{A_n}`\phantom{T(A_{n-1})}`]{600}{1}a
   \putmorphism(300,300)(-1,-1)[`\phantom{A_{n-1}}`]{300}{1}a
   \put(530,530){\line(-1,-1){200}}
   \put(100,250){\makebox(50,50){$\overline{\gamma}_{n-1}$}}
   \put(450,230){\makebox(50,50){$\delta_{n-1}$}}
\end{picture}
\end{center}
The disjoint union of the sets $\{ A_n\}_{n\in \o}$ is the set of opetopes in the sense of T. Leinster.

The proof of the following theorem is uses ordered face structures, cf. \cite{Z}, and will not be given here.

\begin{theorem}
The category of opetopic sets so defined is equivalent to the category of multitopic sets.
\end{theorem}

{\bf Remark} {\em Internal opetopic sets.} Clearly the fibration $!:Set\ra {\bf 1}$ is not the only interesting suitable one to start the process of iteration.
For example,  we can start with $!:\cE\ra {\bf 1}$ where $\cE$ is a sufficiently cocomplete topos.
Thus, we have the category of internal opetopic sets in any Grothendieck topos, even in the category of opetopic sets itself!

\subsection{A tower of fibrations for $n$-categories}\label{tower-n-cat}

If we start with the (fibred) identity monad $1_\cE$ on a fibration $p:\cE\ra \bB$ whose fibres have pullbacks
then the fibration of monoids over $p$, $q_{1_\cE}:Mon(1_\cE,p)\lra \cE$ again has pullbacks in the fibres. Thus we can
iterate this process and get another tower of fibrations based on monoids, this time:
\begin{center} \xext=1600 \yext=2150
\begin{picture}(\xext,\yext)(\xoff,\yoff)

\put(530,2100){\makebox(50,50){$\vdots$}}
\put(1800,2100){\makebox(50,50){$\vdots$}}

\putmorphism(900,1900)(1,0)[\cE_3=Mon(1_{\cE_2},p_{1_{\cE_1}})`Gph(1_{\cE_2},p_{1_{\cE_1}})`\cU_{1_{\cE_2}}]{1000}{1}a
\putmorphism(1800,1920)(-2,-1)[``]{1050}{1}a
\putmorphism(750,1900)(0,-1)[\phantom{\cE_3=Mon(1_{\cE_2},p_{1_{\cE_1}})}`\phantom{\cE_1=Mon(1_{\cE_1},p)}`q_{1_{\cE_2}}]{500}{1}l
\put(1400,1600){\makebox(50,50){$p_{1_{\cE_2}}$}}
\put(250,1900){\oval(60,60)[l]}
\put(250,1930){\vector(1,0){180}}
\put(250,1870){\line(1,0){180}}
\put(130,1870){\makebox(50,50){$1_{\cE_3}$}}


\putmorphism(900,1400)(1,0)[\cE_2=Mon(1_{\cE_1},p_{1_{\cE_0}})`Gph(1_{\cE_1},p_{1_{\cE_0}})`\cU_{1_{\cE_1}}]{1000}{1}a
\putmorphism(1800,1420)(-2,-1)[``]{1050}{1}a
\putmorphism(750,1400)(0,-1)[\phantom{\cE_2=Mon(1_{\cE_1},p_{1_{\cE_0}})}`\phantom{\cE_1=Mon(1_{\cE_0},p)}`q_{1_{\cE_1}}]{500}{1}l
\put(1400,1100){\makebox(50,50){$p_{1_{\cE_1}}$}}
\put(250,1400){\oval(60,60)[l]}
\put(250,1430){\vector(1,0){180}}
\put(250,1370){\line(1,0){180}}
\put(130,1370){\makebox(50,50){$1_{\cE_2}$}}

\putmorphism(900,900)(1,0)[\cE_1=Mon(1_{\cE_0},p)`Gph(1_{\cE_0},p)`\cU_{1_{\cE_0}}]{1000}{1}a
\putmorphism(1800,920)(-2,-1)[``]{1100}{1}a
\putmorphism(750,900)(0,-1)[\phantom{\cE_1=Mon(1_{\cE_0},p)}`\phantom{\cE_0}`q_{1_{\cE_0}}]{500}{1}l
\putmorphism(750,400)(0,-1)[\cE_0`\bB`p]{400}{1}l
\put(1400,600){\makebox(50,50){$p_{1_{\cE_0}}$}}
\put(250,400){\oval(60,60)[l]}
\put(250,430){\vector(1,0){420}}
\put(250,370){\line(1,0){420}}
\put(130,370){\makebox(50,50){$1_{\cE_0}$}}

\put(250,900){\oval(60,60)[l]}
\put(250,930){\vector(1,0){250}}
\put(250,870){\line(1,0){250}}
\put(130,870){\makebox(50,50){$1_{\cE_1}$}}
\end{picture}
\end{center}
If we start with a suitable fibration $p:\cC\ra {\bf 1}$, then after $n$-th iteration
we recover the D. Bourn \cite{Bo} construction of internal $n$-categories in $\cC$.

\section{Amalgamated signatures vs polynomial functors}\label{am-sig}

\subsection{The amalgamated signatures fibration $p_a: Sig_a \ra Set$}

This example is one of the main reasons for considering lax monoidal fibrations in the context of higher category theory at all.
The monoids in this fibration are precisely the (1-level) multicategories
with non-standard amalgamation. They are like the multicategories considered by C. Hermida M.Makkai J. Power in \cite{HMP} to define
the multitopic sets, except that there the 2-level version is used by. This modification will be explained at the end of the
section.

{\em Notation.} Let $[n]=\{0,\ldots,n\}$, $(n]=\{1,\ldots,n\}$,
for $n\in\o$. In particular $[n]=[0]\cup(n]$ and $(0]=\emptyset$. For a set $O$, we put
$O^\dag_n=O^{[n]}$, $O^*_n=O^{(n]}$ and
$O^\dag=\bigcup_{n\in\o}O^{[n]}$, $O^*=\bigcup_{n\in\o}O^{(n]}$.
$S_n$ acts on both $O^\dag_n$ and $O^*_n$ on the right by
composition (i.e. we leave $0$ fixed in the domain of the elements
of $O^\dag_n$). If $d:[n]\ra O$ is a function, then its restriction
to the positive numbers is denoted by $d^+:(n]\ra O$ and to $[0]$
by $d^-:[0]\ra O$. This restrictions establish a bijection $\lk
(-)^-,(-)^+\rk : O^\dag\ra O\times O^*$. Clearly $(-)^\dag : Set \lra Set$ is a functor.

The base category of our fibration is $Set$. The total category $Sig_a$ of our
fibration  has as objects triples, $(A,\partial,O)$ such that $A$ and $O$ are sets and
$\partial:A\ra O^\dag$ is a function.
We write $\partial_a:[n]\ra O$ for the effect of $\partial$ on $a\in
A$, and $n$ in this case will be referred to as $|a|$. A morphism
$(f,\sigma,u):(A,\partial,O)\ra(B,\partial,Q)$ in $Sig_a$ is a pair of
functions $f:A\ra B$ and $u:O\ra Q$, and for any $a\in A$ with $n=|a|$
a permutation $\sigma_a: [n]\ra [n]\in S_n$ (with $\sigma_a(0)=0$) making the
square
\begin{center}
\xext=600 \yext=500
\begin{picture}(\xext,\yext)(\xoff,\yoff)
 \setsqparms[-1`1`1`1;600`400]
 \putsquare(0,20)[[n]`[n]`O`Q;\sigma_a`\partial_a`\partial_{f(a)}`u]
\end{picture}
\end{center}
commute. A morphism $(f,\sigma,u)$ is called {\em strict} if $\sigma_a$ is an identity,  for $a\in A$ .
The projection functor $p:Sig_a\lra Set$ sends the morphism $(f,\sigma,u):(A,\partial,O)\ra(B,\partial,Q)$ to $u:O\ra Q$.

{\bf Remarks}
\begin{enumerate}
  \item The category $Sig_m$ is isomorphic to the full subcategory of $Sig_a$  whose morphisms are strict.
  \item We think of an object $(A,\partial,O)$ of  $Sig_a$ as a signature with $O$ as the set of its types, $A$ the set of its operation symbols,
and $\partial$ the typing function associating arities to function symbols $a: \partial_a(1),\ldots, \partial_a(|a|)\lra \partial_a(0)$, i.e. is
$\partial_a(1),\ldots, \partial_a(|a|)$ are types of the arguments (inputs) of $a$ and $\partial_a(0)$ is the type of values (outputs) of $a$.
\end{enumerate}

\subsection*{The lax monoidal structure on $p_a$}

We have two lax morphisms of fibrations
\begin{center} \xext=1800 \yext=650
\begin{picture}(\xext,\yext)(\xoff,\yoff)
 \setsqparms[1`0`1`0;900`500]
  \putsquare(0,0)[{Sig_a}\times_{Set}{Sig_a}`Sig_a``Set;\otimes``p_a`]
  \setsqparms[-1`0`0`0;900`500]
  \putsquare(900,0)[\phantom{Sig_a}`Set`\phantom{Set}`;I```]
  \put(100,420){\vector(2,-1){740}}
      \put(230,200){$p'_a$}
   \put(1700,430){\vector(-2,-1){760}}
   \put(1450,200){$1_{Set}$}
\end{picture}
\end{center}

Let $(A,\partial,O)$ and $(B,\partial,O)$ be two object in the fibre
over $O$. Their tensor $(A\otimes_O B,\partial^\otimes,O)$ is
defined as follows
\[A\otimes_O B=\{ \lk a,b_i \rk_{i\in (|a|]} : a\in A,\, b_i\in B,\,
\partial_a(i)=\partial_{b_i}(0),\; \mbox{for }  i\in (|a|] \} \]
and for $\lk a,b_i\rk_{i\in (|a|]}\in A\otimes_O B$,
\[ \partial^\otimes_{\lk a,b_i\rk _{i}}=
[ \partial^-_a, \partial^+_{b_i}]_{i} \, :\, [|\lk a,b_i\rk_i|] = [\, \sum_{i=1}^{|a|} |b_i|\, ]  \lra O.\]
Note that just saying that we have a coproduct determines the function  $[ \partial^-_a, \partial^+_{b_i}]_{i}$ only up to a permutation.
In principle we don't need more than that for as far as $[ \partial^-_a, \partial^+_{b_i}]_{i}(0)= \partial^-_a(0)$. But to be on the safe side, we
will always tacitly assume that the domains of $\partial^+_{b_i}$ are placed one after the other.

 For a pair of maps in $Sig_a$
\[ \overline{f}=(f,\sigma,u):(A,\partial,O)\ra(A',\partial,Q),\;\;\;
\overline{g}=(g,\tau, u):(B,\partial,O)\ra(B',\partial,Q)\]
over the same map $u:O\ra Q$ we define the map
\[\overline{f}\otimes_u\overline{g}=(f\otimes_u
g,\sigma\otimes_u\tau,u):(A\otimes_O
B,\partial^\otimes,O)\lra (A'\otimes_Q B',\partial^\otimes,Q) \]
 so that, for $\lk a, b_i\rk_{i\in (|a|]}\in A\otimes_O B$,
\[ f\otimes_u g(\lk a,b_i\rk_{i\in (|a|]}) =
\lk f(a),g(b_{\sigma_a(j)})\rk_{j\in (|f(a)|]}\]
Clearly, $|a|=|f(a)|$ and
\[ n=|\lk a, b_i\rk_{i\in (|a|]}|= \sum_{i\in |a|} |b_i| = \sum_{i\in |f(a)|} |g(b_i)| = |f\otimes_u g(\lk a,b_i\rk_{i\in (|a|]})|. \]
Moreover, we put
 \[ (\sigma\otimes_u\tau)_{\lk a, b_i \rk_i}=
 [\tau^+_{b_{\sigma_a(i)}}]_{i} \]
making the square
\begin{center} \xext=1200 \yext=650
\begin{picture}(\xext,\yext)(\xoff,\yoff)
   \setsqparms[-1`1`1`1;1200`500]
  \putsquare(0,50)[ [n]`[n]`O`Q;(\sigma\otimes_u\tau)_{\lk
a,b_i\rk_i}`\partial^\otimes_{\lk a,b_i\rk_i}`
 \partial^\otimes_{\lk f(a),g(b_i)\rk_i}`u]
\end{picture}
\end{center}
commute. This ends the definition of the tensor $\otimes$.

The unit $I_O$  in the fibre over $O$, is $I_O=(O,\partial^{I_O},O)$
such that for $x\in O$, $\partial^{I_O}_x: [1]\ra O$ is a
constant function equal to $x$. We note for the record

\begin{lemma}The fibration $p_a:Sig_a\ra Set$ with the structure described above is a lax monoidal fibration
whose fibres are strong monoidal categories. $\Box$
\end{lemma}

\subsection*{Pulling back the monoidal structure}
We shall describe how reindexing functors interact with the monoidal structure in the fibration $p_a$.

Any  object $B$ in the fibre over $Q$ of $p_a:Sig_a\ra Set$ can be pulled back along a
function $u:O\ra Q$:
\begin{center} \xext=700 \yext=700
\begin{picture}(\xext,\yext)(\xoff,\yoff)
   \setsqparms[1`1`1`1;700`400]
  \putsquare(0,250)[u^*(B)`B`O^\dag`Q^\dag;\pi_B`\partial`\partial`u^\dag]
  \putmorphism(0,0)(1,0)[O`Q`u]{700}{1}b
\end{picture}
\end{center}
thus
\[ u^*(B)=\{ \lk b,d\rk : b\in B,\, d:[|b|]\ra O,\, \mbox{such that}\; u^\dag(d)=\partial_b \}\]
($u^\dag(d)=u\circ d$) and
\[ \partial_{\lk b,d\rk}=d \]
 We have
\[ u^*(I_{Q})=\{ \lk x,x'\rk\in O^2: u(x)=u(x') \} \]
and
\[ \phi_0: I_O\lra u^*(I_{Q}) \]
\[ x\mapsto \lk x,x\rk \]
Moreover, for objects $A$ and $B$ over $Q$ we have

\[ u^*(A\otimes B)= \{ \lk\lk a,b_i\rk_{i\in (|a|]},d\rk:  \lk a,b_i\rk_{i\in (|a|]}\in A\otimes B,\;   u^\dag(d) =\partial^\otimes_{\lk a,b_i\rk_{i\in (|a|]}} \} \]
and
\[ u^*(A)\otimes u^*(B)= \{ \lk\lk a,d\rk,\lk b_i,d_i\rk\rk_{i\in (|a|]}: a \in A,\; b_i\in B,\; \hskip 5cm \]
\[  \hskip 3cm  u^\dag(d)=\partial_a,\, u^\dag(d_i) =\partial_{b_i},\; d(i)=d_i(0)   ,\; \mbox{for}\; i\in (|a|] \}\]

Thus we have a transformation
 \[ \phi_{2,A,B}: u^*(A)\otimes u^*(B)\lra u^*(A\otimes B)\]
 such that
\[  \lk\lk a,d\rk,\lk b_i,d_i\rk\rk_{i\in (|a|]} \mapsto \lk\lk a,b_i\rk_{i\in (|a|]},[d^-,d^+_i]_{i\in
(|a|]}\rk  \]

All the morphisms defined above $\pi_B$, $\phi_0$, and $\phi_{2,A,B}$ are strict, i.e. with amalgamations being identities.

\begin{lemma} The data $u^*$, $\phi_0$, $\phi_2$ above, make the usual (three)
diagrams for coherence of monoidal functor (not necessarily strong) commute.
\end{lemma}

{\it Proof.}~ Exercise. $~\Box$

Moreover we have

\begin{proposition}  The total category of the fibration $q_a:Mon(Sig_a)\lra Set$ of monoids in
$p_a:Sig_a\lra Set$ is equivalent to the category
of (1-level) multicategories with non-standard amalgamations.
The fibred forgetful functor from the fibration of monoids to the
fibration of amalgamated signatures
${\cal U}: Mon(Sig_a)\lra Sig_a$
\begin{center} \xext=600 \yext=580
\begin{picture}(\xext,\yext)(\xoff,\yoff)
 \settriparms[1`1`1;500] \putVtriangle(0,0)[Sig_a`Mon(Sig_a)`Set;{\cal F} `p_a`q_a]
 \putmorphism(150,420)(1,0)[``{\cal U}]{650}{-1}b
\end{picture}
\end{center}
has a fibred left adjoint $\cF$, the free monoid functor.
\end{proposition}

{\it Proof.}~ Strictly speaking the multicategories with non-standard amalgamations
were defined \cite{HMP} from the single tensor
and additional property, called commutativity there.
But, as it is well known, they can be equivalently defined using the total tensor,
i.e. the one we defined above. For more, see also Subsection \ref{single-tensor}.  $~\Box$

The free functor $\cF$ mentioned in the Proposition above was described in \cite{HMP}.

\subsection{The action of $p_a$ on the basic fibration}

The lax monoidal fibration $p_a$ comes equipped naturally with an action on the basic fibration $cod:Set^\ra\lra Set$
\begin{center} \xext=800 \yext=520
\begin{picture}(\xext,\yext)(\xoff,\yoff)
 \settriparms[1`1`1;400]
  \putVtriangle(0,0)[Sig_a\times_{Set}{Set}^\ra`{Set}^\ra`{Set};\star``cod]
\end{picture}
\end{center}
For $(A,\partial^A,O)$ in $Sig_a$ and $(X,d,O)$ in $Set^\ra$, the set $A\star X$ is defined from the following diagram
\begin{center} \xext=1000 \yext=520
\begin{picture}(\xext,\yext)(\xoff,\yoff)
 \setsqparms[-1`1`1`-1;600`400]
  \putsquare(500,0)[A`A\star X`O^*`X^*;`\partial^{A,+}``d^*]
  \putmorphism(0,400)(1,0)[O`\phantom{A}`\partial^{A,-}]{500}{-1}a
\end{picture}
\end{center}
with the square being  a pullback. $(-)^*$ is the free monoid functor, i.e.
\[ A\star X =\{(a,x_1,\ldots, x_{|a|}) : \partial^A_a(i)=d(x_i),\; i=1,\ldots, |a| \} \]
and $\partial^\star :A\star X \lra O$ is defined by
\[ \partial^\star(a,x_1,\ldots, x_{|a|}) =\partial^A_a(0) \]
Thus it is the composition of the upper horizontal morphism in the above diagram.
On morphisms the action $\star$ is defined in the obvious way.
Thus we have an adjoint morphism of lax monoidal fibrations
\begin{center} \xext=800 \yext=520
\begin{picture}(\xext,\yext)(\xoff,\yoff)
 \settriparms[1`1`1;400]
  \putVtriangle(0,0)[Sig_a`Exp(Set)`{Set};rep_a`p_a`]
\end{picture}
\end{center}
where, as usual, $Exp(Set)$ is the exponent in $\bCat/Set$.

\subsection{Polynomial diagrams and polynomial functors}

In this subsection we collect the definitions and facts concerning polynomial
diagrams and polynomial functors from the literature, that are needed in the following.
For (much!) more, the reader should consult \cite{Ko}, \cite{GK}
and bibliography there. We deal with polynomial functors based on
an arbitrary locally cartesian closed category but with a
special eye on $Set$ and the presheaf category $Set^{\cS_*}$,
where $\cS_*$ is the coproduct  in $Cat$ of the (finite) symmetric groups.
The later category will be important in Section \ref{symm-sig}.
In this section, unless otherwise specified,  $\cE$ is an arbitrary locally
cartesian closed category, and by this we mean that  $\cE$ has the terminal
object, as well.

By a {\em  polynomial diagram (over $O$)} in $\cE$, we mean the
following diagram in $\cE$
\begin{center} \xext=1500 \yext=150
\begin{picture}(\xext,\yext)(\xoff,\yoff)
\putmorphism(500,0)(1,0)[E`A`p]{500}{1}a
\putmorphism(0,0)(1,0)[O`\phantom{E}`s]{500}{-1}a
\putmorphism(1000,0)(1,0)[\phantom{A}`O`t]{500}{1}a
 \end{picture}
\end{center}
The object $O$ is an object of types of the polynomial $(t,p,s)$.
We say that a polynomial diagram $(t,p,s)$ in $Set$ is {\em
finitary} if and only if the function $p$ has finite fibres. A {\em morphism
of polynomial diagrams (over $u:O\ra Q$)} in $\cE$ is a triple
$(f,g,u)$ of morphism making the diagram
\begin{center} \xext=1500 \yext=600
\begin{picture}(\xext,\yext)(\xoff,\yoff)
\setsqparms[1`1`1`1;500`400]
\putsquare(500,50)[E`A`E'`A';p`g`f`p']
\setsqparms[-1`1`0`-1;500`400]
\putsquare(0,50)[O`\phantom{E}`Q`\phantom{E'};s`u``s']
\setsqparms[1`0`1`1;500`400]
\putsquare(1000,50)[\phantom{A}`O`\phantom{A'}`Q;t``u`t']
 \end{picture}
\end{center}
commute, and such that the square in the middle is a pullback.
Morphisms of polynomial diagram compose in the obvious way, by
putting one on top of the other. Let $\cP oly\cD iag(\cE)$ denotes
the category of the polynomial diagrams and morphisms between them.

\vskip 2mm
{\bf Remark} If $\cE$ is the category $Set$, we can think of $A$ as the set of operations, and $E$ as the set of
arguments of all operations in $A$. Thus with this interpretation $p^{-1}(a)$ is
the set of arguments (or arity) of $a$. Then $s(e)$, for $e\in p^{-1}(a)$, can be interpreted
as the type of the argument $e$ of the operation $a$, and $t(a)$ is the type of the
values of the operation $a$.

\vskip 2mm

We have an obvious projection functor
\[ p_{pd,\cE}:\cP oly\cD iag(\cE)\lra \cE\]
sending $(f,g,u)$ to $u$, which is a lax
monoidal fibration. The tensor in fibres is given by
composition of diagrams, cf. \cite{GK} 1.11. Let
\[ p_{pd}:\cP oly\cD iag\ra Set \]
denote the finitary polynomial diagram fibration, the full
subfibration of $p_{pd,Set}$ consisting of finitary polynomial
diagrams. By $p_{poly}:\cP oly\ra Set$ we denote the image of the fibration $p_{pd}$ in $p_{exp,Set}$.
We shall see that $p_{poly}$ is a lax monoidal subfibration of $p_{exp,Set}$.

In our terminology, the connection between polynomial diagrams and polynomial functors can be expressed as the fact that the fibration
of polynomial diagrams comes equipped with a representation morphism into $Exp(\cE)$, i.e. a morphism of lax monoidal fibration,
cf. Proposition \ref{rep-E1} below,
\begin{center} \xext=800 \yext=550
\begin{picture}(\xext,\yext)(\xoff,\yoff)
 \settriparms[1`1`1;450]
  \putVtriangle(0,0)[\cP oly\cD iag(\cE)`Exp(\cE)`\cE;rep_{pd,\cE}`p_{pd,\cE}`p_{exp}]
\end{picture}
\end{center}
whose essential image is, by definition, the (lax monoidal)
fibration of (finitary) polynomial (endo)functors and polynomial
transformations between them. We shall recall this now. Any
morphism $u:O\ra Q$ in a locally cartesian closed category $\cE$
induces three functors
\begin{center} \xext=600 \yext=300
\begin{picture}(\xext,\yext)(\xoff,\yoff)
\putmorphism(0,90)(1,0)[\cE_{/O}`\cE_{/Q}`u^*]{600}{-1}a
\putmorphism(0,30)(1,0)[\phantom{\cE_{/O}}`\phantom{\cE_{/Q}}`u_*]{600}{1}b
\putmorphism(0,200)(1,0)[\phantom{\cE_{/O}}`\phantom{\cE_{/Q}}`{u_!}]{600}{1}a
 \end{picture}
\end{center}
such that $u^*$ is a pullback functor and $u_!\dashv u^*\dashv u_*$.
The unit and counit of the adjunction $u_!\dashv u^*$ will be denoted by $\eta^u$ and $\varepsilon^u$, respectively and
the unit and counit of the adjunction $u^*\dashv u_*$ will be denoted by $\bar{\eta}^u$ and $\bar{\varepsilon}^u$, respectively.

For an object $(t,p,s)$ over $O$, we define a functor
\[ rep_{pd,\cE}(t,p,s)=t_!p_*s^* : \cE_{/O} \lra \cE_{/O}\]
For a morphism of polynomial diagrams $(f,g,u): (t,p,s)\lra
(t',p',s')$, we define a morphism $rep_{pd,\cE}(f,g,u):rep_{pd,\cE}(t,p,s)\ra rep_{pd,\cE}(t',p',s')$ in $Exp(\cE)$ over $u$,
as follows. We have a diagram of categories, functors and natural transformations
\begin{center} \xext=2200 \yext=600
\begin{picture}(\xext,\yext)(\xoff,\yoff)
\setsqparms[1`-1`0`1;600`500]
\putsquare(0,50)[\cE_{/O}`\phantom{\cE_{/E}}`\cE_{/Q}`\phantom{\cE_{/E'}};s^*`u^*``s'^*]
\setsqparms[1`-1`1`1;800`500]
\putsquare(600,50)[\cE_{/E}`\cE_{/A}`\cE_{/E'}`\cE_{/A'};p_*`g^*`f_!`p'_*]
\setsqparms[1`0`1`1;800`500]
\putsquare(1400,50)[\phantom{\cE_{/A}}`\cE_{/O}`\phantom{\cE_{/A'}}`\cE_{/Q};t_!``u_!`t'_!]

\putmorphism(1220,540)(0,1)[``]{500}{-1}a
\put(1120,260){\makebox(50,80){$f^*$}}

\put(1300,300){\makebox(50,80){$\varepsilon^f$}}
\put(1300,200){\makebox(50,80){$\Da$}}

\putmorphism(2120,540)(0,1)[``]{500}{-1}a
\put(2020,260){\makebox(50,80){$u^*$}}

\put(400,400){\makebox(50,80){$\cong$}}
\put(1100,400){\makebox(50,80){$\cong$}}
\put(1950,100){\makebox(50,80){$\cong$}}
 \end{picture}
\end{center}
where $\varepsilon^f : f_!f^*\ra 1_{\cE/A'}$ is the counit of the adjunction $f_!\dashv f^*$.
Thus we have a natural transformation
\[ t'_!(\varepsilon^f){p'_*s'^*}:t'_!f_!f^*p'_*s'^*\lra t'_!p'_*s'^* \]
and passing through the natural isomorphisms indicated in the above diagram we get the corresponding a natural transformation $rep_{pd,\cE}(f,g,u)$  as follows:
\[ u_!t_!f^*p'_*s'^*\lra t'_!p'_*s'^* \]
via right square iso, and this
\[ u_!t_!p_*g^*s'^*\lra t'_!p'_*s'^* \]
via middle square iso, and this
\[ u_!t_!p_*s^*u^*\lra t'_!p'_*s'^* \]
via left square iso. Finally, taking the adjoint ($u_!\dashv u^*$) of this morphism we get
\[ rep_{pd,\cE}(f,g,u) : t_!p_*s^*u^*\lra u^*t'_!p'_*s'^* \]
which is a morphism from $rep_{pd,\cE}(t,p,s)$ to
$rep_{pd,\cE}(t',p',s')$ in $Exp(\cE)$ over $u$. The essential image of
the functor $rep_{pd,\cE}$ is, by definition, the fibration of
polynomial (endo)functors and polynomial transformations $p_{poly,\cE}:\cP oly(\cE)\lra \cE$.

By taking the exponential adjoint to $rep_{pd,\cE}$ we obtain an
action of the fibration $p_{pd,\cE}$ on the basic fibration
$cod:\cE^\ra\lra \cE$
\begin{center} \xext=800 \yext=520
\begin{picture}(\xext,\yext)(\xoff,\yoff)
 \settriparms[1`1`1;420]
  \putVtriangle(0,0)[\cP oly\cD iag(\cE)\times_\cE\cE^\ra`\cE^\ra`\cE;\star^{poly,\cE}``cod]
\end{picture}
\end{center}
even if this point of view is less customary.

\subsection*{The case $\cE=Set$}

We shall make here the above abstract definitions concrete in case
$\cE=Set$. A functor $P:Set_{/Q}\ra Set_{/Q}$  is a {\em
polynomial functor}\footnote{There is an obvious notion of a
polynomial functor $P:Set_{/O}\ra Set_{/Q}$ with $O$ not assumed
to be equal to $Q$, but this can be considered a special case of
the above definition, as such functors are (some) polynomial
functors $Set_{/O+Q}\ra Set_{/O+Q}$. For details see \cite{GK}.}
if and only if it is isomorphic to one of form
\[ \Pi_{(t,p,s)}=t_!p_*s^*:Set_{/Q}\ra Set_{/Q} \] for some polynomial
diagram $(t,p,s)$.  Thus for $d:Y\ra Q$ in $Set_{/Q}$ we have
\[ t_!p_*s^*(Y,d) =\{ \lk b, \vec{y} \rk : b\in B,\; \vec{y}:p^{-1}(b) \ra Y,\; d\circ \vec{y}=s_{\lceil p^{-1}(b)} \} \]
where $s_{\lceil p^{-1}(b)}$ is the restriction of the function $s$ to the fibre of the function $p$ over the element $b\in B$.
It is a routine to verify that a polynomial functor is finitary if and only if it comes from a finitary diagram.

A  morphism of polynomial diagrams $(f,g,1_Q):(t,p,s)\ra
(t',p',s')$ defines a natural transformation $\Pi_{(f,g)}:
\Pi_{(t,p,s)}\lra \Pi_{(t',p',s')}$ so that for $d:Y\ra Q$ in
$Set_{/Q}$, and $\lk b, \vec{y} \rk\in t_!p_*s^*(Y,d)$, we have
\[ \Pi_{(f,g)} (\lk b, \vec{y} \rk) = \lk f(b), \vec{y}\circ (g_{\lceil p^{-1}(b)})^{-1} \rk \]
A natural transformation between polynomial functors is
{\em polynomial} if and only if it is given by the morphism of polynomial diagrams defining them.

The following Theorem is due to many authors. The precise account
of this can be found in  \cite{GK}, 1.18 and 1.19. However the proof, based on ideas of \cite{AV},
seems to be new.

\begin{theorem}  \label{poly-char}
For any set $Q$, the functor \[ \Pi_Q:(\cP oly\cD iag)_Q\lra
Nat(Set_{/Q},Set_{/Q}) \] defined above
is faithful, full on isomorphisms and its essential image consists of finitary functors preserving
wide pullbacks and cartesian natural transformations. 
\end{theorem}

A functor $F:Set_{/Q} \lra Set_{/Q}$ is {\em thin} if there is $q\in Q$ such that $F=\bi_q\circ ev_q \circ F$ and $ev_q\circ F(1)=1$.
The functors $\bi_q$ and  $ev_q$ are inclusion and evaluation functors, respectively, see Subsection \ref{an_char}.

{\it Proof.}~ First note that any functor $P: Set_{/Q}\ra Set_{/Q}$ is a coproduct of thin functors and $P$ preserves wide pullbacks and
filtered colimits if all the thin factors in the coproduct do. As $\Pi_Q$ preserves coproducts we can assume that $P$ is thin say,
$P=\bi_q\circ ev_q\circ P$ for some $q\in Q$. Thus $P_q= ev_q\circ P:Set_{/Q}\ra Set$ preserves all limits and is finitary.
Hence, by the characterization of the representable functors c.f. \cite{CWM} page 130,
it is represented by an object $s:E\ra Q$ in $Set_{/Q}$ with $E$ finite. Now it is a matter of a simple check, that with the diagram
\begin{center} \xext=1500 \yext=0
\begin{picture}(\xext,\yext)(\xoff,\yoff)
\putmorphism(500,0)(1,0)[E`1`p]{500}{1}a
\putmorphism(0,0)(1,0)[Q`\phantom{E}`s]{500}{-1}a
\putmorphism(1000,0)(1,0)[\phantom{1}`Q`t]{500}{1}a
 \end{picture}
\end{center}
in $(\cP oly\cD iag)_Q$, where $t(*)=q$, the functor $\Pi(t,p,s)$ is isomorphic to $P$.

To show that the functor $\Pi_Q$ is full and faithful it is enough to consider morphisms between diagrams with one operation and
cartesian natural transformations between corresponding thin functors, as other cartesian natural transformations between other polynomial functors
must come from those.

So suppose that we have a cartesian natural transformation between two such functors
\[ \tau : P=t_!\circ p_*\circ s^* \lra  P'=t'_!\circ p'_*\circ s'^* \]
where
\begin{center} \xext=1500 \yext=0
\begin{picture}(\xext,\yext)(\xoff,\yoff)
\putmorphism(500,0)(1,0)[E'`1`p']{500}{1}a
\putmorphism(0,0)(1,0)[Q`\phantom{E'}`s']{500}{-1}a
\putmorphism(1000,0)(1,0)[\phantom{1}`Q`t']{500}{1}a
 \end{picture}
\end{center}
As $P$ and $P'$ are thin we must have $t=t'$, say $t(*)=q\in Q$.
Thus $\tau_q=ev_q(\tau): P_q\ra P'_q$ is a cartesian natural transformation between functors that preserves the terminal object.
Hence, as any component of $\tau_q$ is a pullback of $(\tau_q)_{1_Q}$
which is a morphism from the terminal object to itself,  $\tau_q$ is cartesian if and only if it is an isomorphism.
By the Yoneda Lemma, the natural isomorphisms $\tau_q:P_q\ra P'_q$ in $Cat(Set_{/Q},Set)$ between the functors represented by
$s:E\ra Q$ and $s':E'\ra Q'$ correspond to isomorphisms
 \begin{center} \xext=600 \yext=350
\begin{picture}(\xext,\yext)(\xoff,\yoff)
\settriparms[-1`1`1;300]
  \putVtriangle(0,0)[E`E'`Q;g`s`s']
 \end{picture}
\end{center}
in $Set_{/Q}$. But those isomorphisms $g$ are exactly the functions $g$ making the diagram
\begin{center} \xext=1500 \yext=600
\begin{picture}(\xext,\yext)(\xoff,\yoff)
\setsqparms[1`1`1`1;500`400]
\putsquare(500,50)[E`1`E'`1;p`g^{-1}``p']
\setsqparms[-1`1`0`-1;500`400]
\putsquare(0,50)[Q`\phantom{E}`Q`\phantom{E'};s`1_Q``s']
\setsqparms[1`0`1`1;500`400]
\putsquare(1000,50)[\phantom{1}`Q`\phantom{1}`Q;t``1_Q`t]
 \end{picture}
\end{center}
a morphism in $(\cP oly\cD iag)_Q$, i.e.  making the left square commute and the middle square a pullback. The reader may verify that if we take
a morphism of diagrams corresponding to $\tau_q$, then $\Pi_Q$ will send it back to $\tau_q$.
$~\Box$

\vskip 2mm

{\bf Remark} An analog of  Theorem \ref{poly-char} does not hold in all locally cartesian closed categories. Even if $\cE$ is a
presheaf category, then the endofunctors on slices of $\cE$  that are finitary and preserve wide pullbacks do not necessarily
come from polynomial diagrams  in $\cE$. For $\cE=Set^{\ra}$, the functor sending
$(x:X_0\ra X_1)$ to $(\lk 1_{X_0}, x\rk :X_0\ra X_0\times X_1)$ is finitary and preserves all limits but is not polynomial.

\subsection{Some properties of  the representation  $rep_a$}

The main objective of this section is to establish some properties of the representation $rep_a$ and then
show (Corollary \ref{equiv-mona-polymonads}) that the 1-level multicategories with non-standard amalgamations
are the same as the cartesian monads on slices of $Set$ whose functor part is finitary and preserves wide pullbacks.

As some statements concerning polynomial diagrams and functors hold in greater generality, in arbitrary
locally cartesian closed categories, we start in this more general context.

First, we describe  supine and prone morphisms in $p_{pd,\cE}:\cP oly\cD iag(\cE)\ra\cE$. Let  $u:O\ra Q$ be a morphism in $\cE$.
The supine morphism $su_{u,(t,p,s)}:(t,p,s) \ra (u\circ t,p,u\circ s)$ over $u$ with the domain being the polynomial diagram  $(t,p,s)$  in $\cE$ over $O$ is defined by the diagram
\begin{center} \xext=1500 \yext=500
\begin{picture}(\xext,\yext)(\xoff,\yoff)
\setsqparms[-1`1`1`-1;500`400]
\putsquare(0,50)[O`E`Q`E;s`u``u\circ s]
\setsqparms[1`0`1`1;500`400]
\putsquare(500,50)[\phantom{E}`\phantom{B}`\phantom{E}`\phantom{B};p`1_E`1_B`p]
\setsqparms[1`0`1`1;500`400]
\putsquare(1000,50)[B`O`B`Q;t``u`u\circ t]
\end{picture}
\end{center}
i.e. $su_{u,(t,p,s)}=(1_B,1_E,u)$.
The prone morphism $pr_{u,(t,p,s)}:(\tilde{t},\tilde{p},\tilde{s})\ra (t,p,s)$ over $u$ with the domain being the polynomial diagram
$(t,p,s)$  in $\cE$ over $Q$ is defined by the diagram
\begin{center} \xext=2400 \yext=1500
\begin{picture}(\xext,\yext)(\xoff,\yoff)

\setsqparms[-1`1`1`-1;600`400]
\putsquare(0,50)[O`s^*(O)`Q`E;\bar{s}`u``s]
\setsqparms[0`0`1`1;1200`400]
\putsquare(600,50)[\phantom{s^*(O)}`\phantom{p_*s^*(O)}`\phantom{E}`\phantom{B};`\bar{u}`\hat{u}`p]
\setsqparms[0`0`1`1;600`400]
\putsquare(1800,50)[p_*s^*(O)`O`B`Q;``u`t]

\putmorphism(1100,750)(2,-1)[p^*p_*S^*(O)`\phantom{p_*s^*(O)}`]{680}{1}r
 \put(880,680){\vector(-1,-1){160}}
 \put(570,620){\makebox(200,100){$\bar{\varepsilon}^p_{s^*(O)}$}}
  \put(1330,620){\makebox(200,100){$\bar{p}$}}

  \settriparms[1`1`0;650]
 \putbtriangle(1800,450)[\tilde{B}`\phantom{p_*s^*(O)}`\phantom{O};h`\tilde{t}`]

\putmorphism(1100,1400)(0,-1)[\tilde{E}`\phantom{p^*p_*s^*(O)}`\bar{h}]{650}{1}l

\put(1250,1370){\vector(2,-1){480}}
 \put(1400,1280){\makebox(200,100){$\tilde{p}$}}

\put(920,1370){\vector(-1,-1){870}}
 \put(150,750){\makebox(200,100){$\tilde{s}$}}
\put(920,1370){\line(1,0){100}}
\put(1170,1370){\line(1,0){80}}
\end{picture}
\end{center}
i.e. $pr_{u,(t,p,s)}=(\hat{u}\circ h,\bar{u}\circ\bar{\varepsilon}^p_{s^*(O)}\circ\bar{h},u)$.
The above diagram is constructed in the following way.
First we apply the functor $t_!p_*s^*$ to the (left most) morphism $u:O\ra Q$ (in $\cE_{/Q}$) to get $t\circ \hat{u}:p_*s^*(Q)\ra Q$.
$\bar{\varepsilon}^p$ is the counit of the adjunction $p^*\dashv p_*$.
Then we pull it back along $u$ to get $\tilde{t}$ and $h$. Finally, pulling back $h$ along  $\bar{p}$ we get $\tilde{p}$ and $\bar{h}=\bar{p}^*(h)$,  the last part of the polynomial diagram $(\tilde{t},\tilde{p},\tilde{s})$. The morphism $\tilde{s}$ is defined as the composition
$\bar{s}\circ\bar{\varepsilon}^p_{s^*(O)}\circ \bar{p}^*(h)$, and the objects
$\tilde{B}$ and  $\tilde{E}$ are $u^*t_!p_*s^*(O)$ and $\bar{p}_*(\tilde{B})$, respectively.

We note for the record

\begin{proposition} \label{prone-supine}
The prone and supine morphisms in the bifibration \mbox{$p_{pd,\cE}:\cP oly\cD iag(\cE)\ra\cE$} are as described above.
\end{proposition}
{\it Proof.}~  Routine verification.
 $~\Box$

The following Lemma collects some known facts that will be used in Proposition \ref{rep-E1}.

\begin{lemma} \label{folk lemma}
Let $u$, $u'$, $s$ and $p$ be morphisms in a locally cartesian closed category $\cE$, such that $cod(u')=dom(u)$, $cod(u)=cod(s)=dom(p)$. Then
\begin{enumerate}
  \item $s^*(\varepsilon^u)\cong (\varepsilon^{s^*(u)})s^*$;
  \item $p_*(\varepsilon^u)\cong (\varepsilon^{p_*(u)})p_*$;
  \item $\varepsilon^{u\circ u'}\cong \varepsilon^u\circ(u_!(\varepsilon^{u'})u^*)$.$~\Box$
\end{enumerate}
\end{lemma}

\begin{proposition} \label{rep-E1} Let $\cE$ be a locally cartesian closed category.
The morphism
\begin{center} \xext=800 \yext=550
\begin{picture}(\xext,\yext)(\xoff,\yoff)
 \settriparms[1`1`1;450]
  \putVtriangle(0,0)[\cP oly\cD iag(\cE)`Exp(\cE)`\cE;rep_{pd,\cE}`p_{pd,\cE}`p_{exp}]
\end{picture}
\end{center}
defined in the previous section is  a morphism of bifibrations.
\end{proposition}
{\it Proof.}~ We need to show that $rep_{pd,\cE}$ preserves prone and supine morphisms. We show preservation of supine morphisms first.

Let $su_{u,(t,p,s)}=(1_B,1_E,u): (t,p,s) \ra (u\circ t,p,u\circ s)$ be a supine morphism in $\cP oly \cD iag(\cE)$.
The representation of the supine morphism  $rep_{pd,\cE}(su_{u,(t,p,s)})$ is a natural transformation which is isomorphic to the adjoint ($u_!\dashv u^*$) to
 \begin{center} \xext=1600 \yext=100
\begin{picture}(\xext,\yext)(\xoff,\yoff)
\putmorphism(0,0)(1,0)[(ut)_!(1_B)_!1_B^*p_*(us)^*`(ut)_!p_*(us)^*`(ut)_!(\varepsilon^{1_B})p_*(us)^*]{1600}{1}a
\end{picture}
\end{center}
The above morphism is isomorphic to the identity natural transformation
 \begin{center} \xext=1600 \yext=100
\begin{picture}(\xext,\yext)(\xoff,\yoff)
\putmorphism(0,0)(1,0)[u_!t_!p_*s^*u^*`u_!t_!p_*s^*u^*`1_{u_!t_!p_*s^*u^*}]{1600}{1}a
\end{picture}
\end{center}
and the supine morphism with the codomain $rep_{pd,\cE}(t,p,s)$ over $u$ is a natural transformation
 \begin{center} \xext=1600 \yext=100
\begin{picture}(\xext,\yext)(\xoff,\yoff)
\putmorphism(0,0)(1,0)[(t_!p_*s^*)u^*`u^*u_!(t_!p_*s^*)u^*`\eta_{(t_!p_*s^*)u^*}]{1600}{1}a
\end{picture}
\end{center}
adjoint to the above identity morphism. Thus  $rep_{pd,\cE}$ preserves the supine morphisms.

For the prone morphism we use the notation introduced at the beginning of the section. Let
$pr_{u,(t,p,s)}=(\hat{u}\circ h,\bar{u}\circ\bar{\varepsilon}^p_{s^*(O)}\circ\bar{p}^*(h),u):(\tilde{t},\tilde{p},\tilde{s})\ra (t,p,s)$
be a prone morphism.  We have a diagram of categories, functors and natural transformations
\begin{center} \xext=3000 \yext=1650
\begin{picture}(\xext,\yext)(\xoff,\yoff)

\setsqparms[1`0`0`1;800`600]
\putsquare(100,50)[\cE_{/O}`\cE_{/s^*(O)}`\cE_{/Q}`\cE_{/E};\bar{s}^*```s^*]

\putmorphism(0,650)(0,-1)[\phantom{\cE_{/Q}}`\phantom{\cE_{/O}}`u^*]{600}{-1}l
\putmorphism(200,650)(0,-1)[\phantom{\cE_{/Q}}`\phantom{\cE_{/O}}`u_!]{600}{1}r

 \put(0,350){\makebox(200,100){$\varepsilon^u$}}
  \put(0,250){\makebox(200,100){$\Da$}}


\putmorphism(900,50)(1,0)[\phantom{\cE_{/E}}`\phantom{\cE_{/B}}`p_*]{1200}{1}b

\putmorphism(800,650)(0,-1)[\phantom{\cE_{/Q}}`\phantom{\cE_{/O}}`\bar{u}^*]{600}{-1}l
\putmorphism(1000,650)(0,-1)[\phantom{\cE_{/Q}}`\phantom{\cE_{/O}}`\bar{u}_!]{600}{1}r

 \put(800,350){\makebox(200,100){$\varepsilon^{\bar{u}}$}}
  \put(800,250){\makebox(200,100){$\Da$}}

\putmorphism(2000,650)(0,-1)[\phantom{\cE_{/Q}}`\phantom{\cE_{/O}}`\hat{u}^*]{600}{-1}l
\putmorphism(2200,650)(0,-1)[\phantom{\cE_{/Q}}`\phantom{\cE_{/O}}`\hat{u}_!]{600}{1}r

 \put(2000,350){\makebox(200,100){$\varepsilon^{\hat{u}}$}}
  \put(2000,250){\makebox(200,100){$\Da$}}

\setsqparms[0`0`0`-1;800`600]
\putsquare(2100,50)[\cE_{/p_*s^*(O)}`\cE_{/O}`\cE_{/B}`\cE_{/Q};```t_!]

\putmorphism(2800,650)(0,-1)[\phantom{\cE_{/Q}}`\phantom{\cE_{/O}}`u^*]{600}{-1}l
\putmorphism(3000,650)(0,-1)[\phantom{\cE_{/Q}}`\phantom{\cE_{/O}}`u_!]{600}{1}r

 \put(2800,350){\makebox(200,100){$\varepsilon^u$}}
  \put(2800,250){\makebox(200,100){$\Da$}}

\putmorphism(1350,1000)(2,-1)[\cE_{/p^*p_*S^*(O)}`\phantom{\cE_{/p_*s^*(O)}}`]{680}{1}r
 \put(980,710){\vector(1,1){200}}
 \put(800,780){\makebox(200,100){$(\varepsilon^p_{s^*(O)})^*$}}
  \put(1450,720){\makebox(200,100){$\bar{p}_*$}}

\putmorphism(2000,1450)(0,-1)[\phantom{\cE_{/Q}}`\phantom{\cE_{/p_*s^*(O)}}`h^*]{800}{-1}l
\putmorphism(2200,1450)(0,-1)[\phantom{\cE_{/Q}}`\phantom{\cE_{/p_*s^*(O)}}`h_!]{800}{1}r

 \put(2000,1150){\makebox(200,100){$\varepsilon^{h}$}}
  \put(2000,1050){\makebox(200,100){$\Da$}}

\settriparms[0`1`0;850]
\putbtriangle(2100,700)[\cE_{/\tilde{B}}`\phantom{\cE_{/p_*s^*(O)}}`\phantom{\cE_{/O}};`\tilde{t}_!`]

\putmorphism(1330,1520)(0,-1)[\phantom{\cE_{/\tilde{E}}}`\phantom{\cE_{/p^*p_*s^*(O)}}`\bar{h}^*]{550}{-1}l

\putmorphism(1330,1550)(1,0)[\cE_{/\tilde{E}}`\phantom{\cE_{/\tilde{B}}}`\tilde{p}_*]{750}{1}a
\put(150,750){\vector(4,3){1050}}
 \put(570,1150){\makebox(200,100){$\tilde{s}^*$}}
\end{picture}
\end{center}
Intuitively speaking, the adjoint ($u_!\dashv u^*$) natural transformation  to $rep_{pd,\cE}(pr_{u,(t,p,s)})$ is isomorphic to
$t_!(\varepsilon^{\hat{u}h})p_*s^*   =t_!(\varepsilon^{\hat{u}}\circ(\hat{u}_!(\varepsilon^{h})\hat{u}^*))p_*s^*$
i.e. it is defined with the help of the counits $\varepsilon^{\hat{u}}$
and $\varepsilon^h$.  The adjoint to the natural transformation $pr_{u,rep_{pd,\cE}(t,p,s)}$, being the prone morphism in $Exp(Set)$ over $u$
with the codomain $rep_{pd,\cE}(t,p,s)$, is $((\varepsilon^u)t_!p_*s^*)\circ (u_!u^*t_!p_*s^*(\varepsilon^u))$.
To show that these adjoints are isomorphic, we show, using the above Lemma \ref{folk lemma},
that the counit $\varepsilon^{\hat{u}}$ can be 'moved left' to the 'left' counit $\varepsilon^u$ and the counit $\varepsilon^h$ can be 'moved right' to
the 'right' counit $\varepsilon^u$.

In the sequence of morphisms below, we mark on the right side of the line how we pass from a line to another.
Numbers 1. 2. 3. refer
to Lemma \ref{folk lemma}, MEL is the middle exchange law.  We have
\begin{center} \xext=3200 \yext=3100
\begin{picture}(\xext,\yext)(\xoff,\yoff)
\putmorphism(500,3000)(1,0)[ rep_{pd,\cE}(\tilde{t},\tilde{p},\tilde{s})`rep_{pd,\cE}(t,p,s)`rep_{pd,\cE}(pr_{u,(t,p,s)})]{2000}{1}a

\put(0,2930){\line(1,0){3000}}
\put(3020,2890){\makebox(300,100){def pr}}

\putmorphism(500,2750)(1,0)[ \tilde{t}_!\tilde{p}_*\tilde{s}^*u^*`u^*t_!p_*s^*`rep_{pd,\cE}(\hat{u}h,\bar{u}\varepsilon^p_{S^*(O)}\bar{h},u)]{2000}{1}a

\put(0,2660){\line(1,0){3000}}
\put(3020,2620){\makebox(300,100){$u_!\dashv u^*$}}

\putmorphism(500,2500)(1,0)[ u_!\tilde{t}_!\tilde{p}_*\tilde{s}^*u^*`t_!p_*s^*`]{2000}{1}a

\put(0,2410){\line(1,0){3000}}
\put(3020,2370){\makebox(300,100){def rep}}

\putmorphism(500,2250)(1,0)[ t_!(\hat{u}h)_!(\hat{u}h)^* p_*s^*`t_!p_*s^*`t_!(\varepsilon^{\hat{u}h})p_*s^*]{2000}{1}a

\put(0,2160){\line(1,0){3000}}
\put(3020,2120){\makebox(200,100){3.}}

\putmorphism(500,2000)(1,0)[ t_!\hat{u}_!h_!h^*\hat{u}^* p_*s^*`t_!p_*s^*`t_!(\varepsilon^{\hat{u}}\circ(\hat{u}_!(\varepsilon^{h})\hat{u}^*))p_*s^*]{2000}{1}a

\put(0,1910){\line(1,0){3000}}
\put(3020,1870){\makebox(200,100){=}}

\putmorphism(500,1750)(1,0)[t_!\hat{u}_!h_!h^*\hat{u}^*p_*s^*`t_!p_*s^*`(t_!(\varepsilon^{\hat{u}})p_*s^*)\circ(t_!\hat{u}_!(\varepsilon^{h})\hat{u}^*p_*s^*)]{2000}{1}a

\put(0,1660){\line(1,0){3000}}
\put(3020,1620){\makebox(200,100){1.}}

\putmorphism(500,1500)(1,0)[ u_!u^*t_!\hat{u}_!\hat{u}^* p_*s^*`t_!p_*s^*`
(t_!(\varepsilon^{\hat{u}})p_*s^*)\circ((\varepsilon^{u})t_!\hat{u}_!\hat{u}^*p_*s^*)]{2000}{1}a

\put(0,1410){\line(1,0){3000}}
\put(3020,1370){\makebox(200,100){MEL}}

\putmorphism(500,1250)(1,0)[ u_!u^*t_!\hat{u}_!\hat{u}^* p_*s^*`t_!p_*s^*`((\varepsilon^{u})t_!p_*s^*)\circ(u_!u^*t_!(\varepsilon^{\hat{u}})p_*s^*)]{2000}{1}a

\put(0,1160){\line(1,0){3000}}
\put(3020,1120){\makebox(200,100){2.}}

\putmorphism(500,1000)(1,0)[ u_!u^*t_!p_*\bar{u}_!\bar{u}^* s^*`t_!p_*s^*`((\varepsilon^{u})t_!p_*s^*)\circ(u_!u^*t_!p_*(\varepsilon^{\bar{u}})s^*)]{2000}{1}a

\put(0,910){\line(1,0){3000}}
\put(3020,870){\makebox(200,100){1.}}

\putmorphism(500,750)(1,0)[ u_!u^*t_!p_*s^*u_!u^* `t_!p_*s^*`((\varepsilon^{u})t_!p_*s^*)\circ(u_!u^*t_!p_*s^*(\varepsilon^u))]{2000}{1}a

\put(0,660){\line(1,0){3020}}
\put(3030,620){\makebox(300,100){$u_!\dashv u^*$}}

\putmorphism(500,500)(1,0)[ u^*t_!p_*s^*u^*u_! `u^*t_!p_*s^*`u^*t_!p_*s^*(\varepsilon^u)]{2000}{1}a

\put(0,410){\line(1,0){3000}}
\put(3020,370){\makebox(300,100){def rep}}

\putmorphism(500,250)(1,0)[ u^*rep_{pd,\cE}(t,p,s)u_!u^*`u^*rep_{pd,\cE}(t,p,s)`u^*rep_{pd,\cE}(t,p,s)(\varepsilon^u)]{2000}{1}a

\put(0,160){\line(1,0){3000}}
\put(3000,120){\makebox(300,100){def pr}}

\putmorphism(500,0)(1,0)[ u^*rep_{pd,\cE}(t,p,s)u_!`rep_{pd,\cE}(t,p,s)`pr_{u,rep_{pd,\cE}(t,p,s)}]{2000}{1}a
\end{picture}
\end{center}
Thus  $rep_{pd,\cE}$ preserves the prone morphisms, as well.
 $~\Box$

From now on till the end of this subsection we shall consider fibrations over $Set$ only.
The following appears in \cite{GK}. It is an immediate consequence of  Theorem \ref{poly-char} and Proposition \ref{rep-E1}.

\begin{proposition} \label{rep-E2}
The representation morphism $rep_{pd,Set}$ is a morphism of lax monoidal fibrations
which is faithful and full on isomorphisms. $~\Box$
\end{proposition}

\begin{corollary} \label{lmf-for-E}
We have a sequence of morphisms of lax monoidal fibrations
 \begin{center} \xext=2200 \yext=100
\begin{picture}(\xext,\yext)(\xoff,\yoff)
\putmorphism(0,0)(1,0)[\cP oly\cD iag`\cP oly`rep_{pd,Set}]{1000}{1}a
\putmorphism(1000,0)(1,0)[\phantom{\cP oly}`Cart(Set)`]{600}{1}a
\putmorphism(1600,0)(1,0)[\phantom{Cart(Set)}`Exp(Set)`]{600}{1}a
\end{picture}
\end{center}
with the first being an equivalence of bifibrations and the following two being
inclusions full on isomorphisms. The composition of these morphisms is (isomorphic to) $rep_{pd,Set}$.
\end{corollary}
{\it Proof.}~ This is an immediate consequence of Propositions \ref{rep-E1} and \ref{rep-E2}.   $~\Box$

\vskip 2mm

We shall construct a morphism of lax monoidal fibrations
\begin{center} \xext=700 \yext=400
\begin{picture}(\xext,\yext)(\xoff,\yoff)
 \settriparms[1`1`1;350]
  \putVtriangle(0,0)[Sig_a`\cP oly\cD iag`Set;\iota_a`p_a`p_{pd}]
\end{picture}
\end{center}
Let $(A,\partial,O)$ be  a signature in $(Sig_a)_O$. The functor $\iota_a$ sends this signature to a polynomial functor
defined by the following polynomial diagram
\begin{center} \xext=1500 \yext=150
\begin{picture}(\xext,\yext)(\xoff,\yoff)
\putmorphism(500,0)(1,0)[E^A`A`p^A]{500}{1}a
\putmorphism(0,0)(1,0)[O`\phantom{E^A}`s^A]{500}{-1}a
\putmorphism(1000,0)(1,0)[\phantom{A}`O`t^A]{500}{1}a
 \end{picture}
\end{center}
where
\[ E^A=\coprod_{a\in A} (|a|] =\{ \lk a,i \rk : a\in A,\; i\in (|a|]\} \]
$p^A$ is the first projection, $s^A(a,i)=\partial^A_a(i)$, and $t^A(a)=\partial^A_a(0)$, for $a\in A$ and $i\in (|a|]$.

Moreover, for the morphisms of signatures $(f,\sigma,u):(A,\partial^A,O)\lra (B,\partial^B,Q)$ in $Sig_a$ over $u:O\ra Q$
we define a commuting diagram
\begin{center} \xext=1500 \yext=600
\begin{picture}(\xext,\yext)(\xoff,\yoff)
\setsqparms[1`1`1`1;500`400]
\putsquare(500,50)[E^A`A`E^B`B;p^A`g`f`p^B]
\setsqparms[-1`1`0`-1;500`400]
\putsquare(0,50)[O`\phantom{B}`Q`\phantom{B'};s^A`u``s^B]
\setsqparms[1`0`1`1;500`400]
\putsquare(1000,50)[\phantom{B}`O`\phantom{B'}`Q;t^A``u`t^B]
 \end{picture}
\end{center}
where $g(a,i)=\lk f(a),\sigma_a^{-1}(i) \rk$, for $\lk a,i \rk\in E^A$.
The square in the middle is easily seen to be a pullback. Thus the above diagram is a morphism of polynomial diagrams
\[ (f,g,u):(t^A,p^A,s^A)\lra (t^B,p^B,s^B) \]
in $\cP oly\cD iag$.

We have

\begin{proposition} \label{iota-comp}
The morphism $\iota_a$ defined above is a morphism of lax monoidal fibrations, an equivalence of bifibrations,
and it makes the triangle of morphisms of lax monoidal fibrations over $Set$
 \begin{center} \xext=800 \yext=450
\begin{picture}(\xext,\yext)(\xoff,\yoff)
\settriparms[1`1`1;400]
  \putVtriangle(0,0)[Sig_a`\cP oly\cD iag `Exp(Set);\iota_a`rep_a`rep_{pd}]
 \end{picture}
\end{center}
commute, up to a fibred natural isomorphism.
\end{proposition}

{\it Proof.}~ $\iota_a$ is faithful from the construction. Let $(f,g,u):(t^A,p^A,s^A)\ra (t^B,p^B,s^B)$ be a morphism of polynomial diagrams
over $u:O\ra Q$. Then  $(f,\sigma,u):(A,\partial^A,O)\ra (B,\partial^B,Q)$ such that $\sigma_a=(g_{\lc p^{-1}(a)})^{-1}$, for $a\in A$,
is a morphism in $Sig_a$ over $u$. Moreover $\iota_a(f,\sigma,u)=(f,g,u)$. Thus $\iota_a$ is full.

To see that $\iota_a$ is essentially surjective as well, fix a diagram
\begin{center} \xext=1200 \yext=150
\begin{picture}(\xext,\yext)(\xoff,\yoff)
\putmorphism(400,0)(1,0)[E`A`p]{400}{1}a
\putmorphism(0,0)(1,0)[O`\phantom{E}`s]{400}{-1}a
\putmorphism(800,0)(1,0)[\phantom{A}`O`t]{400}{1}a
 \end{picture}
\end{center}
in $\cP oly \cD iag$. For $a\in A$,  choose  bijections $\tau_a:(n_a] \ra p^{-1}(a)$, for some $n_a=|p^{-1}(a)|\in \o$.
Putting
\[ \partial^A_a(i)\left\{ \begin{array}{ll}
                t(a) & \mbox{ if  $i=0$} \\
                s\tau_a(i)    & \mbox{ otherwise. }
                                    \end{array}
                            \right.  \]
we have that $\iota(A,\partial^A,O)$ is isomorphic to $(t,p,s)$, i.e. $\iota_a$ is essentially surjective as well.

The verification that the triangle commutes is also easy and we leave it to the reader.  $~\Box$

\begin{corollary} \label{lmf-for-Set}
We have a sequence of morphisms of lax monoidal fibrations
 \begin{center} \xext=2600 \yext=100
\begin{picture}(\xext,\yext)(\xoff,\yoff)
\putmorphism(0,0)(1,0)[Sig_a`\cP oly\cD iag`\iota_a]{600}{1}a
\putmorphism(600,0)(1,0)[\phantom{\cP oly\cD iag}`\cP oly`rep_{pd}]{800}{1}a
\putmorphism(1400,0)(1,0)[\phantom{\cP oly}`Cart(Set)`]{600}{1}a
\putmorphism(2000,0)(1,0)[\phantom{Cart(Set)}`Exp(Set)`]{600}{1}a
\end{picture}
\end{center}
with the first two being equivalences of bifibrations and the following two being
inclusions of bifibrations full on isomorphisms. The composition of all four morphisms is (isomorphic to) $rep_a$, and hence $rep_a$
is a morphism lax monoidal fibrations, morphism of bifibrations, faithful, and full on isomorphisms.
\end{corollary}
{\it Proof.}~ This is an immediate consequence of  Corollary \ref{lmf-for-E} and Proposition \ref{iota-comp}.  $~\Box$

\vskip 2mm

\begin{theorem}\label{im-rep-a}
The essential image of the morphism of lax monoidal fibrations $rep_a$ in $Exp(Set)\ra Set$
consists of the finitary endofunctors preserving wide pullbacks and cartesian natural transformations.
\end{theorem}

{\it Proof.}~ By Proposition \ref{iota-comp} the image of $rep_a$ in $p_{exp}$ is $p_{poly}:\cP oly\ra Set$, the image of $rep_{pd}$.
By Proposition \ref{rep-E1}, $rep_{pd}$ preserves prone morphisms.
By Corollary \ref{lmf-for-E}, the image of $rep_{pd}$ is contained in the lax monoidal subfibration $Cart(Set)\ra Set$.
Thus, it is enough to verify the statement on fibres. But this is the content of Theorem \ref{poly-char}.
$~\Box$

\vskip 2mm

By Proposition \ref{mor-lmf-induce}, any morphism of lax monoidal fibration induces a morphism of the corresponding fibrations of monoids.
We finish this section by spelling the most important instance of this fact, announced at the beginning of this subsection,
that follows from Corollary \ref{lmf-for-Set}.

\begin{corollary}\label{equiv-mona-polymonads}
The fibration of multicategories with non-standard amalgamation is
equivalent to the fibration that has as objects finitary cartesian
monads on slices of $Set$ whose functor part preserves wide
pullbacks. A morphism in that fibration  between monads over $Set_{/O}$
and $Set_{/Q}$ over a function $u:O\ra Q$ is a cartesian morphism of monads
whose functor part is the pullback functor $u^*:Set_{/Q}\lra Set_{/O}$.
\end{corollary}

\vskip 2mm

{\bf Example}  The following example shows that some polynomial diagrams can be equipped with a monoid structure in the fibration
of polynomial diagrams $\cP oly \cD iag$ but that this monoid structure (unique in this case)  cannot be lifted to the fibration
of monotone diagrams $\cM\cP oly \cD iag$. This is to explain why there are fewer monotone monads than polynomial ones.
The signature has three types $O=\{ x_0,x_1,x_2\}$ and seven operations. Three operations that will serve as units in the monoid
$1_{x_0}: x_0 \ra x_0$, $1_{x_1}: x_1 \ra x_1$, $1_{x_2}: x_2 \ra x_2$ and four other with typing as displayed:
  \begin{center} \xext=1650 \yext=600
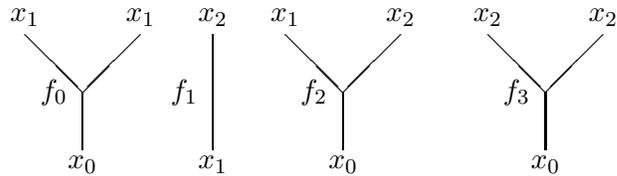
\begin{picture}(\xext,\yext)(\xoff,\yoff)

 \put(0,510){\makebox(100,100){$x_1$}}
  \put(400,510){\makebox(100,100){$x_1$}}
   \put(100,250){\makebox(100,100){$f_0$}}
\put(450,500){\line(-1,-1){200}}
\put(50,500){\line(1,-1){200}}
\put(250,300){\line(0,-1){200}}
 \put(200,0){\makebox(100,100){$x_0$}}

  \put(650,510){\makebox(100,100){$x_2$}}
   \put(550,250){\makebox(100,100){$f_1$}}
\put(700,500){\line(0,-1){400}}
 \put(650,0){\makebox(100,100){$x_1$}}

 \put(900,510){\makebox(100,100){$x_1$}}
  \put(1300,510){\makebox(100,100){$x_2$}}
   \put(1000,250){\makebox(100,100){$f_2$}}
\put(1350,500){\line(-1,-1){200}}
\put(950,500){\line(1,-1){200}}
\put(1150,300){\line(0,-1){200}}
 \put(1100,0){\makebox(100,100){$x_0$}}

 \put(1600,510){\makebox(100,100){$x_2$}}
  \put(2000,510){\makebox(100,100){$x_2$}}
   \put(1700,250){\makebox(100,100){$f_3$}}
\put(2050,500){\line(-1,-1){200}}
\put(1650,500){\line(1,-1){200}}
\put(1850,300){\line(0,-1){200}}
 \put(1800,0){\makebox(100,100){$x_0$}}

 \end{picture}
\end{center}
Then, no matter how we order entries in $f_2$: either $x_1<x_2$ or $x_2<x_1$, we won't be able to define one of the multiplications
 \begin{center} \xext=1650 \yext=800
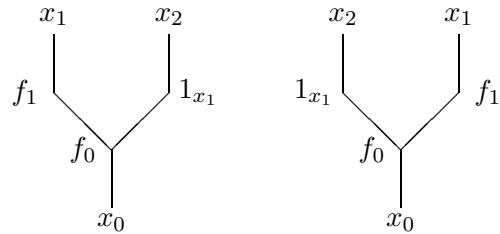
\begin{picture}(\xext,\yext)(\xoff,\yoff)
 \put(100,710){\makebox(100,100){$x_1$}}
  \put(500,710){\makebox(100,100){$x_2$}}
 \put(0,450){\makebox(100,100){$f_1$}}
 \put(600,450){\makebox(100,100){$1_{x_1}$}}
   \put(200,250){\makebox(100,100){$f_0$}}
\put(550,500){\line(-1,-1){200}}
\put(150,500){\line(1,-1){200}}
\put(150,700){\line(0,-1){200}}
\put(550,700){\line(0,-1){200}}
\put(350,300){\line(0,-1){200}}
 \put(300,0){\makebox(100,100){$x_0$}}

 \put(1100,710){\makebox(100,100){$x_2$}}
  \put(1500,710){\makebox(100,100){$x_1$}}
 \put(1000,450){\makebox(100,100){$1_{x_1}$}}
 \put(1600,450){\makebox(100,100){$f_1$}}

   \put(1200,250){\makebox(100,100){$f_0$}}
\put(1550,500){\line(-1,-1){200}}
\put(1150,500){\line(1,-1){200}}
\put(1150,700){\line(0,-1){200}}
\put(1550,700){\line(0,-1){200}}
\put(1350,300){\line(0,-1){200}}
 \put(1300,0){\makebox(100,100){$x_0$}}
 \end{picture}
\end{center}
This problem disappears if we can switch entries in the result, which is possible in amalgamated signatures and polynomial diagrams.

\subsection{The 2-level amalgamated signatures fibration $p_{2a}: Sig_{2a} \ra Set^\ra$}

We describe below a lax monoidal fibration such that
its category of monoids contains as a full subcategory the category
of 2-level multicategories with non-standard amalgamation, c.f. \cite{HMP}. It is fairly clear that this
construction can be farther generalize by making the structure of objects even more
involved but right now we don't see any real applications for such structures.

\subsection*{The fibration $p_{2a}: Sig_{2a} \ra Set^\ra$}

The base category of our fibration is $Set^{\ra}$. A typical object of
$Set^{\ra}$ is a function $\dot{(-)}:O\ra\ddot{O}$, denoted by $\vec{O}$. 
$O$ is referred to as the set of objects of $\vec{O}$, $\ddot{O}$ is referred to as the set of types of $\vec{O}$, and $\dot{(-)}$ is the typing of $\vec{O}$.
A morphism in $Set^{\ra}$, denoted  $\vec{u}=(u,\ddot{u}): \vec{O}\ra \vec{Q}$, is a pair of morphism making the square
\begin{center} \xext=500 \yext=450
\begin{picture}(\xext,\yext)(\xoff,\yoff)
   \setsqparms[1`1`1`1;500`400]
  \putsquare(0,0)[O`Q`\ddot{O}`\ddot{Q};u`\dot{(-)}`\dot{(-)}`\ddot{u}]
\end{picture}
\end{center}
commute. The notation  $\vec{O}$ and $\vec{u}$ will be used exclusively in this subsection.
In spite of the fact that we think of $O$ and $\ddot{O}$ as
disjoint sets, it is convenient to 'test' elements of those sets
for equality in the sense that either they both belong to one set
and are equal or otherwise we move one of the elements to
$\ddot{O}$ and there they are equal. Formally, we define the
`graded equality' $\dot{=}$ so that if $x,y\in O+\ddot{O}$, then
\[ x\dot{=}y\;\; \mbox{iff}\;\;
\left\{ \begin{array}{ll}
                x=y & \mbox{ if  $x,y\in O$ or $x,y\in \ddot{O}$} \\
                x=\dot{y}    & \mbox{ if  $y\in O$ and $x\in \ddot{O}$ }\\
                \dot{x}=y    & \mbox{ if  $x\in O$ and $y\in \ddot{O}$ }\\
                                    \end{array}
                            \right.  \]

By $\vec{O}^\dag$ we denote the sum $\bigcup_{n\in \o} \vec{O}^\dag_n$, where
 \[\vec{O}^\dag_n=\{d:[n]\ra O+\ddot{O} : d((n])\subseteq O    \}\]
is the set of functions from $[n]$ to the disjoint sum $O+\ddot{O}$ such that
positive integers are sent to objects and $0$ is sent either to an object or a type.
Extending the previous conventions we will write
If $\dot{(-)}:O\ra \ddot{O}$ is an identity, we write $O^\dag$ for $\vec{O}^\dag$.
For $d:[n]\ra O+\ddot{O}$ we have restrictions of $d$ to $d^+:(n]^+\ra O+\ddot{O}$ and to $d^-: \{ 0\}\ra O+\ddot{O}$.

The total category of our fibration $p_{2a}: Sig_{2a} \ra Set^\ra$ has as
objects triples $(A,\partial^A,\vec{O})$, such that $A$ is a
set, $\vec{O}$ is an object of $Set^\ra$ and $\partial^A:A\ra
\vec{O}^\dag$ is a function. A morphism
$(f,\sigma,\vec{u}):(A,\partial,\vec{O})\ra(B,\partial,\vec{O})$
is a triple such that $f:A\ra B$ is a function,
$\vec{u}:\vec{O}\ra \vec{Q}$ is a morphism in $Set^{\ra}$
and for any $a\in A$ with $|a|=n$, $\sigma \in S_n$ is a
permutation such that
\begin{center}
\xext=700 \yext=500
\begin{picture}(\xext,\yext)(\xoff,\yoff)
 \setsqparms[-1`1`1`1;700`400]
 \putsquare(0,50)[[n]`[n]`O+\ddot{O}`Q+\ddot{Q};\sigma_a`\partial_a`\partial_{f(a)}`u+\ddot{u}]
\end{picture}
\end{center}
commutes. The projection functor $p:Sig_{2a}\ra Set^{\ra}$ sends
$(f,\sigma,\vec{u}):(A,\partial,\vec{O})\ra(B,\partial,\vec{Q})$
to $\vec{u}:\vec{O}\ra \vec{Q}$.

For $\vec{u}:\vec{O}\ra \vec{Q}$ in $Set^\ra$,  we have a pullback operation
\begin{center} \xext=600 \yext=500
\begin{picture}(\xext,\yext)(\xoff,\yoff)
   \setsqparms[1`1`1`1;600`400]
  \putsquare(0,50)[u^*(B)`B`\vec{O}^\dag`\vec{Q}^\dag;\pi`\partial`\partial`\vec{u}^\dag]
\end{picture}
\end{center}
where $\vec{u}^\dag(d)=(u+\ddot{u})\circ d$, where
$d:[n]\ra O+\ddot{O}$. Then
\[ \vec{u}^*(B) = \{ \lk d,b\rk  :  d:[|b|]\ra O+\ddot{O},\; b\in B,\; \vec{u}^\dag(d)=\partial_b  \} \]

\subsection*{The  monoidal structure in $p_{2a}$}

We have two lax morphisms of fibrations
\begin{center} \xext=1800 \yext=480
\begin{picture}(\xext,\yext)(\xoff,\yoff)
\settriparms[1`1`1;400]
 \putVtriangle(0,0)[Sig_{2a}\times_{Set^\ra}Sig_{2a}`Sig_{2a}`Set^\ra;\otimes`p'`p_{2a}]
\settriparms[1`1`1;400]
 \putVtriangle(1200,0)[Set^\ra`Sig_{2a}`Set^\ra;I`id`p_{2a}]
\end{picture}
\end{center}
The tensor $(A\otimes_{\vec{O}} B,\partial^\otimes,\vec{O})$ of two objects
$(A,\partial,\vec{O})$ and $(B,\partial,\vec{O})$ in a fibre $(Sig_{2a})_{\vec{O}}$  is defined as follows
\[A\otimes_{\vec{O}} B=\{ \lk a,b_i \rk_{i\in (|a|]} : a\in A,\, b_i\in B,\,
\partial_{b_i}(0)\dot{=}\partial_a(i) \}  \]
and for $\lk a,b_i\rk_{i\in (|a|]}\in A\otimes_O B$,
\[ \partial^\otimes(\lk a,b_i\rk _{i})=
[ \dot{\partial}^-_{a}, \partial^+_{b_i}]_{i} \, :\, [\,\sum_{i=1}^{|a|}
|b_i| \,] = [0]+\coprod_{i} (|b_i|] \lra O+\ddot{O}\]
where $\dot{\partial}^-_a(*)=\partial^\otimes(\lk a,b_i\rk_{i})(*)$ is `$\partial_a(*)$ as much as possible', i.e.
\[\dot{\partial}^+_a(*)=
\left\{ \begin{array}{ll}
                \dot{(\partial_a(*))} & \mbox{ if  $\partial_a(*)\in O$ and
                $\exists_{i\in (|a|]}\; \partial_{b_i}(*)\in\ddot{O}$ } \\
               \partial_a(*)    & \mbox{ otherwise.}
                                    \end{array}
                            \right.  \]
{\bf Remark} This definition is slightly more complicated than possible to make sure that composition with identities
(that have objects from $O$, rather than types from $\ddot{O}$, as codomains) is neutral.

\vskip 2mm

The tensor on morphisms is defined as in $Sig_a$. For a pair of maps in $Sig_{2a}$
\[ \overline{f}=(f,\sigma,\vec{u}):(A,\partial,\vec{O})\ra(A',\partial,\vec{Q}),\;\;\; \overline{g}=(g,\tau, \vec{u}):(B,\partial,\vec{O})\ra(B',\partial,\vec{Q})\]
over the same map $\vec{u}:\vec{O}\ra \vec{Q}$ in $Set^\ra$ we define the map
\[\overline{f}\otimes_{\vec{u}}\overline{g}=(f\otimes_{\vec{O}}
g,\sigma\otimes_{\vec{u}}\tau,\vec{u}):(A\otimes_{\vec{O}}
B,\partial^\otimes,\vec{O})\lra (A'\otimes_{\vec{Q}} B',\partial^\otimes,\vec{Q}) \]
 so that, for $\lk a, b_i\rk_{i\in (|a|]}\in A\otimes_{\vec{O}} B$,
\[ f\otimes_{\vec{u}} g(\lk a,b_i\rk_{i\in (|a|]}) =
(\lk f(a),g(b_{\sigma_a(j)})\rk_{j\in (|f(a)|]})\]
Clearly, $|a|=|f(a)|$ and
\[ n=|\lk a, b_i\rk_{i\in (|a|]}|= \sum_{i\in |a|} |b_i| = \sum_{i\in |f(a)|} |g(b_i)| = |f\otimes_{\vec{u}} g(\lk a,b_i\rk_{i\in (|a|]})|. \]
Moreover, we put
 \[ (\sigma\otimes_{\vec{u}}\tau)_{\lk a, b_i \rk_i}=
 [\sigma^-_{a},\tau^+_{b_{\sigma_a(i)}}]_{i} \]
making the triangle
\begin{center} \xext=1200 \yext=650
\begin{picture}(\xext,\yext)(\xoff,\yoff)
   \setsqparms[-1`1`1`1;1200`500]
  \putsquare(0,50)[ [n]`[n]`O`Q;(\sigma\otimes_{\vec{u}}\tau)_{\lk
a,b_i\rk_i}`\partial^\otimes_{\lk a,b_i\rk_i}`
 \partial^\otimes_{\lk f(a),g(b_i)\rk_i}`\vec{u}]
\end{picture}
\end{center}
commute. This ends the definition of the tensor $\otimes$.

The unit $I_{\vec{O}}$ in the fibre $(Sig_{2a})_{\vec{O}}$, is $(O,\partial^{I_O},{\vec{O}})$
where, for $x\in O$, the function $\partial^{I_{\vec{O}}}_x:[1]\ra O$ is constant equal $x$ .

{\bf Remark} As we mentioned earlier, we would like to put
$\partial^{{\vec{O}}}_x(*)=\dot{x}$ to be sure that
the codomains are always the types of $\vec{O}$
 but this would not work as the compatibility of the
identities on the left in 2-level multicategories with
non-standard amalgamation must be on the objects of $\vec{O}$.

\begin{lemma} The fibration $p_{2a}:Sig_{2a}\lra Set^\ra$ together with the above
defined tensor $\otimes$ and unit $I$ is a lax monoidal fibration whose fibres
are strong monoidal categories.
The fibred forgetful functor $\cU_{2a}$
\begin{center} \xext=600 \yext=580
\begin{picture}(\xext,\yext)(\xoff,\yoff)
\settriparms[1`1`1;500] \putVtriangle(0,0)[Sig_{2a}`Mon(Sig_{2a})`Set^\ra;{\cal F}_{2a} `p_{2a}`q_{2a}]
 \putmorphism(150,420)(1,0)[``{\cal U}_{2a}]{650}{-1}b
\end{picture}
\end{center}
has a fibred left adjoint $\cF_{2a}$, the free monoid functor.
\end{lemma}

{\it Proof.}~ The fact that $p_{2a}$ is a lax monoidal fibration is a simple check.
The construction of a free monoid functor $\cF_{2a}$ is essentially the same as the one given in \cite{HMP}, see also \cite{A} \cite{Ke}, \cite{BJT}.
 $~\Box$

\subsection*{The full embeddings}

The 2-level amalgamated signatures fibration $p_{2a}:Sig_{2a}\lra Set^\ra$ contains as a
lax monoidal subfibration the amalgamated signatures fibration $p_{a}:Sig_{a}\lra Set$,
and as a consequence the category of monoids $Mon(Sig_a)$ with respect to $p_a$ has a full
fibred embedding into the category of monoids $Mon(Sig_{2a})$ with respect to $p_{2a}$.
Moreover, the category of 2-level multicategories with non-standard amalgamation of \cite{HMP} is
a full subcategory of $Mon(Sig_{2a})$, as well. Below we describe this in detail.

\vskip 2mm

The first embedding
\begin{center} \xext=1000 \yext=520
\begin{picture}(\xext,\yext)(\xoff,\yoff)
 \setsqparms[1`1`1`1;600`400]
  \putsquare(0,0)[Sig_a`Sig_{2a}`Set`Set^\ra;\iota`p_a`p_{2a}`\delta]
\end{picture}
\end{center}
is given by the diagonal functor $\delta$ and an inclusion $\iota$. The functor $\delta$ has
both adjoints, say $l\dashv \delta \dashv r$, associating codomain and domain, respectively.
The functor $\iota$ has also both fibred adjoints $L\dashv \iota \dashv R$.
The left adjoint $L$ is defined by composition.
For $\xi$ in $Sig_a$ we have
\[ L(\xi : Y\lra \dot{O}\times O^*)= \dot{\xi } \]
where $\dot{\xi}$ is the composition
\begin{center} \xext=1800 \yext=150
\begin{picture}(\xext,\yext)(\xoff,\yoff)
  \putmorphism(0,0)(1,0)[Y`\dot{O}\times O^*`\xi]{800}{1}a
    \putmorphism(800,0)(1,0)[\phantom{\dot{O}\times \dot{O}^*}`\dot{O}\times \dot{O}^*=\dot{O}^\dag`1\times\dot{(-)}^*]{1000}{1}a
\end{picture}
\end{center}
The right adjoint $R$ is defined by the pullback. We have
\[ R(\xi : Y\lra \dot{O}\times O^*)= \tilde{\xi}  \]
where $\tilde{\xi}$ is given by the following pullback
\begin{center} \xext=800 \yext=520
\begin{picture}(\xext,\yext)(\xoff,\yoff)
 \setsqparms[1`1`1`1;800`450]
  \putsquare(0,0)[\tilde{Y}`Y`O^\dag`\dot{O}\times O^*;`\tilde{\xi}`\xi`\dot{(-)}\times 1_{O*}]
\end{picture}
\end{center}

\vskip 2mm

The second embedding is a functor
\[ \Phi : {\rm Multicat} \lra Mon(Sig_{2a}) \]
from the category {\rm Multicat} of 2-level multicategories with
non-standard amalgamation to the category of monoids in the lax monoidal fibration $p_{2a}:Sig_{2a}\ra Set^\ra$.
Let $\bC$ be an object of {\rm Multicat}. Let $A=A(\bC)$,
$O=O(\bC)$, $\ddot{O}=\ddot{O}(\bC)$ denote arrows, objects and
lower level objects (i.e. types) in $\bC$, respectively. The monoid $\Phi(\bC)$ is in the
fibre over $\dot{(-)}:O\ra \ddot{O}$, i.e. $\vec{O}$. The universe of $\Phi(\bC)$ is $A$ and the typing
function $\partial^A: A\ra \vec{O}^\dag$ is defined as follows.
For $a\in A$ we have source of $a$, $s(a):(|a|]\ra O$ and target
of $a$, $t(a)\in \ddot{O}$.  The function $\partial_a :[|a|]\lra
O+\ddot{O}$ is equal $s(a)$ on $(|a|]$ and
\[ \partial_a(0) =\left\{ \begin{array}{ll}
               x    & \mbox{ if  $a=1_x$} \\
               t(a)    & \mbox{ otherwise.}
                                    \end{array}
                            \right.  \]
i.e. it is the object $x$ if $a$ is the identity on $x$
and it is the type $t(a)$ otherwise. The unit map
$(e,\sigma,id_{\vec{O}}): (O,\partial^{I_{\vec{O}}},\vec{O})\ra
(A,\partial^A,\vec{O})$ is sending $x\in O$
to $1_x\in A$, and $\sigma_x=1_{[1]}$.
The multiplication
 \[ (\mu_A,\sigma,id_{\vec{O}}) : (A\otimes_{\vec{O}}A,\partial^\otimes,\vec{O})
 \lra (A,\partial,\vec{O}) \]
is defined with the help of simultaneous composition operation, see \cite{HMP}. We
first describe exactly the tensor
$(A\otimes_{\vec{O}}A,\partial^\otimes,\vec{O})$. We have
\[A\otimes_{\vec{O}} A=\{ \lk a,a_i \rk_{i\in (|a|]} : a_i,a\in A,\,
\partial_{a_i}(0)\dot{=}\partial_a(i) \]
and for $\lk a,a_i\rk_{i\in (|a|]}\in A\otimes_{\vec{O}} B$,
\[ \partial^{\otimes,+}_{\lk a,a_i\rk _{i}}=
\coprod_{i} \partial^+_{a_i} \, :\, (|\lk a_i,b\rk_i|] = \coprod_{i} (|a_i|]  \lra
O+\ddot{O}\]
and
\[\partial^{\otimes,+}_{\lk a,a_i\rk _{i}}(0)=
\left\{ \begin{array}{ll}
                \dot{(\partial_a(0))} & \mbox{ if  $\partial_a(0)\in O$ and
                $\exists_{i\in (|a|]}\; \partial_{a_i}(0)\in\ddot{O}$ } \\
               \partial_a(0)    & \mbox{ otherewise.}
                                    \end{array}
                            \right.  \]

For $\lk a,a_i \rk_{i\in (|a|]}\in A\otimes_{\vec{O}}A$,  we have the
simultaneous composition $b=a(a_i/i : i\in (|a|])$ together with
amalgamating maps $\phi_i:s(a_i)\ra s(b)$ for $i\in (|a|]$ such that
the function
 \[ [\phi_i]_i : \coprod_i s(a_i)\lra s(b) \]
is a bijection. We put $\mu_A(\lk a,a_i \rk_i)$ to be equal $b$, and $\sigma_b$ is such that
$\sigma_b^- = [\phi_i]_i$.

The remaining details, as well as the definition on morphisms, are easy.
The category $Sig_{2a}$ and the tensor $\otimes_{\vec{O}}$ are so defined that
identities must have objects as codomains. As we
mentioned at the beginning, this is so to mimic the behavior of
identities in {\rm Multicat}. However, all the other arrows in the
monoids coming from {\rm Multicat} have types as
codomains. In fact, this characterizes the monoids coming from {\rm
Multicat}. We have

\begin{proposition} The  functor
\[ \Phi : {\rm Multicat} \lra Mon(Sig_{2a}) \]
is full and faithful and its essential image consists of those
monoids in which all arrows but identities have types
as codomains (i.e. either  $\partial_a(0)\in\ddot{O}$ or $a=1_{\partial_a(0)}$ , for $a\in A$).

The muticategories with the object of objects equal to $\dot{(-)}:O\ra
\ddot{O}$ are sent to the monoids in the fibre over $\vec{O}$.
\end{proposition}

{\it Proof.}~ Simple verification.
 $~\Box$

\vskip 2mm

\subsection{Single tensor in the fibration $p_a$}\label{single-tensor}
In this subsection we describe another widely used tensor in the fibration $p_a$.
To distinguish these two tensors, we shall call the one considered so far the {\em total tensor} and denoted it, in this subsection only, by $\otimes^t$.
The tensor we are going to discuss here, the {\em single tensor}, will be denoted by $\otimes^s$.
Both tensors have the same unit.

Let $(A,\partial^A,O)$ and $(B,\partial^B,O)$ be two object in the fibre $(Sig_a)_O$.
The {\em single tensor} $(A\otimes^s B,\partial^{\otimes^s},O)$ is defined as
follows
\[A\otimes^s B=\{ \lk a,i,b\rk : a\in A,\, b\in B,\, i\in (|a|]\;
\partial^A_{a}(i)=\partial^B_b(0)\} \]
and for $\lk a,i,b\rk\in A\otimes^s B$,
\[ \partial^{\otimes^s}(\lk a,i,b\rk)= {\partial^A_a}_{\lceil (|a|]- \{i\})} + \partial^{B,-}_{b}.\]

Note that contrary to the case of the total tensor $\otimes^t$, the coherence morphisms for the associativity $\alpha$ and the right unit $\rho$
are not isomorphisms. This example together with the Burroni fibrations were the main motivation for choosing the directions
of coherence morphisms in the definition of lax monoidal fibrations.

There is a long debate whether single or total tensor is more convenient.  There are arguments for each. The fibration $p_a:Sig_a\ra Set$ equipped with
single tensor also acts on the basic fibration, but this action is not so much in use. The action $\star^s$
\begin{center} \xext=800 \yext=450
\begin{picture}(\xext,\yext)(\xoff,\yoff)
 \settriparms[1`1`1;400]
  \putVtriangle(0,0)[Sig_a\times_{Set}{Set}^\ra`{Set}^\ra`{Set};\star^s``cod]
\end{picture}
\end{center}
is defined as follows. For $(A,\partial^A,O)$ in $Sig_a$ and $(X,d^X)$ in $(Set^\ra)_O$
we put $(A,\partial^A,O)\star^s(X,d^X) = (A\star^s X, d^s)$, where
\[ A\star^s X =\{ \lk a,i,x\rk :  \partial^A_a(i)=d^X(x) \} \]
and $d^s(a,i,x)= \partial^A_a(0)$. One can easily verify that this extends to a definition of an action of the lax monoidal fibration $(Sig_a,p_a,\otimes^s,I)$
on the basic fibration. So, by adjunction we get a morphism of lax monoidal fibrations
\begin{center} \xext=800 \yext=450
\begin{picture}(\xext,\yext)(\xoff,\yoff)
 \settriparms[1`1`1;400]
  \putVtriangle(0,0)[Sig_a`Exp(Set)`{Set};rep^s_a`p_a`p_{exp}]
\end{picture}
\end{center}
We can describe this representation equivalently as a representation of the fibration of polynomial diagrams. That is, if we consider in the fibration of
polynomial diagrams the image under equivalence of categories $\iota_a$ of the tensor $\otimes^s$ in $p_a$, we have a tensor, also denoted $\otimes^s$,
 in $p_{pd}:\cP oly \cD iag \ra Set$.   Then $p_{pd}$ considered with this tensor has a representation $rep^s_{pd}$ corresponding to the representation
 $rep^s_a$, i.e. a morphism of lax monoidal fibrations
 \begin{center} \xext=800 \yext=450
\begin{picture}(\xext,\yext)(\xoff,\yoff)
 \settriparms[1`1`1;400]
  \putVtriangle(0,0)[\cP oly \cD iag`Exp(Set)`{Set};rep^s_{pd}`p_{pd}`p_{exp}]
\end{picture}
\end{center}
that sends the diagram $(t,p,s)$ to the functor $t_!p_!s^*$.

The representation $rep^s_{pd}$ is faithful but not full even on isomorphisms, as the diagrams
\begin{center} \xext=2100 \yext=150
\begin{picture}(\xext,\yext)(\xoff,\yoff)
\putmorphism(300,0)(1,0)[E`A`p]{300}{1}a
\putmorphism(0,0)(1,0)[O`\phantom{E}`s]{300}{-1}a
\putmorphism(600,0)(1,0)[\phantom{O}`O`t]{300}{1}a

\putmorphism(1500,0)(1,0)[E`A'`p']{300}{1}a
\putmorphism(1200,0)(1,0)[O`\phantom{E}`s]{300}{-1}a
\putmorphism(1800,0)(1,0)[\phantom{E}`O`t']{300}{1}a
 \end{picture}
\end{center}
have isomorphic representations if and only if $t\circ p=t'\circ p'$. We want to describe the image of $rep^s_{pd}$.

We shall call a polynomial diagram $(t,p,s)$, a {\em linear diagram}\footnote{This name and notion is taken from \cite{Ko}.} if and only if $p$ is an isomorphism.
We shall denote by $p_{ld}:\cL in\cD iag \ra Set$ the full subfibration of $p_{pd}$ whose objects in the total category are linear diagrams.
The image of $rep^s_{pd}$ is the same as $rep_{ld}$, the image of $rep^s_{pd}$ restricted to $p_{ld}$.
The linear diagrams are closed under both tensors $\otimes^t$ and $\otimes^s$ and both tensors agree on them.
Thus both representations $rep^t_{pd}$ and $rep^s_{pd}$ coincide on linear diagrams. This statement characterizes the linear diagrams.
This is why the representation of linear diagrams is denoted by $rep_{ld}$, with no superscript.

\begin{proposition} \label{im-lin-diag}
The total category of the image of the morphism of lax monoidal fibrations
 \begin{center} \xext=800 \yext=480
\begin{picture}(\xext,\yext)(\xoff,\yoff)
 \settriparms[1`1`1;400]
  \putVtriangle(0,0)[\cL in \cD iag`Exp(Set)`{Set};rep_{ld}`p_{ld}`p_{exp}]
\end{picture}
\end{center}
consists of endofunctors on slices of $Set$ that preserves colimits and wide pullbacks as objects and cartesian natural transformations as morphisms.
\end{proposition}

{\it Proof.}~ From the characterization of the image of $p_{pd}$ in $p_{exp}$ the necessity of the conditions is obvious.
On the other hand, by  Theorem \ref{poly-char}, any endofunctor $P$ on a slices $Set_{/O}$ that preserves colimits and wide pullbacks
is of the form $t_!p_*s^*$ for some polynomial diagram
\begin{center} \xext=900 \yext=100
\begin{picture}(\xext,\yext)(\xoff,\yoff)
\putmorphism(300,0)(1,0)[E`A`p]{300}{1}a
\putmorphism(0,0)(1,0)[O`\phantom{E}`s]{300}{-1}a
\putmorphism(600,0)(1,0)[\phantom{O}`O`t]{300}{1}a
 \end{picture}
\end{center}
Recall that $p$ has finite fibres.
Since $P$ preserves the initial object the fibres of $p$ cannot be empty. Since $p$ preserves binary coproduct the fibres of $p$ cannot have more than one element.
Thus $p$ is an isomorphism, and the diagram representing $P$ is linear. The characterization of natural transformations in the image of $rep_{ld}$ follows directly from Theorem \ref{poly-char}.
 $~\Box$

\vskip 2mm

The polynomial diagrams of form
 \begin{center} \xext=900 \yext=150
\begin{picture}(\xext,\yext)(\xoff,\yoff)
\putmorphism(300,0)(1,0)[E`O`p]{300}{1}a
\putmorphism(0,0)(1,0)[O`\phantom{E}`s]{300}{-1}a
\putmorphism(600,0)(1,0)[\phantom{O}`O`1_O]{300}{1}a
 \end{picture}
\end{center}
are also closed under the tensor $\otimes^t$ in $p_{pd}$, and form a monoidal full subfibration of $p_{pd}$ considered with the total tensor $\otimes$.
Such diagrams correspond to signatures that have exactly one operation of each (output) type. The representations of such diagrams are endofunctors $P$ on $Set_{/O}$
that send a function $d$ to a function $P(d)$ that has as fibres finite products of fibres of $d$. For this reason, we  call such
polynomial diagrams {\em monomial diagrams} and the full subfibration of monomial diagrams will be denoted by $p_{md}:\cM ono \cD iag \ra Set$.
The composition of  $rep_{pd}$ with the inclusion gives a representation morphism $rep_{md}$
   \begin{center} \xext=800 \yext=450
\begin{picture}(\xext,\yext)(\xoff,\yoff)
\settriparms[1`1`1;400]
  \putVtriangle(0,0)[\cM ono\cD iag `Exp(Set)`Set;rep_{md}`p_{md}`rep_{exp}]
 \end{picture}
\end{center}

We have

\begin{proposition} \label{im-mono-diag}
The total category of the image of the morphism $rep_{md}$ of lax monoidal fibrations
 \begin{center} \xext=800 \yext=480
\begin{picture}(\xext,\yext)(\xoff,\yoff)
 \settriparms[1`1`1;400]
  \putVtriangle(0,0)[\cM ono\cD iag`Exp(Set)`{Set};rep_{md}`p_{md}`p_{exp}]
\end{picture}
\end{center}
consists of finitary endofunctors on slices of $Set$ that preserves limits and cartesian natural transformations.
\end{proposition}

{\it Proof.}~This follows from Theorem \ref{poly-char} and the observation, that for any polynomial diagram $(t,p,s)$ the polynomial functor
$t_!p_*s^*$ preserves the terminal object $1$ if and only if $t$ is an isomorphism.  $~\Box$

\vskip 2mm

{\bf Remarks}
\begin{enumerate}
  \item We want to point out that the fibrations of linear diagrams $p_{ld}$ and monomial diagrams $p_{md}$ are certain fibrations of graphs.
The category $Gph$ has as objects parallel pairs of functions $s,t:A\ra O$.  A morphism of graphs $(f,u): (s,t)\ra (s',t')$ is a pair of functions $f:A\ra A'$
$u:O\ra Q$ making the diagram
\begin{center} \xext=1000 \yext=500
\begin{picture}(\xext,\yext)(\xoff,\yoff)
\setsqparms[1`1`1`1;500`400]
\putsquare(500,50)[E`O`E'`Q;p`f`u`p']
\setsqparms[-1`1`0`-1;500`400]
\putsquare(0,50)[O`\phantom{E}`Q`\phantom{E'};s`u``s']
 \end{picture}
\end{center}
commute. The category of cartesian graphs $cGph$ is a subcategory containing the same objects as $Gph$ and a morphism $(f,u)$ in $cGph$
if and only if the right square above is a pullback. These categories are fibred over $Set$ and the projection functors $p_G: Gph\ra Set$, $p_{cG}: Gph\ra Set$
send a morphism $(f,u)$ to $u$. Both fibrations have a lax monoidal structure with tensor given by the obvious pullbacks. Moreover, we have
equivalences of lax monoidal fibrations
 \begin{center} \xext=2000 \yext=400
\begin{picture}(\xext,\yext)(\xoff,\yoff)
 \settriparms[1`1`1;300]
  \putVtriangle(0,0)[Gph`\cL in \cD iag`{Set};`p_{G}`p_{ld}]

  \putVtriangle(1400,0)[cGph`\cM ono\cD iag`{Set};`p_{cG}`p_{md}]
\end{picture}
\end{center}
sending a graph $(t,s:A\ra O)$ to a linear diagram $(t,1_A,s)$ and a cartesian graph $(p,s:A\ra O)$ to a monomial diagram $(t,p,1_O)$.
  \item So far we haven't said anything about monoids in $p_a$ with the single tensor $\otimes^s$.
We can pass from multiplication with respect to the total tensor to
multiplication with respect to the single tensor by putting identities into all places but one.
Thus we have an embedding of the monoids with respect to the total tensor into the the monoids with respect to the single tensor.
  To characterize the image of this embedding
  we shall use a certain natural isomorphism involving $\otimes^s$ and the binary coproduct in fibres $+$.
  Note that for any signatures $A$, $B$, $C$   in the same fibre over $O$ of $Sig_a$, we have an isomorphism:
  \[ ((A \otimes^s B) \otimes^s C) +  (A \otimes^s (C \otimes^s B)) \cong ((A \otimes^s C) \otimes^s B) + (A \otimes^s (B \otimes^s C)) \]
  that 'repairs' the lack of strong associativity for the tensor $\otimes^s$. Intuitively, both sides of the isomorphism contain the part
  \[ (A \otimes^s (C \otimes^s B)) + (A \otimes^s (B \otimes^s C)) \]
  i.e. the $A$' into which we plug either a $B$ with plugged in a $C$ or  a $C$ with plugged in a $B$.
  The remaining part on both  sides of the above isomorphism is $A$'s into which we plug directly  a $B$ and a $C$.
    Clearly, this isomorphism, $\upsilon_{A,B,C}$ is natural in $A$, $B$ and $C$.
    $\upsilon_{A,A,A}$ may look trivial but it is not! Then one can verify that
  a monoid $(M,m,e)$ in $Sig_a$ with respect to $\otimes^s$ comes from a monoid in $Sig_a$ with respect to $\otimes^t$, if the following diagram
  \begin{center} \xext=3000 \yext=650
\begin{picture}(\xext,\yext)(\xoff,\yoff)
\put(1400,600){\makebox(200,100){$((M \otimes^s M) \otimes^s M) +  (M \otimes^s (M \otimes^s M)) \cong ((M \otimes^s M) \otimes^s M) + (M \otimes^s (M \otimes^s M))$}}
\setsqparms[0`1`0`1;1000`550]
\putsquare(500,50)[``(M \otimes^s M)+(M \otimes^s M)`M;`(m\otimes 1_M)+(1_M\otimes m)``{[m,m]}]
\setsqparms[0`0`1`-1;1000`550]
\putsquare(1500,50)[``\phantom{M}`(M \otimes^s M)+(M \otimes^s M);``(m\otimes 1_M)+(1_M\otimes m)`{[m,m]}]
\end{picture}
\end{center}
  commutes, where the unnamed isomorphism is  $\upsilon_{M,M,M}$. This condition corresponds to the commutativity condition in the multicategories with non-standard amalgamations in \cite{HMP}.
\end{enumerate}

\section{Symmetric signatures vs analytic functors}\label{symm-sig}

\subsection{The symmetric signature fibration  $p_s: Sig_s \lra Set$}

\subsection*{The category of symmetric sets}

The category of symmetric sets is equivalent to the category of species, cf. \cite{J1},
however the presentation is slightly different.

A symmetric set $(A,\alpha)$ is a graded set $\{ A_n: n\in\o \}$ with (right)
actions of symmetric groups $\alpha_n: A_n\times S_n\ra A_n$, for
$n\in\o$. We write $a\in A$ to mean that $a\in \coprod_n A_n$ and if $a\in A$ then we write
$|a|=n$ to mean that $a\in A_n$. Thus, for $a\in A$, we have $a\in A_{|a|}$.
In case $\alpha(a,\sigma)$ is defined we usually
write it  as $a\cdot\sigma$ if it does not lead to a confusion. A
morphism of symmetric sets $f:(A,\alpha)\ra (B,\beta)$
is a family of morphisms of actions $f_n:(A_n,\alpha_n)\ra (B_n,\beta_n)$ for $n\in\o$, i.e.
it is a function $f:A\ra B$ commuting with the actions $\alpha$ and $\beta$, in short.
We call such morphisms {\em equivariant}. The category of symmetric sets will be denoted by $\sigma Set$.
$\sigma Set$ is (equivalent to) of the presheaf category $Set^{S_*^{op}}$, where $S_*$ the coproduct of
(finite) symmetric groups in $Cat$.
Clearly the groups $S_*$ act on
$O^\dag$ on the right by composition, leaving $0$ fixed. This symmetric set on $O^\dag$ will be denoted by $O^\ddag$.
Any function $u:O\ra Q$ induces an equivariant map $u^\ddag :O^\ddag\ra Q^\ddag$, so that $u^\ddag(d)=u\circ d$.
Thus we have a functor
\[ (-)^\ddag : Set \lra \sigma Set \]

\subsection*{The operad of symmetries $\bS$}

Recall\footnote{For example from \cite{Le} pp. 51-54.} that the universes of symmetric groups $\lk S_n\rk_{n\in\o}$ form the underlying sets of an operad, called the operad of symmetries $\bS$. The compositions
\[ \ast : S_k\times (S_{n_1}\times\ldots \times S_{n_k}) \lra S_{\sum_{i=1}^k n_i}  \]
\[  (\tau,\sigma_1,\ldots, \sigma_k)\mapsto \tau\ast (\sigma_1,\ldots, \sigma_k) \]
where, for $\tau\in S_k$, $\sigma_1\in S_{n_1}, \ldots, \sigma_k\in S_{n_k}$, $1\leq m_0 \leq k$, $1\leq m_1 \leq n_{k_{m_0}}$
are given by
\[ \tau\ast(\sigma_1,\ldots, \sigma_k)(k_1+\ldots +k_{m_0-1}+m_1)= k_{\tau^{-1}(1)}+\ldots +k_{\tau^{-1}(\tau(m_0)-1)}+\sigma_{m_0}(m_1)  \]

\subsection*{The fibration $p_s$}

Taking the pullback of the basic fibration on the category of symmetric sets $cod: \sigma Set^{\rightarrow} \lra \sigma Set$ along the functor $(-)^\ddag$
\begin{center}
\xext=700 \yext=550
\begin{picture}(\xext,\yext)(\xoff,\yoff)
 \setsqparms[1`1`1`1;700`450]
 \putsquare(0,50)[Sig_s`\sigma Set^{\rightarrow}`Set`\sigma Set;`p_s`cod`(-)^\ddagger]
 \end{picture}
\end{center}
we obtain the symmetric signature fibration $p_s$.
We describe the category $Sig_s$ explicitly.
An object of $Sig_s$ over the set $O$ is a quadruple
$(A,\alpha,\partial^A,O)$ such that $(A,\alpha)$ is a
symmetric set, $\partial^A:(A,\alpha)\ra O^\ddag$ is an
equivariant map called the {\em typing} (or {profile} in \cite{BaezDolan}) map of the
signature. We write $\partial^A_a:[n]\ra O$ for the effect of
$\partial^A$ on $a\in A$, and $n$ in this case can be referred to as
$|a|$. The fact that $\partial^A$ is equivariant means that we have
$\partial^A_{a\cdot\sigma}=\partial^A_a\circ\sigma$, for $a\in A$ and
$\sigma \in S_{|a|}$. A morphism $(f,u): (A,\alpha,\partial^A,O)\lra
(B,\beta,\partial^B,Q)$ in $Sig_s$ over a function $u:O\ra Q$
is a commuting square of equivariant maps:
\begin{center}
\xext=700 \yext=600
\begin{picture}(\xext,\yext)(\xoff,\yoff)
 \setsqparms[1`1`1`1;700`450]
 \putsquare(0,50)[(A,\alpha)`(B,\beta)`O^\ddag`Q^\ddag;f`\partial^A`\partial^B`u^\ddag]
 \end{picture}
\end{center}

\subsection*{The monoidal structure in the fibres of $p_s$}

We define two lax morphisms of fibrations
\begin{center} \xext=1800 \yext=480
\begin{picture}(\xext,\yext)(\xoff,\yoff)
\settriparms[1`1`1;400]
 \putVtriangle(0,0)[{Sig_s}\times_{Set}{Sig_s}`{Sig_s}`Set;\otimes`p'_s`p_s]
\settriparms[1`1`1;400]
 \putVtriangle(1200,0)[Set`{Sig_s}`Set;I`id`p_s]
\end{picture}
\end{center}

Let $(A,\alpha,\partial^A,O)$, $(B,\beta,\partial^B,O)$ be two
objects in the fibre over $O$ in the fibration $p_s$. The {\em tensor product}
\[ (A,\alpha,\partial^A,O)\otimes_O(B,\beta,\partial^B,O)=
(A\otimes_OB,\alpha\otimes_O\beta,\partial^\otimes,O)  \]
is defined as follows
\begin{flushleft}
$(A\otimes_OB)_n =$
\end{flushleft}
\[ = \{ \lk a,\lk b_i\rk_{i\in (|a|]},\sigma \rk :\sum_i|b_i|=n,\, b_i\in B,\,a\in A,
\partial_{b_i}(0)=\partial_a(i),\, {\rm for\;}i\in (|a|],\, \sigma\in S_n \}_{/\sim}\]
where the equivalence relation $\sim$ is defined as follows:
\[ \lk a\cdot\tau;\lk b_{\tau(i)}\cdot\sigma_{\tau(i)}\rk_i;\sigma\rk\sim
\lk a,\lk b_{i}\rk_i,\tau\ast (\sigma_{\tau(1)}, \ldots,\sigma_{\tau(|a|)}) \circ\sigma\rk \]
where $\tau\in S_{|a|}$, $\sigma_i\in S_{|b_i|}$, $\sigma\in
S_{\sum_i|b_i|}$, $*$ is the composition in the operad of symmetries, and $\circ$ is the usual composition of permutations.
The equivalence class of the element $\lk a,\lk b_i\rk_{i\in (|a|]},\sigma \rk$ with be denoted by $[\lk a,\lk b_i\rk_{i\in (|a|]},\sigma \rk ]_\sim$.

Let $\kappa_i : (|b_i|]\ra (\sum_i|b_i|]$ be the $i$-th inclusion into the coproduct, for $i=1,\ldots, |a|$.
Clearly, there are many such inclusions that make  $(\sum_i|b_i|]$ into a coproduct of $(|b_i|]$'s (in $Set$) but we will always mean the simplest,
that is embedding blocks  $(|b_i|]$ one after the other into $(\sum_i|b_i|]$
(i.e. $\kappa_i(j)=j+\sum_{k=1}^{i-1}|b_k|$ for $j\in (|b_i|])$.

We define
\[ \partial^\otimes([\lk a,\lk b_i\rk_{i\in (|a|]},\sigma \rk ]_\sim): [\sum_i|b_i|] \lra O \]
so that
 \[ \partial^\otimes([\lk a,\lk b_i\rk_{i\in (|a|]},\sigma \rk ]_\sim)(0)=\partial_a(0) \]
and the squares
\begin{center} \xext=850 \yext=530
\begin{picture}(\xext,\yext)(\xoff,\yoff)
   \setsqparms[1`-1`1`1;850`450]
  \putsquare(0,30)[(\sum_i|b_i|]`(\sum_i|b_i|]`(|b_i|]`O;
  \sigma^{-1}`\kappa_i`\partial^\otimes_{[\lk a,\lk b_i\rk_{i},\sigma \rk ]_\sim}`\partial_{b_i}]
\end{picture}
\end{center}
commute, for all $i\in  (|a|]$. So the type of the codomain of the 'operation' $[\lk a,b_i,\sigma \rk_{i\in (|a|]}]_\sim$
in $A\otimes_OB$ is the same as the type of the codomain of $a$ in $A$ and the types of the domain  of
the 'operation'  $[\lk a,\lk b_i\rk_{i\in (|a|]},\sigma \rk ]_\sim$ in $A\otimes_OB$ are the types of the domains
of $b_i$'s in $B$ put one next to the other and permuted by $\sigma$.

The action of $S_n$ on $(A\otimes_OB)_n $
is defined so that
 \[ [\lk a,\lk b_i\rk_{i\in (|a|]},\sigma \rk ]_\sim \cdot\sigma'=[\lk a,\lk b_i\rk_{i\in (|a|]},\sigma\circ\sigma' \rk ]_\sim  \]
for $\sigma'\in S_{\sum_i|b_i|}$.

For two morphisms
 \[ (f,u):(A,\alpha,\partial,O)\lra
 (B,\beta,\partial,Q),\hskip 1cm (f',u):(A',\alpha',\partial,O)\lra (B',\beta',\partial,Q) \]
 over $u$, we define their tensor to be
\[ (f\otimes_u f',u):(A\otimes_O A',\alpha\otimes_O\alpha',\partial^\otimes,O)
\lra (B\otimes_Q B',\beta\otimes_Q\beta',\partial^\otimes,Q)\]
in the following way. For $\lk a,\lk a'_i\rk_{i\in (|a|]},\sigma \rk \in A\otimes_O A'$, we put
 \[ f\otimes_u f'([\lk a,\lk a'_i\rk_{i\in (|a|]},\sigma \rk ]_\sim)=[\lk f(a),\lk f'(a'_{i})\rk_{i\in (|a|]},\sigma \rk ]_\sim\]
Note that $|a|=|f(a)|$ and we have
\[ \partial^B_{f(a)}(i)=u\circ\partial^A_a(i)=u\circ\partial^{A'}_{a'_{i}}(0)=
\partial^{B'}_{f'(a'_{i})}(0) \] so $[\lk
f(a),\lk f'(a'_{i})\rk_{i\in(|f(a)|]},\sigma \rk]_\sim$ belongs to
$B\otimes_Q B'$ indeed.
This ends the definition of the tensor product functor $\otimes$ in $p_s$.

The unit $I_O=(O,1,\partial,O)$  for the tensor $\otimes_O$ in the fibre $(Sig_s)_O$
is defined as follows. For $x\in O$,
$\partial_x:[1]\ra O$ is a function  such that $\partial_x(0)=\partial_x(1)=x$. So
only the group $S_1$ acts on $O$ and it acts trivially. The
association $O\mapsto I_O$ is clearly the object part of a lax morphism of fibrations, as it should be.

\begin{lemma} The functors $\otimes$ and $I$ together with obvious associativity, left unit, and right unit isomorphisms $\alpha$, $\lambda$, $\rho$
make the fibres of $p_s$ into (strong) monoidal categories.
\end{lemma}

\subsection*{Pulling back the monoidal structure in $p_s$}

Any  object $\lk B,\beta,\partial,Q\rk$ in the fibre $(Sig_s)_Q$ can be pulled back
along a function $u:O\ra Q$
\begin{center} \xext=500 \yext=450
\begin{picture}(\xext,\yext)(\xoff,\yoff)
   \setsqparms[1`1`1`1;700`400]
  \putsquare(0,0)[u^*(B)`B`O^\ddag`Q^\ddag;\pi_B`\partial^{u^*(B)}`\partial^B`u^\ddag]
\end{picture}
\end{center}
where $u^\ddag(d)=u\circ d$. We have
\[ u^*(B)=\{ \lk b,d\rk : b\in B,\, d:[|b|]\ra O,\, {\mathrm such\, that}\; u\circ d=\partial_b \}\]
and
\[ \partial^{u^*(B)}_{\lk b,d\rk}=d \]
The action in $u^*(B)$ applies the permutation to both arguments, i.e.
\[ \lk b,d \rk \cdot \sigma = \lk b\cdot \sigma,d\circ \sigma \rk \]

Let, for $x\in O$, $d_x [1]\ra O$ be the function such that $d_x(0)=d_x(1)=x$.   We have
\[ u^*(I_{Q})=\{ \lk x,d_x\rk :  x\in O \} \]
and
\[ \phi_0: I_O\lra f^*(I_{O'}) \]
\[ x\mapsto \lk u(x),d_x\rk \]
Moreover, for $\lk A,\alpha\rk$ and $\lk
B,\beta\rk$ in $(Sig_s)_Q$, we have
\[ u^*(A\otimes B)= \{ \lk\lk a,\lk b_i\rk_{i\in (|a|]},\sigma\rk,d^+_a+\coprod_{i\in
(|a|]}d^-_{b_i}\rk:
 u^\ddag(d^+_a+\coprod_{i\in (|a|]}d^-_{b_i})
=\partial^+_a+\coprod_{i\in (|a|]}\partial^-_{b_i}\}\]
($d_a^+:[0]:\ra O$, $d_{b_i}^-:(|b_i|]:\ra O$, for $i\in (|a|]$) and
\[ u^*(A)\otimes u^*(B)= \{ \lk\lk a,d\rk ,\lk \lk b_i,d_i\rk\rk_{i\in (|a|]},\sigma\rk:
 u^\ddag(d)=\partial_a,\,
u^\ddag(d_i) =\partial_{a_i},\; \mathrm{for}\; i\in
(|a|] \}\] ($d:[|a|]:\ra O$, $d_i:[|b_i|]:\ra O$, for $i\in (|a|]$).
 Thus we have a transformation
 \[ \phi_{2,A,B}: u^*(A)\otimes u^*(B)\lra u^*(A\otimes B)\]
 such that
\[  \lk\lk a,d\rk,\lk\lk b_i,d_i\rk\rk_{i\in (|a|]},\sigma\rk \mapsto
\lk\lk a,\lk b_i\rk_{i\in (|a|]},\sigma\rk,d^+_a+\coprod_{i\in
(|a|]}d^-_{b_i}\rk  \]

\begin{lemma} The map $(\pi_B,u):u^*(B)\ra B$ is a prone arrow over $u$.
The data $u^*$, $\phi_0$, $\phi_2$ above make the usual (three) diagrams
of a (lax) monoidal functor commute, i.e. $p_s$ equipped with $\otimes$, $I$,
$\alpha$, $\lambda$, $\rho$ is a lax monoidal fibration. $\Box$
\end{lemma}

Moreover, we have

\begin{proposition}  The total category of the fibration of monoids
$q_s:Mon(Sig_s)\lra Set$ is equivalent to the category
of symmetric multicategories.
The fibred forgetful functor from the fibration of monoids to the
fibration of symmetric signatures
${\cal U}: Mon(Sig_s)\lra Sig_s$
\begin{center} \xext=600 \yext=580
\begin{picture}(\xext,\yext)(\xoff,\yoff)
 \settriparms[1`1`1;500] \putVtriangle(0,0)[Sig_s`Mon(Sig_s)`Set;{\cal F}_s`p_s`q_s]
 \putmorphism(150,420)(1,0)[``{\cal U}_s]{650}{-1}b
\end{picture}
\end{center}
is a morphism of fibrations and has a left adjoint $\cF_s$, the free monoid functor,
which is a lax morphism of fibrations. $\Box$
\end{proposition}

{\bf Remark} The $p_s$ is a bifibration as, for an object $(A,\alpha,\partial^A,O)$
and a function $u:O\ra Q$ the morphism
\[  (1_A,u): (A,\alpha,\partial^A,O)\lra (A,\alpha,u^\ddag\circ\partial^A,Q)\]
is a supine morphism.

\subsection{The action of $p_s$ on the basic fibration and analytic functors}\label{act-sigs}

The fibration $p_s$ acts on the basic fibration $cod:Set^\ra\lra
Set$ as follows
\begin{center} \xext=1000 \yext=480
\begin{picture}(\xext,\yext)(\xoff,\yoff)
\settriparms[1`1`1;400]
 \putVtriangle(0,0)[{Sig_s}\times_{Set}{Set^\ra}`{Set^\ra}`Set;\star``cod]
\end{picture}
\end{center}
 In the following, we often denote an object $(A,\alpha,\partial^A,O)$ in $Sig_s$
 as $A$ and an object  $d^X:X\ra O$ in $Set^\ra$ as $X$, when it does not lead to confusion.

The object $d^\star : A\star X\lra O$ is defined as the quotient of the set
\[ \{ (a,\vec{x}): a\in A,\, \vec{x}:(|a|]\ra X,\; \partial^{A,+}_a=d^X\circ \vec{x} \}\]
 by an equivalence $\sim$ so that
 \[ (a,\vec{x})\sim (a\cdot\sigma,\vec{x}\circ\sigma)\]
 for $a\in A$, $\vec{x}:(|a|]\ra  X$, and $\sigma\in S_{|a|}$.
 The function $d^\star : A\star
 X\ra O$ is defined as
 \[ d^\star([ a,\vec{x}]_\sim)=\partial_a(0) \]

The action $\star$ is defined on morphisms as follows. For maps
\[ (f,u):(A,\alpha,\partial^A,O)\lra(B,\beta,\partial^B,Q) \]
in $Sig_s$,
\[ (g,u):(X,d^X)\ra (Y,d^Y)\]
in $Set^\ra$ over $u:O\ra Q$ and for an element $[ a,\vec{x}]_\sim\in A\star X$ we put
\[ f\star_u g([a,\vec{x}]_\sim)=[ f(a),g\circ\vec{x}]_\sim \]
so that the square
\begin{center} \xext=1000 \yext=600
\begin{picture}(\xext,\yext)(\xoff,\yoff)
   \setsqparms[1`1`1`1;1000`500]
  \putsquare(0,20)[A\star X`B\star Y`O`Q;f\star_u g`d^\star`d^\star`u]
\end{picture}
\end{center}
commutes.

We have an adjoint morphism of lax monoidal fibrations
\begin{center} \xext=800 \yext=520
\begin{picture}(\xext,\yext)(\xoff,\yoff)
 \settriparms[1`1`1;400]
  \putVtriangle(0,0)[Sig_s`Exp(Set)`{Set};rep_s`p_s`]
\end{picture}
\end{center}
For an object $A$ in $(Sig_s)_O$, we have a functor
\[ rep_s(A)=A\star(-): Set_{/O}\lra Set_{/O} \]
and for a morphism $(f,u):A\ra B$ in $Sig_s$ over $u:O\ra Q$ we have a natural transformation
\[ rep_s(f,u) : A\star u^*(-) \lra u^*(B\star (-)) \]
so that for a morphism $g: Y \ra Y'$ in $(Set^\ra)_Q$  we have a commuting square
 \begin{center} \xext=1200 \yext=780
\begin{picture}(\xext,\yext)(\xoff,\yoff)
\setsqparms[1`1`1`1;1200`550]
\putsquare(0,80)[A\star u^*(Y)`u^*(B\star Y)`A\star u^*(Y')`u^*(B\star Y');rep_s(f,u)_Y`1_A\star u^*(g)`u^*(1_B\star g)`rep_s(f,u)_{Y'}]
 \end{picture}
\end{center}
We note for the record, that for $[a,\vec{x}]_\sim\in A\star u^*(Y)$ so that $a\in A$ and $\vec{x}:(|a|]\ra u^*(Y)$ we have
\[ rep_s(f,u)_Y([a,\vec{x}]_\sim) = \lk \partial^A_a(0),[f(a),u^Y\circ\vec{x}]_\sim\rk \]
where
 \begin{center} \xext=500 \yext=500
\begin{picture}(\xext,\yext)(\xoff,\yoff)
\setsqparms[1`1`1`1;500`400]
\putsquare(0,50)[u^*(Y)`Y`O`Q;u^Y`u^*(d^Y)`d^Y`u]
 \end{picture}
\end{center}
is a pullback in $Set$. We have

\begin{proposition} \label{reps-is-bifib}
The morphism of lax monoidal fibrations $rep_s$ defined above is a morphism of bifibrations and it preserves coproducts in the fibres.
\end{proposition}

{\it Proof.}~ The proof of this proposition can be made more abstract but we prefer it to be concrete. Preservation of coproduct is trivial.

First, we shall show that $rep_s$ is a morphism of fibrations. Let $B=(B,\beta,\partial^B,Q)$ be a symmetric signature.
The prone morphism over $u:O\ra Q$ in the fibration $p_s$ is defined via pullback in the category of symmetric sets
 \begin{center} \xext=800 \yext=500
\begin{picture}(\xext,\yext)(\xoff,\yoff)
\setsqparms[1`1`1`1;800`400]
\putsquare(0,50)[u^*(B,\beta)`(B,\beta)`O^\ddag`Q^\ddag;u^B`\partial^{u^*(B)}`\partial^B`u^\ddag]
 \end{picture}
\end{center}
We write $u^*(B)$ for $(u^*(B,\beta),\partial^{u^*(B)},O)$. Then the prone morphism is $pr_{u,B}=(u^B,u) : u^*(B) \ra B$.
We will use the usual representation of this pullback in $Set$ i.e.
\[ u^*(B)=\{\lk b,d\rk : b\in B,\; d:[|b|]\ra O,\; u\circ d =\partial^B_b  \} \]
and
\[ \partial^{u^*(B)}(\lk b,d\rk)= d, \hskip 5mm u^B(\lk b,d\rk)= b \]
for $\lk b,d\rk \in u^*(B)$.

We will show that the following natural transformations
 \begin{center} \xext=1600 \yext=150
\begin{picture}(\xext,\yext)(\xoff,\yoff)
\putmorphism(0,30)(1,0)[ u^*(B)\star u^*(-) `u^*(B\star (-)) `rep_s (pr_{u,B}) = rep_s (u^B,u) ]{2000}{1}a
\end{picture}
\end{center}
and
 \begin{center} \xext=1600 \yext=150
\begin{picture}(\xext,\yext)(\xoff,\yoff)
\putmorphism(0,30)(1,0)[ u^*(B\star u_!u^*(-))` u^*(B\star (-)) `pr_{u, rep_s(B)} = u^*(B\star \varepsilon^u_{(-)})]{2000}{1}a
\end{picture}
\end{center}
in $Cat(Set_{/Q},Set_{/O})$ are isomorphic as objects of $Cat(Set_{/Q},Set_{/O})_{/u^*(B\star (-))}$.
To this end we shall define a natural isomorphism
 \begin{center} \xext=1200 \yext=100
\begin{picture}(\xext,\yext)(\xoff,\yoff)
\putmorphism(0,0)(1,0)[u^*(B)\star(-)`u^*(B\star u_!u^*(-))`\xi]{1200}{1}a
\end{picture}
\end{center}
(and its inverse) so that $rep_s(u^B,u)=u^*(B\star \varepsilon^u_{(-)})\circ \xi$, i.e. for any $d^Y:Y\ra Q$ the triangle
 \begin{center} \xext=1000 \yext=800
\begin{picture}(\xext,\yext)(\xoff,\yoff)
\putmorphism(0,700)(0,-1)[u^*(B)\star(d^Y)`u^*(B\star u_!u^*(d^Y))`\xi_{d^Y}]{600}{1}l

 \put(1400,350){\makebox(200,100){$u^*(B\star d^Y)$}}
 \put(350,150){\vector(4,1){850}}
  \put(350,650){\vector(4,-1){850}}

  \put(700,600){\makebox(300,100){$rep_s(u^B,u)_{d^Y}$}}
   \put(700,110){\makebox(300,100){$u^*(B\star \varepsilon^u_{d^Y})$}}
       \put(2000,400){\makebox(300,100){$(1)$}}
\end{picture}
\end{center}
commutes. We fix $d^Y:Y\ra Q\in Set_{/Q}$. The following diagram
 \begin{center} \xext=1000 \yext=950
\begin{picture}(\xext,\yext)(\xoff,\yoff)
\setsqparms[1`1`1`1;600`400]
\putsquare(0,250)[(|b|]`u^*(Y)`[|b|]`O;\vec{x}```d]

\setsqparms[1`0`1`1;600`400]
\putsquare(600,250)[\phantom{u^*(Y)}`Y`\phantom{O}`Q;u^Y`u^*(d^Y)`d^Y`u]

\put(600,30){\makebox(200,100){$\partial^B_b$}}

 \put(1200,0){\vector(0,1){190}}
  \put(0,0){\line(0,1){200}}
  \put(0,0){\line(1,0){1200}}

\put(600,800){\makebox(200,100){$\vec{y}$}}

 \put(1200,900){\vector(0,-1){190}}
  \put(0,900){\line(0,-1){190}}
  \put(0,900){\line(1,0){1200}}
    \put(2000,400){\makebox(300,100){$(2)$}}
 \end{picture}
\end{center}
where the right hand square is a pullback is to fix the notation. We do not assume now that other morphisms exist, but if they do they have domains
and codomains as displayed. Similarly, other figures in this diagram are not assumed to commute unless we explicitly say so. We will refer often
to this diagram in the rest of the proof.

We note that, for $b\in B$,
\[ \lk o,[b,\vec{y}]_\sim \rk\in u^*(B\star d^Y) \hskip 5mm {\rm iff} \hskip 5mm o\in O,\; u(o)=\partial^B_b(0),\;  \partial^{B,+}_b = d^Y\circ \vec{y}\]

\[   \lk o,[b,\vec{x}]_\sim \rk\in u^*(B\star u_!u^*(d^Y))   \hskip 5mm {\rm iff} \hskip 5mm  o\in O,\; u(o)=\partial^B_b(0),\; \partial^{B,+}_b=  u\circ u^*(d^Y)\circ\vec{x}\]

\[   [ \lk b,d \rk \vec{x}]_\sim \in u^*(B)\star u^*(d^Y)   \hskip 5mm {\rm iff} \hskip 5mm   u\circ \partial^B_b=d,\; d^+ =  u^*(d^Y)\circ\vec{x}\]

With the above notation we spell out the three functions occurring in $(1)$:
 \begin{center} \xext=1500 \yext=450
\begin{picture}(\xext,\yext)(\xoff,\yoff)
\putmorphism(0,300)(1,0)[u^*(B)\star u^*(d^Y)`u^*(B\star d^Y)`rep_s(u^B,u)_{d^Y}]{1500}{1}a
\putmorphism(200,40)(1,0)[[\lk b,d \rk,\vec{x} ]_\sim`\lk d(0),[b,u^Y\circ \vec{x}]_\sim \rk `]{1100}{1}a
 \put(450,0){\line(0,1){80}}
\end{picture}
\end{center}
and
 \begin{center} \xext=1500 \yext=450
\begin{picture}(\xext,\yext)(\xoff,\yoff)
\putmorphism(0,300)(1,0)[u^*(B\star u_!u^*(d^Y))`u^*(B\star d^Y)`u^*(B\star \varepsilon^u_{d^Y})]{1500}{1}a
\putmorphism(200,40)(1,0)[\lk o,[b,\vec{x}]_\sim \rk`\lk o,[b,u^Y\circ u^Y\circ \vec{x}]_\sim \rk `]{1100}{1}a
 \put(450,0){\line(0,1){80}}
\end{picture}
\end{center}
and
 \begin{center} \xext=1500 \yext=650
\begin{picture}(\xext,\yext)(\xoff,\yoff)
\putmorphism(0,500)(1,0)[u^*(B)\star u^*(d^Y)`u^*(B\star u_!u^*(d^Y))`\xi_{d^Y}]{1500}{1}a
\putmorphism(200,240)(1,0)[ [\lk b,d \rk,\vec{x} ]_\sim`\lk d(0),[b,\vec{x}]_\sim \rk `]{1100}{1}a
 \put(450,200){\line(0,1){80}}

\putmorphism(200,40)(1,0)[ [\lk b,\bar{d} \rk,\vec{x} ]_\sim`\lk o,[b,\vec{x}]_\sim \rk `]{1100}{-1}a
 \put(1050,0){\line(0,1){80}}
 \end{picture}
\end{center}
where $\bar{d}:[|b|]\ra O$ is so defined that $\bar{d}(0)=o$ and $\bar{d}^+=u^*(d^Y)\circ \vec{x}$.
Now a simple check shows that $(1)$ commutes, i.e. $rep_s$ preserves prone morphisms.

Now we shall show that $rep_s$ preserves supine morphisms, i.e. it is a morphism of opfibrations.
Let $(A,\alpha,\partial^A,O)$ be a symmetric signature. The supine morphism $su_{u,A}$ in $p_s$ over $u:O\ra Q$ with domain $A$
is defined from the square
 \begin{center} \xext=600 \yext=500
\begin{picture}(\xext,\yext)(\xoff,\yoff)
\setsqparms[1`1`1`1;600`400]
\putsquare(0,50)[(A,\alpha)`(A,\alpha)`O^\ddag`Q^\ddag;1_A`\partial^A`u^\ddag\circ \partial^A`u^\ddag]
 \end{picture}
\end{center}
we write $A_!$ for $(A\alpha,u^\ddag \circ\partial^A,Q)$ and $su_{u,A}=(1_A,u):A\ra A_!$.

We shall show that the natural transformations
 \begin{center} \xext=1600 \yext=200
\begin{picture}(\xext,\yext)(\xoff,\yoff)
\putmorphism(0,30)(1,0)[ A\star u^*(-)`u^*(A_!\star (-))`rep_s(su_{u,A})=rep_s(1_A,u)]{2000}{1}a
\end{picture}
\end{center}
and
 \begin{center} \xext=1600 \yext=150
\begin{picture}(\xext,\yext)(\xoff,\yoff)
\putmorphism(0,30)(1,0)[ A\star u^*(-)`u^*u_!A\star u^*(-)`su_{u,rep_s(A)} =\eta^u_{A\star u^*(-)}]{2000}{1}a
\end{picture}
\end{center}
are isomorphic in $_{A\star u^*(-) \backslash}Cat(Set_{/Q},Set_{/O})$.

We shall define a natural isomorphism
 \begin{center} \xext=1200 \yext=100
\begin{picture}(\xext,\yext)(\xoff,\yoff)
\putmorphism(0,0)(1,0)[u^*u_!(A\star u^*(-))`u^*(A\star (-))`\zeta]{1200}{1}a
\end{picture}
\end{center}
(and its inverse) so that $rep_s(1_A,u)=\zeta \circ \eta^u_{A\star u^*(-)}$, i.e. for any $d^Y:Y\ra Q$ the triangle
 \begin{center} \xext=1600 \yext=900
\begin{picture}(\xext,\yext)(\xoff,\yoff)
\putmorphism(1600,700)(0,-1)[u^*(A_!\star d^Y)`u^*u_!(A\star u^*(d^Y))`\zeta_{d^Y}]{600}{-1}r

\put(0,350){\makebox(200,100){$A\star u^*(d^Y)$}}
 \put(350,450){\vector(4,1){900}}
  \put(350,350){\vector(4,-1){850}}

  \put(680,650){\makebox(300,100){$rep_s(1_A,u)_{d^Y}$}}
   \put(660,110){\makebox(300,100){$\eta^u_{A\star u^*(d^Y)}$}}

   \put(2000,400){\makebox(300,100){$(3)$}}
\end{picture}
\end{center}
commutes. Using the notation from diagram $(1)$ we note that, for $a\in A$, we have
\[ [a,\vec{x}]_\sim \in A\star u^*(d^Y) \hskip 5mm {\rm iff} \hskip 5mm
\partial^{A,+}_a = u(d^Y)\circ \vec{x}\]

\[ \lk o,[a,\vec{y}]_\sim \rk\in u^*(A_!\star d^Y) \hskip 5mm {\rm iff} \hskip 5mm
o\in O,\; o=\partial^A_a(0),\;  u\circ\partial^{A,+}_a = d^Y\circ \vec{y}\]

\[ \lk o,[a,\vec{x}]_\sim \rk\in u^*u_!(A\star u^*(d^Y)) \hskip 5mm {\rm iff} \hskip 5mm
o\in O,\; u(o)=u\circ\partial^A_a(0),\;  u\circ\partial^{A,+}_a = d^Y\circ \vec{x}\]

Now we spell out explicitly the function occurring in $(3)$.
 \begin{center} \xext=1500 \yext=450
\begin{picture}(\xext,\yext)(\xoff,\yoff)
\putmorphism(0,300)(1,0)[A\star u^*(d^Y))`u^*(A_!\star d^Y)`rep_s(1_A,u)_{d^Y}]{1500}{1}a
\putmorphism(200,40)(1,0)[ [a,\vec{x}]`\lk \partial^A_a(0),[a,u^Y\circ \vec{x}]_\sim \rk `]{1100}{1}a
 \put(335,0){\line(0,1){80}}
\end{picture}
\end{center}
and
 \begin{center} \xext=1500 \yext=450
\begin{picture}(\xext,\yext)(\xoff,\yoff)
\putmorphism(0,300)(1,0)[A\star u^*(d^Y))`u^*u_!(A\star u^*(d^Y))`\eta^u_{A\star u^*(d^Y)}]{1500}{1}a
\putmorphism(200,40)(1,0)[ [a,\vec{x}]`\lk \partial^A_a(0),[a,\vec{x}]_\sim \rk `]{1100}{1}a
 \put(335,0){\line(0,1){80}}
\end{picture}
\end{center}
and
 \begin{center} \xext=1500 \yext=650
\begin{picture}(\xext,\yext)(\xoff,\yoff)
\putmorphism(0,500)(1,0)[ u^*u_!(A\star u^*(d^Y))`u^*(A_!\star d^Y)`\zeta_{d^Y}]{1500}{1}a
\putmorphism(200,240)(1,0)[ \lk o,[a,\vec{x}]_\sim \rk`\lk o,[a,u^Y\circ\vec{x}]_\sim \rk `]{1100}{1}a
 \put(450,200){\line(0,1){80}}

\putmorphism(200,40)(1,0)[ [ \lk o,[a,\vec{x}]_\sim \rk`\lk o,[a,\vec{y}]_\sim \rk `]{1100}{-1}a
 \put(1050,0){\line(0,1){80}}
 \end{picture}
\end{center}
In the last correspondence, $\vec{y}\mapsto \vec{x}:(|a|]\ra u^*(Y)$ is defined using the fact that right square in $(2)$ is a pullback and $d^Y\circ \vec{y}=u\circ\partial^{A,+}_a$.

Again, a simple check shows that $(3)$ commutes, i.e. $rep_s$ is a morphism of opfibrations, as well. $~\Box$
\vskip 2mm

Later we will show that $rep_s$ is faithful and full on isomorphisms.

The fibration that is the essential image of the representation  $rep_s$ will be denoted by $p_{an}:An\ra Set$,
we take it as the definition of the fibration of (multivariable) analytic (endo)functors and analytic transformations between them.
Thus by an {\em analytic functor} on $Set_{/O}$, where $O$ is a set, we understand a functor (isomorphic to one) of the form $A\star(-):Set_{/O}\ra Set_{/O}$ for a symmetric signature
$A=(A,\alpha,\partial^A,O)$. Moreover, by an {\em analytic transformation} over a function $u:O\ra Q$ between two analytic functors   $A\star(-):Set_{/O}\ra Set_{/O}$
and $B\star(-):Set_{/Q}\ra Set_{/Q}$ we mean a natural transformation of the form $rep_s(f,u) : A\star u^*(-)\lra  u^*(B\star(-))$ for a morphism
symmetric signatures $(f,u):(A,\alpha,\partial^A,O)\lra  (B,\beta,\partial^B,Q)$.

{\bf Remarks}
\begin{enumerate}
  \item Note that for the one element set, say $[0]$, we have $Set_{/[0]}\cong Set$
  and the fibre of $p_{an}$ over $[0]$ is (isomorphic to) the category of usual
(one-variable) analytic functors that was characterized in \cite{J2}, see also
\cite{AV}, as the category of finitary endofunctors on $Set$ that weakly
preserve pullbacks and
weakly cartesian natural transformations between them. We will see in  the next
section how this characterization extends from this fibre to the whole fibration.
  \item In \cite{J2} multivariable analytic functors were defined as certain functors
  $Set^I\ra Set$, were $I$ is a finite set.
Ignoring the `size problems', this definition can be extended to infinite sets $Q$
 by saying that the class(!) of multivariable analytic functors
$Set_{/Q}\ra Set$ is a cofiltered 'limit' of the classes(!) of analytic functors
$Set^I\ra Set$ where $I\subseteq Q$ and $I$ is finite. In these terms, what we call
an analytic functor on $Set_{/Q}$, is just a $Q$-tuple of multivariable analytic
functors $Set_{/Q}\ra Set$.
\item The multivariable analytic functors  $Set_{/Q}\ra Set$ can be also described more
explicitly avoiding all the `size problems'. Let $Q$ be a set and $q\in Q$. We have
the evaluation functor and the inclusion of a fibre\label{bi-ev-functors}
 \begin{center} \xext=1200 \yext=100
\begin{picture}(\xext,\yext)(\xoff,\yoff)
\putmorphism(0,0)(1,0)[Set_{/Q}`Set`ev_q]{600}{1}a
\putmorphism(600,0)(1,0)[\phantom{Set}`Set_{/Q}`\bi_q]{600}{1}a
\end{picture}
\end{center}
such that $ev_{q}(X,d)=d^{-1}(z)$  and $\bi_q(B) : B\ra Q$ is the function defined by $ \bi_q(B)(b)=z$ for $b\in B$.
A {\em multivariable analytic functor} from $Set_{/Q}$ to $Set$ is a functor of the form
 \begin{center} \xext=1600 \yext=100
\begin{picture}(\xext,\yext)(\xoff,\yoff)
\putmorphism(0,0)(1,0)[Set_{/Q}`Set_{/Q}`B\star(-)]{800}{1}a
\putmorphism(800,0)(1,0)[\phantom{Set_{/Q}}`Set`ev_q]{800}{1}a
\end{picture}
\end{center}
where $(B,\beta,\partial^B,Q)$ is a symmetric signature,  $z\in Q$ and  $ev_z$
is the evaluation functor such that $ev_z(d:X\ra Q)=d^{-1}(z)$.
If we restrict symmetric signatures to those for which $\partial^B_b(0)=z$
then such functors, as we shall see, determine the signatures up to an isomorphism.
\end{enumerate}

\subsection{A characterization of the fibration of analytic functors}\label{an_char}

The following extends the characterization of analytic functors and
analytic transformations, c.f. \cite{J2}, from the fibre over $[0]$
to the whole fibration of analytic functors $p_{an}$.

\begin{theorem} \label{analytic char}
The  lax monoidal fibration $p_{an}:An\lra Set$ has as its objects finitary
endofunctors on categories $Set_{/Q}$ that weakly preserve wide pullbacks
and  weakly cartesian natural transformations as morphisms between them.
\end{theorem}

Before we prove a series of lemmas needed to establish the above theorem,
we shall immediately present the following obvious Corollary that is even more interesting.

\begin{corollary} \label{analytic monads}
The fibration  of symmetric multicategories is equivalent to the fibration of
weakly cartesian analytic monads and weakly cartesian morphisms of monads whose
functor parts are pullback functors between them. Under this correspondence
the free symmetric multicategories correspond to the free analytic monads. $\Box$
\end{corollary}

The following definition is an extension of a notion from \cite{AV}. Let $Q$ be a set.
The functor $F:Set_{/Q}\ra Set$ is {\em superfinitary} if and only if there
is an object $d:I\ra Q$ in $Set_{/Q}$ with $I$ finite such that, for any $d^X:X\ra Q$ in $Set_{/Q}$
\[ F(X,d^X)= \bigcup_{f:(I,d)\ra (X,d^X)} F(f)(F(I,d)) \]
i.e. the elements of $F(I,d)$ generates the whole functor $F$.
The following two Lemmas and their proofs are 'colored' versions of Theorem 2.6 and Corollary 2.7  and their proofs from \cite{AV}.

\begin{lemma}\label{superfini}
Let $F:Set_{/Q}\ra Set$ be a superfinitary functor. Then $F$ is a multivariable analytic functor if and only if $F$ weakly preserves pullbacks.
\end{lemma}

{\it Proof.}~ By observations of V. Trnkov\'a, cf. \cite{T}, any functor $F:Set_{/Q}\ra Set$ is a coproduct (indexed by $F(1_Q)$)
of functors that preserves the terminal object. If $F=\coprod_i F_i$, then $F$ weakly preserves pullbacks (is analytic) if and only if
all $F_i$'s do (are). Thus it is enough to prove the lemma for a superfinitary functor $F:Set_{/Q}\ra Set$ such that $F(1_Q)=1$.

So suppose that $F$ weakly preserves pullbacks. Fix a minimal object $d^F:(n]\ra Q$ such that for any $d^X:X\ra Q$
 \[ F(X,d^X)= \bigcup_{f:((n],d^F)\ra (X,d^X)} F(f)(F((n],d^F)) \]
By 'minimal' we mean that there is no proper subobject of $((n],d^F)$ in $Set_{/Q}$ with this property.
Thus there is an element $g^F\in F(d^F)$ that it is not in the image of any proper inclusion into
$d^F$. The pair $(g^F,d^F)$ or just the element $g^F$ if $d^F$ is understood, will be called
{\em generic}, cf. \cite{J2}, \cite{AV}. Thus, if we have a morphism
\begin{center} \xext=800 \yext=440
\begin{picture}(\xext,\yext)(\xoff,\yoff)
 \settriparms[1`1`1;400]
  \putVtriangle(0,0)[X`(n]`Q;f`d^X`d^F]
\end{picture}
\end{center}
in $Set_{/Q}$ such that $g^F\in F(f)(X,d^X)$ then $f$ is onto.
Therefore, any endomorphism of $((n],d^F)$ in $Set_{/Q}$ leaving  $g^F$ fixed, is a bijection.
We can define a subgroup of $S_n$ as follows
\[ \cG^F =\{ \sigma: d^F\ra d^F \in Set_{/Q} :\, F(\sigma )(g^F)=g^F  \} \subseteq S_n \]
Let $F^o$ be ${S_n}_{/\cG^F}$, the set of right cosets of $S_n$ over $\cG^F$. The class of $\tau$ in  ${S_n}_{/\cG^F}$ will be denoted by
$[\tau ]_{\sim F}$. We have a right action, say $\phi$, of $S_n$ on  $F^o$ acting by composition on the right.
We define $\partial^{F^o}: F^o\ra Q^\ddag$ so that, for $\tau\in S_n$,
\[ \partial^{F^o,+}_{[\tau ]_{\sim F}}=d^F\circ \tau:(n]\ra Q. \]
and $\partial^{F^o}_{[\tau ]_{\sim F}}(0)=z$
where $z$ is any element of $Q$ (if $Q$ is empty there is nothing to prove).

The functor $ev_z\circ (F^o,\phi,\partial^{F^o},Q)\star(-)$ will be denoted
from now on simply as $F^o\star(-)$. For $d^X:X\ra Q$ in $Set_{/Q}$ we put
 \begin{center} \xext=1600 \yext=400
\begin{picture}(\xext,\yext)(\xoff,\yoff)
\putmorphism(0,300)(1,0)[F^o\star(d^X)`F(d^X)`\kappa^F_{(d^X)}]{1600}{1}a
\putmorphism(300,50)(1,0)[{[[\tau]_{\sim F},u:d^F\ra d^X ]}`F(u)(g^F)`]{1100}{1}a
   \put(765,0){\line(0,1){100}}
\end{picture}
\end{center}
We shall show that $\kappa^F:F^o\star(-)\lra F$ is a natural isomorphism.
To this aim it is enough to show
\begin{description}
  \item[$(a)$] for every $d^X:X\ra Q$ and $x\in F(X,d^X)$ there is $u:d^F \ra d^X$ in $Set_{/Q}$ such that $x=F(u)(g^F)$;
  \item[$(b)$] for $u,v:d^F\ra d^X$ we have, $F(u)(g^F)= F(v)(g^F)$ if and only if there is $\sigma\in \cG^F$ such that $u=v\circ \sigma$;
  \item[$(c)$] $\kappa^F$ is natural.
\end{description}

We  establish $(a)$ and $(b)$. Then, as $F$ is a functor, $(c)$ will be obvious.

Ad (a). It is enough to show $(a)$ for elements $x\in F((n],d^F)$. As $F$ weakly preserves pullbacks and $F(1)=1$,
$F$ weakly preserves binary products. We have a binary product
\begin{center} \xext=1000 \yext=600
\begin{picture}(\xext,\yext)(\xoff,\yoff)
 \settriparms[-1`1`1;500]
  \putqtriangle(0,0)[(n]`(n]\times_Q(n]`Q;\pi_1`d^F`d^\pi]
 \settriparms[1`0`1;500]
  \putptriangle(500,0)[\phantom{(n]\times_Q(n]}`(n]`\phantom{Q};\pi_2``d^F]
\end{picture}
\end{center}
of $d^F$ with itself in $Set_{/Q}$ and hence a weak product
 \begin{center} \xext=1600 \yext=100
\begin{picture}(\xext,\yext)(\xoff,\yoff)
\putmorphism(0,0)(1,0)[F(d^F)`F(d^\pi)`F(\pi_1)]{800}{-1}a
\putmorphism(800,0)(1,0)[\phantom{F(d^\pi)}`F(d^F)`F(\pi_2)]{800}{1}a
  \end{picture}
\end{center}
in $Set$. Hence there is $p\in F(d^\pi)$ such that
\[ F(\pi_1)(p)=g^F,\hskip 5mm  F(\pi_2)(p)=x \]
Since $F$ is superfinitary and by assumption on $d^F$, there are a morphism $f:d^F\ra d^\pi$ and $y\in F(d^F)$ such that
$F(f)(y)=p$. Thus $g^F=F(\pi_1)(p)=F(\pi_1\circ f)(y)$, and hence $\pi_1\circ f :d^F\ra d^F$ is epi and then iso.
Putting $u=\pi_2\circ f\circ (\pi_1\circ f)^{-1}: d^F\lra d^F$, we have
\[ x=F(\pi_2)=F(\pi_2\circ f)(y) = F(\pi_2\circ f\circ (\pi_1\circ f)^{-1})(g^F)=F(u)(g^F) \]
as needed.

Ad (b). If for some $\sigma\in \cG^F$ we have $u=v\circ \sigma$, then
\[ F(u)(g^F)= F(v\circ \sigma)(g^F)=  F(v)(F( \sigma)(g^F))= F(v)(g^F) \]

Now assume that $F(u)(g^F) = F(v)(g^F)$.  We form a pullback in $Set_{/Q}$
 \begin{center} \xext=700 \yext=500
\begin{picture}(\xext,\yext)(\xoff,\yoff)
\setsqparms[1`1`1`1;700`400]
\putsquare(0,50)[(P,d^P)`((n],d^F)`((n],d^F)`(X,d^X);\bar{v}`\bar{u}`u`v]
 \end{picture}
\end{center}
As $F$ weakly preserves pullbacks, there is $p\in F(d^P)$ such that $F(\bar{u})(p)=g^F=F(\bar{v})(p)$.
Using $(a)$ we get $f:((n],d^F) \ra (P,d^P)$ such that $F(f)(g^F)=p$. Thus $\bar{u}\circ f,\,\bar{v}\circ f \in \cG^F$. We have
\[ u = u\circ (\bar{v}\circ f)\circ (\bar{v}\circ f)^{-1} = v\circ (\bar{u}\circ f)\circ (\bar{v}\circ f)^{-1} \]
and $(\bar{u}\circ f)\circ (\bar{v}\circ f)^{-1}\in \cG^F$.

Now suppose that $F$ is a composition of functors
\begin{center} \xext=1600 \yext=100
\begin{picture}(\xext,\yext)(\xoff,\yoff)
\putmorphism(0,0)(1,0)[Set_{/Q}`Set_{/Q}`B\star(-)]{800}{1}a
\putmorphism(800,0)(1,0)[\phantom{B\star(-)}`Set`ev_z]{800}{1}a
  \end{picture}
\end{center}
for a symmetric signature $(B,\beta,\partial^B,Q)$ and $z\in Q$. Since $F(1_Q)=1$, $(B,\beta)$ has just one orbit.
We can assume that $\partial^B_b(0)=z$ for $b\in B$. Thus by a slight abuse we shall identify $B\star(-)$ with $F$.
Let
 \begin{center} \xext=700 \yext=500
\begin{picture}(\xext,\yext)(\xoff,\yoff)
\setsqparms[1`1`1`1;700`400]
\putsquare(0,50)[(P,d^P)`(Y,d^Y)`(X,d^X)`(Z,d^Z);\pi_Y`\pi_X`g`f]
 \end{picture}
\end{center}
be a pullback in $Set_{/Q}$. We need to show that
 \begin{center} \xext=1000 \yext=550
\begin{picture}(\xext,\yext)(\xoff,\yoff)
\setsqparms[1`1`1`1;1000`400]
\putsquare(0,50)[B\star P`B\star Y`B\star X`B\star Z;B\star \pi_Y`B\star \pi_X`B\star g`B\star f]
 \end{picture}
\end{center}
is a weak pullback. Fix $b\in B$ and let $n=|b|$. Suppose for
\begin{center} \xext=2000 \yext=440
\begin{picture}(\xext,\yext)(\xoff,\yoff)
 \settriparms[1`1`1;400]
  \putVtriangle(0,0)[(n]`X`Q;u`\partial^{B,+}_b`d^X]
  \putVtriangle(1200,0)[(n]`Y`Q;v`\partial^{B,+}_b`d^Y]
\end{picture}
\end{center}
we have
\[  [b,f\circ u]_\sim =(B\star f)([b,u]_\sim) = (B\star g)([b,v]_\sim) = [b,g\circ v]_\sim \]
i.e. there is $\sigma\in S_n$ such that
\[ b\cdot\sigma =b,\hskip 5mm  f\circ u \circ \sigma = g\circ v   \]
Using the property of the above pullback, we get $w:(n]\ra P$ such that
\[ \pi_X \circ w=u\circ \sigma ,\hskip 5mm \pi_Y \circ w=v  \]
Then
\[ B\star\pi_Y([b,w]_\sim) = [b,\pi_Y\circ w]_\sim =  [b,v]_\sim \]
and
\[ B\star\pi_X([b,w]_\sim) = [b,\pi_X\circ w]_\sim =  [b\cdot\sigma,u\circ \sigma]_\sim =  [b,u]_\sim ~~~\Box \]

Recall that the functor $F:Set_{/Q} \lra Set_{/Q}$ is {\em thin} if there is $z\in Q$ such that $F=\bi_q\circ ev_q \circ F$ and $ev_q\circ F(1)=1$.

Clearly, every functor $F:Set_{/Q}\lra Set_{/Q}$ is a coproduct of thin functors indexed by the domain of $F(1_Q)$ and every natural transformation
between such functors is a coproduct of transformations between thin functors. We have

\begin{lemma}\label{fini}
Let $F:Set_{/Q}\ra Set_{/Q}$ be a finitary functor. The following are equivalent
\begin{enumerate}
  \item $F$ is a multivariable analytic functor;
  \item $F$ weakly preserves wide pullbacks of power $\leq \aleph_0\, +\, |Q|$;
  \item $F$ weakly preserves wide pullbacks.
\end{enumerate}
\end{lemma}

{\it Proof.}~ As every functor $Set_{/Q}\ra Set_{/Q}$ is a coproduct of thin functors, we can assume that $F$ is thin. Thus we consider $F$ as a functor
$Set_{/Q}\ra Set$ such that $F(1)=1$.

$3. \Ra 2.$ is obvious and $1. \Ra 3.$ can be proved as in Lemma \ref{superfini}.

In order to show $2. \Ra 1.$ we shall show that $F$ is superfinitary and use Lemma \ref{superfini}. Suppose on the contrary that for each $n\in \o$ and $f:(n]\ra Q$ there is $d^f:X^f\ra Q$ and $x^f\in F(d^f)$ such that $x^f\not\in \bigcup_{h:f\ra d^f} F(h)(f)$.
Since $F$ weakly preserves pullbacks of power $\leq \aleph_0\, +\, |Q|$,  there is $p\in F(\prod_{n\in\o,\, f:(n]\ra Q } (X^f,d^f))$ such that $F(\pi^f)(p)=x^f$ for $n\in\o$ and $f:(n]\ra Q$. Since $F$ is finitary, there are
\begin{center} \xext=800 \yext=460
\begin{picture}(\xext,\yext)(\xoff,\yoff)
 \settriparms[1`1`1;400]
  \putVtriangle(0,0)[(m]`\prod_{n\in\o,\, f:(n]\ra Q } (X^f,d^f)`Q;g`f_0`]
\end{picture}
\end{center}
and $y\in F((m],f_0)$ such that $F(g)(y)=p$. Then  $\pi^{f_0}\circ g : ((m],f_0)\ra (X^{f_0},d^{f_0})$
and
\[ F(\pi^{f_0}\circ g)(y)=F(\pi^{f_0}(p))=x^{f_0} \]
contrary to the assumption.
 $~\Box$

The following fact is not needed for the proof of Theorem \ref{analytic char} but it follows easily from the proofs of the above Lemmas and
puts some light on the correspondence between orbits of $(B,\beta)$ and elements of $B\star(1_Q)$.

\begin{scholium}\label{one-orbit}
Let $F:Set_{/Q}\ra Set_{/Q}$ be a finitary functor that preserves wide pullbacks.
Then $F$ is isomorphic to a functor
\[ (B,\beta,\partial^B,Q)\star(-)  :Set_{/Q}\ra Set_{/Q} \]
for a symmetric signature $(B,\beta,\partial^B,Q)$ such that $(B,\beta)$
contains as many orbits as the cardinality of the domain $F(1_Q)$ orbits.
In particular, if $F$ is thin iff $(B,\beta)$ has exactly one orbit.
$~\Box$
\end{scholium}

{\bf Remark} In the proof of the above Lemma \ref{superfini} we introduced the notions of a  minimal object, a generic element and a generic  pair.
The notion of a generic element is a variant of a notion introduced in \cite{J2}, see also \cite{AV}. For the later use, we describe below such generic pairs
for the functors of form $B\star(-)$, where $(B,\beta)$ has one orbit.

\begin{lemma}\label{generic elt}
Let $(B,\beta,\partial^B,Q)$  be a symmetric signature such that $(B,\beta)$ has a single orbit, $b\in B$ and let $n=|b|$. Then the positive typing for $b$, i.e. $\partial^{B,+}_b :(n]\ra Q$ is the minimal object for functor $B\star(-)  :Set_{/Q}\ra Set_{/Q}$ and $[b,1_{(n]}]_\sim\in B\star((n],\partial^{B,+}_b)$ is a generic
element for $B\star(-)$. More generally, let $d^Y:Y\ra Q$, $\vec{y}:Y\ra Q$ and $b\in B$. Then $\lk b,\vec{y}\rk$ represents a generic element
$[b,\vec{y}]_\sim\in B\star d^Y$  for $B\star(-)$ if and only if $\vec{y}$ is a bijection and $\partial^{B,+}_b=d^Y\circ \vec{y}$.
\end{lemma}

{\it Proof.}~ Exercise. $~\Box$

\begin{lemma}\label{wpres-pb}
Let $(f,u):(A,\alpha,\partial^A,O)\lra (B,\beta,\partial^B,Q)$ be a morphism of symmetric signatures in $Sig_s$ over a function $u:O\ra Q$. Then the natural transformation
in $Cat(Set_{/Q}, Set_{/O})$ representing $(f,u)$
\[ rep_s(f,u): rep_s(A)\circ u^* \lra u^*\circ rep_s(B) \]
is weakly cartesian.
\end{lemma}

{\it Proof.}~ Thus we have a commuting square
 \begin{center} \xext=500 \yext=500
\begin{picture}(\xext,\yext)(\xoff,\yoff)
\setsqparms[1`1`1`1;500`400]
\putsquare(0,50)[A`B`O^\ddag`Q^\ddag;f`\partial^A`\partial^B`u^\ddag]
 \end{picture}
\end{center}
in the category of symmetric sets. Let
\begin{center} \xext=800 \yext=400
\begin{picture}(\xext,\yext)(\xoff,\yoff)
 \settriparms[1`1`1;350]
  \putVtriangle(0,0)[Y`Y'`Q;g`d^Y`d^{Y'}]
\end{picture}
\end{center}
be a morphism in $Set^\ra_{/Q}$ i.e. in $Set_{/Q}$. Then we pullback this morphism along $u$, and we get a diagram
 \begin{center} \xext=1500 \yext=900
\begin{picture}(\xext,\yext)(\xoff,\yoff)
\setsqparms[1`1`1`1;500`400]
\putsquare(550,450)[u^*(Y)`Y`u^*(Y')`Y';u^Y`u^*(g)`g`]
\setsqparms[0`1`1`1;500`400]
\putsquare(550,50)[\phantom{u^*(Y')}`\phantom{Y}`O`Q;u^{Y'}`u^*(d^{Y'})`d^{Y'}`u]

  \put(-30,380){\makebox(50,80){$u^*(d^Y)$}}
     \put(150,50){\vector(1,0){320}}
   \put(150,50){\line(0,1){800}}
   \put(150,850){\line(1,0){200}}

    \put(1520,380){\makebox(50,80){$d^Y$}}
     \put(1450,50){\vector(-1,0){320}}
   \put(1450,50){\line(0,1){800}}
   \put(1450,850){\line(-1,0){300}}
 \end{picture}
\end{center}
in which three squares are pullbacks. We need to show that the square
 \begin{center} \xext=1200 \yext=750
\begin{picture}(\xext,\yext)(\xoff,\yoff)
\setsqparms[1`1`1`1;1200`500]
\putsquare(0,100)[A\star u^*(Y)`u^*(B\star Y)`A\star u^*(Y')`u^*(B\star Y');
rep_s(f,u)_{d^Y}`1_A\star u^*(g)`u^*(1_B\star g)`{rep_s(f,u)}_{d^{Y'}}]
 \end{picture}
\end{center}
is a weak pullback.
So let $[a,\vec{x}]_\sim \in A\star u^*(Y')$ i.e.  $a\in A$, $\vec{x}:(|a|]\ra u^*(Y)$, so that $\partial^{A,+}_a=u^*(d^Y)\circ\vec{x}$
and let $\lk o,[b,\vec{y}]_\sim \rk \in u^*(B\star Y)$ i.e. $o\in O$, $b\in B$ and $\vec{y}:(|b|]\ra Q$, so that
$u(o)=\partial^B_b(0)$ and $\partial^{B,+}_ b=d^Y\circ \vec{y}$. Moreover, assume that
\[ {rep_s(f,u)}_{d^{Y'}}([a,\vec{x}]_\sim) =  \lk \partial^A_a(0),[f(a),u^{Y'}\circ \vec{x}]_\sim \rk =\]
\[=\lk o, [b,g\circ \vec{y}]_\sim \rk = u^*(1_B\star y')(\lk o,[b,\vec{y}]_\sim \rk )   \]
i.e. $\partial^A_a(0)=o$, and there is $\sigma\in S_{|a|}$ such that $f(a)=b\cdot \sigma$ and $u^{Y'}\circ\vec{x}=g\circ\vec{y}\circ \sigma$.
Thus, using the upper pullback above, we get a function $\vec{z}:(|a|]\ra u^*(Y)$ such that
\[ u^Y\circ \vec{z}=\vec{y}\circ \sigma,\hskip 15mm u^*(g)\circ \vec{z}=\vec{x}. \]
Then
\[ 1_A\star u^*(g)([a,\vec{z}]_\sim )= [a,u^*(g)\circ\vec{z}]_\sim =[a,\vec{x}]_\sim. \]
Moreover
\[ o=\partial^A_a(0),\;\;\; {\rm and}\;\;\;\;\;
\lk f(a),u^Y\circ\vec{z}\rk =  \lk b\cdot\sigma , \vec{y}\circ\sigma \rk \sim \lk b,\vec{y}\rk  \]
i.e. $rep_s(f,u)_{d^Y}([a,\vec{z}]_\sim )= \lk o,[b,\vec{y}]_\sim \rk $.  Thus the above square
is a weak pullback, as required.
 $~\Box$

\begin{lemma} \label{full-in-fib}
Let $(A,\alpha,\partial^A,O)$ and $(B,\beta,\partial^B,O)$ be two symmetric signatures in $(Sig_s)_O$. If $\xi : A\star(-)\lra B\star(-)$
 is a weakly cartesian natural transformation then there is a unique morphism of symmetric signatures $(f,1_O):(A,\alpha,\partial^A,O)\lra (B,\beta,\partial^B,O)$
 in $(Sig_s)_O$ such that $rep_s(f,1_O)=\xi$.
\end{lemma}

{\it Proof.}~ Let $(A,\alpha,\partial^A,O)$, $(B,\beta,\partial^B,O)$ be two symmetric signatures in $(Sig_s)_O$ and let $\xi : A\star(-)\lra B\star(-)$
  be a weakly cartesian natural transformation. By remark after Scholium \ref{one-orbit}, we can assume that both $(A,\alpha)$ and $(B,\beta)$ have one orbit. Then the existence of
 $\xi$ as above implies that for some $z\in O$ we have $\partial^A_a(0)=\partial^B_b(0)$ for all $a\in A$ and $b\in B$. Hence we shall not consider these values any more.
Fix $d^A:(n^A]\ra O$ and $x^A\in A\star d^A$ so that $(d^A,x^A)$ is a generic pair for the functor $A\star(-)$, and
$d^B:(n^B]\ra O$ and $x^B\in B\star d^B$ so that $(d^B,x^B)$ is a generic pair for the functor $B\star(-)$. Thus there is a morphism
\begin{center} \xext=800 \yext=400
\begin{picture}(\xext,\yext)(\xoff,\yoff)
 \settriparms[1`1`1;350]
  \putVtriangle(0,0)[(n^B]`(n^A]`O;u`d^B`d^A]
\end{picture}
\end{center}
in $Set_{/O}$, such that $B\star u (x^B)=\xi_{d^A}(x^A)$. Since $\xi$ is weakly cartesian, the square
 \begin{center} \xext=1000 \yext=500
\begin{picture}(\xext,\yext)(\xoff,\yoff)
\setsqparms[1`1`1`1;1000`400]
\putsquare(0,50)[A\star d^B`B\star d^B`A\star d^A`B\star d^A;\xi_{d^B}`A\star u `B\star u`\xi_{d^A}]
 \end{picture}
\end{center}
is a weak pullback and there is $x\in A\star d^B$, such that
\[ \xi_{d^B}(x)=x^B, \;\;\;\; {\rm and}\;\;\;\;  A\star u(x)=x^A. \]
 Since $x^A$ is generic for $A\star (-)$ there is a morphism
 \begin{center} \xext=800 \yext=400
\begin{picture}(\xext,\yext)(\xoff,\yoff)
 \settriparms[1`1`1;350]
  \putVtriangle(0,0)[(n^A]`(n^B]`O;v`d^A`d^B]
\end{picture}
\end{center}
in $Set_{/O}$, such that $A\star v (x^A)=x$. Thus $A\star (u\circ v)(x^A) = x^A$ and $u\circ v$ is iso. By naturality of $\xi$
\[ B\star v(\xi_{d^A}(x^A))=x^B \]
and $B\star (v\circ u)(x^B) = x^B$ and $v\circ u$ is iso.  Therefore, both $u$ and $v$ are bijections, and $n^A=n^B=n$.
Thus we can assume that $u$ is an identity, $d^A=d^B=d$ and $x^B=\xi_{d}(x^A)$.
Moreover, by Lemma \ref{generic elt}, we can assume that $x^A=[a,1_{(n]}]_\sim$ for some $a\in A$, $n=|a|$ and $d=\partial^{A,+}_a$.
Furthermore, again by Lemma \ref{generic elt}, we can assume that $\xi_{d}(x^A)=x^B=[b,\vec{x}]_\sim$ for some
$b\in B$ and a bijection $\vec{x}:(n]\ra (n]$ such that $\partial^{B,+}_b=\partial^{A,+}_a\circ \vec{x}$.

Thus $\partial^{B,+}_{b\cdot (\vec{x})^{-1}}= \partial^{B,+}_b\circ (\vec{x})^{-1}=\partial^{A,+}_a$. Hence the association $a\mapsto b\cdot (\vec{x})^{-1}$
extends to a morphism of symmetric signatures $(f,1_O):(A,\alpha,\partial^A,O)\ra (B,\beta,\partial^B,O)$ such that
$f(a\cdot\sigma)= (b\cdot (\vec{x})^{-1})\cdot\sigma$. We shall show that $rep_s(f,1_O)=\xi$.

First note that
\[   rep_s(f,1_O)_d(x^A)=rep_s(f,1_O)([a,1_{(n]}]_\sim )=[f(a),1_{(n]}]_\sim = \]
\[ = [b\cdot (\vec{x})^{-1},1_{(n]}]_\sim =[b,\vec{x}]_\sim = x^B = \xi_{d}(x^A) \]
i.e.  $rep_s(f,1_O)$ and $\xi$ agree on $x^A$. Now let $d^X:X\ra O$ and $x\in A\star d^X$  be arbitrary. Since $x^A$ is generic we have a morphism $w:d\ra d^X$
such that $A\star w(x^A) = x$. Now using the naturality of $\xi$ and $ rep_s(f,1_O)$ on $w$, i.e. serial commutativity of the diagram
 \begin{center} \xext=1000 \yext=720
\begin{picture}(\xext,\yext)(\xoff,\yoff)
\setsqparms[0`1`1`0;1200`500]
\putsquare(0,130)[A\star d`B\star d`A\star d^X`B\star d^X;`A\star w `B\star w`]

\putmorphism(0,160)(1,0)[\phantom{A\star d^X}`\phantom{B\star d^X}`\xi_{d^X}]{1200}{1}a
\putmorphism(0,100)(1,0)[\phantom{A\star d^X}`\phantom{B\star d^X}`rep_s(f,1_O)_{d^X}]{1200}{1}b

\putmorphism(0,660)(1,0)[\phantom{A\star d}`\phantom{B\star d}`\xi_{d}]{1200}{1}a
\putmorphism(0,600)(1,0)[\phantom{A\star d}`\phantom{B\star d}`rep_s(f,1_O)_{d}]{1200}{1}b
\end{picture}
\end{center}
we get that $rep_s(f,1_O)_{d^X}(x)  =\xi_{d^X}(x)$ and hence $rep_s(f,1_O)=\xi$.

If $(g,1_O):(A,\alpha,\partial^A,O)\ra (B,\beta,\partial^B,O)$ is another morphism of symmetric signatures such that $rep_s(g,1_O)=\xi$,
then in particular
\[ [g(a),1_{(n]}]_\sim =rep_s(g,1_O)_d(x^A) = \xi_{d}(x^A)= rep_s(f,1_O)_d(x^A) =[f(a),1_{(n]}]_\sim \]
This implies that $g(a)=f(a)$ and, since $(A,\alpha)$ has one orbit, $(g,1_O)=(f,1_O)$, as required.
$~\Box$

\vskip 2mm

{\it Proof of Theorem \ref{analytic char}.}~ From Lemma \ref{fini} we know that the objects in the essential image of the representations
$rep_s$ are finitary functors that weakly preserve wide pullbacks. From Lemma \ref{wpres-pb}
we know that morphisms in the essential image are weakly cartesian natural transformations.
Let $\xi:A\star(-) \lra B\star(-)$ be a morphism in $Exp(Set)$ over $u:O\ra Q$ which is a weakly cartesian natural transformation.
By Proposition \ref{reps-is-bifib} $rep_s$ is a prone morphisms of fibrations.
Hence $\xi$ can be factored in $Exp(Set)$, in an essentially unique way,  via a prone morphism  $rep_s(pr_{u,B}):u^*(B)\star(-)\ra B\star(-)$ in $p_{exp}$ and vertical morphisms $\xi':A\star(-)\ra u^*(B)\star(-)$ in the fibre over $O$, so that
$\xi=rep_s(pr_{u,B}) \circ \xi'$. From Proposition \ref{subfib of exp} both morphisms $pr_{u,B\star(-)}$, $\xi'$ are weakly cartesian. As
again by Lemma \ref{full-in-fib} $rep_s$ is faithful and full on weakly cartesian arrows in fibres, we obtain that $rep_s$
is faithful and full on weakly cartesian arrows in the whole fibration.
$~\Box$

\subsection{The analytic diagrams vs analytic functors}\label{an-vs-poly}

In section 6, we have shown that the concepts of amalgamated signature, polynomial diagram and polynomial functor are
equivalent when organized into lax monoidal fibrations. The signatures are the most explicit and the functors are the most abstract among these
concepts. The diagrams constitute a useful and important link between them.
In section \ref{act-sigs}, we have described the direct connection between lax monoidal fibrations of symmetric signatures and
of analytic functors. We provide here the missing link in this approach, the analytic diagrams. They correspond to analytic functors in much the same way
as polynomial diagrams correspond to polynomial functors.  In fact, these representing diagrams will constitute a full subcategory of the category of polynomial
diagrams in the category of symmetric sets $\sigma Set$. However the monoidal structure is {\em not} inherited from
$p_{poly,\sigma Set}:\cP oly\cD iag(\sigma Set)\ra \sigma Set$.
In the remaining part of the paper we want to indicate the relevant definitions and the obvious facts
leaving a more comprehensive study of the analytic diagrams to another paper.

Recall that the diagonal functor $\delta : Set\ra \sigma Set$ induced by the unique functor $\cS_*\ra 1$, has both adjoints
$orb\dashv\delta\dashv fix$. This adjunction can be also sliced. For details concerning all such functors that we shall consider in the following see the Appendix.

By a {\em pseudo-analytic diagram (over a set $O$)} in $Set$ we mean a diagram in $\sigma Set$
\begin{center} \xext=1800 \yext=150
\begin{picture}(\xext,\yext)(\xoff,\yoff)
\putmorphism(600,0)(1,0)[(E,\varepsilon)`(A,\alpha)`p]{600}{1}a
\putmorphism(0,0)(1,0)[\delta(O)`\phantom{(E,\varepsilon)}`s]{600}{-1}a
\putmorphism(1200,0)(1,0)[\phantom{(A,\alpha)}`\delta(O)`t]{600}{1}a
 \end{picture}
\end{center}
such that the fibres of $p$ are finite.

The object $O$ is an object of types of the polynomial $(t,p,s)$.
A {\em morphism of pseudo-analytic diagrams (over a function $u:O\ra Q$)} is a triple $(f,g,u)$,
with $f$ and $g$ morphisms in $\sigma Set$, making the diagram
\begin{center} \xext=1800 \yext=600
\begin{picture}(\xext,\yext)(\xoff,\yoff)
\setsqparms[1`1`1`1;600`400]
\putsquare(600,50)[(E,\varepsilon)`(A,\alpha)`(E',\varepsilon')`(A',\alpha');p`g`f`p']
\setsqparms[-1`1`0`-1;600`400]
\putsquare(0,50)[\delta(O)`\phantom{(E,\varepsilon)}`\delta(Q)`\phantom{(E',\varepsilon')};s`\delta(u)``s']
\setsqparms[1`0`1`1;600`400]
\putsquare(1200,50)[\phantom{(A,\alpha)}`\delta(O)`\phantom{(A',\alpha')}`\delta(Q);t``\delta(u)`t']
 \end{picture}
\end{center}
commute and the square in the middle is a pullback. Morphisms of pseudo-analytic diagrams compose in the obvious way, by putting one on top of the other.

Let $\cP\cA n\cD iag$  denotes the category of pseudo-analytic diagrams and morphisms between them. The category of analytic diagrams
$\cA n\cD iag$ is the slice $\cP\cA n\cD iag$ over the pseudo-analytic diagram $\cR$
\begin{center} \xext=1800 \yext=100
\begin{picture}(\xext,\yext)(\xoff,\yoff)
\putmorphism(0,0)(1,0)[\delta(1)`(R,\rho)`]{600}{-1}a
\putmorphism(600,0)(1,0)[\phantom{(R,\rho)}`\delta(1)`]{600}{1}a
\putmorphism(1200,0)(1,0)[\phantom{\delta(1)}`\delta(1)`]{600}{1}a
\end{picture}
\end{center}
where  $R_n =\{n\}\times (n]$ with action $\rho_n(\lk n,i\rk,\tau)=\lk n,\tau^{-1}(i)\rk$ for $\tau\in S_n$.
As $\delta(1)$ is the terminal symmetric set all the
morphisms in the above diagram are uniquely determined.
Thus an analytic diagram is a pseudo-analytic diagram
\begin{center} \xext=1800 \yext=500
\begin{picture}(\xext,\yext)(\xoff,\yoff)
\setsqparms[1`1`1`1;600`400]
\putsquare(600,50)[(E,\varepsilon)`(A,\alpha)`(R,\rho)`\delta(1);p```]
\setsqparms[-1`0`0`0;600`400]
\putsquare(0,50)[\delta(O)`\phantom{(E,\varepsilon)}``\phantom{(R,\rho)};```]
\setsqparms[1`0`0`0;600`400]
\putsquare(1200,50)[\phantom{(A,\alpha)}`\delta(O)`\phantom{\delta(1)}`;```]
 \end{picture}
\end{center}
so that $p$ is a pullback of $(R,\rho)\lra \delta(1)$ along $(A,\alpha)\lra \delta(1)$. Thus the fibre of $p$ over $a\in A$
can and will be identified as
$\{ a\}\times {(|a|]}$ with the action of $S_{|a|}$ so that $\lk a, i\rk \cdot \tau = \lk a\cdot \tau, \tau^{-1}(i)\rk$,
for $\tau \in S_{|a|}$. With such an identification if $(f,g,u): (t,p,s)\ra (t',p',s')$ is a morphism of analytic diagrams
then $g(a,i) =\lk f(a), i\rk$, for $\lk a,i\rk \in E$. Hence we shall not specify $g$ in the morphism of analytic diagrams anymore
and we shall denote it when necessary as $\bar{f}$.

We have an obvious projection functor
\[ p_{and}:\cA n\cD iag\lra Set,\]
sending $(f,u)$ to $u$, which is a lax monoidal fibration. However the tensor in fibres is {\em not} the one induced
by the tensor of those diagrams as if they were
polynomial diagrams in $\cP oly\cD iag(\sigma Set)\ra \sigma Set$. Its will be described below in an indirect way.

As in the case of the polynomial diagram fibration, the fibration of analytic diagrams comes equipped
with a representation morphism into the exponential
fibration $Exp(Set)\ra Set$. The representation functor
 \begin{center} \xext=600 \yext=520
\begin{picture}(\xext,\yext)(\xoff,\yoff)
\settriparms[1`1`1;450]
  \putVtriangle(0,0)[\cA n \cD iag`Exp(Set)`Set;rep_{and}`p_{and}`]
 \end{picture}
\end{center}
is defined as follows. For an analytic diagram $(t,p,s)$ over $O$ as displayed above,
we define a functor $rep_{and}(t,p,s)$  from $Set_{/O}$ to $Set_{/O}$
as the composition of five functors
\begin{center} \xext=3500 \yext=200
\begin{picture}(\xext,\yext)(\xoff,\yoff)
\putmorphism(0,0)(1,0)[Set_{/O}`\phantom{\sigma Set_{/\delta(O)}}`\delta_{/O}]{700}{1}a
\putmorphism(700,0)(1,0)[\sigma Set_{/\delta(O)}`\phantom{\sigma Set_{/(E,\varepsilon)}}`s^*]{700}{1}a
\putmorphism(1400,0)(1,0)[\sigma Set_{/(E,\varepsilon)}`\phantom{\sigma Set_{/(A,\alpha)}}`p_*]{700}{1}a
\putmorphism(2100,0)(1,0)[\sigma Set_{/(A,\alpha)}`\phantom{\sigma Set_{/\delta(O)}}`t_!]{700}{1}a
\putmorphism(2800,0)(1,0)[\sigma Set_{/\delta(O)}`Set_{/O}`orb_{/O}]{700}{1}a
 \end{picture}
\end{center}
i.e. we take a diagonal functor $\delta_{/O}$ to move from the 'set context' to the 'symmetric set context',
then we apply the usual polynomial functor over
the category of symmetric sets and we come back again to the set context via $orb_{/O}$,
by taking orbits of whatever we collected on the way.

We describe it with more details in three steps. For an object $(X,d^X:X\ra O)$ in $Set_{/O}$ the domain of  $p_* s^* (\delta_{/O})(X,d^X)$
in $\sigma Set_{/(A,\alpha)}$ is the set
\[  \{ \lk a,\vec{x}\rk : a\in A,\; \vec{x}:p^{-1}(a)\ra X,\; d^X\circ \vec{x}= \pi_O\circ s_{\lc p^{-1}(a)} \} \]
where $\pi_O:\delta(O)=\o\times O\ra O$ is the obvious projection, equipped with an action $\xi$ acting by conjugation, i.e. for $\sigma\in S_{|a|}$
\[ \xi(\lk a,\vec{x}\rk, \sigma)=\lk a,\vec{x}\rk\cdot \sigma = \lk a\cdot \sigma,\vec{x}((-)\cdot\sigma^{-1})\cdot\sigma\rk \]
The typing function sends $\lk a,\vec{x}\rk$ to $a\in A$.
The functor $t_!$ changes only the typing, i.e. $t_!p_* s^* (\delta_{/O})(X,d^X)$
in $\sigma Set_{/\delta(O)}$ has typing sending $\lk a,\vec{x}\rk$ to $\lk |a|,t(a) \rk\in \delta(O)$.
Finally, $orb_{/O}$ associates the orbits to what we've got so far, i.e.  the domain of $(orb_{/O})t_!p_* s^* (\delta_{/O})(X,d^X)$
in $Set_{/O}$ is the set of equivalence classes $[a,\vec{x}]_\sim$ of pairs $\lk a,\vec{x}\rk$ as above, divided by the action $\xi$.
The value of this morphism on the class $[a,\vec{x}]_\sim$ is sent to $\pi_Ot(a)\in O$.

For a morphism of analytic diagrams  $(f,u): (t,p,s)\lra (t',p',s')$ over $u$ as defined above, we define a morphism in $Exp(Set)$ over $u$, i.e.
a natural transformation
\[ rep_{and}(f,u) : (orb_{/O}) t_! p_* s^* (\delta_{/O}) u^* \lra  u^* (orb_{/Q}) t'_! p'_* s'^* (\delta_{/Q})   \]
using the diagram
\begin{center} \xext=3500 \yext=700
\begin{picture}(\xext,\yext)(\xoff,\yoff)
\putmorphism(0,600)(1,0)[Set_{/O}`\phantom{\sigma Set_{/\delta(O)}}`\delta_{/O}]{700}{1}a
\putmorphism(700,600)(1,0)[\sigma Set_{/\delta(O)}`\phantom{\sigma Set_{/(E,\varepsilon)}}`s^*]{700}{1}a
\putmorphism(1400,600)(1,0)[\sigma Set_{/(E,\varepsilon)}`\phantom{\sigma Set_{/(A,\alpha)}}`p_*]{700}{1}a
\putmorphism(2100,600)(1,0)[\sigma Set_{/(A,\alpha)}`\phantom{\sigma Set_{/\delta(O)}}`t_!]{700}{1}a
\putmorphism(2800,600)(1,0)[\sigma Set_{/\delta(O)}`Set_{/O}`orb_{/O}]{700}{1}a

\putmorphism(0,600)(0,-1)[\phantom{\sigma Set_{/\delta(O)}}`\phantom{\sigma Set_{/\delta(O)}}`u^*]{600}{-1}l
\putmorphism(700,600)(0,-1)[\phantom{\sigma Set_{/\delta(O)}}`\phantom{\sigma Set_{/\delta(O)}}`\delta(u)^*]{600}{-1}l
\putmorphism(1400,600)(0,-1)[\phantom{\sigma Set_{/\delta(O)}}`\phantom{\sigma Set_{/\delta(O)}}`\bar{f}^*]{600}{-1}l
\putmorphism(2000,600)(0,-1)[\phantom{\sigma Set_{/\delta(O)}}`\phantom{\sigma Set_{/\delta(O)}}`f^*]{600}{-1}l
\putmorphism(2200,600)(0,-1)[\phantom{\sigma Set_{/\delta(O)}}`\phantom{\sigma Set_{/\delta(O)}}`f_!]{600}{1}r
\putmorphism(2800,600)(0,-1)[\phantom{\sigma Set_{/\delta(O)}}`\phantom{\sigma Set_{/\delta(O)}}`\delta(u)_!]{600}{1}r
\putmorphism(3400,600)(0,-1)[\phantom{\sigma Set_{/\delta(O)}}`\phantom{\sigma Set_{/\delta(O)}}`u_!]{600}{1}l
\putmorphism(3500,600)(0,-1)[\phantom{\sigma Set_{/\delta(O)}}`\phantom{\sigma Set_{/\delta(O)}}`u^*]{600}{-1}r

\putmorphism(0,0)(1,0)[Set_{/Q}`\phantom{\sigma Set_{/\delta(Q)}}`\delta_{/Q}]{700}{1}a
\putmorphism(700,0)(1,0)[\sigma Set_{/\delta(Q)}`\phantom{\sigma Set_{/(E',\varepsilon')}}`s'^*]{700}{1}a
\putmorphism(1400,0)(1,0)[\sigma Set_{/(E',\varepsilon')}`\phantom{\sigma Set_{/(A',\alpha')}}`p'_*]{700}{1}a
\putmorphism(2100,0)(1,0)[\sigma Set_{/(A',\alpha')}`\phantom{\sigma Set_{/\delta(Q)}}`t'_!]{700}{1}a
\putmorphism(2800,0)(1,0)[\sigma Set_{/\delta(Q)}`Set_{/Q}`orb_{/Q}]{700}{1}a

\put(2050,350){$\varepsilon^f$}
\put(2050,260){$\Da$}
\end{picture}
\end{center} as follows.
By  adjunction $u_!\dashv u^*$ it is enough to define a natural transformation between these functors
\[  u_!(orb_{/O}) t_! p_* s^* (\delta_{/O}) u^* \lra (orb_{/Q}) t'_! p'_* s'^* (\delta_{/Q})   \]
and using the commutativity of some squares (including Beck-Chevalley condition) we define a natural transformation between functors isomorphic to those above
\[ (orb_{/Q}) t'_! (\varepsilon^f)_{p'_* s'^* (\delta_{/Q})}   : (orb_{/Q}) t'_! f_!f^*p'_* s'^* (\delta_{/Q})  \lra (orb_{/Q}) t'_! p'_* s'^* (\delta_{/Q})  \]
Tracing this definition back through the adjunctions we find that the so defined natural transformation $rep_{and}(f,g,u)$ applied to an object
$(Y,d^Y:Y\ra Q)$ in $Set_{/Q}$ sends the element
\[ [a,\vec{y}:p^{-1}(a)\ra u^*(Y)]_\sim \]
in the domain of $(orb_{/O}) t_! p_* s^* (\delta_{/O}) u^*(Y,d^Y)$    to the element
\[ \lk t(a),[f(a),u^Y\circ \vec{y}]_\sim \rk \]
in the domain of $u^* (orb_{/Q}) t'_! p'_* s'^* (\delta_{/Q})(Y,d^Y)$, where the notation is as in the following diagram
 \begin{center} \xext=1750 \yext=500
\begin{picture}(\xext,\yext)(\xoff,\yoff)
\setsqparms[-1`1`1`-1;450`390]
\putsquare(0,50)[Y`u^*(Y)`Q`O;u^Y`d^Y`u^*(d^Y)`u]
\putmorphism(450,50)(1,0)[\phantom{O}`\o\times O`\pi_O]{550}{-1}b
\putmorphism(1000,50)(1,0)[\phantom{\o\times O`}`p^{-1}(a)\subseteq E`s_{\lc p^{-1}(a)}]{800}{-1}b
 \put(1650,120){\line(-2,1){200}}
 \put(1050,320){\makebox(200,100){$\vec{y}$}}
 \put(1450,220){\vector(-4,1){850}}
 \end{picture}
\end{center}
We leave for the reader the verification that the so defined $rep_{and}$ is a lax morphism of fibrations.

In order to show that the essential image of $rep_{and}$ is the fibration of analytic functors,
we shall define first a morphism (in fact an equivalence) of fibrations
 \begin{center} \xext=600 \yext=520
\begin{picture}(\xext,\yext)(\xoff,\yoff)
\settriparms[1`1`1;450]
  \putVtriangle(0,0)[Sig_s`\cA n \cD iag`Set;\iota_s`p_s`p_{and}]
 \end{picture}
\end{center}
To an object $(A,\alpha,\partial^A,O)$ in $Sig_s$, $\iota_s$ assigns an analytic diagram as follows
\begin{center} \xext=1900 \yext=150
\begin{picture}(\xext,\yext)(\xoff,\yoff)
\putmorphism(0,0)(1,0)[\delta(O)`\phantom{(E^A,\overline{\alpha})}`s^A]{600}{-1}a
\putmorphism(600,0)(1,0)[(E^A,\overline{\alpha})`(A,\alpha)`p^A]{700}{1}a
\putmorphism(1300,0)(1,0)[\phantom{(A,\alpha)}`\delta(O)`t^A]{600}{1}a
 \end{picture}
\end{center}
where
\[ E^A=\{ \lk a,i \rk : a\in A,\; i\in (|a|] \} \]
and
\[ \overline{\alpha}(\lk a,i \rk,\sigma)=\lk a,i \rk\cdot \sigma = \lk a\cdot \sigma, \cdot \sigma^{-1}(i) \rk \]
Moreover
\[ s^A(a,i)=\lk |a|,\partial^A_a(i) \rk,\hskip 1cm  p^A(a,i)=a, \hskip 1cm  t^A(a)=\lk |a|,\partial^A_a(0)\rk   \]
for $\lk a,i \rk\in E^A$, $a\in A$.

If $(f,u): (A,\alpha,\partial^A,O)\lra (A',\alpha',\partial^{A'},Q)$ is a morphism in $Sig_s$ over $u:O\ra Q$, then $\iota_s$ assigns to it
the following morphism of diagrams
\begin{center} \xext=1900 \yext=550
\begin{picture}(\xext,\yext)(\xoff,\yoff)
\setsqparms[1`1`1`1;700`400]
\putsquare(600,50)[(E^A,\overline{\alpha})`(A,\alpha)`(E'^{A'},\overline{\alpha'})`(A',\alpha');p^A`\bar{f}`f`p'^{A'}]
\setsqparms[-1`1`0`-1;600`400]
\putsquare(0,50)[\delta(O)`\phantom{(E^A,\overline{\alpha})}`\delta(Q)`\phantom{(E^{A'},\overline{\alpha'})};s^A`\delta(u)``s'^{A'}]
\setsqparms[1`0`1`1;600`400]
\putsquare(1300,50)[\phantom{(A,\alpha)}`\delta(O)`\phantom{(A',\alpha')}`\delta(Q);t^A``\delta(u)`t'^{A'}]
 \end{picture}
\end{center}
so that $\bar{f}(a,i) = \lk f(a),i\rk$ for $ \lk f(a),i\rk \in E^A$, as before.

\begin{proposition} \label{iota-equi}
The association $\iota_s$ is an equivalence of fibrations.
\end{proposition}

{\it Proof.}~ $\iota_s$ is in fact an isomorphism if we restrict only to those analytic diagrams
that
\begin{center} \xext=1800 \yext=150
\begin{picture}(\xext,\yext)(\xoff,\yoff)
\putmorphism(600,0)(1,0)[(E,\varepsilon)`(A,\alpha)`p]{600}{1}a
\putmorphism(0,0)(1,0)[\delta(O)`\phantom{(E,\varepsilon)}`s]{600}{-1}a
\putmorphism(1200,0)(1,0)[\phantom{(A,\alpha)}`\delta(O)`t]{600}{1}a
 \end{picture}
\end{center}
for which $(E,\varepsilon)$ is identified with $(A,\alpha)\times_{\delta(O)}(R,\rho)$. Thus
it is an equivalence indeed.
$~\Box$

\vskip 2mm

\begin{proposition} \label{iota-rep-comm}
The following triangle of morphisms of fibrations
 \begin{center} \xext=800 \yext=450
\begin{picture}(\xext,\yext)(\xoff,\yoff)
\settriparms[1`1`1;400]
  \putVtriangle(0,0)[Sig_s`\cA n \cD iag`Exp(Set);\iota_s`rep_s`rep_{and}]
 \end{picture}
\end{center}
commutes up to an isomorphism.
\end{proposition}

{\it Proof.}~ All the necessary items were defined. We shall check that the values of both functors on objects
(that are functors on slices of $Set$) agree on objects. The remaining details are left for the reader.

Let $A=(A,\alpha,\partial^A,O)$ be a symmetric signature and $X=(X,d^X)$ an object in $Set_{/O}$.
Both values $rep_s(A)(X)$ and  $rep_{and}\circ \iota_s(A)(X)$ are functions into the set $O$ whose
domains are (= can be identified) with the set
\[ \{ \lk a,\vec{x} \rk : a\in A,\; \vec{x}:(|a|]\ra X,\; d^X\circ \vec{x}=\partial^{A,+}_a \} \]
divided by an equivalence relation. In the former case the relation identifies the pair
$\lk a,\vec{x} \rk$ with $\lk a\cdot \sigma,\vec{x}\circ\sigma \rk$ for $\sigma\in S_{|a|}$.
In the second case the relation identifies the elements of the same orbit of the action such that for $\sigma\in S_{|a|}$
\[  \lk a,\vec{x} \rk \cdot \sigma = \lk a\cdot \sigma,\vec{x}((-)\cdot \sigma^{-1})\cdot \sigma \rk  \]
For $i\in (|a|]$, we have
\[ (\vec{x}((-)\cdot \sigma^{-1})\cdot \sigma )(i) = \vec{x}(i\cdot \sigma^{-1})\cdot \sigma =
\vec{x}(\sigma(i))\cdot \sigma = \vec{x}\circ\sigma(i) \]
The last equality follows from the fact that the action in $X$ (in fact $\delta(X)$) is constant.
Thus these equivalence relations are the same and hence the whole morphisms into $O$ sending
the equivalence class of $[a,\vec{x}]_\sim$ to $\partial^{A,+}_a(0)$ are the same.
$~\Box$

\vskip 2mm

As $\iota_s$ is an equivalence of fibrations by Proposition \ref{iota-equi} and the essential
image of $rep_s$ is (by definition) the fibration of analytic functors,
we get from the above Proposition \ref{iota-rep-comm}

\begin{corollary} \label{rep-and-image}
The essential image of the representation functor
 \begin{center} \xext=800 \yext=450
\begin{picture}(\xext,\yext)(\xoff,\yoff)
\settriparms[1`1`1;400]
  \putVtriangle(0,0)[\cA n \cD iag`Exp(Set)`Set;rep_{and}`p_{and}`]
 \end{picture}
\end{center}
is the fibration of analytic functors.
\end{corollary}

\subsection{Comparing the polynomial and the analytic approaches}

In Section \ref{am-sig} we have shown that the lax monoidal fibrations of amalgamated signatures, polynomial diagrams and polynomial functors are
equivalent. In the previous Subsection \ref{an-vs-poly} we have introduced the notion of an analytic diagram and we have shown that the lax monoidal fibrations of symmetric signatures, analytic diagrams, and analytic functors are equivalent. Thus in each case, we have three different ways of presenting essentially the same notion.
Below we compare these notion at all three levels, i.e. we shall define the missing functors $K_{sig}$, $K_{diag}$ and  natural transformations $\Phi$, $\Psi$ in
the following diagram
 \begin{center} \xext=800 \yext=1680
\begin{picture}(\xext,\yext)(\xoff,\yoff)
\setsqparms[1`1`1`0;800`400]
\putsquare(0,1200)[ Sig_a`Sig_s`\phantom{\cP oly\cD iag}`\phantom{\cA n\cD iag};K_{sig}`\iota_{a}`\iota_{s}`]
\setsqparms[1`1`1`0;800`400]
\putsquare(0,800)[\cP oly\cD iag`\cA n\cD iag`\phantom{\cP oly}`\phantom{\cA n};K_{diag}`rep_{pd}`rep_{and}`]
\setsqparms[1`1`1`0;800`400]
\putsquare(0,400)[\cP oly`\cA n`\phantom{Cart(Set)}`\phantom{wCart(Set)};K_{fu}```]
\settriparms[1`1`1;400]
  \putVtriangle(0,0)[Cart(Set)`wCart(Set)`Exp(Set);``]
\put(540,930){\makebox(200,100){$\Psi$}}
\put(560,850){\makebox(200,100){$\Ra$}}
\put(540,1330){\makebox(200,100){$\Phi$}}
\put(560,1250){\makebox(200,100){$\Ra$}}
 \end{picture}
\end{center}
All the arrows are morphisms of lax monoidal fibrations and of bifibrations over $Set$. The four named vertical arrows are equivalence of lax monoidal fibrations.
The five unnamed arrows are inclusions. The three named horizontal arrows are morphisms comparing signatures, diagrams, and functors, respectively.

So we begin by describing the functor
 \begin{center} \xext=800 \yext=450
\begin{picture}(\xext,\yext)(\xoff,\yoff)
 \settriparms[1`1`1;400]
  \putVtriangle(0,0)[Sig_a`Sig_s`Set;K_{sig}`p_a`p_s]
\end{picture}
\end{center}
The morphism $(f,\sigma,u) : (A,\partial^A:A\ra O^\dag,Q)\lra  (B,\partial^B,Q)$ in $Sig_a$ over $u$ is sent to a morphism
\[ (s(f,\sigma),u): (s(A),\alpha,\partial^{s(A)}:s(A)\ra O^\ddag,O)\lra (s(B),\beta,\partial^{s(B)},Q) \]
so that
\[  s(A)=\{ \lk a,\tau \rk : n\in\o,\; a\in A_n,\; \tau\in S_n \}  \]
and, for  $\lk a,\tau \rk \in s(A)$
\[ \partial^{s(A)}(a,\tau)= \partial^A_a\circ \tau: [|a|]\ra O,\hskip1cm     s(f,\sigma)(a,\tau)  = \lk f(a),\sigma_a^{-1}\circ \tau \rk \]
We have, for$ \lk a,\tau \rk\in s(A)$,
\[ \partial^{s(B)}\circ s(f,\sigma)(a,\tau) = \partial^{s(B)}(f(a),\sigma_a^{-1}\circ \tau)= \]
\[  = \partial_{f(b)}^B\circ \sigma_a^{-1}\circ \tau = u\circ \partial^A_{a}\circ \tau =u^\ddag \circ \partial^A (a,\tau)    \]
i.e. the square
 \begin{center} \xext=1000 \yext=520
\begin{picture}(\xext,\yext)(\xoff,\yoff)
\setsqparms[1`1`1`1;1000`400]
\putsquare(0,50)[(s(A),\alpha)`(s(B),\beta)`O^\ddag`Q^\ddag;s(f,\sigma)`\partial^{s(A)}`\partial^{s(B)}`u^\ddag]
 \end{picture}
\end{center}
commutes and $K_{sig}$ is a well defined functor. We note for the record

\vskip 2mm

\begin{proposition} \label{Ksig}
The functor $K_{sig}$ is full, faithful, and its essential image consists of those symmetric signatures that have free actions.
\end{proposition}

{\it Proof.}~ Simple check. $~\Box$

\vskip 2mm

Next we define the functor
 \begin{center} \xext=800 \yext=500
\begin{picture}(\xext,\yext)(\xoff,\yoff)
 \settriparms[1`1`1;400]
  \putVtriangle(0,0)[\cP oly \cD iag`\cA n \cD iag`Set;K_{diag}`p_{pd}`p_{and}]
\end{picture}
\end{center}
To a polynomial diagram
\begin{center} \xext=1500 \yext=200
\begin{picture}(\xext,\yext)(\xoff,\yoff)
  \putmorphism(0,50)(1,0)[O`E`s]{500}{-1}a
  \putmorphism(500,50)(1,0)[\phantom{E}`A`p]{500}{1}a
  \putmorphism(1000,50)(1,0)[\phantom{A}`O`t]{500}{1}a
\end{picture}
\end{center}
$K_{diag}$ associates an analytic diagram
\begin{center} \xext=1800 \yext=150
\begin{picture}(\xext,\yext)(\xoff,\yoff)
\putmorphism(600,0)(1,0)[(\tilde{E},\tilde{\varepsilon})`(\tilde{A},\tilde{\alpha})`\tilde{p}]{600}{1}a
\putmorphism(0,0)(1,0)[\delta(O)`\phantom{(\tilde{E},\tilde{\varepsilon})}`\tilde{s}]{600}{-1}a
\putmorphism(1200,0)(1,0)[\phantom{(\tilde{A},\tilde{\alpha})}`\delta(O)`\tilde{t}]{600}{1}a
 \end{picture}
\end{center}
so that
\[  \tilde{A} = \{ \lk a,h \rk : a\in A,\; h:(|a|]\stackrel{\cong}{\lra} p^{-1}(a)  \} \]
\[  \tilde{E} = \{ \lk a,h,i \rk : a\in A,\; h:(|a|]\stackrel{\cong}{\lra} p^{-1}(a),\; i\in (|a|]  \} \]
where $|a|$ is the number of elements of $p^{-1}(a)$. For $\lk a,h \rk\in \tilde{A}$ we have
\[ \tilde{\alpha}(\lk a,h\rk,\tau)=\lk a,h\circ\tau \rk,\hskip 1cm \tilde{t}(a,h)=\lk |a|,t(a)\rk \]
and for  $\lk a,h,i \rk\in \tilde{E}$,
\[ \tilde{\varepsilon}(\lk a,h,i\rk,\tau)=\lk a,h\circ\tau,\tau^{-1}(i) \rk,\hskip 5mm \tilde{p}(a,h,i)=\lk a, h\rk, \hskip 5mm \tilde{s}(a,h,i)=\lk |a|, s\circ h(i)\rk. \]
To a morphism of polynomial diagrams
\begin{center} \xext=1500 \yext=500
\begin{picture}(\xext,\yext)(\xoff,\yoff)
\setsqparms[1`1`1`1;500`400]
\putsquare(500,50)[E`A`E'`A';p`g`f`p']
\setsqparms[-1`1`0`-1;500`400]
\putsquare(0,50)[O`\phantom{A}`Q`\phantom{E'};s`u``s']
\setsqparms[1`0`1`1;500`400]
\putsquare(1000,50)[\phantom{B}`O`\phantom{'}`Q;t``u`t']
 \end{picture}
\end{center}
$K_{diag}$ associates a morphism of analytic diagrams
\begin{center} \xext=1800 \yext=550
\begin{picture}(\xext,\yext)(\xoff,\yoff)
\setsqparms[-1`1`0`-1;600`400]
\putsquare(0,50)[\delta(O)`\phantom{(\tilde{E},\tilde{\varepsilon})}`\delta(O)`\phantom{(\tilde{E'},\tilde{\varepsilon'})};
\tilde{s}`\delta(u)``\tilde{s'}]
\setsqparms[1`1`1`1;800`400]
\putsquare(600,50)[(\tilde{E},\tilde{\varepsilon})`(\tilde{A},\tilde{\alpha})`(\tilde{E'},\tilde{\varepsilon'})`(\tilde{A'},\tilde{\alpha'});
\tilde{p}`\tilde{g}`\tilde{f}`\tilde{p'}]
\setsqparms[1`0`1`1;600`400]
\putsquare(1400,50)[\phantom{(\tilde{A},\tilde{\alpha})}`\delta(O)`\phantom{(\tilde{A'},\tilde{\alpha'})}`\delta(Q);
\tilde{t}``\delta(u)`\tilde{t'}]
 \end{picture}
\end{center}
so that, for $\lk a,h\rk\in \tilde{A}$ and $\lk a,h,i\rk\in \tilde{E}$
\[ \tilde{f}(a,h)= \lk f(a),g_{\lc p^{-1}(a)}\circ h\rk,\hskip 5mm \tilde{g}(a,h,i)= \lk f(a),g_{\lc p^{-1}(a)}\circ h, i \rk. \]
This ends the definition of $K_{diag}$.

Now we shall define the natural isomorphism $\Phi$. Fix $(A,\partial^A,O)$ in $Sig_a$. We need to define a morphism
\[ \Phi_{(A,\partial^A,O)} :  K_{diag}\circ \iota_a (A,\partial^A,O)\lra \iota_s\circ K_{sig}(A,\partial^A,O) \]
in $\cA n \cD iag$ in the fibre over $O$, i.e. a morphism of analytic diagrams
\begin{center} \xext=2800 \yext=550
\begin{picture}(\xext,\yext)(\xoff,\yoff)
\setsqparms[-1`1`0`-1;600`400]
\putsquare(800,50)[\delta(O)`\phantom{(\widetilde{E^A},\widetilde{\varepsilon})}`\delta(O)`\phantom{(E^{s(A)},\tilde{\varepsilon})};
\widetilde{s^A}`\delta(1_O)``s^{s(A)}]
\setsqparms[1`1`1`1;800`400]
\putsquare(1400,50)[(\widetilde{E^A},\widetilde{\varepsilon})`(\tilde{A},\tilde{\alpha})`(E^{s(A)},\overline{\alpha})`(s(A),\alpha);
\widetilde{p^A}`\Phi_1`\Phi_0`p^{s(A)}]
\setsqparms[1`0`1`1;600`400]
\putsquare(2200,50)[\phantom{(\tilde{A},\tilde{\alpha})}`\delta(O)`\phantom{(s(A),\alpha)}`\delta(O);
\widetilde{t^A}``\delta(1_O)`t^{s(A)}]
\putmorphism(0,450)(0,-1)[ K_{diag}\circ \iota_a (A,\partial^A,O)=`\iota_s\circ K_{sig}(A,\partial^A,O)=`\Phi_{(A,\partial^A,O)} ]{400}{1}l
 \end{picture}
\end{center}
An element of $\widetilde{A}$ is a pair $\lk a,h\rk$ such that $a\in A$ and $h:(|a|]\ra p^{-1}(a)=\{\lk a,i\rk : i\in (|a|] \}$ is a bijection.
An element of $(s(A)$ is a pair $\lk a,\tau\rk$ so that $a\in A$ and  $\tau\in S_{|a|}$. Thus we can put
\[ \Phi_0(a,h) = \lk a,\pi_2\circ h\rk \]
with $\pi_2(a,i)=i$, for $i\in (|a|]$.  Clearly,  $\Phi_0$ is a bijection.
An element of $\widetilde{E^A}$ is a triple $\lk a,h,i\rk$ so that $\lk a,h\rk\in \widetilde{A}$ and $i\in (|a|]$.
An element of $E^{s(A)}$ is a triple $\lk a,\tau,i\rk$ such that $\lk a,\tau\rk\in s(A)$ and $i\in (|a|]$.
Clearly, we put  $\Phi_1(a,h) = \lk a,\pi_2\circ h,i\rk$, and  $\Phi_1$ is a bijection as well.

We have

 \vskip 2mm

\begin{proposition} \label{Phi-nt}
The transformation $\Phi :  K_{diag}\circ \iota_a\ra \iota_s\circ K_{sig}$ defined above is a natural isomorphism.
\end{proposition}

{\it Proof.}~ We have already seen that the components of $\Phi$ are isomorphisms. The verification that $\Phi$ is natural
is left for the reader. $~\Box$

\vskip 2mm

Finally, we define the natural isomorphism $\Psi$. We fix a polynomial diagram
\begin{center} \xext=1500 \yext=200
\begin{picture}(\xext,\yext)(\xoff,\yoff)
  \putmorphism(0,50)(1,0)[O`E`s]{500}{-1}a
  \putmorphism(500,50)(1,0)[\phantom{E}`A`p]{500}{1}a
  \putmorphism(1000,50)(1,0)[\phantom{A}`O`t]{500}{1}a
\end{picture}
\end{center}
in $\cP oly\cD iag_O$. We need to define a morphism
\[ \Psi_{(t,p,s)} :  K_{fu}\circ rep_{pd} (t,p,s)\lra rep_{and}\circ K_{diag}(t,p,s)\]
in $\cA n$ in the fibre over $O$, i.e. a natural transformation
\[ \Psi_{(t,p,s)} : t_!p_*s^* \lra (orb_{/O})\tilde{t}_!\tilde{p}_*\tilde{s}^*(\delta_{/O})  \]
between endofunctors on $Set_{/O}$. To this end, we need to define its components
\[ (\Psi_{(t,p,s)})_{(X,d^X)} : t_!p_*s^*(X,d^X) \lra (orb_{/O})\tilde{t}_!\tilde{p}_*\tilde{s}^*(\delta_{/O})(X,d^X)  \]
for any object $(X,d^X:X\ra O)$ in $Set_{/O}$.  Fix $(X,d^X)$ in $Set_{/O}$.
An element of $t_!p_*s^*(X,d^X)$ is a pair $\lk a, \vec{x}\rk$ so that $a\in A$, $\vec{x}: (|a|]\ra X$ is a function such that
$d^X\circ \vec{x}= s_{\lc p^{-1}(a)}$.
An element $[\lk a,h,\vec{x} \rk ]_\sim$ of $(orb_{/O})\tilde{t}_!\tilde{p}_*\tilde{s}^*(\delta_{/O})(X,d^X)$ is an equivalence class of triples
$\lk a,h,\vec{x} \rk$ so that $\lk a,\vec{x} \rk$ is an element of $t_!p_*s^*(X,d^X)$ and $h:(|a|]\ra p^{-1}(a)$ is a bijection. The action of $S_{|a|}$ is defined so
that $\lk a,h,\vec{x} \rk\cdot \tau =\lk a,h\circ\tau,\vec{x} \rk$. Thus any two  triples $\lk a,h,\vec{x} \rk$ and $\lk a',h',\vec{x'} \rk$ are identified
if and only if $a=a'$ and $\vec{x}=\vec{x'}$.  Thus we can put
\[  (\Psi_{(t,p,s)})_{(X,d^X)}(a,\vec{x}) =\{ \lk a,h,\vec{x} \rk : h:(|a|]\stackrel{\cong}{\lra} p^{-1}(a)  \} \]
i.e.  we associate to  $\lk a,\vec{x} \rk$ the equivalence class of all triples whose second component is a bijection of  $(|a|]$ and $p^{-1}(a)$.
As these sets have, by definition, the same number of elements, $\Psi$ is well defined.

We have

 \vskip 2mm

\begin{proposition} \label{Psi-nt}
The transformation $\Psi :   K_{fu}\circ rep_{pd} \lra rep_{and}\circ K_{diag}$ defined above is a natural isomorphism.
\end{proposition}

{\it Proof.}~ As before, the verification that $\Psi$ is natural is left for the reader.
From the considerations above it should be clear that for any polynomial diagram $(t,p,s)$ and $(X,d^X)$ in $Set_{/O}$
$(\Psi_{(t,p,s)})_{(X,d^X)}$ is a bijection. So $(\Psi_{(t,p,s)})$ is a natural isomorphism and hence $\Psi$ is an isomorphism, as well.
$~\Box$

In that way we have completed the description of the diagram of categories, functors and natural transformations from the beginning of this subsection.
Thus, we know that the whole diagram commutes (at least up to an equivalence), moreover the named horizontal
functors are equivalences of categories. As we note, Proposition \ref{Ksig}, $K_{sig}$ is full and faithful. Therefore both
$K_{diag}$ and $K_{fu}$ are full and faithful, as well.  Hence using the characterizations of fibrations of polynomial and analytic functors, Proposition \ref{im-rep-a}, Theorem \ref{analytic char}, we obtain a statement, a bit surprising at first sight.

\begin{corollary} \label{emb-full}
Any weakly cartesian natural transformation between polynomial functors is cartesian.
\end{corollary}

\section{Appendix}
We spell below in detail some well known definitions of various adjoint functors between slices of $Set$ and $\sigma Set$.

First, recall that the unique functor $\cS_*\ra 1$  induces by composition the diagonal functor $\delta$ that has both adjoints
 \begin{center} \xext=800 \yext=300
\begin{picture}(\xext,\yext)(\xoff,\yoff)
\putmorphism(0,120)(1,0)[Set`\sigma Set`\delta]{800}{1}a
\putmorphism(0,50)(1,0)[\phantom{Set}`\phantom{\sigma Set}`fix]{800}{-1}b
\putmorphism(0,240)(1,0)[\phantom{Set}`\phantom{\sigma Set}`orb]{800}{-1}a
\end{picture}
\end{center}
$orb\dashv\delta\dashv fix$. The functor $\delta$ sends set $X$ to $\o\times X$,
i.e. to $\o$ copies of $X$, with $n$-th copy of $X$ equipped with a trivial action of $S_n$.
The functor $orb$ sends a symmetric set $(A,\alpha)$ to the set of its orbits with respect to all actions $A_{/\alpha}$.
The functor $fix$ sends a symmetric set $(A,\alpha)$ to the product over $\o$  of the sets of fix points with respect to each action $S_n$, i.e.
\[ fix(A,\alpha) =\prod_{n\in \o} fix_n(A_n,\alpha_n) \]
where $fix_n(A_n,\alpha_n)= \{ a\in A_n : a\cdot \sigma =a\;\; {\rm for}\;\; \sigma\in S_n \}$.
The functor $fix$ is not used directly but its existence shows that
$\delta$ preserves colimits.

For any set $O$ the above adjunction can be sliced, i.e. we have functors
 \begin{center} \xext=1000 \yext=400
\begin{picture}(\xext,\yext)(\xoff,\yoff)
\putmorphism(0,140)(1,0)[Set_{/O}`\sigma Set_{/\delta(O)}`\delta_{/O}]{1000}{1}a
\putmorphism(0,50)(1,0)[\phantom{Set_{/O}}`\phantom{\sigma Set_{/\delta(O)}}`fix_{/O}]{1000}{-1}b
\putmorphism(0,290)(1,0)[\phantom{Set_{/O}}`\phantom{\sigma Set_{/\delta(O)}}`orb_{/O}]{1000}{-1}a
\end{picture}
\end{center}
such that  $orb_{/O}\dashv \delta_{/O}\dashv fix_{/O}$.
For $d^X:X\ra O$ in $Set_{/O}$
\[ \delta_{/O}(X,d^X)= \delta(d^X):\delta(X)\ra \delta(O) \]
is the sliced diagonal functor. Moreover, for  $d^Y:(Y,\zeta)\ra \delta(O)$ in $\sigma Set_{/\delta(O)}$, we have
\[ orb_{/O}((Y,\zeta),d^Y) : orb(Y,\zeta) \lra O \]
so that $orb_{/O}((Y,\zeta),d^Y)([y]_\sim)=o$ if $d^Y(y)=\lk n,o\rk$ for some $n\in\o$.
Finally, for $((Y,\zeta),d^Y)$ as above
\[ fix_{/O}((Y,\zeta),d^Y):\coprod_{o\in O}\prod_{n\in \o} fix_{n,o}((Y,\zeta),d^Y) \lra O \]
is the obvious projection function, where \[ fix_{n,o}((Y,\zeta),d^Y)= \{ y\in Y_n : d^Y(y)=o,\, y\cdot \sigma =y\;\; {\rm for}\;\; \sigma\in S_n \}\]

Next we recall the pullback functor and its adjoint in the category of symmetric sets $\sigma Set$.
Any morphism $p:(E,\varepsilon)\ra (A,\alpha)$ in $\sigma Set$ induces three functors
 \begin{center} \xext=1100 \yext=300
\begin{picture}(\xext,\yext)(\xoff,\yoff)
\putmorphism(0,120)(1,0)[\sigma Set _{/(E,\varepsilon)}`\sigma Set_{/(A,\alpha)}`p^*]{1100}{-1}a
\putmorphism(0,50)(1,0)[\phantom{\sigma Set_{/(E,\varepsilon)}}`\phantom{\sigma Set_{/(A,\alpha)}}`p_*]{1100}{1}b
\putmorphism(0,240)(1,0)[\phantom{\sigma Set_{/(E,\varepsilon)}}`\phantom{\sigma Set_{/(A,\alpha)}}`p_!]{1100}{1}a
\end{picture}
\end{center}
so that $p_!\dashv p^* \dashv p_*$.  $p^*$ is defined by pulling back along $p$,  $p_!$ is defined by composing with $p$. The actions are
defined in the obvious way. For $d^X:(X,\xi)\ra (E,\varepsilon)$ the universe of  $p_*((X,\xi),d^X)$ is
\[ \{ \lk a,\vec{x}: E_a\ra (X,\xi)\rk : a\in A,\, d^X\circ \vec{x}=i_a \} \]
where $E_a$ and $i_a$ are defined from the following pullback in $\sigma Set$
\begin{center} \xext=600 \yext=450
\begin{picture}(\xext,\yext)(\xoff,\yoff)
 \setsqparms[1`-1`-1`1;600`400]
 \putsquare(0,0)[(E,\varepsilon)`(A,\alpha)`E_a`S_{|a|};p`i_a`\bar{a}`]
 \end{picture}
\end{center}
and $\bar{a}:S_{|a|}\lra (A,\alpha)$ is the morphism from the symmetric set $S_{|a|}$ (with action of $S_{|a|}$ on the right) sending identity
on $(|a|]$ to $a$.

The action in $p_*((X,\xi),d^X)$ is defined by conjugation
\[ \lk a,\vec{x} \rk \cdot \sigma = \lk a\cdot\sigma, \vec{x}((-)\cdot \sigma^{-1})\cdot \sigma \rk \]
and the typing sends $\lk a,\vec{x} \rk$ to $a$.

Thus we can draw a diagram of categories and functors
\begin{center} \xext=1200 \yext=1250
\begin{picture}(\xext,\yext)(\xoff,\yoff)
\setsqparms[-1`-1`-1`-1;1200`800]
\putsquare(230,200)[\sigma Set_{/\delta(O)}`\sigma Set_{/\delta(Q)}`Set_{/O}`Set_{/Q};
\delta(u)^*`\delta_{/O}`\delta_{/Q}`u^*]

\setsqparms[0`1`1`0;1660`800]
\putsquare(0,200)[\phantom{\sigma Set_{/\delta(O)}}`\phantom{\sigma Set_{/\delta(Q)}}`\phantom{\sigma Set_{O}}`\phantom{\sigma Set_{Q}};
`orb_{/O}`fix_{/Q}`]

\setsqparms[1`0`0`1;1200`1100]
\putsquare(230,50)[\phantom{\sigma Set_{/\delta(O)}}`\phantom{\sigma Set_{/\delta(Q)}}`\phantom{\sigma Set_{O}}`\phantom{\sigma Set_{Q}};
\delta(u)_!```u_*]

\putmorphism(350,1000)(0,-1)[\phantom{\sigma Set_{/\delta(O)}}`\phantom{\sigma Set_{O}}`fix_{/O}]{800}{1}r
\putmorphism(1300,1000)(0,-1)[\phantom{\sigma Set_{/\delta(O)}}`\phantom{\sigma Set_{O}}`orb_{/Q}]{800}{1}l

\putmorphism(230,900)(1,0)[\phantom{\sigma Set_{/\delta(O)}}`\phantom{\sigma Set_{/\delta(Q)}}`\delta(u)_*]{1200}{1}b
\putmorphism(230,300)(1,0)[\phantom{\sigma Set_{/O}}`\phantom{\sigma Set_{/Q}}`u_!]{1200}{1}a
 \end{picture}
\end{center}
in which we have a natural isomorphism of functors \[ \delta(u)^*\circ \delta_{/Q} \cong \delta_{/O}\circ u^*\] and hence of their left
\[ orb_{/Q}\circ \delta(u)_! \cong u_! \circ orb_{/O}\]
and right adjoints
\[ fix_{/Q}\circ \delta(u)_! \cong u_! \circ fix_{/O}\]


\end{document}